\documentclass[reqno,11pt]{amsart}

\usepackage{amsthm, amsmath, amssymb}
\usepackage[utf8]{inputenc}
\usepackage[T1]{fontenc}
\usepackage{tipa}
\usepackage{stmaryrd}
\usepackage{thmtools}

\usepackage[unicode, hidelinks, hypertexnames=false]{hyperref}
\usepackage{microtype}

\usepackage[margin=1in]{geometry}
\usepackage[shortlabels]{enumitem}
\usepackage[normalem]{ulem}

\usepackage[noabbrev,capitalize]{cleveref}
\crefname{equation}{Eq.}{Eqs.}

\renewcommand{\thepart}{\Roman{part}}

\let\oldtocsection=\tocsection
\let\oldtocsubsection=\tocsubsection
\let\oldtocsubsubsection=\tocsubsubsection

\renewcommand{\tocsection}[2]{\hspace{0em}\oldtocsection{#1}{#2}}
\renewcommand{\tocsubsection}[2]{\hspace{1em}\oldtocsubsection{#1}{#2}}
\renewcommand{\tocsubsubsection}[2]{\hspace{2em}\oldtocsubsubsection{#1}{#2}}

\newtheorem{scthm}{Theorem}[section]

\newtheorem{theorem}[scthm]{Theorem}
\newtheorem{proposition}[scthm]{Proposition}
\newtheorem{lemma}[scthm]{Lemma}
\newtheorem{claim}[scthm]{Claim}
\newtheorem{corollary}[scthm]{Corollary}
\newtheorem{conjecture}[scthm]{Conjecture}

\theoremstyle{definition}
\newtheorem{definition}[scthm]{Definition}

\newtheorem{example}[scthm]{Example}

\theoremstyle{remark}
\newtheorem*{remark}{Remark}

\crefname{claim}{Claim}{Claims}

\newcommand{\abs}[1]{\left\lvert#1\right\rvert}
\newcommand{\norm}[1]{\left\lVert#1\right\rVert}

\newcommand{\sang}[1]{\langle #1 \rangle}

\newcommand{\paren}[1]{\left( #1 \right)}

\newcommand{\set}[1]{\left\{ #1 \right\}}

\newcommand{\ol}{\overline}

\DeclareMathOperator{\spn}{span}
\DeclareMathOperator{\proj}{proj}

\DeclareMathOperator{\dist}{dist}
\DeclareMathOperator{\rank}{rank}
\DeclareMathOperator{\codim}{codim}
\DeclareMathOperator{\Flec}{Flec}

\DeclareMathOperator{\tr}{tr}
\DeclareMathOperator{\Id}{Id}

\DeclareMathOperator{\len}{length}
\DeclareMathOperator{\Phy}{Phy}
\DeclareMathOperator{\lcm}{lcm}
\DeclareMathOperator{\LT}{HT}
\DeclareMathOperator{\LM}{HM}
\DeclareMathOperator{\LC}{HC}

\newcommand{\aff}{\mathrm{aff}}

\newcommand{\CC}{\mathbb{C}}

\newcommand{\FF}{\mathbb{F}}
\newcommand{\RR}{\mathbb{R}}

\newcommand{\OO}{\mathbb{O}}
\newcommand{\ZZ}{\mathbb{Z}}
\newcommand{\PP}{\mathbb{P}}

\renewcommand{\SS}{\mathbb{S}}

\newcommand{\F}{\mathbb{F}}
\newcommand{\G}{\mathbb{G}}
\renewcommand{\P}{\mathbb{P}}
\newcommand{\R}{\mathbb{R}}
\renewcommand{\S}{\mathbb{S}}

\newcommand{\Z}{\mathbb{Z}}

\newcommand{\fm}{\mathfrak{m}}
\newcommand{\fp}{\mathfrak{p}}
\newcommand{\fq}{\mathfrak{q}}
\newcommand{\fr}{\mathfrak{r}}

\newcommand{\fF}{\mathfrak{F}}
\newcommand{\fL}{\mathfrak{L}}
\newcommand{\fP}{\mathfrak{P}}
\newcommand{\fR}{\mathfrak{R}}

\newcommand{\cC}{\mathcal C}
\newcommand{\cD}{\mathcal D}
\newcommand{\cE}{\mathcal E}
\newcommand{\cF}{\mathcal F}
\newcommand{\cG}{\mathcal G}
\newcommand{\cI}{\mathcal I}

\newcommand{\cL}{\mathcal L}
\newcommand{\cM}{\mathcal M}
\newcommand{\cO}{\mathcal O}
\newcommand{\cP}{\mathcal P}
\newcommand{\cQ}{\mathcal Q}
\newcommand{\cR}{\mathcal R}
\newcommand{\cS}{\mathcal S}
\newcommand{\cT}{\mathcal T}

\newcommand{\cV}{\mathcal V}
\newcommand{\cX}{\mathcal X}
\newcommand{\cZ}{\mathcal{Z}}

\newcommand{\cFb}{\overline{\mathcal F}}

\newcommand{\tj}{\tilde\jmath}

\renewcommand{\ss}{\subseteq}
\newcommand{\ns}{\not\subseteq}

\newcommand{\new}{\mathrm{new}}
\renewcommand{\epsilon}{\varepsilon}

\def\mc{\mathcal }

\makeatletter
\renewcommand\part{%
   \if@noskipsec \leavevmode \fi
   \par
   \addvspace{4ex}%
   \@afterindentfalse
   \secdef\@part\@spart}

\def\@part[#1]#2{%
    \ifnum \c@secnumdepth >\m@ne
      \refstepcounter{part}%
      \addcontentsline{toc}{part}{\thepart\hspace{1em}#1}%
    \else
      \addcontentsline{toc}{part}{#1}%
    \fi
       \begin{center}\Large \partname\nobreakspace\thepart.\nobreakspace\scshape #2\end{center}%
    \nobreak
    \vskip 1.5ex
    \@afterheading}
\def\@spart#1{%
    {\parindent \z@ \raggedright
     \interlinepenalty \@M
     \normalfont
     \huge \bfseries #1\par}%
     \nobreak
     \vskip 3ex
     \@afterheading}
\makeatother

\title{The Erd\H{o}s Distinct Distances Problem in \texorpdfstring{$\R^3$}{R3}}

\author[Tidor]{Jonathan Tidor}
\author[Yu]{Hung-Hsun Hans Yu}
\author[Zakharov]{Dmitrii Zakharov}

\address{Department of Mathematics, Princeton University, Princeton, NJ 08544, USA}
\email{\{jtidor,hansonyu\}@princeton.edu}

\address{Department of Mathematics, Massachusetts Institute of Technology, Cambridge, MA 02139, USA}
\email{zakhdm@mit.edu}

\begin{document}

\begin{abstract}
We prove that $N$ points in $\R^3$ determine at least $N^{2/3-o(1)}$ distinct distances.
\end{abstract}

\maketitle

\tableofcontents

\section{Introduction}
\label{sec:intro}

One of the most well-known problems in discrete geometry is the Erd\H{o}s distinct distances problem which asks for the minimum number of distinct distances determined by a set of $N$ points in $\R^d$. That is, given a collection of $N$ points in $\R^d$, measure the Euclidean distance between each pair of points; what is the smallest possible size that the resulting set of distances can attain? In 1946, Erd\H{o}s observed that an $N^{1/d}\times\cdots\times N^{1/d}$ section of the integer lattice $\Z^d$ determines on the order of $N/\sqrt{\log N}$ distinct distances for $d=2$ and on the order of $N^{2/d}$ distinct distances for $d\geq 3$ \cite{Erd46}. This construction is conjectured to be asymptotically optimal.

The Erd\H{o}s distinct distances problem is one of many problems in discrete geometry where one aims to extremize a local quantity over discrete geometric configurations. For many of these problems it is believed that the optimal examples have the global structure of a lattice, and moreover, that the integer lattice is an asymptotically optimal example. Such settings include the Szemer\'edi--Trotter theorem on point-line incidences, the Erd\H{o}s unit distance problem in $\R^2$ and $\R^3$, the Erd\H{o}s distinct distances problem in $\R^d$, and the sum-product problem over $\R$. 

The Szemer\'edi--Trotter theorem confirms this belief for point-line incidences \cite{ST83}.\footnote{Note however, there are examples based on other lattices which perform as well as the integer lattice \cite{GS25,Cur25}. Furthermore, it is not known that lattice structure is necessary in asymptotically optimal examples.}
For the distinct distances problem in $\R^2$, Erd\H{o}s originally showed that $N$ points determine at least $N^{1/2}$ distinct distances \cite{Erd46} and this exponent was improved via a long series of works \cite{Mos52,Chu84,CST92,Sze97,ST01,Tar03,KT04}. Finally, breakthrough work of Guth and Katz showed that $N$ points in $\R^2$ determine at least on the order of $N/\log N$ distinct distances \cite{GK15}. This implies that the integer lattice is within a $\sqrt{\log N}$ factor of the optimal construction.

In contrast to these results, very recent work disproves this conjecture both for the unit distance problem in $\R^2$ (see \cite{ABGLS26,Saw26}) and for the sum-product problem \cite{BSSZ26}, showing that other lattices can outperform the integer lattice by a polynomial factor.

The distinct distances problem in dimensions $d\geq 3$ has resisted improvement. Prior to this work, the best bounds on this problem came from a result of Solymosi--Vu which bootstraps a distinct distances bound in lower dimension to one in higher dimension \cite{SV08}. Applying this result with the Guth--Katz bound implies that $N$ points in $\R^3$ determine at least on the order of $N^{3/5}/\log^{2/5} N$ distinct distances.

In this paper, we resolve the distinct distances problem in $\R^3$ up to an $N^{o(1)}$ factor, proving that the integer lattice determines, up to subpolynomial factors, the minimal number of distinct distances in $\R^3$.

\begin{theorem}
\label{thm:main}
There exists a function $\epsilon=\epsilon(N)$ satisfying $\epsilon(N)\lesssim \sqrt{\log \log N/\log N}$ so that any set of $N$ points in $\R^3$ determines at least $N^{2/3-\epsilon}$ distinct distances.
\end{theorem}

To prove this result, we build upon the polynomial method techniques introduced by Guth and Katz in their solution to the 2-dimensional problem as well as later extensions of these techniques due to Walsh \cite{Wal20,Wal23}. We develop new polynomial method tools in order to study incidence geometry problems involving higher-dimensional objects. These include the resolution of an algebro-geometric conjecture of Walsh \cite[Conjecture A]{Walsh22talk} which he highlighted as an important step towards attacking the distinct distances problem. Using these tools, we then prove a near-optimal point-flat incidence bound.

\subsection{The Guth--Katz argument in \texorpdfstring{$\mathbb{R}^2$}{R2}}
The breakthrough work of Guth and Katz on the distinct distances conjecture in $\R^2$ relies on an elegant reduction, due to Elekes and Sharir, to a point-line incidence problem in $\R^3$ \cite{ES11}. Given a set of points $\cP\subset\R^2$, we lower bound the number of distinct distances determined by $\cP$ by upper bounding the \emph{distance energy}:
\[\abs{\set{(p,q,p',q')\in\cP^4:\norm{p-p'}=\norm{q-q'}}}.\]
The key observation is that $(p,q,p',q')$ is a \emph{distance quadruple}, i.e., $\norm{p-p'}=\norm{q-q'}$, if and only if there exists a (unique) orientation-preserving rigid motion of the plane which sends $p$ to $q$ and $p'$ to $q'$. The group $E^+(2)$ of orientation-preserving rigid motions of the plane is 3-dimensional\footnote{One can parametrize $E^+(2)$ by $(x,y,\theta)$, representing a rotation counterclockwise by $\theta$ around the point $(x,y)$.} and thus each distance quadruple corresponds to a point in this 3-dimensional space. Furthermore, the set of rigid motions which send $p$ to $q$ forms a subset $\ell_{pq}\subset E^+(2)$ which is 1-dimensional.\footnote{Any such motion can be written as the translation sending $p$ to $q$ composed with a rotation around $q$ by an arbitrary angle $\theta$.} Now $(p,q,p',q')$ is a distance quadruple if and only if the sets $\ell_{pq}$ and $\ell_{p'q'}$ intersect. Guth and Katz showed how to choose coordinates identifying a dense open subset of $E^+(2)$ with $\R^3$ in such a way that the sets $\ell_{pq}$ are identified with lines in $\R^3$. In this way, the distinct distance problem in the plane reduces to bounding the number of incident pairs of lines among the set $\{\ell_{pq}\}_{p,q\in\cP}$ of $|\cP|^2$ lines in $\R^3$.

Given a set $\cL$ of lines, we use $\cP_k(\cL)$ to denote the set of \emph{$k$-rich points}: the set of points lying on at least $k$ lines from $\cL$. The classical Szemer\'edi--Trotter theorem \cite{ST83} gives an asymptotically tight bound on the number of $k$-rich points determined by a set of $L$ lines in $\R^2$. For $2\leq k\leq L^{1/2}$, the number of $k$-rich points is $O(L^2/k^3)$ and there are examples matching this bound. Without further restrictions, this bound is also tight in $\R^3$ as any example in $\R^2$ can be placed inside $\R^3$. To prove their result, Guth and Katz showed that if the collection of lines is ``truly 3-dimensional'' then the number of $k$-rich points is much smaller: $O(L^{3/2}/k^2)$ instead of $O(L^2/k^3)$. Here ``truly 3-dimensional'' means a bound on the \emph{algebraic concentration} of the lines: in this case, the hypothesis that at most $O(L^{1/2})$ lines lie in any plane or regulus.\footnote{A regulus is an algebraic surface of degree 2 (specifically a one-sheeted hyperbola or a hyperbolic paraboloid) that is doubly ruled by lines.} To prove this incidence bound for lines in $\R^3$, Guth and Katz developed several new tools, including a polynomial method result called real polynomial partitioning. To complete the proof, Guth and Katz showed that the set of lines $\{\ell_{pq}\}_{p,q\in\cP}$ satisfies the appropriate bound on the algebraic concentration.

\subsection{First attempts to extend to \texorpdfstring{$\R^3$}{R3}}
We now discuss the most natural extension of this proof strategy to the distinct distances problem in $\R^3$ and the issues that arise. The group $E^+(3)$ of orientation-preserving rigid motions of $\R^3$ is 6-dimensional. If $(p,q,p',q')$ is a distance quadruple with $(p,q)\neq (p',q')$, then there is a 1-dimensional family of rigid motions which send $p$ to $q$ and $p'$ to $q'$. These correspond to taking a fixed map with this property and composing with a rotation about the line $\aff\{q,q'\}$ through $q$ and $q'$. For each pair $p,q$, the set of rigid motions which send $p$ to $q$ form a 3-dimensional subset $F_{pq}\subset E^+(3)$. Bardwell-Evans and Sheffer showed how to choose coordinates identifying a dense open subset of $E^+(3)$ with $\R^6$ so that each $F_{pq}$ is identified with a 3-flat \cite{BES19}. Thus the problem reduces to taking a collection of 3-flats in $\R^6$ and bounding the number of pairs which intersect at a line. One way to study this problem (though not the one we will eventually use) is to intersect with a generic hyperplane. This results in a point-plane incidence bound in $\R^5$.

Let us focus on bounding the number of 2-rich points determined by this set. Given a collection $\cF$ of planes (i.e., 2-flats) in $\R^n$, write $\cP_2(\cF)$ for the set of 2-rich points. One immediate problem is that $\cP_2(\cF)$ may be infinite: if any two of the planes intersect in a line, then all the points on that line are 2-rich for $\cF$. Fortunately this issue does not occur in our setting: the collection of 3-flats $\{F_{pq}\}_{p,q\in\cP}$ is \emph{1-transverse}, meaning that all pairwise intersections have dimension at most 1. Thus the collection of planes obtained by intersecting with a generic hyperplane is necessarily transverse (i.e., any two distinct planes in the collection are either disjoint or intersect in a single point). Now given a collection $\cF$ of transverse planes in $\R^n$, let us try to bound the size of $\cP_2(\cF)$. The trivial upper bound of $\binom {|\cF|}2$ is achieved for a generic collection of planes in $\R^4$. One might hope that the number of 2-rich points is significantly smaller for a collection of planes in $\R^5$ with small algebraic concentration.

For a set of flats $\cF$ in $\R^n$, define the \emph{$k$-th algebraic concentration}
\begin{equation}
\label{eq:D_k}
\cD_k(\cF)=\max_{V\subset\R^n:\dim V=k}\frac{\abs{\set{F\in\cF:F\subseteq V}}}{\deg V},
\end{equation}
where the maximum is taken over irreducible $k$-dimensional varieties $V\subseteq\R^n$. These quantities were introduced by Walsh in a work which generalized the Guth--Katz incidence bound to a point-line incidence bound in $\R^n$ for all $n\geq 2$, proving the following result \cite{Wal23}. (For simplicity, we state only one special case of the result here.)

\begin{theorem}[{\cite[Theorem 1.3]{Wal23}}]
Let $\cL$ be a set of lines in $\F^n$. Then
\[|\cP_2(\cL)|\lesssim_n\sum_{k=2}^n |\cL|\cD_k(\cL)^{\tfrac1{k-1}}.\]
\end{theorem}

Walsh's techniques are one of very few techniques which give tight incidence bounds in high dimensions. Note that to apply the Guth--Katz incidence bound, one must study how the lines cluster in surfaces of degree at most 2; in contrast, to apply Walsh's incidence bound, one must study how the lines cluster in varieties of any degree, weighted by the degree of the variety.

Given $\cF$, a set of transverse planes in $\R^5$, we have
\[1=\cD_2(\cF)\leq\cD_3(\cF)\leq\cD_4(\cF)\leq\cD_5(\cF)=|\cF|.\]
One way for $\cF$ to have a large number of 2-rich points is for $\cD_4(\cF)$ to be large. For example, if all the planes lie in a 4-flat, we have $\cD_4(\cF)=|\cF|$ and there can be $\binom {|\cF|}2$ 2-rich points. There is a similar potential issue if $\cD_3(\cF)$ is large.\footnote{It is harder to construct an example in this case as the transverse assumption implies that each intersection occurs at a singular point of the variety, forcing any such example to live in a very high degree and very singular variety.} Our goal will be to prove that $|\cP_2(\cF)|\lesssim |\cF|^{5/3}$ under the assumption that $\cD_3(\cF)\lesssim |\cF|^{1/3}$ and $\cD_4(\cF)\lesssim |\cF|^{2/3}$. (Using parameter counting, one can show that these bounds are the smallest values $\cD_3(\cF)$ and $\cD_4(\cF)$ can possibly take.) However, these two concentration bounds are not sufficient. As another example, consider a set of $|\cF|$ lines in a plane in $\R^5$ with $\binom {|\cF|}2$ 2-rich points. Generically extending each of these lines to a plane in $\R^5$, we obtain a set of transverse planes in $\R^5$ with $\cD_3(\cF),\cD_4(\cF)$ small, yet there are $\binom {|\cF|}2$ 2-rich points.

To rule out this example, we also need to consider the \emph{$k$-th codimension 1 algebraic concentration}
\begin{equation}
\label{eq:D1_k}
\cD_k^1(\cF)=\max_{V\subset\R^n:\dim V=k}\frac{\sum_{F\in\cF}\sum_{\gamma\subseteq F\cap V:\dim \gamma=1}\deg\gamma}{\deg V}
\end{equation}
where the second sum is over the irreducible 1-dimensional components $\gamma$ of $F\cap V$. We have the bounds
\[1=\cD_1^1(\cF)\leq\cD_2^1(\cF)\leq\cD_3^1(\cF)\leq\cD_4^1(\cF)=|\cF|.\]

One of the new contributions of this paper is an incidence bound which implies that if all of the algebraic concentrations are bounded, then $\cF$ determines few 2-rich points. In particular, using our methods, we can show that if $\cF$ is a collection of transverse planes in $\R^5$ satisfying $\cD_2^1(\cF),\cD_3(\cF)\lesssim |\cF|^{1/3}$ and $\cD_3^1(\cF),\cD_4(\cF)\lesssim |\cF|^{2/3}$, then $|\cP_2(\cF)|\lesssim |\cF|^{5/3}\log^2 |\cF|$. This incidence bound is tight up to the polylogarithmic factor. It is one of very few results which give near-tight bounds for incidence problems involving higher dimensional objects (see \cite{ST12,TYZ22} for two other examples). 

With this incidence bound, one potential strategy to proving that every set $\cP$ of $N$ points in $\R^3$ determines $N^{2/3-o(1)}$ distinct distances is to show that the distance energy of $\cP$ is bounded by $N^{10/3+o(1)}$. Letting $\cF$ be the set of $N^2$ planes in $\R^5$ formed by intersecting $\{F_{pq}\}_{p,q\in\cP}$ by a generic hyperplane, we would need to show that $|\cP_2(\cF)|\leq N^{10/3+o(1)}$. This follows if $\cD_2^1(\cF),\cD_3(\cF)\lesssim N^{2/3}$ and $\cD_3^1(\cF),\cD_4(\cF)\lesssim N^{4/3}$. We will prove the codimension 0 concentration bounds:  $\cD_3(\cF)\lesssim N^{2/3}$ and $\cD_4(\cF)\lesssim N^{4/3}$. However, we run into issues when studying the codimension 1 concentration. Our first issue is quite technical in nature: while studying $\cD_3^1(\cF)$, we can bound the codimension 1 concentration in 3-dimensional varieties of small degree, but are unable to handle varieties of very large degree. The second issue is more major and proves fatal to this proof strategy: the desired bound $\cD_2^1(\cF)\lesssim N^{2/3}$ is simply false; in fact, we can have $\cD_2^1(\cF)=N$.

The counterexample is geometric in nature and is easier to describe in $\R^6$, before intersecting with a generic hyperplane. Let $\tau$ be an orientation-\emph{reversing} rigid motion of $\R^3$. It turns out that $F_\tau$, defined to be the set of orientation-preserving rigid motions which agree with $\tau$ on a plane, is a 3-flat.\footnote{To see that this set should be 3-dimensional, note that any such rigid motion can be written as the composition of $\tau$ with the reflection across a plane in $\R^3$. As the space of planes in $\R^3$ is 3-dimensional, so is this set.} Furthermore, if $\tau(p)=q$, then $F_{pq}$ is a 3-flat such that $F_{pq}\cap F_\tau$ is a 2-flat.\footnote{To see the dimension, note that any element of $F_{pq}\cap F_\tau$ can be written as the composition of $\tau$ with a reflection across a plane containing $q$; this is a 2-dimensional set.} Now intersecting with a generic hyperplane $H$, we have the plane $H\cap F_\tau$ whose intersection with many of the planes $H\cap F_{pq}$ is a line. In particular, if $\tau$ is $k$-rich, meaning that there are at least $k$ pairs $(p,q)\in\cP^2$ with $\tau(p)=q$, we have produced a plane which intersects $k$ of the planes in $\cF$ in a line, showing that $\cD_2^1(\cF)\geq k$. Since it is possible that there exists a single rigid motion $\tau$ which is $N$-rich, we may be in the situation where $\cD_2^1(\cF)=N$.

\subsection{New ideas in \texorpdfstring{$\R^3$}{R3}}
We make several changes to the above proof strategy in order to resolve these two issues. First we work in $\R^6$ instead of intersecting with a generic hyperplane. This may not appear to materially change the problem, but it preserves the underlying geometric structure of the problem in a useful way. Let $\cF$ be the set of $N^2$ 3-flats in $\R^6$. Write $\cL_k(\cF)$ for the set of $k$-rich lines: those which lie on at least $k$ of the 3-flats of $\cF$. As remarked previously, these 3-flats are 1-transverse so $\cL_k(\cF)$ is finite for each $k\geq 2$. We introduce a trick based on Beck's theorem which reduces the problem to bounding $|\cL_2(\cF)|$ instead of $|\cL_k(\cF)|$ for all $k\geq 2$.

In \cref{part:i} of the argument, we prove an incidence bound which, from $\cF$, produces a set of surfaces $\cS$ of controlled total degree, so that
\[|\cL_2(\cF)|\lesssim N^{10/3}\log^2N +|\cL_2(\cS)|.\]
This bound holds for any set $\cF$ of $N^2$ 3-flats subject to the codimension 0 concentration bounds $\cD_4(\cF)\lesssim N^{2/3}$ and $\cD_5(\cF)\lesssim N^{4/3}$ without any hypothesis about the codimension 1 concentration bounds. (Note that we are dealing with 3-flats instead of planes here so we have replaced $\cD_3,\cD_4$ with $\cD_4,\cD_5$.) To prove this result, we extend the techniques of Walsh \cite{Wal23} which he used to prove tight point-curve incidence bounds to apply to higher-dimensional objects. We then show the codimension 0 concentration bounds $\cD_4(\cF)\lesssim N^{2/3}$ and $\cD_5(\cF)\lesssim N^{4/3}$ subject to some mild conditions on $\cP$. This part of the argument involves developing new tools for studying flats lying in high degree varieties, including a generalization of the Cayley--Monge--Salmon theorem for flats and new tools in algebraic geometry about approximate complete intersections which substantially generalize prior results of Walsh \cite{Wal20}. We discuss these tools more in \cref{ssec:outline-i}.

In \cref{part:ii} of the argument, we need to bound $|\cL_2(\cS)|$ where $\cS$ is a collection of surfaces in $\R^6$. These surfaces are produced in \cref{part:i} of the argument by intersecting the 3-flats $\cF$ with various polynomials of degree $N^{2/3}$. The main contribution to $|\cL_2(\cS)|$ comes from surfaces which are (a) singly ruled or doubly ruled or (b) planes. In \cref{part:ii} of the argument, we bound the contribution to $|\cL_2(\cS)|$ where $\cS$ is a collection of singly or doubly ruled surfaces. Such surfaces contain a 1-dimensional family of lines. A standard way of controlling $|\cL_2(\cS)|$ would be to intersect with a generic hyperplane; this produces a collection of curves in $\R^5$ for which we want to bound the number of 2-rich points. This quantity can be bounded by a result of Walsh \cite{Wal23} in terms of the algebraic concentration of the collection of curves. However, bounding these algebraic concentrations is essentially equivalent to bounding the codimension 1 concentrations $\cD_3^1(\cF),\cD_4^1(\cF)$ which is exactly the original issue that we ran into. 

Here is the first place where staying in $\R^6$, instead of intersecting with a generic hyperplane, is crucially used. Indeed, knowing that a line $\ell$ lies in two ruled surfaces $S,S'$ is significantly more information than seeing the same configuration intersected by a generic hyperplane. Each of $S,S'$ contains a 1-dimensional family of lines and $\ell$ must be a member of both of these families. Using this geometric information, we transfer the problem to a different parameter space where the problem has more structure. To study the resulting incidence problem, we develop a fairly novel technique that combines standard polynomial method techniques with more sophisticated algebraic geometry, including the so-called ``trisecant lemma.''

In \cref{part:iii}, it remains to bound $|\cL_2(\cS)|$ where $\cS$ is a collection of planes. Here the underlying geometry of the problem again gives us significantly more structure. Each of the planes in $\cS$ can be written in the form $F_{pq}\cap F_\tau$ where $\tau$ is a $k$-rich orientation-reversing rigid motion with $k\gg N^{2/3}$. To complete the proof, we must show that there cannot be too many $\tau$ which produce too many 2-rich lines $\cL_2(\cS)$. 

First we note that for $k\gg N^{2/3}$, we can turn the problem of bounding the number of $k$-rich $\tau$ into a 2-flat--3-flat incidence problem in $\R^6$. Indeed, as we discussed previously, $\tau$ is $k$-rich if there exist $k$ 3-flats $F_{pq}$, each of which intersects $F_\tau$ in a 2-flat. The polynomial method tools which we develop in \cref{part:i} also apply to this problem, showing that there are at most $k^{-3/2}N^3\log^{3/2}N$ such $\tau$. Unfortunately this bound is too weak for us to apply outright -- it implies that the contribution to $\cL_2(\cS)$ from $k$-rich orientation-reversing rigid motions is at most $k^{1/2}N^3\log^{3/2}N$. For every $k\gg N^{2/3}$, this is not strong enough to attain the bound of $N^{10/3+o(1)}$. We believe that this bound is not tight and later state a conjectured bound which would be strong enough to complete the proof (\cref{conj:2-flat--3-flat}). Unfortunately, this bound seems out of reach of these techniques.

Instead, we combine the weak bound with a very involved combinatorial and geometric argument. Assuming that there is a collection of very rich orientation-reversing rigid motions $\tau$ which cause $|\cL_2(\cS)|\gg N^{10/3}$, we find increasingly more geometric structure in the problem, eventually reaching a contradiction.

Let us point out that as we originally described the reduction, we work in $E^+(3)$, the 6-dimensional space of orientation-preserving rigid motions. However, we could symmetrically have worked in $E^-(3)$, the 6-dimensional space of orientation-\emph{reversing} rigid motions. Define $F^-_{pq}\subset E^-(3)$ to be the 3-flat consisting of orientation-reversing rigid motions which map $p$ to $q$. The arguments in \cref{part:i,part:ii} work equally well in $E^+(3)$ and $E^-(3)$. Working in $E^+(3)$, \cref{part:i,part:ii} successfully bound the contribution to $|\cL_2(\cF)|$ except from the contribution coming from exceptional 3-flats $F^+_\tau\subset E^+(3)$: the set of orientation-preserving rigid motions of $\R^3$ which agree with the orientation-reversing rigid motion $\tau$ on a plane. Working in $E^-(3)$, we can instead bound the contribution to $|\cL_2(\cF)|$ except from the contribution coming from exceptional 3-flats $F^-_\rho\subset E^-(3)$: the set of orientation-reversing rigid motions of $\R^3$ which agree with the orientation-preserving rigid motion $\rho$ on a plane. Thus in \cref{part:iii}, it suffices to bound the contribution to $|\cL_2(\cF)|$ which is captured by both an exceptional $F^+_\tau$ and an exceptional $F^-_\rho$. This extra information will be crucial for the argument.

In \cref{sec:outline}, we give a more thorough overview of the proof, discuss the new techniques we introduce, and remark on higher-dimensional generalizations. The proof formally begins in \cref{sec:setup} where we set up the notation we use and give the reduction of the distinct distance problem to an incidence problem in $\R^6$.

\vspace{7pt}\noindent\textbf{Acknowledgments.}
We are indebted to Yufei Zhao for his guidance early in this project's development. We thank Miguel Walsh, Joshua Zahl, and Ruixiang Zhang for several enlightening discussions. We thank Larry Guth for helpful comments on the manuscript.

Tidor was supported by NSF Award DMS-2554092 and partially supported by a Stanford Science Fellowship. 
Yu was partially supported by NSF Award DMS-2246682 and a Jane Street Graduate Fellowship.
He thanks Shuang-Yen Lee, Junyao Peng, and Yi-Tsung Wang for helpful conversations on algebraic geometry.
Zakharov was partially supported by a Jane Street Graduate Fellowship and a Simons Dissertation Fellowship. He thanks Ivan Karpov for instructive conversations on intersection theory.

Parts of this work took place at the Simons Laufer Mathematical Sciences Institute during the Spring 2025 semester, supported by NSF Award DMS-1928930; the conference ``Graph expansion and its applications,'' organized by Matija Buci\'c; the workshop on ``High-dimensional phenomena in discrete analysis,'' hosted by the American Institute of Mathematics; and the workshop ``Structural Results,'' hosted by the Simons Institute for the Theory of Computing. We thank SLMath, Matija Buci\'c, AIM, and the Simons Institute for their hospitality.

\vspace{7pt}\noindent\textbf{AI use statement.}
The mathematical content of this paper and all of the writing were generated solely by the authors. ChatGPT-5.6 was used during the editorial process to identify mathematical and typographical errors.
\section{Proof overview}
\label{sec:outline}

\subsection{Choosing coordinates}

In \cref{sec:setup}, we choose coordinates on $E^+$ and $E^-$, the spaces of orientation-preserving and orientation-reversing rigid motions of $\R^3$, respectively. We refer to these two spaces as \emph{positive space} and \emph{negative space}. We will show how to choose coordinates so that all the sets $F^+_{pq},F^+_\tau,F^-_{pq},F^-_\rho$ defined in \cref{sec:intro} are 3-flats in $E^+$ or $E^-$. (Starting in \cref{sec:setup} we will also distinguish between an affine open set $E^+_o\cong\R^6$ and $E^+$, which will be a carefully chosen compactification of $E^+_o$ containing all the orientation-preserving rigid motions and additional points at infinity.)
For the purpose of this proof overview, it suffices to treat $E^+,E^-$ as copies of $\R^6$ and to know that there is some choice of coordinates so that all the objects in question are 3-flats.

There is a 3-fold symmetry between $E^+,E^-$, and a third space $E^{\Phy}$ called \emph{physical space}, which is essentially the parameter space of pairs $(p,q)\in\R^3\times\R^3$.
To make the symmetry between the spaces more apparent, the pair $(p,q)$ is represented by the point $((p+q)/2,(p-q)/2)\in E^{\Phy}$. For the rest of the paper we use the notation $\phi$ for the map 
\[\phi(p,q)=((p+q)/2,(p-q)/2).\]

We add two more collections of 3-flats, $F^{\Phy}_\rho, F^{\Phy}_\tau$ to the previously defined ones. We defined $F^+_{pq}$, a \emph{physical 3-flat} in positive space, to be the collection of $\rho$ so that $\rho(p)=q$. Also, $F^+_\tau$, a \emph{negative 3-flat} in positive space, was defined to be the collection of $\rho$ which agree with $\tau$ on a plane. To complete the set, we define $F^{\Phy}_\rho$, a \emph{positive 3-flat} in physical space, to be essentially the graph of the function $\rho$: it is the collection of $\phi(p,q)$ where $\rho(p)=q$. The other 3-flats, $F^-_{pq},F^-_\rho,F^{\Phy}_\tau$ are defined symmetrically.

There is a strong symmetry between these collections of 3-flats: for example, $\rho\in F^+_{pq}$ if and only if $\phi(p,q)\in F^{\Phy}_\rho$. These 3-flats have a number of nice properties which can be checked geometrically (and will be checked rigorously in \cref{sec:setup}). For example, two 3-flats from the same collection have intersection either a line or empty, while two 3-flats from different collections have intersection either a plane or a point.

One final object that plays a central role in the proof is a special collection of lines called \emph{isotropic lines}. We will give an intrinsic definition in \cref{sec:setup}, but for now we give an equivalent geometric definition. A line $\ell^{\Phy}$ in $E^{\Phy}$ connecting the points $\phi(p,q),\phi(p',q')\in\R^3\times\R^3$ is isotropic if $\norm{p-p'} = \norm{q-q'}$. A line $\ell^+$ in $E^+$ is isotropic if it is the 1-dimensional family of rigid motions $\rho$ such that $\rho(p)=q$ and $\rho(p')=q'$ for some fixed quadruple of points $(p,p',q,q')$ with $\norm{p-p'} = \norm{q-q'}$ and similarly for a line $\ell^-$ in $E^-$. From this definition there is a \emph{geometric triality}: a 1-to-1-to-1 correspondence between the isotropic lines in each of the three spaces. We always use the naming convention $(\ell^{\Phy},\ell^+,\ell^-)$ for a triple of corresponding lines. See \cref{tab:geometry} for a summary of these relations.

\renewcommand{\arraystretch}{1.4}
\begin{table}[t]
    \centering
    $\begin{array}{c | c | c}
         E^{\Phy} & E^+ & E^-  
         \\\hline\hline
         \dim F_{\rho}^{\Phy} \cap F_{\tau}^{\Phy} =2 & \rho \in F_\tau^+ & \tau \in F_\rho^- 
         \\\hline
         \phi(p,q) \in F_\tau^{\Phy} &  \dim F_{pq}^+ \cap F_\tau^+ =2 & \tau \in F_{pq}^-
         \\\hline
         \phi(p,q) \in F_\rho^{\Phy} & \rho \in F_{pq}^+ & \dim F_{pq}^- \cap F_\rho^- =2 
         \\\hline
         \phi(p,q), \phi(p',q') \in \ell^{\Phy} & \ell^+ \subset F_{pq}^+ \cap F_{p'q'}^{+} & \ell^{-} \subset F_{pq}^{-} \cap F_{p'q'}^{-} 
         \\\hline
         \ell^{\Phy} \subset F_\rho^{\Phy} & \rho \in \ell^+ &  \ell^{-} \subset F_{\rho}^-  
         \\\hline
         \ell^{\Phy} \subset F_\tau^{\Phy} &  \ell^{+} \subset F_{\tau}^+   & \tau \in \ell^-
    \end{array}$

    \vspace{7pt}
    \caption{Some geometric relations represented in physical, positive, and negative space.}
    \label{tab:geometry}
\end{table}
\renewcommand{\arraystretch}{1}

\subsection{Initial reductions}

Given a set $\cP\subset\R^3$ of $N$ points, our goal is to prove that it determines at least $N^{2/3-o(1)}$ distinct distances. To do this, it suffices to show that $\cP$ determines at most $N^{10/3+o(1)}$ distance quadruples. Now any distance quadruple $(p,p',q,q')\in\cP^4$ with $\norm{p-p'}=\norm{q-q'}$ and $(p,q)\neq (p',q')$ corresponds to the isotropic line $\ell^{\Phy}$ through $\phi(p,q),\phi(p',q')$ as well as the isotropic line $\ell^+$ which is $F^+_{pq}\cap F^+_{p'q'}$ and the isotropic line $\ell^-$ which is $F^-_{pq}\cap F^-_{p'q'}$. Working in positive space, our goal is to show that the collection $\cF=\{F^+_{pq}\}_{p,q\in \cP}$ of $N^2$ physical 3-flats determines few 2-rich lines. In the following three parts of the paper we will carry out this strategy.

To make the argument in \cref{part:iii} work, we work with a slight generalization of this problem. Given $H\subseteq\cP\times\cP$, let $\cF$ be either $\{F^+_{pq}\}_{(p,q)\in H}$ or $\{F^-_{pq}\}_{(p,q)\in H}$, a set of $|H|$ physical 3-flats. Our goal will be to prove a bound of the form
\[|\cL_{[2,C]}(\cF)|\lesssim N^{4/3+o(1)}|H|,\]
where $\cL_{[2,C]}(\cF)$ is the set of lines contained in at least 2 and at most $C$ elements of $\cF$ for some absolute constant $C$. 
The added flexibility of working both in positive and negative space as well as with arbitrary subsets $H\subseteq\cP\times\cP$ will be used heavily in \cref{part:iii}.
The reason we can restrict attention to lines which are at most $C$-rich is due to a trick based on Beck's theorem (\cref{thm:beck-trick}) that we introduce during the reduction to a bound on the number of distances quadruples.

We also place three assumptions on the set $\cP\subset\R^3$:
\begin{enumerate}[(i)]
    \item $|\cP\cap\gamma|\leq N^{2/3}\deg\gamma$ for every curve $\gamma\subset\CC^3$;
    \item $|\cP\cap F|\leq N^{2/3}$ for every plane $F\subset\R^3$;
    \item $|\cP\cap\S|\leq N^{2/3}$ for every sphere $\S\subset\R^3$.
\end{enumerate}
These assumptions are harmless in proving \cref{thm:main}: an easy application of B\'ezout's theorem shows that if $\cP$ fails (i), then the set $\cP\cap\gamma$ already determines $\gtrsim N^{2/3}$ distinct distances, while if $\cP$ fails (ii) or (iii), then the Guth--Katz theorem (and its generalization to the sphere) shows that $\cP\cap F$ or $\cP\cap\S$ determines $\gtrsim N^{2/3}/\log N$ distinct distances \cite{GK15,Tao11blogpost}. In \cref{part:i}, assumption (i) will be crucially used to bound the algebraic concentration of $\cF$. In \cref{part:ii}, both assumptions (i)(ii) will be used, while in \cref{part:iii}, we will need all three assumptions on $\cP$.

\subsection{Part I: line--3-flat incidences}
\label{ssec:outline-i}

In this part we show how to replace each 3-flat $F\in\cF$ with a collection of surfaces so that most of the lines which are 2-rich for $\cF$ are 2-rich for the collection of surfaces. The main goal of \cref{part:i} is to produce a set of surfaces $\cS$ with the property that
\[|\cL_{[2,C]}(\cF)|\lesssim |\cL_{[2,C]}(\cS)|+N^{4/3}|H|\log^2N\]
and so that
\[\deg\cS\lesssim N^{2/3}|H|\log N.\]
(Here and throughout the paper, we use $\deg\cS$ to refer to the sum of the degrees of the surfaces $S\in\cS$.) In fact we have slightly more control over these two quantities, see \cref{thm:codim-0-incidence-bound}.

\subsubsection{Incidence bounds}
In \cref{sec:codim-0-incidences}, we simplify and substantially generalize the techniques of Walsh to prove incidence bounds involving flats of arbitrary dimension. In \cite{Wal23}, Walsh proved tight incidence bounds for point-line incidences in $\R^n$ and $\F^n$ for every dimension $n$. These bounds are in terms of the algebraic concentrations of the set of lines. (These also generalize to $k$-flat--$(k+1)$-flat incidences and from lines to arbitrary curves.) We simplify his methods, which allows us to generalize them to incidence problems of arbitrary codimension. This simplification comes with the drawback that the results are no longer tight, but rather incur polylogarithmic loss.

As one example of our techniques, we can prove point--$r$-flat incidence bounds in $\R^n$ and $\F^n$ for $n\geq 2r$. These bounds are in terms of the algebraic concentrations of the set of flats. Generalizing \cref{eq:D_k,eq:D1_k}, one needs to be able to bound the total degree of the $(r-i)$-dimensional components of the intersections of the flats with a variety $V$ of arbitrary degree and dimension. These results also generalize to $s$-flat--$(r+s)$-flat incidences and from $r$-flats to arbitrary varieties. We do not pursue the full extent of these techniques in this paper, instead, we defer these to the companion paper \cite{TYZ26b}. One example of the type of results that we can prove is the following.

\begin{theorem}[{\cite{TYZ26b}}]
\label{thm:2-rich-companion-paper}
Let $\cF$ be a family $(r+s)$-flats in $\F^n$ whose pairwise intersection have dimension at most $s$. Then writing $\cL_2(\cF)$ for the set of $s$-flats contained in at least two elements of $\cF$, we have
\[|\cL_2(\cF)|\lesssim_n\sum_{k=s+2}^{n-r+1}|\cF|(\log|\cF|)^{r}\cD^{r-1}_k(\cF)^{\tfrac{r}{k-s-1}}.\]
\end{theorem}

In this result, $\cD^i_k(\cF)$ is the $k$-th codimension $i$ algebraic concentration, defined by
\begin{equation*}
\cD_k^i(\cF)=\max_{V\subset\F^n:\dim V=k}\frac{\sum_{F\in\cF}\sum_{U\subseteq F\cap V:\dim U=\dim F-i}\deg U}{\deg V}
\end{equation*}
where the second sum is over the irreducible $(\dim F-i)$-dimensional components $U$ of $F\cap V$. In many regimes of the parameters, \cref{thm:2-rich-companion-paper} is tight up to the polylogarithmic factor. 

In this paper, we focus on applying our new techniques to the line--3-flat incidence problem in $\F^6$. We note that in general, one can prove stronger incidence bounds over $\R^6$ than over arbitrary fields. For example, for $2\leq k\leq n^{1/2}$, a set $\cL$ of $n$ lines in $\F^2$ satisfies $|\cP_k(\cL)|\lesssim n^2/k^2$, while if the lines lie in $\R^2$, they satisfy $|\cP_k(\cL)|\lesssim n^2/k^{3}$. Note that these bounds are of the same strength for $k=O(1)$, but the bound over the reals is stronger for large values of $k$. For the Guth--Katz distinct distances bound in $\R^2$ it was necessary to bound the number of $k$-rich points for all $k\geq 2$, thus it was necessary to study the incidence problem in $\R^3$ instead of $\F^3$. However, due to the trick using Beck's theorem in the initial reduction, we only need to study 2-rich lines; thus it suffices to work in $\F^6$ instead of $\R^6$. We note that this Beck trick appears to apply to the distinct distance problem in $\R^d$ for all $d\geq 3$, but not for $d=2$.

Now our techniques allow us to bound $|\cL_2(\cF)|$ in terms of the codimension 0 algebraic concentrations $\cD_4(\cF)$ and $\cD_5(\cF)$, see \cref{eq:D_k}, and the codimension 1 algebraic concentrations $\cD_3^1(\cF)$ and $\cD_4^1(\cF)$, see \cref{eq:D1_k}. We are able to prove the bounds $\cD_4(\cF)\lesssim N^{2/3}$ and $\cD_5(\cF)\lesssim N^{4/3}$. It turns out that $\cD_4(\cF)\gtrsim|\cF|^{1/3}$ and $\cD_5(\cF)\gtrsim|\cF|^{2/3}$ for any set of 3-flats $\cF$. As $|\cF|=|H|\leq N^2$, these bounds are asymptotically tight, at least in the case $H=\cP\times\cP$.

However, we run into serious technical difficulties when trying to bound $\cD_4^1(\cF)$ and, as mentioned in \cref{sec:intro}, the desired bound on $\cD_3^1(\cF)$ is simply false due to the presence of rich rigid motions of the opposite orientation.\footnote{To be precise, the desired bound on $\cD_4^1(\cF)$ is also false: one can take a 1-dimensional family $\cT$ of rigid motions of which contains many rich members. The union of $F_\tau^+$ over $\tau\in \cT$ is a 4-dimensional variety which shows that $\cD_4^1(\cF)$ is large. Even isolating this type of 4-dimensional variety, the aforementioned technical issues still prevent us from bounding $\cD_4^1(\cF)$ over high degree 4-dimensional varieties not of this form.} Our strategy is to instead only bound the codimension 0 algebraic concentrations. Our techniques allow us to use this information to reduce the line--3-flat incidence problem to a line-surface incidence problem.

Let us compare the hypothesis of Walsh's result to other incidence bounds known in high dimensions. The Guth--Katz incidence bound gives a tight bound on the number of 2-rich points determined by a set of lines in $\R^3$, requiring only that the lines do not concentrate on planes and reguli \cite{GK15}. However, this proof does not seem to generalize to higher dimension. A different polynomial method technique known as \emph{bounded-degree polynomial partitioning} was used by Guth and Guth--Zahl to solve the same problem (with subpolynomial loss) first in $\R^3$ and then in $\R^4$, now requiring non-concentration on varieties of bounded degree \cite{Gut15lowdeg,GZ17}. These techniques also run into difficulties in generalizing to higher dimensions. In contrast, Walsh's techniques give tight bounds on the number of 2-rich points in any dimension, but require working inside varieties of high degree. As a result, we are forced to control how our flats concentrate inside varieties of arbitrary degree, which poses significant difficulties.

\subsubsection{Algebraic geometry tools}
Walsh's techniques require working inside varieties of arbitrary degree and require precise quantitative control over these varieties. This control is provided by new results in algebraic geometry that he proved earlier \cite{Wal20}. To summarize, the key result is that for an irreducible variety $V\subset\F^n$ of codimension $r$, one can express $V$ as an \emph{approximate complete intersection}. That is, there exist polynomials $f_1,\ldots,f_r\in\F[x_1,\ldots,x_n]$ so that $V$ is an irreducible component of $Z(f_1,\ldots,f_r)$ and so that $\deg V\sim \deg f_1\deg f_2\cdots\deg f_r$. This result will also play an important role in this paper; we give a self-contained exposition in \cref{sec:Walsh}.

To control the algebraic concentrations $\cD_4(\cF)$ and $\cD_5(\cF)$, we require a significant strengthening of this result. We wish to understand how many of the 3-flats $ F\in \cF$ can lie in a variety $V$; to do this, we need to understand the tangent spaces $T_pV$. The tangent space $T_pV$ can be computed as the kernel of the matrix whose columns are $\vec\nabla f(p)$ for all $f\in I(V)$.\footnote{We use standard algebraic geometry notation: $I(V)$, the \emph{ideal} of $V$, is the ring of polynomials $f\in\F[x_1,\ldots,x_n]$ which vanish identically on $V$.} This characterization is too unwieldy for our purposes. We would like to write $V$ as an approximate complete intersection $Z(f_1,\ldots,f_r)$ with the properties above and the additional property that the kernel of the \emph{Jacobian matrix} whose columns are $\vec\nabla f_1(p),\ldots,\vec\nabla f_r(p)$ agrees with $T_pV$ for almost all $p\in V$. We prove this result in \cref{thm:approx-complete-intersection}, resolving a conjecture of Walsh. It is an essential component of our strategy to bound the algebraic concentrations $\cD_4(\cF)$ and $\cD_5(\cF)$.

\subsubsection{Ruled surface theory}

Given a variety $V\subseteq E^+$, we study its \emph{Fano variety}
\[\Lambda=\{\phi(p,q)\in\F^6:F^+_{pq}\subseteq V\}\subseteq E^{\Phy},\]
in other words, the collection of physical 3-flats which it contains. We say that a variety $V$ is \emph{ruled} (by physical 3-flats) if for a generic point $\rho\in V$, there exists some $\phi(p,q)\in\Lambda$ with $\rho\in F_{pq}^+\subseteq V$.

Our goal is to prove the bounds $\cD_4(\cF)\lesssim N^{2/3}$ and $\cD_5(\cF)\lesssim N^{4/3}$. For this discussion let $\cF=\{F^+_{pq}\}_{p,q\in\cP}$. (The distinction between $H$ and $\cP\times\cP$ will not matter for this part of the argument. We also want to apply this argument both in $E^+$ and $E^-$, but the two are symmetric.)
Thus our goal is to bound
\[\abs{{\Lambda\cap\phi(\cP\times\cP)}}\lesssim N^{2/3}\deg V\]
for every irreducible 4-dimensional variety $V$, and 
\[\abs{{\Lambda\cap\phi(\cP\times\cP)}}\lesssim N^{4/3}\deg V\]
for every irreducible 5-dimensional variety $V$, where $\Lambda$ is as defined above.

We prove these two bounds in \cref{sec:concentration}. We break into cases depending on if $V$ is ruled. In general, if $V$ is ruled, then $\dim \Lambda$ may be large, but it is also quite structured. We classify the possibilities for $\dim\Lambda$ and, in each case, bound $\deg\Lambda$ in terms of $\deg V$. With this classification, it remains to bound $\abs{\Lambda\cap\phi(\cP\times\cP)}$. However, using the fact that $\cP\times\cP$ is a grid and assumption (i) gives that no curve contains too many points of $\cP$, we are able to prove the desired bounds on the algebraic concentration in the ruled case.

For the unruled case we need to develop some theory about varieties ruled by families of flats. This generalizes classical work on ruled surfaces in $\CC^3$. In particular, the Cayley--Monge--Salmon theorem implies that for an irreducible surface $V\subset\CC^3$, there is a polynomial $\Flec V$ called the \emph{flecnodal polynomial} such that
\begin{itemize}
    \item $\deg\Flec V\leq 11\deg V$;
    \item every line in $V$ lies in $V\cap Z(\Flec V)$; and
    \item if $\Flec V$ vanishes identically on $V$, then $V$ is ruled by lines.
\end{itemize}
As a consequence, if $V$ is an unruled variety then every line in $V$ also lies in the proper subvariety $V\cap Z(\Flec V)$. This has dimension 1, and by B\'ezout's theorem has degree at most $11(\deg V)^2$. Thus the Cayley--Monge--Salmon theorem implies that every unruled surface $V$ contains at most $11(\deg V)^2$ lines.

In \cref{sec:csm-for-flats}, we prove a wide-reaching generalization of the Cayley--Monge--Salmon theorem.\footnote{See \cite{GZ18} for another generalization of the Cayley--Monge--Salmon theorem, this one in the setting of varieties ruled by a family of bounded-degree curves.} As one consequence, we show that if $V\subseteq E^+$ is not ruled by physical 3-flats, there exists a polynomial $g$ which does not vanish identically on $V$, yet vanishes identically on all physical 3-flats lying on $V$. Controlling $\deg g$ will be crucial, yet the bound $\deg g\lesssim\deg V$ will not be strong enough for our application. Instead we need to bound $\deg g$ by a smaller parameter known as the \emph{partial degree} of $V$. This will be defined precisely later (\cref{def:partial-degree}) but informally is equal to $\max\{\deg f_1,\ldots,\deg f_r\}$ where $Z(f_1,\ldots,f_r)$ is a representation of $V$ as an approximate complete intersection. (Recall that $\deg f_1\deg f_2\cdots\deg f_r\sim \deg V$, so for varieties of codimension at least 2, the partial degree can be much smaller than the degree.) To prove our generalization of the Cayley--Monge--Salmon theorem with this precise control on $\deg g$, we crucially use the results of \cref{sec:approx-complete-intersection}, namely the representation of an arbitrary variety $V$ as an approximate complete intersection whose Jacobian matrix defines the tangent space almost everywhere.

\subsection{Part II: line-surface incidences}
\label{ssec:outline-ii}

Following \cref{part:i}, the remaining problem is to bound $|\cL_{[2,C]}(\cS)|\lesssim N^{4/3+o(1)}|H|$ where $\cS$ is a set of surfaces satisfying $\deg\cS\lesssim N^{2/3}|H|\log N$. These surfaces are produced by intersecting the 3-flats $\cF$ with polynomials of degree $\lesssim N^{2/3}$; in particular, each surface has degree $\lesssim N^{2/3}$ and is contained in a physical 3-flat $F^+_{pq}$ for some $(p,q)\in H$. (We also perform the symmetric argument where each surface is contained in a physical 3-flat $\cF^-_{pq}$ in negative space. The argument is identical, so for this overview we focus on the argument in positive space.)

By classical ruled surface theory, each surface $S\in\cS$ is either unruled, singly ruled, doubly ruled, or infinitely ruled. This means that a generic point of $S$ is contained in zero, one, two, or infinitely many lines that lie in $S$, respectively. An infinitely ruled surface is a plane; we bound their contributions in \cref{part:iii}. The contribution to $|\cL_{[2,C]}(\cS)|$ from unruled surfaces is easy to bound: by the Cayley--Monge--Salmon theorem, such a surface $S$ contains at most $11(\deg S)^2\lesssim N^{2/3}\deg S$ lines, so the total contribution to $|\cL_{[2,C]}(\cS)|$ is bounded by $\lesssim N^{2/3}\deg\cS\lesssim N^{4/3}|H|\log N$.

Defining $\cS'\subseteq\cS$ to be the set of singly or doubly ruled surfaces, the main result of \cref{part:ii} (\cref{thm:ruled-ruled}) is the bound
\[|\cL_{[2,C]}(\cS')|\lesssim N^{2/3}\deg\cS'\log N.\]
As $\deg\cS'\lesssim N^{2/3}|H|\log N$, this gives an overall bound of $N^{4/3}|H|\log^2N$. We prove this result in \cref{sec:ruled-ruled} after collecting some standard result from classical ruled surface theory in \cref{sec:ruled-surface-classical} and some tools from algebraic geometry in \cref{sec:multiplicity}.

We use a novel strategy for studying this incidence problem between lines and ruled surfaces. A standard way to study an incidence problem between lines and ruled surfaces in $\F^6$ would be to intersect with a generic hyperplane, producing a point-curve incidence problem in $\F^5$. However, this loses a significant amount of information. As each of our surfaces is ruled, they contain a 1-dimensional family of lines. Moving to the parameter space of isotropic lines in $E^+$, each line maps to a point in this space, and each ruled surface maps to a curve.
One could view this as a point-curve incidence problem in the parameter space of isotropic lines in $E^+$ (a 9-dimensional space). Instead of studying this 9-dimensional problem directly, we use geometric triality to transfer the line-surface incidence problem to a more structured incidence problem in $E^{\Phy}$.

Each ruled surface $S\subset F^+_{pq}\subset E^+$ produces a 1-dimensional family of isotropic lines $\ell^+$ that it contains. Using geometric triality, such an $S$ also produces the corresponding 1-dimensional family of isotropic lines $\ell^{\Phy}$ in $E^{\Phy}$. As all the $\ell^+$ lie in $F^+_{pq}$, all the $\ell^{\Phy}$ will pass through the point $\phi(p,q)$. We transform each ruled surface $S$ to the union of these lines, which is a cone $C_S$ in physical space with apex $\phi(p,q)$.

In this way, we turn the incidence problem between lines and ruled surfaces in positive space to an incidence problem between points and cones in physical space. Indeed, we are interested in lines in physical space that lie in two cones $C_S,C_{S'}$ with apexes $\phi(p,q)$ and $\phi(p',q')$. For our argument, we simplify the setting by forgetting the line $\ell^{\Phy}$ connecting $\phi(p,q),\phi(p',q')$ and only recording the incidence between the point $\phi(p',q')$ and the cone $C_S$. Note that since in we are studying $|\cL_{[2,C]}(\cS)|$ for some absolute constant $C$, all the lines in question will contain at most $C$ points. Thus this simplification is harmless as it loses at most a factor of $C$.

So far we have reduced the problem to a point-cone incidence problem in $E^{\Phy}$ between $\phi(H)\subseteq\phi(\cP\times\cP)$ and the collection of cones $\{C_S\}_{S\in\cS'}$. Next we apply $\phi^{-1}$ and project the problem onto either the first three or last three coordinates, producing a point-cone incidence problem in $\CC^3$ where the collection of points is $\cP$. 

We study this point-cone incidence problem in $\CC^3$ using a combination of standard and novel tools. First we use the polynomial method (parameter counting) to find a collection of curves $\Gamma$ with $\deg\Gamma\lesssim N^{2/3}$ such that the curves contain all the points $\cP$. The main term in the point-cone incidence problem is a curve-cone incidence problem, i.e., the number of pairs of curve and cone where the curve is contained within the cone. We bound this problem using the so-called ``trisecant lemma'' in algebraic geometry. Morally speaking, this result says that it is unlikely for three distinct curves to be contained in the same cone. More precisely, suppose that the curves and cones all have degree $d$. Then the incidence graph between the curves and cones will be $K_{3,d^2+1}$-free.\footnote{To be more precise, the result as stated only holds with the further assumption that no curve passes through the apex of another cone.} This allows us to control the number of curve-cone incidences, producing the desired bound on the original incidence problem between lines and ruled surfaces.

\subsection{Part III: very rich 3-flats}
\label{ssec:outline-iii}

Combining the results of \cref{part:i,part:ii}, it remains to bound the contribution to $\cL_{[2,C]}(\cS)$ where at least one of the surfaces involved is a plane. For any plane $S\in\cS$, we know that $S\subset F^+_{pq}$ for some $(p,q)\in H$. It turns out that any such $S$ satisfies $S=F_{pq}^+\cap F_\tau^+$ for some $\tau\in E^-$.\footnote{Intuitively, one should expect this to be true since for every $\tau\in E^-$ such that $\tau(p)=q$, the intersection $F_{pq}^+ \cap F_\tau^+$ is a plane. There is a 3-dimensional family of such $\tau$ and a 3-dimensional family of planes $S \subset F_{pq}^+$.} Now if $\tau$ is exactly $k$-rich, meaning that $\tau(p)=q$ for exactly $k$ choices of $(p,q)\in H$, then $S=F^+_{pq}\cap F^+_\tau$ contributes at most $k$ to $\cL_{[2,C]}(\cS)$. In particular, it contributes the lines $F_{pq}^+ \cap F_{p'q'}^{+}$ where $\tau(p')=q'$. Thus the contribution from $\tau$ which are at most $N^{2/3+o(1)}$-rich is bounded by $N^{2/3+o(1)}|\cS|\leq N^{4/3+o(1)}|H|$.

Moving from positive space back to physical space, it remains to bound the number of isotropic lines $\ell^{\Phy}$ which contain at least 2 and at most $C$ points of $\phi(H)$ and are contained in $F_\tau^{\Phy}$ for some $\tau\in E^-$ that is (at least) $N^{2/3+o(1)}$-rich. Since the arguments in \cref{part:i,part:ii} treat positive space and negative space symmetrically, we can further restrict our attention to those isotropic lines which lie in both $F_\rho^{\Phy}$ and $F_\tau^{\Phy}$ for some $N^{2/3+o(1)}$-rich $\rho\in E^+$ and some $N^{2/3+o(1)}$-rich $\tau\in E^-$. Write $\cL$ for this set of lines. 

For simplicity, suppose that we have sets $\cR\subset E^+$ and $\cT\subset E^-$ where each $\rho\in\cR$ is exactly $k$-rich for $\phi(H)$ and each $\tau\in\cT$ is exactly $k'$-rich for $\phi(H)$ for some $k,k'\gg N^{2/3}$. Furthermore, suppose that each line $\ell^{\Phy}\in\cL$ is contained in $F_\rho^{\Phy}$ for a unique $\rho\in\cR$ and in $F_\tau^{\Phy}$ for a unique $\tau\in\cT$. Finally, suppose $F_\rho^{\Phy}$ contains $\gtrsim k^2$ lines $\ell^{\Phy}\in\cL$ for each $\rho\in\cR$ and $F_\tau^{\Phy}$ contains $\gtrsim k'^2$ lines for each $\tau\in\cT$. We will show that it is always possible to pass to a subset for which these regularity properties approximately hold, which will be crucial to carry out the argument in \cref{part:iii}. Unfortunately this regularization is the most quantitatively costly part of the argument: it is the reason why we lose a factor of $N^{O(\sqrt{\log\log N/\log N})}$ instead of a polylogarithmic factor. In addition, this is the reason why we work with a subset $H\subseteq\cP\times\cP$, as we cannot guarantee that the regularization holds with $H=\cP\times\cP$. 

Our goal is to show that such a configuration can only determine many 2-rich lines if it has a large amount of underlying geometric structure. We will find increasingly more structure in the problem, eventually reaching a contradiction.

First consider a point $(p,q)\in H$. We transfer the portion of the setup which interacts with $(p,q)$ to a point-plane incidence problem in the 3-flat $F_{pq}^+$. Consider the rigid motions $\cR_{pq}\subseteq\cR$ and $\cT_{pq}\subseteq\cT$ which map $p$ to $q$ as well as the lines $\cL_{pq}\subseteq\cL$ which pass through $\phi(p,q)$. We can view this setup as a collection of points, planes, and lines in the 3-flat $F_{pq}^+$. Indeed, each $\rho\in\cR_{pq}$ is a point in $F_{pq}^+$, for each line $\ell^{\Phy}\in\cL_{pq}$, the corresponding line $\ell^+$ lies in $F_{pq}^+$, and for each $\tau\in\cT_{pq}$, we have $F_{pq}^+\cap F_\tau^+$ is a plane in $F_{pq}^+$. As each line $\ell^{\Phy}\in\cL_{pq}$ is captured by a unique pair $(\rho,\tau)\in\cR_{pq}\times\cT_{pq}$, each line in the collection lies in a unique point-plane incidence in this setup.

Here we use a result of de Zeeuw \cite{deZ16} to decompose this point-plane incidence problem in $F_{pq}^+\cong\R^3$ into a structured part and a small exceptional part. One way to have many point-plane incidences is for all of the points to lie on a special line $\ell^+_\ast$ and have all of the planes contain $\ell^+_\ast$. The result of de Zeeuw shows that we must have a small number of configurations of this form and a small number of exceptional incidences. 

To study the structured contribution, we move from positive space back to physical space. If a line $\ell^{\Phy}$ through $\phi(p,q),\phi(p',q')$ was captured by the structured part, this means that there is some line $\ell^+_\ast\subset F_{pq}^+$ so that $\rho\in \ell^+_\ast\subset F_\tau^+\cap F_{pq}^+$. However, we can do the same argument in $F_{p'q'}^+$. If $\ell^{\Phy}$ is captured by the structured part here, we have some line $(\ell')^+_\ast\subset F_{p'q'}^+$ so that $\rho\in (\ell')^+_\ast\subset F_\tau^+\cap F_{p'q'}^+$. Passing back to physical space, we have two lines $\ell^{\Phy}_\ast$ and $(\ell')^{\Phy}_\ast$ which both lie in the 2-flat $F_\rho^{\Phy}\cap F_\tau^{\Phy}$. In particular, these two lines intersect (possibly at infinity). Now the regularization implies that each $F_\rho^{\Phy}$ is $\sim k$-rich and also contains $\gtrsim k^2$ distance quadruples. Thus if the structured case dominates, we get $\sim k$ lines $\ell^{\Phy}_\ast$ in $F_\rho^{\Phy}$, yet $\gtrsim k^2$ pairs of them pairwise intersect.

To study this situation, we prove a structural result about a collection of lines in $\P_{\RR}^3$ where a constant fraction pairwise intersect. We show that this can only occur if many of the lines coincide, many of the lines are concurrent, many of the lines are coplanar, or many of the lines lie in a common regulus. In each of these four cases, the additional structure is enough to prove the desired bound on the number of distance quadruples. We study these structured cases in \cref{sec:degenerate}. This handles the case where the structured contribution dominates.

Now if the exceptional incidences dominate, this means that we have a small number of point-plane incidences which are capturing a large number of distance quadruples. Again passing from positive space to physical space, fix a rich rigid motion $\rho\in\cR$ and consider the part of the problem living in $F_\rho^{\Phy}$. There are $\sim k$ points of $\phi(H)$ in $F_\rho^{\Phy}$ and $\gtrsim k^2$ of these correspond to a distance quadruple, meaning that $\phi(p,q),\phi(p',q')$ both lie in some plane $F_\rho^{\Phy}\cap F_\tau^{\Phy}$. Using more incidence geometry in the 3-flat $F_\rho^{\Phy}$, we can show that this can only happen if most of the configuration is captured by rich lines. This puts us in the first structured setting that we studied above, completing the proof.

The new tools used to study incidence geometry tools in $\R^3$ are developed in \cref{sec:incidence-geometry-r3}. Then the above argument is carried out in full in \cref{sec:main-argument}. We note that this argument relies heavily on the regularity assumptions. We show how to pass to $H\subseteq\cP\times\cP$ which regularizes the setup in \cref{ssec:regularization}. As mentioned previously, this is the only part of the argument which loses more than a polylogarithmic factor.

\subsection{Remarks on higher dimensions}
\label{ssec:higher-dimen}

To study the distinct distances problem in higher dimensions, consider a set of $N$ points $\cP\subset\R^d$. Defining $F_{pq}\subset E^+(d)$ to be the set of orientation-preserving rigid motions which send $p$ to $q$, we have that $E^+(d)$ is $\binom{d+1}2$-dimensional and $F_{pq}$ is $\binom d2$-dimensional. If $(p,q,p',q')$ is a distance quadruple, then $F_{pq}\cap F_{p'q'}$ is $\binom{d-1}2$-dimensional, otherwise it is empty. Bardwell-Evans and Sheffer showed how to identify an affine open subset of $E^+(d)$ with $\R^{\binom{d+1}2}$ so that each $F_{pq}$ is identified with a $\binom d2$-flat \cite{BES19}. For $d>3$ there is no triality between the three spaces (since physical space has dimension $2d$ while positive and negative space have dimension $\binom{d+1}2$) but one could still attempt a similar proof strategy to study this incidence problem between $\binom{d-1}2$-flats and $\binom d2$-flats in positive space.

The same trick using Beck's theorem applies in all dimension $d\geq 3$, so it suffices to bound the number of 2-rich $\binom{d-1}2$-flats determined by the set $\cF=\{F_{pq}\}$ of $N^2$ $\binom d2$-flats.

The main complication in dimension $d=3$ is the presence of the 3-flats $F_\tau^+$ which cause $\cD_3^1(\cF)$ to be too large. However, it seems plausible that this issue is an artifact of working in dimension $d=3$. Indeed, for $d\geq 4$ it may be possible to perform a computation similar to the argument in \cref{sec:concentration} to show that for every $i$ with $0\leq i\leq d-2$ and $\binom d2\leq i+k\leq \binom{d+1}2$ and every $k$-dimensional variety $V\subseteq \R^{\binom{d+1}2}$, the number of $\binom d2$-flats in $\cF$ which have a $(\binom d2-i)$-dimensional intersection with $V$ is bounded by $O_{\deg V}(N^{2(k+i-\binom d2)/d})$. (In fact, proving this for $i=d-2$ implies it for all smaller $i$.) If true, this bound would be best-possible. Unfortunately, it is out of reach of our current techniques to prove this bound with a linear dependence on $\deg V$ for any $i>0$. This is the technical difficulty that we run into in our argument which means that we cannot bound $\cD_4^1(\cF)$ even ignoring the issue with the 3-flats $F_\tau^+$ that cause $\cD_3^1(\cF)$ to be too large. If it were possible to overcome this difficulty, then one could give the following optimal bounds on the algebraic concentrations.

\begin{conjecture}
\label{conj:concentration-conj-higher-dim}
For $d\geq 4$, let $\cP\subset\R^d$ be a set of points such that $|\cP\cap V|\leq N^{2k/d}\deg V$ for every irreducible $k$-dimensional variety $V\subset\CC^d$. Then defining $\cF=\{F_{pq}\}_{p,q\in\cP}$, a collection of $N^2$ $\binom d2$-flats in $\R^{\binom{d+1}2}$ whose pairwise intersections have dimension at most $\binom{d-1}2$, the bound
\[\cD_k^{d-2}(\cF)\lesssim N^{2(k-\binom {d-1}2-1)/d}\]
holds for all $\binom {d-1}2+2\leq k\leq \binom{d}2+2$.
\end{conjecture}

While \cref{conj:concentration-conj-higher-dim} is false for $d=3$, it is plausibly true for each $d\geq 4$. Combining \cref{thm:2-rich-companion-paper,conj:concentration-conj-higher-dim} would imply the distinct distance conjecture in $\R^d$ (up to $\log^{d-1}N$ error) for each $d\geq 4$.
\section{The space of rigid motions}
\label{sec:setup}

A rigid motion of $\R^3$ is a map $\rho\colon\R^3\to\R^3$ that preserves distances. The space of rigid motions is a 6-dimensional space with two connected components: the orientation-preserving rigid motions and the orientation-reversing rigid motions. Given a pair of points $p,q\in\R^3$, the set of rigid motions which map $p$ to $q$ is a 3-dimensional space. We would like to choose coordinates so that one copy of $\R^6$ parametrizes most of the orientation-preserving rigid motions and another copy of $\R^6$ parametrizes most of the orientation-reversing rigid motions. We would like to do this in a way such that the set of orientation-preserving rigid motions which map $p$ to $q$ correspond to an affine 3-dimensional subspace of $\R^6$ (a 3-flat) under this parametrization. We will refer to this 3-flat as $F^+_{pq}$. Similarly, we also want the orientation-reversing rigid motions which map $p$ to $q$ to correspond to a 3-flat $F^-_{pq}$.

Now the pair $(p,q)$ lives in $\R^3\times\R^3$, so there are three copies of $\R^6$ that we are dealing with: physical space, corresponding to pairs of points $(p,q)$; positive space, corresponding to orientation-preserving rigid motions; and negative space, corresponding to orientation-reversing rigid motions. There is a deep symmetry between these three spaces which will become more apparent if we apply the change of coordinates $\phi\colon\R^3\times\R^3\to\R^6$ defined by 
\begin{equation}
\label{eq:phi-defn}
\phi(p,q)=\paren{\frac{p+q}2,\frac{p-q}2}
\end{equation}
to physical space.

One of the reasons we care about this symmetry is that we need to consider more 3-flats than just $F^+_{pq}$ and $F^-
_{pq}$. For an orientation-preserving rigid motion $\rho\colon\R^3\to\R^3$, define its graph $\Gamma_\rho\subset\R^3\times\R^3$ to be the usual $\Gamma_{\rho}=\{(p,\rho(p)):p\in\R^3\}$. This is a 3-flat in $\R^3\times\R^3$. Under the change of coordinates $\phi$, set $F^{\Phy}_\rho=\phi(\Gamma_\rho)$, a 3-flat in physical space. Similarly, define $F^{\Phy}_\tau=\phi(\Gamma_\tau)$ for any orientation-reversing rigid motion $\tau$. To complete the picture, we define $F^-_{\rho}$ for each orientation-preserving rigid motion $\rho$ and $F^+_\tau$ for each orientation-reversing rigid motion $\tau$. The former, $F^-_{\rho}$, is a 3-flat in negative space, consisting of the orientation-reversing rigid motions $\tau$ such that $\rho,\tau$ agree on a 2-flat. (This is 3-dimensional since there is a 3-dimensional space of 2-flats in $\R^3$ on which the motions can agree. Once $\rho$ and the 2-flat are fixed, $\tau$ is determined.) The latter, $F^+_\tau$, is a 3-flat in positive space, defined analogously. 

In the following sections, we define the parameterization of orientation-preserving and orientation-reversing rigid motions by $\R^6$ such that all these sets are indeed 3-flats. We need to determine which rigid motions were missed by our parametrization. We also need to understand the intersection between these 3-flats -- in what circumstances does the intersection take on which dimensions? In addition, for some of our arguments we want to work in $\R^6$; however, some work better over an algebraically closed field and so we complexify to $\CC^6$. Other times we want to work over a projective space; it will turn out that compactifying $\R^6$ into $\P_\R^6$ is not the right thing to do -- the symmetry between the three spaces becomes the most apparent when we compactify $\R^6$ into a certain quadric hypersurface in $\P_\R^7$. 

We note that prior work of Bardwell-Evans--Sheffer also shows how to choose coordinates on positive space so that the $F_{pq}^+$ are 3-flats \cite{BES19}. Their construction works for all dimensions $d\geq 2$; specializing to $d=3$, their construction agrees with ours, up to signs and a permutation of the coordinates.
We choose to perform the reduction in our own way since we will need to use the symmetry between the three spaces. We also point out that this symmetry is a special property of the dimension $d=3$. In general, physical space has dimension $2d$, while positive and negative space have dimension $\binom{d+1}2$; these only agree for $d=3$. 

\subsection{Graphs of rigid motions and compactification}
\label{sec:graph-and-compactification}
Given an affine linear map $\rho\colon\R^3\to\R^3$, its graph $\Gamma_{\rho}=\{(p,\rho(p)):p\in\RR^3\}\subset \RR^3\times\RR^3$ is a 3-flat. Now $\rho$ is a rigid motion if and only if it preserves distances, i.e., if $\norm{p-p'}=\norm{q-q'}$ for all $(p,q),(p',q')\in \Gamma_{\rho}$. This can be rewritten as 
\[0=\sum_{i=1}^3\paren{p_i-p_i'}^2-\paren{q_i-q_i'}^2=\sum_{i=1}^{3}\left((p_i+q_i)-(p_i'+q_i')\right)\left((p_i-q_i)-(p_i'-q_i')\right).\]
Under the change of coordinates $\phi\colon\R^3\times\R^3\to\R^6$ defined in \cref{eq:phi-defn}, write $\phi(p,q)=x=(\vec{x},\vec{x}_*)$. Then the above condition can be rewritten as $(\vec{x}-\vec{x}')\cdot (\vec{x}_*-\vec{x}_*')=0$ for every $x,x'\in \phi(\Gamma_{\rho})$. (Here $\cdot$ denotes the usual dot product on $\R^3$.)

\begin{definition}
\label{defn:isotropic}
A line $\ell\subset\R^6$ is \emph{isotropic} if for every $x=(\vec x,\vec x_*)$ and $x'=(\vec x',\vec x'_*)$ on $\ell$, we have
\[(\vec{x}-\vec{x}')\cdot (\vec{x}_*-\vec{x}_*')=0.\]
Note that if this condition holds a single pair of distinct points on $\ell$, then it holds for every pair of points on $\ell$.
\end{definition}

The argument above shows that $\rho$ is a rigid motion if and only if every line contained in the 3-flat $\phi(\Gamma_{\rho})$ is isotropic. In addition, $(p,q,p',q')$ is a distance quadruple (meaning $\norm{p-p'}=\norm{q-q'}$) if and only if the line $\aff\{\phi(p,q),\phi(p',q')\}$ is isotropic.

This condition can be nicely reformulated in $\PP^7$. We use coordinates $[x_0:x_1:x_2:x_3:x_{0^*}:x_{1^*}:x_{2^*}:x_{3^*}]$ for $\PP^7$, and we write $\vec x=(x_1,x_2,x_3)$ and $\vec x_*=(x_{1^*},x_{2^*},x_{3^*})$.
Define the embedding $\iota\colon \RR^6\hookrightarrow\PP^7$ by
\begin{equation}
\label{eq:iota-phys-defn}
\iota(\vec{x},\vec{x}_*)=[1:\vec{x}:(\vec{x}\cdot \vec{x}_*):\vec{x}_*]
\end{equation}
Then $\iota$ embeds $\R^6$ into the zero set of the homogeneous quadratic polynomial $N$, defined by
\begin{equation}
\label{eq:N-defn}
N(x) = x_0x_{0^*}-\vec{x}\cdot \vec{x}_*.
\end{equation}
Each line in $\R^6$ gets embedded into a curve in $Z(N)$; it turns out that the image is a line if and only if the initial line is isotropic. Indeed, for every $(\vec{x},\vec{x}_*),(\vec{x}',\vec{x}'_*)$ spanning an isotropic line in $\RR^6$, we have
\begin{align*}
(\lambda \vec{x}+(1-\lambda)\vec{x}')\cdot (\lambda \vec{x}_*+(1-\lambda)\vec{x}_*')
&=\lambda \vec{x}\cdot \vec{x}_*+(1-\lambda)\vec{x}'\cdot \vec{x}_*'-\lambda(1-\lambda)(\vec{x}-\vec{x}')\cdot (\vec{x}_*-\vec{x}_*')\\
&=\lambda \vec{x}\cdot \vec{x}_*+(1-\lambda)\vec{x}'\cdot \vec{x}_*',
\end{align*}
showing that $\iota(\lambda x+(1-\lambda)x')$ lies on the line connecting $\iota(x)$ and $\iota(x')$. Conversely, if $\iota(\lambda x+(1-\lambda)x')=\lambda\iota(x)+(1-\lambda)\iota(x')$, then the above calculation shows that $\lambda(1-\lambda)(\vec{x}-\vec{x}')\cdot (\vec{x}_*-\vec{x}_*')=0$, implying that $(\vec{x},\vec{x}_*),(\vec{x}',\vec{x}'_*)$ span an isotropic line in $\RR^6$.

Now if $\rho$ is a rigid motion, then every line in $\phi(\Gamma_\rho)$ is isotropic, which implies that $\iota(\phi(\Gamma_\rho))$ is a (dense open) subset of a $3$-plane completely lying in $Z(N)$. 

This motivates us to extend the space of pairs of points from $\RR^6$ to $Z(N)\subset\PP^7$. We use the term \emph{physical space}, denoted $E^{\Phy}$, for a copy of $Z(N)\subset\PP^7$. By the discussion above, it makes sense to think of rigid motions as $3$-planes in $E^{\Phy}$. This is not literally true since there are a few additional 3-planes ``at infinity'' which do not correspond to a geometric rigid motion, but this will not affect the argument. We will see next that the 3-planes in $Z(N)$ are exactly parametrized by two copies of $Z(N)$, corresponding to the orientation-preserving and orientation-reversing rigid motions.

\subsection{Split octonions, totally isotropic subspaces, and positive and negative space}
We wish to understand the 3-planes in $E^{\Phy}=Z(N)\subset\P^7$. Equivalently, we want to characterize the 4-dimensional linear subspaces of $\R^8$ on which $N(x)$ vanishes identically. The quadratic function $N(x)$ is well-studied; it can be viewed as the norm on the split octonions. Specifically, the coordinates used in \cref{eq:N-defn} when we represent the split octonions using Zorn's vector-matrix algebra \cite{Zorn33}.

\begin{definition}
\label{defn:zorn-split-octonion}
Let $\F=\R$ or $\CC$. Define the \emph{split octonions}, denoted $\OO_s(\F)$\footnote{We drop the $\F$ when the field is clear from context.}, to be the following 8-dimensional composition algebra over $\F$. For a point $x=(x_0,x_1,x_2,x_3,x_{0^*},x_{1^*},x_{2^*},x_{3^*})\in\OO_s$, write $\vec x=(x_1,x_2,x_3)$ and $\vec x_*=(x_{1^*},x_{2^*},x_{3^*})$. Then arrange the coordinates into a ``matrix''
\[x=\begin{bmatrix}x_0&\vec{x}\\\vec{x}_*&x_{0^*}\end{bmatrix},\]
where some of the entries are scalars and some are vectors. Then $N(x)$, defined by \cref{eq:N-defn}, is a nondegenerate quadratic form on $\OO_s$. Multiplication on $\OO_s$ is defined by
\[x\star y = \begin{bmatrix}x_0&\vec{x}\\\vec{x}_*&x_{0^*}\end{bmatrix}\star \begin{bmatrix}y_0&\vec{y}\\\vec{y}_*&y_{0^*}\end{bmatrix}=\begin{bmatrix}
    x_0y_0+\vec{x}\cdot \vec{y}_*&x_0\vec{y}+y_{0^*}\vec{x}+\vec{x}_*\times \vec{y}_*\\
    y_0\vec{x}_*+x_{0^*}\vec{y}_*-\vec{x}\times\vec{y}& x_{0^*}y_{0^*}+\vec{x}_*\cdot \vec{y}
\end{bmatrix}.\]
(Here $\times$ refers to the usual vector cross product in $\F^3$.) Define the bilinear form on $\OO_s$ by
\[\langle x,y\rangle = N(x+y)-N(x)-N(y) = x_0y_{0^*}-\vec{x}\cdot \vec{y}_*+x_{0^*}y_0-\vec{x}_*\cdot\vec{y}.\]
Write $1\in\OO_s$ for the element $(1,0,0,0,1,0,0,0)\in\OO_s$, corresponding to the ``identity matrix''
\[1=\begin{bmatrix}
    1 & \vec 0\\ \vec 0 & 1.
\end{bmatrix}\]
Finally, define the conjugation operation on $\OO_s$ by 
\[\bar{x} = \begin{bmatrix}
    x_{0^*}&-\vec{x}\\ -\vec{x}_*&x_0
\end{bmatrix}.\]
\end{definition}

This defines an algebra as the multiplication operation distributes over vector addition and is compatible with scalar multiplication, though the multiplication operation is nonassociative. The element 1 is the multiplicative identity and it is known that $N(x\star y)=N(x)N(y)$. (See \cite{Zorn33}, though all of these properties can also be checked easily by hand.) This means that $\OO_s$ is a composition algebra. Several simple properties of $\OO_s$ include the identities
\begin{align*}
x\star (\bar{x}\star y)=(y\star \bar{x})\star x&=N(x)y,\\
\overline{x\star y}&=\bar{y}\star \bar{x},\\
N(x)&=N(\bar x),
\end{align*}
which can all be checked easily (see also \cite[Eq. (1.2)]{vdBSpr60}).

We are interested in the set $E^{\Phy}=\P(Z(N))$. A point $x\in\OO_s$ is \emph{isotropic} if $N(x)=0$ and a linear subspace of $\OO_s$ is \emph{totally isotropic} if all of its elements are isotropic. We are interested in 3-planes in $E^{\Phy}$, i.e., 4-dimensional totally isotropic subspaces of $\OO_s$.

This problem was studied in 1960 by van der Blij and Springer \cite{vdBSpr60}. We will use the results from \cite{vdBSpr60} freely, though each statement we use can be checked tediously by hand as well. The composition algebra $\OO_s$ is split, meaning that there exists a nonzero isotropic vector. It is well-known that the maximal totally isotropic subspaces of a split 8-dimensional composition algebra have dimension 4 (see \cite[Section 2]{vdBSpr60}). Note that since $N$ respects the multiplication operation, for every $x\in Z(N)\subset \OO_s$, we have that $N$ vanishes identically on both $x\star \OO_s$ and $\OO_s\star x$, both of which are subspaces of $\OO_s$. It turns out that these are 4-dimensional and are exactly the maximal totally isotropic subspaces of $\OO_s$.

\begin{theorem}[{\cite[Theorems 3 and 4]{vdBSpr60}}]
\label{thm:4-dim-isotropic-characterization}
Every totally isotropic subspace of $\OO_s$ is contained in a 4-dimensional totally isotropic subspace which has the form $x\star \OO_s$ or $\OO_s\star x$ for some nonzero isotropic $x\in \OO_s$. Moreover, up to scaling $x$, this representation is unique.
\end{theorem}

We wish to parameterize the 3-planes in $E^{\Phy}$, equivalently the 4-dimensional totally isotropic subspaces of $\OO_s$, by two copies of $\P(Z(N))$, which we call positive space and negative space, $E^+$ and $E^-$, respectively. For our convenience, we use a slightly different parameterization than the one given by \cref{thm:4-dim-isotropic-characterization}.

\begin{definition}
$E^{\Phy},E^+,E^-$, refer to three copies of $\P(Z(N))$, called \emph{physical space}, \emph{positive space}, and \emph{negative space}, respectively. For each $\rho\in E^+$, define the \emph{positive 3-plane in physical space}
\[F_\rho^{\Phy}=\{x\in E^{\Phy}:\rho\star x=0\}\subset E^{\Phy},\]
and for each $\tau\in E^-$, define the \emph{negative 3-plane in physical space}
\[F_\tau^{\Phy}=\{x\in E^{\Phy}:x\star \tau=0\}\subset E^{\Phy}.\]    
\end{definition}

We next show that this definition agrees with \cref{thm:4-dim-isotropic-characterization} and with our earlier discussion about the graphs of rigid motions.

\begin{proposition}\label{prop:all-3-flats}
For each $\rho=[y]\in E^+$, we have $F_\rho^{\Phy}=\P(\bar y\star\OO_s)$. For each $\tau=[y]\in E^-$, we have $F_\tau^{\Phy}=\P(\OO_s\star \bar y)$. Thus the set of $3$-planes in $E^{\Phy}$ is precisely $\{F_{\rho}^{\Phy}\}_{\rho\in E^+\sqcup E^-}$. 
Moreover, for any rigid motion $\rho\colon\R^3\to\R^3$, its graph satisfies $\iota(\phi(\Gamma_\rho))\subset F_{\tilde\rho}^{\Phy}$ for a unique $\tilde\rho\in E^+\sqcup E^-$.
\end{proposition}

\begin{proof}
For $y,z\in Z(N)$, we have $y\star(\bar y\star z)=N(y)z=0$. Thus $\bar y\star\OO_s\subseteq \{z\in\OO_s: y\star z=0\}$. The left-hand side is a 4-dimensional subspace by \cref{thm:4-dim-isotropic-characterization}. Now $z\mapsto y\star z$ is a linear operator on $\OO_s$, so the rank-nullity theorem implies that $\dim \{z\in\OO_s: y\star z=0\}+\dim(y\star \OO_s)=8$. The second summand is 4, again by \cref{thm:4-dim-isotropic-characterization}, implying that $\{z\in\OO_s: y\star z=0\}$ has dimension 4 and thus is equal to $\bar y\star\OO_s$. Thus if $\rho=[y]\in E^+$, we have $F_{\rho}^{\Phy}=\P(\bar y\star \OO_s)$.

Similarly, we have $(z\star\bar y)\star y=z N(y)=0$, so $\OO_s\star \bar y\subseteq \{z\in\OO_s: z\star y=0\}$. Then the same proof shows that if $\tau=[y]\in E^-$, we have $F_{\tau}^{\Phy}=\P(\OO_s\star\bar y)$.

By \cref{thm:4-dim-isotropic-characterization}, each $3$-plane in $E^{\Phy}$ is $\P(y\star \OO_s)$ or $\P(\OO_s\star y)$ for some nonzero $y\in Z(N)$ where $y$ is unique up to scaling. By what we have shown so far, each of these 3-planes is $F_{\rho}^{\Phy}$ for a unique $\rho\in E^+\sqcup E^-$.

For the last part, recall that we showed \cref{sec:graph-and-compactification}, for every rigid motion $\rho\colon\R^3\to\R^3$, its graph is such that $\iota(\phi(\Gamma_\rho))$ is embedded into a 3-plane in $E^{\Phy}$. By what we have shown so far, this means that its closure is $F_{\tilde\rho}^{\Phy}$ for a unique $\tilde\rho\in E^+\sqcup E^-$.
\end{proof}

So far we have shown that $E^+\sqcup E^-$ contains the set of rigid motions of $\R^3$. We will shortly show that the orientation-preserving rigid motions lie in $E^+$ and the orientation-reversing in $E^-$. From now on we will identify rigid motions $\rho\colon\R^3\to\R^3$ with the unique point $\tilde\rho\in E^+\sqcup E^-$ provided by \cref{prop:all-3-flats}. We also use the word rigid motion to refer to an element of $E^+\sqcup E^-$. The ones that correspond to an actual rigid motion of $\R^3$ will be referred to as \emph{genuine rigid motions}, denoted $E^+_g\subset E^+$ and $E^-_g\subset E^-$. We will later show that the genuine rigid motions form a dense open subset of $E^+\sqcup E^-$, so $E^+$ is a compactification of the space of orientation-preserving rigid motions and similarly with $E^-$. Note that cutting away infinity, i.e., passing back to a copy of $\R^6$ embedded in $E^+$ or $E^-$ cuts away some genuine rigid motions. Despite this, this is an operation that we sometimes wish to perform. Later in this section we characterize which genuine rigid motions we miss in this way.

\subsection{Geometric triality}
\label{ssec:triality}
Analogous to $F_\rho^{\Phy}$ and $F_\tau^{\Phy}$, we define two families of $3$-planes in each of $E^+$ and $E^-$.

\begin{definition}
For each $x\in E^{\Phy}$, define the \emph{physical 3-plane in positive space}
\[F_x^{+}=\{\rho\in E^{+}:\rho\star x=0\}\subset E^{+},\]
and for each $\tau\in E^-$, define the \emph{negative 3-plane in positive space}
\[F_\tau^{+}=\{\rho\in E^{+}:\tau\star \rho=0\}\subset E^{+}.\]    
Similarly, for each $x\in E^{\Phy}$, define the \emph{physical 3-plane in negative space}
\[F_x^{-}=\{\tau\in E^{-}:x\star \tau=0\}\subset E^{-},\]
and for each $\rho\in E^+$, define the \emph{positive 3-plane in negative space}
\[F_\rho^{-}=\{\tau\in E^{-}:\tau\star \rho=0\}\subset E^{-}.\]    
\end{definition}

Note that these six families of $3$-planes exhibit many symmetries. For example, it follows immediately from the definitions that for $x\in E^{\Phy}$ and $\rho\in E^+$, we have $x\in F_\rho^{\Phy}$ if and only if $\rho\in F_x^+$. A similar result holds for any pair of the three spaces. We now explore more of these symmetries and some useful geometric consequences of them.
This phenomenon is known as \emph{geometric triality}. (For a more comprehensive exposition on this topic, see \cite[Chapter IV]{Che54}.)

The properties listed below will be apparent geometrically for genuine pairs of points $(p,q)\in\R^3\times\R^3$ and genuine rigid motions $\rho,\tau$.
However, to show that the properties hold in our compactified spaces and over both $\R$ and $\CC$, we would need to either compute equations explicitly or use properties of $\OO_s$. We will prove the following results rigorously using the theory of the split octonions and also mention the corresponding geometric analogue for genuine pairs of points and genuine rigid motions which should be self-evident.

\begin{proposition}
\label{prop:intersection-from-diff-family}
For each $\rho\in E^+$ and $\tau\in E^-$, the intersection $F_{\rho}^{\Phy}\cap F_{\tau}^{\Phy}$ is either a single point or a $2$-plane.
Moreover, the following are equivalent:
    \begin{itemize}
        \item $\rho\in F_{\tau}^+$;
        \item $\tau\in F_{\rho}^-$;
        \item $F_{\rho}^{\Phy}\cap F_{\tau}^{\Phy}$ is a $2$-plane.
    \end{itemize}
The analogous statement holds for $x\in E^{\Phy}, \rho\in E^+$ and also $\tau\in E^-, x\in E^{\Phy}$.
\end{proposition}

\begin{proof}
By \cref{prop:all-3-flats}, we have that $F_{\rho}^{\Phy}=\PP(\bar{\rho}\star \OO_s)$ and $F_{\tau}^{\Phy}=\PP(\OO_s\star \bar{\tau})$. By \cite[Theorem 5]{vdBSpr60}, the intersection $(\bar{\rho}\star \OO_s)\cap (\OO_s\star \bar{\tau})$ has dimension $1$ or $3$, with the dimension $3$ case occurring if and only if $\bar{\rho}\star \bar{\tau}=0$. This shows the first part of the result.

Moreover, $F_{\rho}^{\Phy}\cap F_{\tau}^{\Phy}$ is a $2$-plane if and only if $\bar{\rho}\star \bar{\tau}=0$, which occurs if and only if $\tau \star \rho = 0$, and thus is equivalent to $\rho\in F_{\tau}^+$ and also $\tau\in F_{\rho}^-$.

The analogous statements follow via similar proofs.
\end{proof}

Now if $x\in E^{\Phy}$ corresponds to a genuine pair of points $(p,q)$, (i.e., $x=\iota(\phi(p,q))$) and $\rho\in E^+$, $\tau\in E^-$ are genuine rigid motions, then we have the following geometric interpretations for all six families:
\begin{itemize}
    \item $x\in F_{\rho}^{\Phy}$ if and only if $\rho\in F_x^+$ if and only if $\rho(p)=q$;
    \item $x\in F_{\tau}^{\Phy}$ if and only if $\tau\in F_x^-$ if and only if $\tau(p)=q$;
    \item $\rho\in F_{\tau}^+$ if and only if $\tau\in F_{\rho}^-$ if and only if $\rho,\tau$ agree on a plane.
\end{itemize}

This follows by \cref{prop:all-3-flats,prop:intersection-from-diff-family} since $(p,q)\in\Gamma_\rho$ if and only if $\rho(p)=q$ and $\rho,\tau$ agree on a plane if and only if $\Gamma_\rho\cap\Gamma_\tau$ is a 2-flat.

Having dealt with intersections of $3$-planes from different families, we now deal with intersections of $3$-planes from the same family.

\begin{proposition}\label{prop:intersection-from-same-family}
Given distinct $x, x'\in E^{\Phy}$, the intersection $F_x^+\cap  F_{x'}^+$ is either a line or empty. Moreover, the following are equivalent:
    \begin{itemize}
        \item $F_x^+\cap F_{x'}^+$ is a line;
        \item $\langle x,x'\rangle=0$;
        \item the line connecting $x,x'$ lies completely in $E^{\Phy}$.
    \end{itemize}
The analogous statement holds for any two of $E^{\Phy}, E^+, E^-$.
\end{proposition}

\begin{proof}
By \cref{prop:all-3-flats}, we have that $F_x^+ = \PP(\OO_s\star\bar{x})$ and $F_{x'}^+=\PP(\OO_s\star\bar{x}')$. By \cite[Theorem 4]{vdBSpr60}, the intersection $F_x^+\cap F_{x'}^+$ is either empty or a line, and the latter occurs if and only if $\langle \bar{x}, \bar{x}'\rangle=0$.
Note that 
    \[\langle x,x'\rangle = N(x+x')-N(x)-N(x') = N(\overline{x+x'})-N(\bar{x})-N(\bar{x}') = \langle \bar{x},\bar{x}'\rangle,\]
showing that it is equivalent to $\langle x,x'\rangle = 0$ as well.

Now the line connecting $x,x'$ lies completely in $E^{\Phy}$ if and only if $N(\lambda x+\mu x')=\tfrac12\langle \lambda x+\mu x',\lambda x+\mu x'\rangle=0$ for all $\lambda,\mu$. Since $N(x)=\tfrac12\langle x,x\rangle=N(x')=\tfrac12\langle x',x'\rangle=0$, this condition is equivalent to $\langle x,x'\rangle=0$, as desired.

The analogous statements hold via similar proofs.
\end{proof}

The corresponding geometric interpretation is as follows.
If $x=\iota(\phi(p,q))$ and $x'=\iota(\phi(p',q'))$ for genuine pairs of points $(p,q),(p',q')\in \RR^6$, then there is some orientation-preserving rigid motion $\rho$ with $\rho(p)=q$ and $\rho(p')=q'$ if and only if $\norm{p-p'}=\norm{q-q'}$ which occurs if and only if the line spanned by $\phi(p,q),\phi(p',q')$ is isotropic. As we showed in \cref{sec:graph-and-compactification}, the embedding $\iota$ sends isotropic lines to lines that lie entirely in $E^{\Phy}$. Moreover, if such $\rho$ exists, then we may precompose any rotation around the line connecting $p,p'$ as well, getting a one-dimensional family of such $\rho$ -- this family is $F_x^+\cap F_{x'}^+$.

We are now in a position to rigorously show that orientation-preserving rigid motions correspond to points of $E^+$ and orientation-reversing to points of $E^-$.

\begin{proposition}
\label{prop:realizable-rigid-motions}
Let $\rho\colon\R^3\to\R^3$ be an orientation-preserving rigid motion. There exists a unique $\tilde\rho\in E^+$ such that $\iota(\phi(\Gamma_\rho))\subset F_{\tilde\rho}^{\Phy}$. Similarly if $\tau\colon\R^3\to\R^3$ is an orientation-reversing rigid motion, there exists a unique $\tilde\tau\in E^-$ such that $\iota(\phi(\Gamma_\tau))\subset F_{\tilde\tau}^{\Phy}$.
\end{proposition}

\begin{proof}
By \cref{prop:all-3-flats}, for each rigid motion $\rho\colon\R^3\to\R^3$, there is a unique $\tilde\rho\in E^+\sqcup E^-$ such that $\iota(\phi(\Gamma_\rho))\subset F_{\tilde\rho}^{\Phy}$. We need to check that the orientation-preserving rigid motions correspond with points of $E^+$ and the orientation-reversing rigid motions correspond with points of $E^-$.

Let $\Id\colon\R^3\to\R^3$ be the identity rigid motion. Then $\phi(\Gamma_{\Id})=\R^3\times\{0\}$ and $\iota(\phi(\Gamma_{\Id}))=\{[1:\vec x:0:\vec 0]:\vec x\in\R^3\}$. Now define $\tilde\Id=[0:\vec 0:1:\vec 0]\in E^+$. We compute
\[F_{\tilde\Id}^{\Phy} = \{x:\tilde\Id\star x=0\} =([x_0:\vec{x}:0:\vec 0])_{(x_0,\vec{x})\neq (0,\vec 0)}.\]
This set contains $\iota(\phi(\Gamma_{\Id}))$, showing that $\tilde\Id$ is indeed the point corresponding to the rigid motion $\Id$.

If $\rho,\tau\colon\R^3\to\R^3$ are rigid motions which agree on a plane, then $\iota(\phi(\Gamma_\rho\cap\Gamma_\tau))$ is a 2-flat which is contained in $F_{\tilde\rho}^{\Phy}\cap F_{\tilde\tau}^{\Phy}$. Now if $\tilde\rho,\tilde\tau$ both lay in $E^+$ or both lay in $E^-$, this would contradict \cref{prop:intersection-from-same-family}. Thus, if $\tau=r\circ\rho$ for some reflection $r$, then one of $\tilde\rho,\tilde\tau$ lies in $E^+$ and the other in $E^-$. From this we conclude the desired result, as every orientation-reversing rigid motion can be written as the composition of 1 or 3 reflections and every orientation-preserving rigid motion can be written as the composition of 2 or 4 reflections.
\end{proof}

Using \cref{prop:intersection-from-same-family,prop:intersection-from-diff-family}, we may now characterize the $2$-planes in the three spaces and also show a bijection between the lines of the three spaces.

\begin{proposition}\label{prop:all-2-flats}
Any $2$-plane in $E^{\Phy}$ lies in a unique positive $3$-plane and a unique negative $3$-plane. The analogous statements holds for $E^+,E^-$.
\end{proposition}

\begin{proof}
We will first show that any $3$-dimensional totally isotropic subspace $W$ in $\OO_s$ is contained in precisely two maximal isotropic subspaces.

Write $e_0,e_1,e_2,e_3,e_{0^*},e_{1^*},e_{2^*},e_{3^*}$ for the standard basis vectors of $\OO_s$. Note that $\spn\{e_1,e_2,e_3\}$ is a totally isotropic subspace, as each element of $\spn\{e_1,e_2,e_3\}$ is of the form $(0,\vec x,0,\vec 0)$ and clearly satisfies $N(0,\vec x,0,\vec 0)=0$. Let $\{v_1,v_2,v_3\}$ be an arbitrary basis of $W$, and consider the linear isomorphism $L\colon W\to\spn\{e_1,e_2,e_3\}$ which maps $v_i$ to $e_i$. This map is an isometry (meaning that it is an injective linear map which preserves the bilinear form), since the bilinear form vanishes identically on both $W$ and $\spn\{e_1,e_2,e_3\}$. Thus by Witt's theorem \cite{Wit37}, the isometry $L$ extends to an isometry $\tilde L\colon \OO_s\to\OO_s$. Thus $\tilde L$ gives a 1-to-1 correspondence between the isotropic subspaces of $\OO_s$ containing $W$ and those containing $\spn\{e_1,e_2,e_3\}$.

Now it is not hard to compute that the only 4-dimensional totally isotropic subspaces containing $\spn\{e_1,e_2,e_3\}$ are $\spn\{ e_0,e_1,e_2,e_3\}$ and $\spn\{ e_{0^*},e_1,e_2,e_3\}$. To see this, note that any such subspace must contain a nonzero vector $x=(x_0,\vec 0,x_{0^\ast},\vec x_\ast)$. Since $0=N(x)=x_0x_{0^\ast}-\vec 0\cdot\vec x_\ast$, we see either $x_0=0$ or $x_{0^\ast}=0$. Then since $\spn\{x,e_1,e_2,e_3\}$ is totally isotropic, we see that 
\[\sang{x,(0,\vec y,0,\vec 0)}=\vec x_\ast\cdot \vec y=0\]
for all $\vec y\in\F^3$, showing that $x=e_0$ or $x=e_{0^\ast}$.

    This shows that every $2$-plane in $E^{\Phy}$ lies in exactly two $3$-planes in $E^{\Phy}$.
    By \cref{prop:all-3-flats}, each of them is either a positive $3$-plane or a negative $3$-plane.
    If they are both positive or both negative, this would contradict \cref{prop:intersection-from-same-family}.
    Therefore one must be positive and one negative.
\end{proof}

The geometric intuition of this proposition is that if we prescribe how a rigid motion maps points on a given plane, then there are two ways to extend this to the entire $\RR^3$ -- one orientation-preserving, and one orientation-reversing.

\begin{proposition}\label{prop:all-lines}
Let $\cL^{\Phy}, \cL^+$ and $\cL^-$ be the set of lines in $E^{\Phy}, E^+$ and $E^-$, respectively.
There is a unique 1-to-1-to-1 mapping between lines $\ell^{\Phy}\in\cL^{\Phy}$, lines $\ell^+\in\cL^+$, and lines $\ell^-\in\cL^-$ with the following property: for every $x\in \ell^{\Phy}$ and $\rho\in \ell^+$ and $\tau\in \ell^-$, we have
\[x\in F^{\Phy}_{\rho}\cap F^{\Phy}_{\tau}\quad\text{and}\quad\rho\in F^+_x\cap F^+_{\tau}\quad\text{and}\quad\tau\in F^-_x\cap F^-_{\rho}.\]
\end{proposition}
\begin{proof}
For any line $\ell^{\Phy}\subset E^{\Phy}$, we first show that there is only one candidate for the corresponding $\ell^+$.
We know that $\ell^+$ must satisfy $\ell^+\subset F^+_x$ for every $x\in \ell^{\Phy}$.
Fix distinct $x',x''\in\ell^{\Phy}$. We claim that
\begin{equation}
\label{eq:two-intersection-all-intersection}
F_{x'}^+\cap F_{x''}^+=\bigcap_{x\in\ell^{\Phy}} F_x^+.
\end{equation}
First note that since $x,x'\in\ell^{\Phy}$ (which is a line contained in $E^{\Phy}$), by \cref{prop:intersection-from-same-family}, we have that the left-hand side is a line.

Now pick two distinct $\rho,\rho'\in F_{x'}^+\cap F_{x''}^+$. By definition, we also have $x',x''\in F_\rho^{\Phy}\cap F_{\rho'}^{\Phy}$. The latter set is a line by \cref{prop:intersection-from-same-family} (since it is not empty), so it must be equal to $\ell^{\Phy}$. Thus for any $x\in \ell^{\Phy}$, we have $x\in F_\rho^{\Phy}\cap F_{\rho'}^{\Phy}$, implying that $\rho\in F_x^+$. However $\rho$ was an arbitrary element of $F_{x'}^+\cap F_{x''}^+$ and $x$ was an arbitrary element of $\ell^{\Phy}$. Thus we have shown \cref{eq:two-intersection-all-intersection}; furthermore, this quantity is the only possibility for $\ell^+$.
We may similarly define $\ell^- = \bigcap_{x\in \ell^{\Phy}} F_x^-$.

With these definitions, it immediately follows that for all $x\in\ell^{\Phy}$, we have $\ell^+\subset F_x^+$ and $\ell^-\subset F_x^-$. Furthermore, we showed above that for all $\rho\in\ell^+$, we have $\ell^{\Phy}\subset F_\rho^{\Phy}$ and for all $\tau\in\ell^-$, we have $\ell^{\Phy}\subset F_{\tau}^{\Phy}$. To complete the proof, it remains to show that for all $\rho\in\ell^+$, we have $\ell^-\subset F_\rho^-$ and for all $\tau\in\ell^-$, we have $\ell^+\subset F_\tau^+$.

To prove these, pick arbitrary $\rho\in\ell^+$ and $\tau\in\ell^-$. It suffices to show that $\tau\in F_\rho^-$ and $\rho\in F_\tau^+$. To see this, note that $F_{\rho}^{\Phy}\cap F_{\tau}^{\Phy}$ contains $\ell^{\Phy}$, and thus it is a $2$-flat by \cref{prop:intersection-from-diff-family}. Moreover, \cref{prop:intersection-from-diff-family} implies that $\rho\in F_{\tau}^+$ and $\tau\in F_{\rho}^-$, as desired.
\end{proof}

In the rest of the paper, given a line $\ell^{\Phy}$ in $E^{\Phy}$, we will use $\ell^+$ and $\ell^-$ to denote the corresponding lines given by \cref{prop:all-lines} and similarly if we start with a line in the other two spaces.

Geometrically, for any two lines $\ell,\ell'$ in $\RR^3$ and any isometry $i\colon \ell\to \ell'$, the set of pairs $\{\phi(p,q):(p,q)\in \ell\times \ell', i(p)=q\}$ is a line $\ell^{\Phy}$ in $E^{\Phy}$.
The set of orientation-preserving rigid motions $\rho$ satisfying $\rho|_{\ell}=i$ corresponds to the line $\ell^+\subset E^+$, and the set of orientation-reversing rigid motions $\tau$ satisfying $\tau|_{\ell}=i$ corresponds to the line $\ell^-\subset E^-$.

\begin{corollary}\label{cor:all-lines-through-a-point}
For a point $x\in E^{\Phy}$, define $\cL^{\Phy}(x)$ to be the set of lines in $E^{\Phy}$ through $x$ and $\cL^+(x)$ to be the set of lines contained in $F^+_x$. Then $\ell^{\Phy}\mapsto \ell^+$ is a bijection from $\cL^{\Phy}(x)$ to $\cL^+(x)$.
\end{corollary}
\begin{proof}
\cref{prop:all-lines} implies that if $x\in \ell^{\Phy}$, then $\ell^+\subset F_x^+$. Therefore if $\ell^{\Phy}\in \cL^{\Phy}(x)$, then $\ell^+\in \cL^+(x)$.
Conversely, if $\ell^+\in \cL^+(x)$, let $\rho$ and $\rho'$ be two distinct points on $\ell^+$. By \cref{prop:all-lines}, we know that $\ell^{\Phy}\subseteq F^{\Phy}_\rho\cap F^{\Phy}_{\rho'}$. By \cref{prop:intersection-from-same-family}, this containment must be an equality. Finally $\rho,\rho'\in\ell^+\subset F_x^+$ implies that $x\in F^{\Phy}_{\rho}\cap F^{\Phy}_{\rho'} = \ell^{\Phy}$. This shows that $\ell^{\Phy}\in\cL^{\Phy}(x)$, as desired.
\end{proof}

Since $F_x^+\cong \PP^3$, the set $\cL^+(x)$ -- and thus also $\cL^{\Phy}(x)$ -- can be identified with the Grassmanian of lines in $\PP^3$, denoted $\G(1,3)$. As we will show in \cref{ssec:plucker}, $\cL^{\Phy}(x)$ inherits a coordinate system from $E^{\Phy}$ that is (up to a linear transformation) the same as the usual \emph{Pl\"ucker embedding} of $\G(1,3)$ in $\PP^5$.

\subsection{Affine coordinates}
To apply various tools from incidence geometry, there are several times throughout the argument where it will be more convenient to work in $\FF^6$ instead of $E^{\Phy},E^+,E^-$. Just as we embedded $\RR^6$ into $E^{\Phy}$, we will choose a distinguished dense open subset of $E^+, E^-$ that we identify with $\F^6$.

\begin{definition}
Let $\FF=\RR$ or $\CC$. The notation $E^{\Phy}_o(\FF),E^+_o(\FF),E^-_o(\FF)$ refers to three copies of $\FF^6$. These are identified with dense open subsets of $E^{\Phy}(\FF),E^+(\FF),E^-(\FF)$, respectively, via the embedding $\iota\colon E^{\Phy}_o(\FF)\to E^{\Phy}(\FF)$ defined in \cref{eq:iota-phys-defn}, the embedding $\iota_+\colon E^{+}_o(\FF)\to E^{+}(\FF)$, and the embedding $\iota_-\colon E^{-}_o(\FF)\to E^{-}(\FF)$ where the latter two are defined by the same formula
\begin{align*}
\iota_+(\vec{\rho},\vec{\rho}_*)&=[\vec{\rho}\cdot \vec{\rho}_*:\vec{\rho}:1:\vec{\rho}_*],\\
\iota_-(\vec{\tau},\vec{\tau}_*)&=[\vec{\tau}\cdot \vec{\tau}_*:\vec{\tau}:1:\vec{\tau}_*].
\end{align*}
\end{definition}

Note that while $\iota_+,\iota_-$ are defined by the same formula, $\iota$ is defined by a different formula. Now $\iota$ identifies $E^{\Phy}_o$ with the affine open subset $\{x_0\neq 0\}\subset E^{\Phy}$, while $\iota_+$ identifies $E^+_o$ with $\{\rho_{0^\ast}\neq 0\}\subset E^+$, and $\iota_-$ identifies $E^-_o$ with $\{\tau_{0^\ast}\neq 0\}\subset E^-$.

We point out one potential confusion with our terminology: for two points $x,x'\in E^{\Phy}_o$, they span a line $\aff\{x,x'\}$ in $E^{\Phy}_o$. However, the corresponding points $\iota(x),\iota(x')$ span a line $\aff\{\iota(x),\iota(x')\}$ in $\PP^7$ which may or may not lie in $E^{\Phy}$. Indeed, we showed earlier that $\aff\{\iota(x),\iota(x')\}$ lies in $E^{\Phy}$ if and only if $\aff\{x,x'\}$ is isotropic. We believe that this will not cause too much confusion in the argument; indeed we will never consider a non-isotropic line in $E^{\Phy}_o$.

Now consider $x=(\vec x,\vec x_*)\in E^{\Phy}_o$ and $\rho=(\vec\rho,\vec\rho_*)\in E^+_o$. We want to understand when $\rho\in F_x^+$. (Technically, this should say $\iota_+(\rho)\in F_{\iota(x)}^+$, but we will abuse notation by dropping the $\iota$'s here and throughout the paper.) This occurs if and only if 
\[\begin{bmatrix}
\vec\rho\cdot\vec\rho_* & \vec\rho \\ \vec \rho_*&1
\end{bmatrix}
\star
\begin{bmatrix}
1 & \vec x \\ \vec x_*&\vec x\cdot \vec x_*
\end{bmatrix}
=0.\]
Expanding the definition and simplifying, this is equivalent to
\begin{equation}
\label{eq:f-x-positive}
\vec{\rho}_*=-\vec{x}\times \vec{\rho}-\vec{x}_*
\end{equation}
A similar computation shows that $\tau=(\vec{\tau},\vec{\tau}_*)\in F_x^-$ if and only if
\begin{equation}
\label{eq:f-x-negative}
\vec{\tau} = -\vec{x}_*\times \vec{\tau}_*-\vec{x}.
\end{equation}

From \cref{eq:f-x-positive}, it is clear that for each $x\in E^{\Phy}_o$, the set $\iota_+^{-1}(F_x^+\cap E^+_o)$ is a 3-flat in $E^+_o$. In other words, if $x\in E^{\Phy}_o$ then $F_x^+\not\subseteq E^+\setminus E^+_o$. We will abuse notation and use $F_x^+$ to refer to both the 3-plane in $E^+$ and the corresponding 3-flat in $E^+_o$ (and similarly for the other types of 3-planes). In addition, for any $p,q\in \RR^3$, we write $F_{pq}^{+}$ for $F_x^{+}$ where $x = \phi(p,q)$ and similarly we write $F_{pq}^-$ for $F_x^-$.

Next we study the objects that we miss when we restrict to $E^{\Phy}_o, E^+_o, E^-_o$.

\begin{proposition}\label{prop:actual-rigid-motion}
A point $\rho\in E^+(\RR)$ is a genuine orientation-preserving rigid motion if and only if $F_{\rho}^{\Phy}\not\subseteq E^{\Phy}\setminus E^{\Phy}_o$. The analogous statement holds for $\tau\in E^-(\RR)$.
\end{proposition}
\begin{proof}
If $\rho$ is a genuine orientation-preserving rigid motion, then its graph satisfies $\phi(\Gamma_\rho)=\iota^{-1}(F_{\rho}^{\Phy})$, which shows that $F_{\rho}^{\Phy}\cap E^{\Phy}_o$ is non-empty. On the other hand, suppose that $\iota^{-1}(F_{\rho}^{\Phy})$ is non-empty for some $\rho\in E^+(\R)$. We can write $F_{\rho}^{\Phy}=\P(V)$ for some 4-dimensional linear subspace $V\subseteq Z(N)\subseteq\OO_s$. Then $\iota^{-1}(F_{\rho}^{\Phy})$ is the intersection of $V$ with the affine hyperplane $\{x_0=1\}$, projected onto the $(\vec x,\vec x_{\ast})$ coordinates. Since the intersection is non-empty, it must be a 3-flat, and as the projection is a linear map that is 1-to-1 on $Z(N)\cap\{x_0=1\}$, we see that $\iota^{-1}(F_{\rho}^{\Phy})$ is a 3-flat.
    
Now for every $(p,q),(p',q')\in \phi^{-1}(\iota^{-1}(F_{\rho}^{\Phy}))$, we know that the line between $\iota(\phi(p,q))$ and $\iota(\phi(p',q'))$ lies in $F_\rho^{\Phy}\subset E^{\Phy}$. As we showed earlier, this implies that the line between $\phi(p,q)$ and $\phi(p',q')$ is isotropic, so $\norm{p-p'}^2=\norm{q-q'}^2$.
    
As we are working over $\R$, this first implies that there do not exist $(p,q),(p,q')\in\phi^{-1}(\iota^{-1}(F_{\rho}^{\Phy}))$ with $q\neq q'$ (since we would have $\norm{q-q'}^2=\norm{p-p}^2=0$). Thus this set is the graph of some function $\tilde\rho\colon\R^3\to\R^3$. Furthermore, we see that $\norm{p-p'}=\norm{\tilde\rho(p)-\tilde\rho(p')}$ for all $p,p'\in\R^3$, so $\tilde\rho$ is a rigid motion. Finally, by \cref{prop:realizable-rigid-motions}, we see that $\tilde\rho$ must be orientation-preserving.
\end{proof}

A consequence is that if we are dealing with $\rho\in E^+(\RR)$ and $\tau\in E^-(\RR)$ that interact with a genuine pair of points, then those $\rho,\tau$ are genuine rigid motions.
Therefore in most scenarios (when working over $\R$), we can assume that we are only dealing with genuine rigid motions and not the additional points that we add in $E^+$ and $E^-$.

However, the converse to this idea is not true: there are some genuine rigid motions that do not lie in $E_o^+(\R)$ and $E^-_o(\R)$. Fortunately, there are few of these and they are easy to characterize.

\begin{proposition}\label{prop:missing-rigid-motion}
If $\rho\in E^+(\R)\setminus E^+_o(\R)$ is a genuine orientation-preserving rigid motion, then $\rho$ is a screw displacement with angle $\pi$.

If $\tau\in E^-(\R)\setminus E^-_o(\R)$ is a genuine orientation-reversing rigid motion, then $\tau$ is a glide reflection.
\end{proposition}

\begin{proof}
For such a $\rho$ we have $\rho=[\rho_0:\vec \rho:0:\vec\rho_*]$ where $\vec\rho\cdot\vec\rho_*=0$. For $(p,q)$ with $q=\rho(p)$, write $x=\iota(\phi(p,q))\in E^{\Phy}$, i.e., $x=[1:\vec x: (\vec x\cdot \vec x_*):\vec x_*]$ where $\vec x=(p+q)/2$ and $\vec x_*=(p-q)/2$. Now $x\in F_\rho^{\Phy}$ which implies $\rho\star x = 0$, i.e.,
\[\begin{bmatrix}
\rho_0 & \vec\rho \\ \vec \rho_*&0
\end{bmatrix}
\star
\begin{bmatrix}
1 & \vec x \\ \vec x_*&\vec x\cdot \vec x_*
\end{bmatrix}
=0.\]
Expanding the definition and simplifying, this is equivalent to
\begin{align*}
\rho_0+\vec{\rho}\cdot \vec x_*&=0\\
\vec{\rho}_*-\vec{\rho}\times \vec{x}&=0.
\end{align*}

The first equation implies that implies that $\vec{\rho}, \vec \rho_*$ cannot be simultaneously zero. Using $\vec x=(p+q)/2$, we see
\[\vec{\rho}\times \frac{p+q}{2} = \vec{\rho}_*\]
holds for all $p,q\in\R^3$ with $\rho(p)=q$. For fixed vectors $\vec \rho,\vec \rho_*$, not simultaneously zero, the set of solutions to $\vec{\rho}\times v=\vec\rho_*$ is a fixed line (or empty). This shows that the midpoint of $p$ and $\rho(p)$ lies on a fixed line for all $p \in\R^3$. Every genuine orientation-preserving rigid motion of $\R^3$ is a screw displacement by some angle $\theta$. It is easy to see that the only ones with this property are those where $\theta=\pi$.
    
Similarly, consider $\tau=[\tau_0:\vec \tau:0:\vec\tau_*]$. For $(p,q)$ with $q=\tau(p)$, write $x=\iota(\phi(p,q))$ as before. Then $x\star \tau=0$, i.e.,
\[
\begin{bmatrix}
1 & \vec x \\ \vec x_*&\vec x\cdot \vec x_*
\end{bmatrix}
\star
\begin{bmatrix}
\tau_0 & \vec\tau \\ \vec \tau_*&0
\end{bmatrix}
=0,\]
which is equivalent to
\begin{align*}
\tau_0+\vec{x}\cdot \vec \tau_*&=0\\
\vec{\tau}+\vec{x}_*\times \vec{\tau}_*&=0.
\end{align*}

Once again, the first equation implies that that $\vec{\tau}, \vec{\tau}_*$ cannot be simultaneously zero. Thus
\[\vec{\tau}_*\times \left(\frac{p-q}{2}\right)=\vec{\tau}\]
holds for every pair $p,q\in\R^3$ with $\tau(p)=q$. As before, this implies that $p-\tau(p)$ lies on a fixed line for all $p\in\R^3$. Every genuine orientation-reversing rigid motion of $\R^3$ is either a glide reflection or a rotoreflection.\footnote{A reflection through a plane composed with rotation by some angle $\theta$ around a line perpendicular to the plane.} It is easy to see that the only rigid motions with the desired property are the glide reflections (and the rotoreflections with angle $\theta=0$ which also are a degenerate glide reflection).
\end{proof}

Now we characterize the set of isotropic lines that we miss, i.e., those $\ell^{\Phy}\subset E^{\Phy}_o(\R)$ with $\ell^+\subset E^+(\R)\setminus E^+_o(\R)$ and $\ell^-\subset E^-(\R)\setminus E^-_o(\R)$.

\begin{proposition}
\label{prop:missing-parallelogram}
Let $p,p',q,q'\in\R^3$ be such that \[\phi(\aff\{(p,q),(p',q')\})=\ell^{\Phy}\subset E^{\Phy}_o\]
is an isotropic line.
If $\ell^+\subset E^+(\R)\setminus E^+_o(\R)$, then $pp'qq'$ is a parallelogram (meaning that $p-p'=-(q-q')$).
If $\ell^-\subset E^-(\R)\setminus E^-_o(\R)$, then $pp'q'q$ is a parallelogram (meaning that $p-p'=q-q'$).
\end{proposition}

\begin{proof}
Since $\ell^{\Phy}$ is isotropic, its image $\iota(\ell^{\Phy})$ in $E^{\Phy}$ is a dense open subset of a line which we also call $\ell^{\Phy}$. Now consider the corresponding lines $\ell^+,\ell^-$ produced by \cref{prop:all-lines}. As $\phi(p,q),\phi(p',q')\in\ell^{\Phy}$, we know that for each $\rho\in\ell^+$, we have $\phi(p,q),\phi(p',q')\in F_\rho^{\Phy}$. As these are genuine pairs of points, \cref{prop:actual-rigid-motion} implies that each such $\rho$ is a genuine rigid motion.

Thus we have found a 1-dimensional family of genuine rigid motions $\rho$ for which $\rho(p)=q$ and $\rho(p')=q'$. Suppose for the sake of contradiction that $\ell^+\subset E^+(\R)\setminus E^+_o(\R)$ but $pp'qq'$ is not a parallelogram, i.e., that $(p-p')+(q-q')\neq 0$. Then \cref{prop:missing-rigid-motion} implies that the 1-dimensional family of rigid motions $\rho$ we found previously are all screw displacements with angle $\pi$. We will reach a contradiction by showing that there is in fact at most one such $\rho$.

As $\rho(p)=q$ where $\rho$ is a screw displacement with angle $\pi$, we know that $(p+q)/2$ lies on the axis of the screw displacement. Thus the axis of the screw displacement must pass through $(p+q)/2$ and $(p'+q')/2$. Our assumption implies that these two points are distinct, so they uniquely determine the axis. Now we know the axis of $\rho$ and that its angle is $\pi$, so the only undetermined quantity is the magnitude of the displacement. However, the equation $\rho(p)=q$ determines this quantity, uniquely determining $\rho$.

Suppose that $\ell^-\subset E^-(\R)\setminus E^-_o(\R)$, but $(p-p')- (q-q')\neq 0$. The same argument as before shows that there is a 1-dimensional family of glide reflections $\tau$ for which $\tau(p)=q$ and $\tau(p')=q'$.

Suppose $\tau$ is a glide reflection across a plane with normal vector $\hat v$ and with displacement $\vec d$ where $\hat v\cdot\vec d=0$. Since $\tau(p)=q$ we see that $p-q-\vec d\in\spn\{\hat v\}$. The same holds for $p'-q'-\vec d$, so we conclude that $(p-p')-(q-q')\in\spn\{\hat v\}$ is parallel to $\hat v$. By hypothesis, this vector is nonzero, which uniquely determines $\pm\hat v$. Furthermore, the reflection plane must pass through $(p+q)/2$ which uniquely determines the plane. Finally, the equation $\tau(p)=q$ uniquely determines $\vec d$, uniquely determining $\tau$.
\end{proof}

\subsection{The \texorpdfstring{Pl\"ucker}{Plucker} embedding}
\label{ssec:plucker}
In this section, we expand on the relationship between the Pl\"ucker embedding and the set of lines through a single point in $E_o^{\Phy}$. From \cref{defn:isotropic}, it follows that a line $\ell^{\Phy}\subset E_o^{\Phy}\cong \F^6$ is isotropic if it has direction $[\vec d:\vec d_\ast]$ where $\vec d\cdot\vec d_{\ast}=0$. These are the lines which map to a line under the embedding $\iota\colon E^{\Phy}_o\to E^{\Phy}$.

For a point $x\in E^{\Phy}_o$, note that the set $\cL^{\Phy}(x)$, which we previously defined as the set of lines through $x$ which are contained in $E^{\Phy}$, can be equivalently thought of as the set of isotropic lines through $x$ in $E^{\Phy}_o\cong \F^6$. With the latter viewpoint in mind, we define homogeneous coordinates on $\cL^{\Phy}(x)$ as follows.

\begin{definition}
For $x\in E^{\Phy}_o$, define $\cL^{\Phy}(x)=\{\ell^{\Phy}\subset E^{\Phy}: x\in \ell^{\Phy}\}$ and let $\cL^+(x)=\{\ell^+\subset F_x^+\subset E^+\}$. 
Define
\[\psi_x\colon \cL^{\Phy}(x)\hookrightarrow \{[\vec d:\vec d_*]\in\P^5:\vec d\cdot \vec d_*=0\}\subset\P^5\]
by
\[\psi_x(\ell^{\Phy})=[\vec y-\vec x:\vec y_\ast-\vec x_*]\]
for any $y\in(\ell^{\Phy}\setminus \{x\})\cap E^{\Phy}_o$.
\end{definition}

Clearly $\psi_x\colon \cL^{\Phy}(x)\hookrightarrow\PP^5$ is a well-defined embedding whose image is the set of isotropic directions $\{[\vec{d}:\vec{d}_\ast]:\vec{d}\cdot \vec{d}_\ast=0\}$. 
By \cref{cor:all-lines-through-a-point}, we know that $\cL^{\Phy}(x)$ and $\cL^+(x)$ are isomorphic, so $\psi_x$ also defines an embedding $\cL^+(x)\hookrightarrow\PP^5$ with the same image.

On the other hand, let $\cL$ be the set of lines in $\PP^3$, where we use the coordinates $[x_0:\vec{x}] = [x_0:x_1:x_2:x_3]$. These lines are also parametrized by the classical Pl\"ucker embedding.

\begin{definition}
Write $\G(1,3)$ for the Grassmanian of lines in $\P^3$. Define the \emph{Pl\"ucker embedding}
\[\psi_{\textup{Pl\"ucker}}\colon\G(1,3)\hookrightarrow \{[\vec{d}:\vec{d}_\ast]:\vec{d}\cdot \vec{d}_\ast=0\}\subset\PP^5\]
by
\[\psi_{\textup{Pl\"ucker}}(\ell) = [x_0\vec{y}-y_0\vec{x}:\vec{x}\times \vec{y}]\in \PP^5\]
where $[x_0:\vec x],[y_0:\vec y]\in\ell$ are two distinct points.
\end{definition}

Note that the coordinates of $\psi_{\textup{Pl\"ucker}}$ can be viewed as the determinants of the $2\times 2$ minors of the $2\times 4$ matrix with columns $x,y$.
It is a standard fact (see, e.g., \cite[pg. 102]{SS18}) that $\psi_{\textup{Pl\"ucker}}$ is a well-defined embedding with the given image.

To apply various results from the literature on ruled surface theory, it will be helpful to know that the natural coordinate system in our setting agrees with the Pl\"ucker coordinates.

\begin{proposition}
\label{thm:plucker-embeding}
For any $x\in E^{\Phy}_o$, there is a linear isomorphism $i\colon F^+_x\to \PP^3$ so that for any line $\ell^+\subset F^+_x$, we have
\[ \psi_x (\ell^{\Phy}) = \psi_{\textup{Pl\"ucker}}\circ i(\ell^+).\]
\end{proposition}
\begin{proof}
Recall that $F^+_x\subset E^+$ consists of the points $\rho=[\rho_0:\vec\rho:\rho_{0^\ast}:\vec\rho_\ast]\in E^+$ with $\rho \star \iota(x)=0$.
Expanding the definition, we find
\begin{align*}
\vec{\rho}_* &= \vec{\rho}\times \vec{x}-\rho_{0^*}\vec x_*,\\
\rho_0 &= -\vec{\rho}\cdot \vec{x}_*.
\end{align*}
(This is similar to the computation used to derive \cref{eq:f-x-positive}, but these equations also hold for $\rho\in E^+\setminus E^+_o$.)
Defining the linear map $i\colon F^+_x\to \PP^3$ by $i(\rho) = [\rho_{0^*}:\vec{\rho}]$, the above equations show that this is an isomorphism.

Now for any $\ell^+\subset F^+_x$, we know that $\ell^{\Phy}\cap E^{\Phy}_o$ admits a parametrization $(\vec{x}+\vec{y}t, \vec x_*+\vec y_*t)$ where $\vec{y}\cdot \vec{y}_\ast=0$.
With this parametrization, we have $\psi_x(\ell^{\Phy}) = [\vec{y}:\vec{y}_*]$. Moreover, as $x+yt\in \ell^{\Phy}$ for all $t\in\F$, \cref{prop:all-lines} implies that $\ell^+\subset F^+_{x+yt}$. Then \cref{prop:intersection-from-same-family} shows that $\ell^+=F^+_x\cap F^+_{x+yt}$ for any $t\neq 0$. Expanding the equations for $F_x^+$ defined above and examining the coefficient of $t$, we see that $\ell^+$ consists of $\rho\in F^+_x$ with 
\begin{align*}
\vec{\rho}\cdot \vec{y}_*&=0,\\
\rho_{0^*}\vec{y}_*-\vec{\rho}\times \vec{y}&=0.
\end{align*}
    
We now compute $\psi_{\textup{Pl\"ucker}}\circ i(\ell^+)$. First, if $\vec{y}\neq \vec{0}$, then there exists $\vec{d}\in \FF^3$ such that $\vec{y}_* = \vec{d}\times \vec{y}$. We claim that $[-\vec d\cdot\vec x_\ast:\vec d:1:\vec d\times\vec x-\vec x_\ast]$ and $[-\vec y\cdot \vec x_\ast:\vec y:0:\vec y\times\vec x]$ both lie on $\ell^+$. This holds as long as these points satisfy all four equations that we have derived above. The first two equations are immediate, while the last two equations hold since $\vec d\cdot\vec y_\ast=\vec d\cdot(\vec d\times\vec y)=0$ and $\vec y_\ast -\vec d\times\vec y=\vec0$ for the first point and $\vec y\cdot\vec y_\ast=0$ and $-\vec y\times\vec y =\vec0$ for the second point. Applying $i$, we see that $[1:\vec{d}]$ and $[0:\vec{y}]$ are two (distinct) points on $i(\ell^+)$. Thus
\[\psi_{\textup{Pl\"ucker}}\circ i(\ell^+) = [\vec{y}:\vec{d}\times \vec{y}]=[\vec{y}:\vec{y}_*].\]

On the other hand, if $\vec{y}=\vec{0}$, pick $\vec{d},\vec{d}_\ast\in \FF^3$ such that $\vec{y}_*=\vec{d}\times \vec{d}_\ast$ and $\vec d\cdot\vec d_\ast=0$. (We can do this by picking $\vec d,\vec d_\ast$ to be two perpendicular vectors in the orthogonal complement of $\vec y_*$ of appropriate lengths.)
In this case, we claim that $[-\vec d\cdot\vec x_\ast:\vec d: 0:\vec d\times\vec x]$ and $[-\vec d_\ast\cdot\vec x_\ast:\vec d_\ast: 0:\vec d_\ast\times\vec x]$ both lie on $\ell^+$. As before, the first two equations are obvious, while the last two hold as $\vec d\cdot\vec y_\ast=\vec d\cdot(\vec d\times\vec d_\ast)=0$ and $-\vec d\times \vec 0=\vec 0$ for the first and $\vec d_\ast\cdot\vec y_\ast=\vec d_\ast\cdot(\vec d\times\vec d_\ast)=0$ and $-\vec d_\ast\times \vec 0=\vec0$ for the second. Applying $i$, we see that $[0:\vec d],[0:\vec d_*]$ are two distinct points on $i(\ell^+)$. Thus
\[\psi_{\textup{Pl\"ucker}}\circ i(\ell^+) = [\vec{0}:\vec{d}\times\vec{d}_\ast]=[\vec{y}:\vec{y}_*].\]

In either case we have shown that  $\psi_x(\ell^{\Phy}) =[y]= \psi_{\textup{Pl\"ucker}}\circ i(\ell^+)$.
\end{proof}

\subsection{Setting up distinct distances}

In this section, we will first apply the distance energy argument together with Beck's theorem to reduce the distinct distances problem in $\RR^3$ to an upper bound on the number of 2-rich isotropic lines in $\phi(\cP\times \cP)\subset E^{\Phy}_o(\R)$.
Then we will use what we have set up in this section to move the problem into positive space.

We first recall Beck's theorem.

\begin{theorem}[Beck \cite{Bec83}]\label{thm:beck}
    There are constants $c, C>0$ such that the following holds for all $d\ge 2$ and $N\ge 2$. 
    Let $\cP \subset \R^d$ be a set of size $N$ such that $|\ell \cap \cP| \le N/2$ for any line $\ell$. Then the number of lines $\ell$ such that $|\ell \cap \cP| \in [2, C]$ is at least $cN^2$.
\end{theorem}

\begin{proposition}\label{thm:beck-trick}
    There exists a constant $C\geq 2$ so that the following holds.
    Let $\cP\subset \RR^3$ be any set of points.
    Let $\cL$ be the set of isotropic lines $\ell^{\Phy}\subset E^{\Phy}_o(\R)$ so that
    \[\abs{\left\{\ell^{\Phy}\cap \phi(\cP\times \cP)\right\}}\in [2,C].\]
    Then $\cP$ spans $\gtrsim\min\set{N^4/|\cL|,N}$ distinct distances.
\end{proposition}
\begin{proof}
Note that if there exists $N/2$ points on a line in $\cP$, then these points alone already span at least $N/2-1$ distinct distances. The statement holds in this case, so from now on we assume that no line contains more than $N/2$ points of $\cP$.

By Beck's theorem (\cref{thm:beck}) there are constants $c,C$ so that there are at least $cN^2$ lines containing at least 2 points of $\cP$ and at most $C$ points of $\cP$. Let $B\subseteq \cP\times \cP$ be the set of pairs $(p,p')\in \cP$ so that $p\neq p'$ and $\abs{\aff\{p,p'\}\cap \cP}\in [2,C]$. Then $\abs{B}\geq cN^2$.

Consider the following set of distance quadruples:
    \[\cE_{\mathrm{Beck}} = \left\{(p,q,p',q')\in \cP^4: (p,p'),(q,q')\in B, \, \norm{p-p'}=\norm{q-q'}\right\}.\]
Suppose that $\cP$ spans $t$ distinct distances $\delta_1,\ldots,\delta_t$. Define $r_i$ to be the number of pairs $(p,p')\in B$ so that $\norm{p-p'}=\delta_i$. Then
\[\sum_{i=1}^tr_i=|B|\geq cN^2\qquad\text{ and }\qquad\sum_{i=1}^tr_i^2=\abs{\cE_{\mathrm{Beck}}}.\]
By the Cauchy--Schwarz inequality, we conclude that 
\[\abs{\cE_{\mathrm{Beck}}}\geq \frac{\abs{B}^2}{t} \geq \frac{c^2N^4}{t}.\]

Next we claim that $\abs{\cE_{\mathrm{Beck}}}\leq C^2\abs{\cL}$, implying that, in this case, $t\gtrsim N^4/|\cL|$. To see this, first note that for $(p,q,p',q')\in\cE_{\mathrm{Beck}}$, we have $(p,p'),(q,q')\in B$, implying $p\neq p'$ and $q\neq q'$. Now map the quadruple $(p,q,p',q')$ to the line $\ell^{\Phy}=\aff\{\phi(p,q),\phi(p',q')\}$ in $E^{\Phy}_o$. As $\norm{p-p'}=\norm{q-q'}$, this line is isotropic. Moreover, as $p\neq p'$, all pairs $(p'',q'')\in \cP\times\cP$ with $\phi(p'',q'')\in \ell^{\Phy}$ have distinct $p$-coordinates, and so
\[\abs{\aff\{\phi(p,q),\phi(p',q')\}\cap \phi(\cP\times \cP)}\leq \abs{\aff\{p,p'\}\cap \cP}\leq C.\]
Thus $\ell^{\Phy}\in\cL$. Furthermore, there are at most $C^2$ quadruples in $\cE_{\mathrm{Beck}}$ which correspond to the same line $\ell^{\Phy}$, as each quadruple corresponds to a pair of elements of $\ell^{\Phy}\cap\phi(\cP\times\cP)$. This implies the desired bound $\abs{\cE_{\textup{Beck}}}\leq C^2|\cL|$.
\end{proof}

This result shows that to prove \cref{thm:main}, it suffices to show that the number of $2$-rich isotropic lines for $\phi(\cP\times \cP)\subset E^{\Phy}_o$ is upper bounded by $N^{10/3+\varepsilon}$. Furthermore, we can restrict our attention to the lines which are at most $C$-rich for some absolute constant $C$.

With the theory we have set up in this section, we can now move this problem from physical space to either positive or negative space.

\begin{lemma}\label{thm:analog-elekes-sharir}
Let $\cP\subset \RR^3$ be a set of $N$ points, and let $H\subseteq \cP\times \cP$ be any subset. Define $\cL_H$ to be the set of isotropic lines $\ell^{\Phy}\subset E_o^{\Phy}(\R)$ that are 2-rich for $\phi(H)$. Consider the family of $3$-flats $\cF=\set{F_{pq}^+}_{(p,q)\in H}$ in $E^+_o(\RR)$. Write $\cL_2(\cF)$ for the set of 2-rich lines $\ell^+\subset E_o^+$, i.e., those that lie in at least two of the 3-flats $F\in\cF$. Then for all but at most $N|H|$ lines $\ell^{\Phy}\in\cL_H$, the corresponding line $\ell^+$ lies in $\cL_2(\cF)$.

Analogously, consider the family of $3$-flats $\cF=\set{F_{pq}^-}_{(p,q)\in H}$ in $E^-_o(\RR)$. Write $\cL_2(\cF)$ for the set of 2-rich lines $\ell^-\subset E_o^-$. Then for all but at most $N|H|$ lines $\ell^{\Phy}\in\cL_H$, the corresponding line $\ell^-$ lies in $\cL_2(\cF)$.
\end{lemma}

\begin{proof}
By definition, for any $\ell^{\Phy}\in\cL_H$, there exist $(p,q),(p',q')\in H$ so that $\phi(p,q),\phi(p',q')\in \ell^{\Phy}$. By \cref{prop:all-lines}, we see that $\ell^+\subset F^+_{pq},F^+_{p',q'}$. This means that $\ell^+\in\cL_2(\cF)$ unless $\ell^+$ is a line at infinity, i.e., lies entirely in $E^+\setminus E^+_o$. By \cref{prop:missing-parallelogram}, this implies that $pp'qq'$ is a parallelogram. Thus $q'$ is determined by $(p,q)\in H$ and $p'\in\cP$, so there are at most $N|H|$ of these quadruples which are missed.

For the analogous statement with $F_{pq}^-$, \cref{prop:missing-parallelogram} now implies that $pp'q'q$ is a parallelogram, and the rest of the proof is unchanged.
\end{proof}

\newpage

\part{Line--3-flat incidences}
\label{part:i}
Towards the proof of the main theorem, the main result of \cref{part:i} is \cref{thm:codim-0-incidence-bound}, which reduces the line--3-flat incidence problem to a line-surface incidence problem. The second key result proved in this part is \cref{thm:very-rich-bound-main}, which bounds the number of very rich rigid motions determined by a set of points in $\R^3$ by viewing this problem as a 3-flat--3-flat incidence problem in $\R^6$.

Both of these incidence bounds are proved in \cref{sec:codim-0-incidences}, using the polynomial method and algebraic geometry machinery that we develop in \cref{sec:approx-complete-intersection}. These two bounds both require control on the codimension 0 concentration of the 3-flats $\{F_{pq}^+\}_{(p,q)\in H}$ and $\{F_{pq}^-\}_{(p,q)\in H}$. We prove this in \cref{sec:concentration} (\cref{thm:concentration-codim-0}) using a generalized ruled surface theory that we develop in \cref{sec:csm-for-flats}. This ruled surface theory (specifically \cref{thm:general-CSM}, a Cayley--Monge--Salmon theorem for flats) is built upon novel algebraic geometry results proved in \cref{sec:approx-complete-intersection}, specifically \cref{thm:approx-complete-intersection}.

\section{Approximate complete intersections}\label{sec:approx-complete-intersection}
In this section, we first motivate and reprove some results of Walsh \cite{Wal20} that we require. Then we prove a novel result that allows one to approximate a variety and its tangent spaces as an approximate complete intersection. The techniques that we will use to study incidence problems in high-dimensional space require us to study these incidence problems inside varieties of high-dimension and high degree. That is why the results of this section -- specifically their good quantitative dependence on degree -- are so necessary.

The nicest algebraic varieties are \emph{complete intersections}: a variety $V$ of codimension $r$ in $\CC^n$ is a complete intersection if there exist polynomials $f_1,\ldots,f_r$ so that $V=Z(f_1,\ldots,f_r)$ and $\deg V=\deg(f_1)\deg(f_2)\cdots\deg(f_r)$.

Unfortunately, not all irreducible varieties are complete intersections. The most well-known example of such is the twisted cubic curve $\gamma=\{(t,t^2,t^3):t\in\CC\}$. It is easy to verify that $\gamma = Z(z-xy, y-x^2, xz-y^2)$. However, if we choose any two of these polynomials and consider the variety they cut out, it turns out that there will always be another component. For example $Z(z-xy,xz-y^2)=\gamma\cup\ell$, where $\ell$ is the line $y=z=0$.

Despite this, it is still possible to write any irreducible variety as an approximate complete intersection. This means that for a variety $V$ of codimension $r$, there exist polynomials $f_1,\ldots,f_r$ so that $V$ is an irreducible component of $Z(f_1,\ldots,f_r)$ in such a way that $\deg(V)\sim \deg(f_1)\deg(f_2)\cdots\deg(f_r)$ \cite[Theorem 1.3]{Wal20}. However, this alone is not enough for our purposes. We will need to study how flats can cluster in varieties; to do so, we study the tangent space $T_pV$. Since $V\subseteq Z(f_1,\ldots,f_r)$, the tangent space is always contained in the kernel of the Jacobian matrix $J_p(f_1,\ldots,f_r)$ with entries $(\partial f_i/\partial x_j)(p)$. Ideally the Jacobian would have full rank at a generic point of $V$ and so $T_pV$ would coincide with the kernel of the Jacobian. Unfortunately, this is not always the case. For example, suppose $V=Z(f,g)$. Then, as sets, $V=Z(f^2+g,g)$ as well. However, the Jacobian evaluates to
\[\begin{bmatrix}
| & | \\
2f\vec\nabla f+\vec\nabla g &
\vec\nabla g\\
| & |
\end{bmatrix}=\begin{bmatrix}
| & | \\
\vec\nabla g &
\vec\nabla g \\
| & |
\end{bmatrix}\]
on $V$ and thus always has rank at most 1. In particular, the kernel of the Jacobian always has codimension at most 1, while the tangent space $T_pV$ has codimension 2 for a generic point $p\in V$.

Despite such examples, we will show that one can always write $V$ as an approximate complete intersection whose Jacobian has full rank at a generic point of $V$.

\begin{theorem}\label{thm:approx-complete-intersection}
Let $V\subseteq\CC^n$ be a codimension-$r$ irreducible variety.
Then there exist polynomials $f_1,\ldots, f_r\in\CC[x_1,\ldots,x_n]$ such that:
    \begin{enumerate}
        \item $V$ is an irreducible component of $Z(f_1,\ldots, f_r)$;
        \item $\prod_{i=1}^{r}\deg f_i \lesssim_n \deg V$; and
        \item the Jacobian of $f_1,\ldots, f_r$ has full rank (that is, rank $r$) on a Zariski-dense open subset of $V$.
    \end{enumerate}
\end{theorem}

Without part (3), this is \cite[Theorem 1.3]{Wal20}. We reprove this result in \cref{sec:Walsh}. In \cref{sec:Jacobian} we prove (3) for varieties of codimension at most 2. The full proof is more challenging technically so we defer it to a companion paper \cite{TYZ26b}. In this paper we only use \cref{thm:approx-complete-intersection} for varieties of codimension at most 2.

\begin{remark}
\cref{thm:approx-complete-intersection}(1,2) (and all results appearing in \cref{sec:Walsh}) holds over any algebraically closed field. The full strength of \cref{thm:approx-complete-intersection} (and all results appearing in \cref{sec:Jacobian}) holds over any algebraically closed field of characteristic zero.
\end{remark}

\subsection{Parameter counting}
\label{sec:Walsh}
For completeness, we first give a self-contained proof of the results of Walsh which we need \cite{Wal20}.
In this section $\F$ will refer to an arbitrary algebraically closed field.

A key definition for stating these results is the notion of partial degrees. A codimension $r$ variety $V$ is associated with a list of partial degrees $\delta_1(V)\leq\delta_2(V)\leq\cdots\leq\delta_r(V)$. In \cref{thm:approx-complete-intersection}, the polynomials produced will satisfy $\deg f_i\sim_n \delta_i(V)$ so we will also have $\delta_1(V)\delta_2(V)\cdots\delta_r(V)\sim_n \deg(V)$. We now give an intrinsic definition of these quantities.

\subsubsection{Partial degrees}
To prove \cref{thm:approx-complete-intersection}, a natural approach is to find $f_1,\ldots, f_r$ greedily.
More precisely, let $f_1$ be a polynomial with minimal degree that vanishes identically on $V$.
Then $Z(f_1)$ is an irreducible hypersurface, so let  $f_2$ be a polynomial with minimal degree that vanishes identically on $V$, but does not vanish identically on $Z(f_1)$.
The resulting $Z(f_1,f_2)$ is a variety of codimension $2$ containing $V$, but it may no longer be irreducible (the twisted cubic curve is one such example).
We may nonetheless carry on, defining $f_3$ to be a polynomial of minimal degree vanishing identically on $V$ while not vanishing identically on any irreducible component of $Z(f_1,f_2)$ that contains $V$.

Motivated by this idea, the partial degrees of $V$ are defined as follows.

\begin{definition}\label{def:partial-degree}
Let $V$ be an irreducible variety of codimension $r$ in $\F^n$. 
For $1\leq s\leq r$, define the \emph{$s$-th partial degree} of $V$, denoted $\delta_s(V)$, to be the smallest integer such that there exist polynomials $f_1,\ldots, f_s\in \F[x_1,\ldots, x_n]$ with $\deg f_i\leq \delta_s(V)$ for all $i\in[s]$ such that each irreducible component of $Z(f_1,\ldots, f_s)$ which contains $V$ has codimension $s$.
    
We also define $\delta_0(V)=0$ and $\delta_s(V) = \infty$ if $s>r$. We write $\delta(V) = \delta_r(V)$, the \emph{partial degree} of $V$.
\end{definition}

Note that this is not exactly the same definition as \cite[Definition 2.2]{Wal20}, though we will shortly show that the definitions are equivalent. 

The partial degree, $\delta(V)$, is the smallest integer for which there exist polynomials $f_1,\ldots, f_r$ of degree at most $\delta(V)$ so that $V$ is an irreducible component of $Z(f_1,\ldots,f_r)$. Instead one might want to cut out $V$ exactly, possibly by using more polynomials. That is, what is the degree necessary so that there exist polynomials $g_1,\ldots,g_t$ of at most that degree so that $V=Z(g_1,\ldots,g_t)$? We will give an example showing that this quantity can be much larger than $\delta(V)$.

\begin{example}
    Let $V\subset \CC^4$ be the Zariski closure of
    \[\left\{\left(u,v,u^{-2A^2}v^{2A^2-2A+1},u^{-(2A^2+2A+1)}v^{2A^2}\right):u,v\in\CC^{\times}\right\}\] where $A$ is some large positive integer.
    Using $w,x,y,z$ as the coordinates on $\CC^4$, then $I(V)$ is generated, as a $\CC$-vector space, by the polynomials $w^px^qy^rz^s-w^{p'}x^{q'}y^{r'}z^{s'}$ where 
    \[(p-p')-2A^2(r-r')-(2A^2+2A+1)(s-s')=0\]
    and
    \[(q-q')+(2A^2-2A+1)(r-r')+2A^2(s-s')=0.\]
    For example, the polynomials $w^Az^A-x^{A-1}y^{A+1}$ and $w^{A+1}y^A-x^Az^{A-1}$ are both in $I(V)$.

    One can check that $Z(w^Az^A-x^{A-1}y^{A+1},w^{A+1}y^A-x^Az^{A-1})$ has three irreducible components: the variety $V$, the plane $Z(w,x)$, and the plane $Z(y,z)$. Since $V$ is an irreducible component, this shows that $\delta(V)\leq A+1$. To cut out $V$ without including the plane $Z(w,x)$, a polynomial
    of the form $w^px^qy^rz^s-w^{p'}x^{q'}y^{r'}z^{s'}$ needs to satisfy $p=q=0$ or $p'=q'=0$. Without loss of generality, assume that $p'=q'=0$.
    Write $\Delta_r = r-r'$ and $\Delta_s = s-s'$: integers such that $(\Delta_r,\Delta_s)\neq (0,0)$.
    Then we know that
    \[2A^2\Delta_r+(2A^2+2A+1)\Delta_s=p\geq 0\]
    and
    \[(2A^2-2A+1)\Delta_r+2A^2\Delta_s=-q\leq 0.\]
    As a consequence, $(2A^2)/(2A^2+2A+1)\leq \Delta_s/(-\Delta_r)\leq (2A^2-2A+1)/(2A^2)$.
    However, as $(2A^2-2A+1)(2A^2+2A+1) = (2A^2)^2+1$, it follows from basic properties of the Farey sequence that $\Delta_s\geq 2A^2-2A+1$ and $-\Delta_r\geq 2A^2$.
    Thus if the polynomial $w^px^qy^rz^s-w^{p'}x^{q'}y^{r'}z^{s'}$ does not vanish identically on $Z(w,x)$, both monomials in it have degree at least $2A^2-2A+1$.
    Hence, if we were to find a list of polynomials that cut out precisely $V$, then some polynomial has to have degree at least $2A^2-2A+1$, which is far larger than $\delta(V)$.
\end{example}

Because of this example, it is quantitatively prohibitive to cut out varieties exactly. This is why we always only require $V$ to be an irreducible component of $Z(f_1,\ldots,f_r)$. It is also the reason why it is important that our results are in terms of the partial degrees of $V$ instead of simply $\deg(V)$.

It is not immediately obvious from our definition that there exist $f_1,\ldots,f_r$ with $\deg f_i=\delta_i(V)$ so that $V$ is an irreducible component of $Z(f_1,\ldots,f_r)$. We now resolve this issue and also show that \cref{def:partial-degree} agrees with \cite[Definition 2.2]{Wal20}.

\begin{lemma}[{cf. \cite[Lemma 5.3]{Wal20}}]\label{lem:same-partial-degree}
Let $V$ be an irreducible variety of codimension $r$ in $\F^n$. 
For every positive integer $D$, let $I_{\leq D}(V)$ be the vector space of polynomials vanishing on $V$ with degree at most $D$.
Then for every $s\in[r]$, the $s$-th partial degree $\delta_s(V)$ is the smallest integer $D$ such that all irreducible components $Z(I_{\leq D}(V))$ containing $V$ have codimension at least $s$.
Moreover, there exist polynomials $f_1,\ldots, f_r$ such that $\deg f_i = \delta_i(V)$ for each $i\in[r]$ and such that for all $s\in[r]$, all components of $Z(f_1,\ldots, f_s)$ containing $V$ have codimension precisely $s$.
\end{lemma}

The proof of this result uses the following simple linear-algebraic fact.

\begin{lemma}\label{lem:lin-alg}
Let $V_1,\ldots,V_k\subseteq\F^n$ be irreducible varieties. If $g_1,\ldots,g_k\in\F[x_1,\ldots,x_n]$ are such that $g_i$ does not vanish identically on $V_i$, then a generic linear combination of the $g_i$ does not vanish identically on any of $V_1,\ldots,V_k$.
\end{lemma}

\begin{proof}
Pick points $x_i\in V_i$ so that $g_i(x_i)\neq 0$. Define $g=c_1g_1+\cdots+c_kg_k$ for a generic choice of $c\in\F^k$. Let $H_i\subset\F^k$ be the subspace defined by $H_i=\{c\in\F^k:g(x_i)=0\}$. As $g_i(x_i)\neq 0$, this is a proper subspace. As $\F$ is an infinite field, $\F^k$ cannot be written as the union of finitely many proper subspaces. In particular, a generic choice of $c\in\F^k$ does not lie in $H_1\cup\cdots\cup H_k$, proving the desired result.
\end{proof}

\begin{proof}[Proof of \cref{lem:same-partial-degree}]
We construct the polynomials $f_1,\ldots,f_r$ one-by-one and prove the first part by induction. The case of $s=1$ is immediate.

Now suppose $s\geq 2$, the statement holds for $s-1$, and we have constructed $f_1,\ldots, f_{s-1}$ be polynomials satisfying the desired condition.
Let $D$ be the smallest integer such that all components of $Z(I_{\leq D}(V))$ containing $V$ have codimension at least $s$.
Then from the definition of partial degree, we immediately have $\delta_s(V)\geq D$.
Moreover, suppose that $U_1,\ldots, U_k$ are the components of $Z(f_1,\ldots, f_{s-1})$ containing $V$.
By the inductive hypothesis, we know that each of them has codimension $s-1$.

We claim that for each $i\in [k]$, there exists a polynomial $g_i\in I_{\leq D}(V)$ that does not vanish identically on $U_i$. If this were not the case, then $U_i$ would be contained in $Z(I_{\leq D}(V))$, which contradicts our definition of $D$.

Now by \cref{lem:lin-alg}, taking $f_s\in I_{\leq D}(V)$ to be a generic linear combination of the $g_i$, we see that $f_s$ is not identically zero on any of the $U_i$.
As a consequence, all components of $Z(f_1,\ldots, f_s)$ containing $V$ have codimension strictly larger than $s-1$; by Krull's Hauptidealsatz (see, e.g., \cite[Chapter I, Proposition 7.1]{Har77}) they have codimension exactly $s$. This also shows that $\delta_s(V)=D$, completing the induction.
\end{proof}

Note that an immediate corollary of \cref{lem:same-partial-degree} is that $\delta_i(V)\leq \delta_j(V)$ for all $i\leq j$.
This is because each component of $Z(I_{\leq \delta_j(V)}(V))$ that contains $V$ has codimension at least $j\geq i$, showing that $\delta_i(V)\leq \delta_j(V)$ by the minimality of $\delta_i(V)$.

\subsubsection{Bounds on the partial degrees}
For any irreducible variety $V$ of codimension $r$, we have shown that there exist polynomials $f_1,\ldots, f_r$ such that $V$ is an irreducible component of $Z(f_1,\ldots,f_r)$ and such that $\deg f_i=\delta_i(V)$. By B\'ezout's theorem, we have $\deg(V)\leq\deg(f_1)\cdots\deg(f_r)=\delta_1(V)\cdots\delta_r(V)$. To prove \cref{thm:approx-complete-intersection}(2), it suffices to show that $\delta_1(V)\cdots \delta_r(V) = O_n(\deg V)$.

To do this we want to upper bound the partial degrees in terms of $\deg(V)$, i.e., for each $s$ we need to show that there exists a low-degree polynomial $f_s$ vanishing identically on $V$ without vanishing identically on any component of $Z(f_1,\ldots, f_{s-1})$ which contains $V$.

This is accomplished through parameter counting. If $U$ is a component of $Z(f_1,\ldots, f_{s-1})$ which contains $V$, we want to consider the space of functions $U\to\F$ which are the restriction of a polynomial of degree at most $D$. If this $\F$-vector space has large dimension, we can find the desired polynomial $f_s$.

For a variety $U$ defined over a field $\F$, its (affine) \emph{Hilbert function}, denoted $H_U(D)$ is defined to be $\dim_\F\paren{\F_{\leq D}[x_1,\ldots,x_n]/I_{\leq D}(U)}$. That is, the dimension of the space of polynomials of degree at most $D$ quotiented out by those which vanish identically on $U$.

It is a well-known fact that for sufficiently large $D$, we have $H_U(D)=p_U(D)$ where $p_U$ is a polynomial (the Hilbert polynomial) satisfying $p_U(D) = \deg U\binom{D}{n-r}+O(D^{n-r-1})$ (see, e.g., \cite[Chapter I, Section 7]{Har77}).
The threshold above which the identity starts to hold is the \emph{Castelnuovo--Mumford regularity} of $U$ which in the worst case can be on the order of $\deg(U)$. We need the following lower bound on the Hilbert function which applies when $D$ is much smaller -- on the order of $\delta(U)$.

\begin{lemma}[{\cite[Corollaire 3]{CP99}}]\label{lem:Hilb-func-lower-bound}
Given polynomials $f_1,\ldots, f_m\in \FF[x_1,\ldots, x_n]$ with degrees $d_1\leq \cdots\leq d_m$, let $U$ be an irreducible component of $Z(f_1,\ldots, f_m)$ of codimension $r$.
Then the Hilbert function of $U$ satisfies
    \[H_U(D)\geq \binom{D+n-(d_{m-r+1}+\cdots+d_m)}{n-r}\deg U\]
for every positive integer $D>d_{m-r+1}+\cdots+d_m-r$.
\end{lemma}

\begin{remark}
We remark that the original statement was formulated using homogeneous ideals.
As we are working with affine varieties in $\F^n$, homogeneous polynomials are replaced by polynomials, and we also formulate the statement using varieties instead of their ideals.
\end{remark}

Note that as an immediate corollary, for every irreducible variety $V$ of codimension $r$, we see that $H_V(D)\gtrsim_n \deg V \cdot D^{n-r}$ for every $D>2r\delta(V)$.
We can extend this further to smaller $D$'s.

\begin{corollary}[{\cite[Lemma 3.5]{Wal20}}]\label{cor:Hilb-func-lower-bound}
Let $V$ be an irreducible variety of codimension $r$ in $\F^n$.
Then for every $0\leq s\leq r$ and every positive integer $D$ with $2s\delta_s(V)<D<\delta_{s+1}(V)$, we have 
\[H_V(D)\gtrsim_n \frac{\deg V}{\prod_{i=s+1}^{r}\delta_i(V)}\cdot D^{n-s}.\]
\end{corollary}
\begin{proof}
Let $f_1,\ldots,f_r$ be the polynomials given by \cref{lem:same-partial-degree}.
Define $U_1,\ldots, U_k$ to be the irreducible components of $Z(f_1,\ldots, f_s)$ containing $V$. By \cref{lem:same-partial-degree}, these components all have codimension $s$. (\cref{lem:same-partial-degree} only applies for $s\geq 1$, but the $s=0$ case is trivially true.)
    
For every $j\in[k]$, we know that $V$ is an irreducible component of $U_j\cap Z(f_{s+1},\ldots,f_r)$. Thus by B\'ezout's theorem, $\deg V \leq \delta_{s+1}(V)\cdots \delta_r(V)\deg U_j$, implying that $\deg U_j\geq \deg V/\prod_{i=s+1}^{r}\delta_i(V)$.
As a consequence, \cref{lem:Hilb-func-lower-bound} gives
\[H_{U_j}(D)\gtrsim_n \frac{\deg V}{\prod_{i=s+1}^{r}\delta_i(V)}\cdot D^{n-s}\]
for all $D>2s\delta_s(V)$.
    
Given $D$ satisfying $2s\delta_s(V)<D<\delta_{s+1}(V)$, it suffices to show that $H_{U_j}(D) \leq H_V(D)$ for some $j\in[k]$. If this were not the case, then $I_{\leq D}(U_j)\supsetneq I_{\leq D}(V)$ for all $j$. In other words, for each $j\in[k]$ there would exist a polynomial $g_j\in I_{\leq D}(V)$ that does not vanish identically on $U_j$.
By \cref{lem:lin-alg}, letting $g\in I_{\leq D}(V)$ be a generic linear combination of the $g_j$, we have that $g$ does not vanish identically on any of the $U_j$. However, then every irreducible component of $Z(f_1,\ldots,f_s,g)$ which contains $V$ has codimension at least $s+1$, implying that $\delta_{s+1}(V)\leq D$. This contradicts our assumption on $D$, and thus proves the desired result.
\end{proof}

We have the complementary upper bound $H_V(D)\lesssim_n \deg(V)D^{n-r}$ for any variety $V$ of pure codimension $r$ \cite[Proposition 10]{Som97}. Putting those together, we can now do parameter counting to find a polynomial of small degree vanishing on some small variety without vanishing identically on some larger variety.

\begin{lemma}[{\cite[Theorem 4.6]{Wal20}}]\label{lem:Walsh-param-counting}
Let $V$ be an irreducible variety of codimension $r$ in $\F^n$, and let $T$ be a variety of pure dimension $\ell$ with $\ell<n-r$.
Then there exists a polynomial $f$ that vanishes identically on $T$ but does not vanish identically on $V$ with
\[\deg f\lesssim_n\max_{0\leq s\leq r}\left(\frac{\deg T\prod_{i={s+1}}^{r}\delta_i(V)}{\deg V}\right)^{\frac{1}{n-s-\ell}}.\]
\end{lemma}

\begin{proof}
Let $C=C(n)$ be a large constant.
Define $0\leq s_*\leq r$ to be the value of $s$ which maximizes the quantity
\[\left(\frac{C\deg T\prod_{i={s+1}}^{r}\delta_i(V)}{\deg V}\right)^{\frac{1}{n-s-\ell}},\]
and define $D_*$ to be the maximum. If there are multiple maximizers, let $s_*$ be the largest among them.
    
We will first show that $D_*<\delta_{s_*+1}(V)$.
When $s_*=r$ this is clear, since $\delta_{r+1}(V)=\infty$.
When $s_*<r$, we see that
\[\left[\left(\frac{C\deg T\prod_{i={s_*+2}}^{r}\delta_i(V)}{\deg V}\right)^{\frac{1}{n-s_*-1-\ell}}\right]^{n-s_*-1-\ell}\delta_{s_*+1}(V) = D_*^{n-s_*-\ell},\]
and by the maximality of $D_*$, we get that $D_*< \delta_{s_*+1}(V)$.
A similar argument shows that $D_*\geq \delta_{s_*}(V)$.

Now let $D$ be the minimum integer such that $D\geq D_*$ and $2s\delta_s(V)<D<\delta_{s+1}(V)$ for some $0\leq s\leq r$.
Note that this $D$ exists as $\delta_{r+1}(V)=\infty$ and the corresponding $s$ satisfies $s\geq s_*$. By the definition of $D,s$, we have the bound
\[2s'\delta_{s'}(V)+1\geq \delta_{s'+1}(V)\]
for each $s_*\leq s'<s$, implying that $\delta_s(V)\sim_n \delta_{s'}(V)$ for any $s'$ between $s_*\leq s'<s$.
Moreover, we see that $D=\max\{D_*,2s\delta_s(V)+1\}$, so $D\sim_n D_*$ and if $s>s_*$ we also have $D\sim_n\delta_s(V)$.
By \cref{cor:Hilb-func-lower-bound}, we get that
\[H_V(D)\gtrsim_n\frac{\deg V}{\prod_{i=s+1}^r\delta_i(V)}D^{n-s}\sim_n\frac{\deg V}{\prod_{i=s_*+1}^{r}\delta_i(V)}D_*^{s-s_*}D_*^{n-s}  = C\deg T\cdot D_*^{\ell}.\]

On the other hand, standard upper bounds on the Hilbert function \cite[Proposition 10]{Som97} imply $H_T(D)\lesssim_n \deg(T) D^{\ell}\sim_n \deg(T)D_*^{\ell}$.
By choosing $C$ large enough, we can guarantee that $H_V(D)>H_T(D)$.
Thus $I_{\leq D}(V)\subsetneq I_{\leq D}(T)$, meaning that there exists a polynomial $f$ of degree at most $D$ that vanishes identically on $T$ but not on $V$.
By the definition of $D_*$ and the fact that $D\sim_n D_*$, this proves the desired result.
\end{proof}

In the above result the variety $V$ must be irreducible, but we can easily extend it to varieties of pure codimension $r$ where each of the partial degrees of the components are roughly equal.

\begin{corollary}\label{cor:Walsh-param-counting}
Given positive integers $d_1,\ldots, d_r$, let $V$ be a pure variety of codimension $r$ such that for each irreducible component $V_i$ of $V$, we have $\delta_s(V_i)\leq d_s$ for each $s\in [r]$.
Let $T$ be a variety of pure dimension $\ell$ with $\ell <n-r$.
Then there exists a polynomial $f$ that vanishes identically on $T$ but not identically on any irreducible component of $V$ with
\[\deg f\lesssim_{n}\max_{0\leq s\leq r}\left(\frac{\deg T\prod_{i=s+1}^r d_i}{\min_j\deg V_j}\right)^{\frac{1}{n-s-\ell}}.\]
\end{corollary}

\begin{proof}
For every irreducible component $V_i$ of $V$, by \cref{lem:Walsh-param-counting} we may find a polynomial $f_i$ with
\[\deg f_i\lesssim_n \max_{0\leq s\leq r}\paren{\frac{\deg T\prod_{j=s+1}^r \delta_j(V_i)}{\deg V_i}}^{\frac1{n-s-\ell}}\leq \max_{0\leq s\leq r}\paren{\frac{\deg T\prod_{j=s+1}^r d_j}{\deg V_i}}^{\frac1{n-s-\ell}}.\]
that vanishes identically on $T$ but not identically on $V_i$.
By \cref{lem:lin-alg}, we may find $f$, a generic linear combination of the $f_i$, that vanishes identically on $T$ without vanishing identically on any irreducible component of $V$.
\end{proof}

Now we are ready to show that $\delta_1(V)\cdots \delta_r(V)\sim_n \deg V$ for any irreducible variety $V$ of codimension $r$.
To prove this we will induct on a slightly stronger statement.

\begin{lemma}[{\cite[Theorem 5.5]{Wal20}}]\label{lem:partial-deg-of-levels-above}
Given $C\geq 1$, let $V$ be an irreducible variety in $\F^n$ of codimension $r$.
Let $f_1,\ldots, f_r$ be polynomials such that $V$ is an irreducible component of $Z(f_1,\ldots, f_r)$ and $ \deg f_i\leq C\delta_i(V)$ for each $i\in [r]$.
Then for any $s\in[r]$ and any irreducible component $U$ of $Z(f_1,\ldots, f_s)$ containing $V$, we have that $\dim U=n-s$ and $\deg U \sim_{C,n} \delta_1(V)\cdots \delta_s(V)$ and $\delta_i(U)\sim_{C,n} \delta_i(V)$ for each $i\in [s]$.
\end{lemma}
\begin{proof}
We prove the result by induction on $s$. The case $s=1$ is immediate as $\delta_1(V)\leq \delta_1(U)=\deg U\leq \deg f_1\lesssim_{C} \delta_1(V)$.

Now suppose that the statement holds for $s-1$ where $s\geq 2$. Let $U$ be an irreducible component of $Z(f_1,\ldots,f_s)$ containing $V$. Then let $W$ be the union of the irreducible components of of $Z(f_1,\ldots,f_{s-1})$ which contain $V$. By the inductive hypothesis, we see that each irreducible component $W_i$ of $W$ satisfies $\dim W_i=n-s+1$ and $\deg W_i\sim_{C,n} \delta_1(V)\cdots\delta_{s-1}(V)$ and $\delta_j(W_i)\sim_{C,n}\delta_j(V)$ for each $i\in[s-1]$.

Next we apply \cref{cor:Walsh-param-counting} to $U\subset W$, producing a polynomial $g$ such that 
    \[\deg g\lesssim_{C,n} \max_{0\leq s'\leq s-1}\left(\frac{\deg U\prod_{i=s'+1}^{s-1}\delta_i(V)}{\prod_{i=1}^{s-1}\delta_i(V)}\right)^{\frac{1}{s-s'}}=\max_{0\leq s'\leq s-1}\left(\frac{\deg U}{\prod_{i=1}^{s'}\delta_i(V)}\right)^{\frac{1}{s-s'}}\]
which vanishes identically on $U$ without vanishing identically on any irreducible component of $W$.

Now since $W$ was defined to be the union of the irreducible components of of $Z(f_1,\ldots,f_{s-1})$ which contain $V$, we know that each irreducible component of  $Z(f_1,\ldots, f_{s-1},g)$ which contains $V$ has strictly smaller dimension than $W$. By Krull's Hauptidealsatz, each such component has dimension $\dim W-1=n-s$. Furthermore, for each $s'\leq s$, we claim that each irreducible component of $Z(f_1,\ldots,f_{s'-1},g)$ which contains $V$ has dimension $n-s'$. Indeed, if $X$ were some component with $\dim X>n-s'$, then Krull's Hauptidealsatz implies that every component of $X\cap Z(f_{s'},\ldots,f_{s-1})$ has dimension at least $\dim X-(s-s')>n-s$. As $X$ contains $V$, then at least one of these components also contains $V$ and lies in $Z(f_1,\ldots,f_{s-1},g)$, yet has too large dimension.

Our next goal is to show that $\deg g\gtrsim_{C,n}\delta_s(V)$. Let $1\leq s'\leq s$ maximal so that $\delta_{s'}(V)> C\delta_{s'-1}(V)$ holds. (Such an $s'$ exists as $\delta_0(V)=0$.) Note that by maximality, $\delta_s(V)\leq C^{s-s'}\delta_{s'}(V)$. Thus it suffices to show that $\deg g\geq \delta_{s'}(V)$. Suppose for contradiction that $\deg g<\delta_{s'}(V)$. Note that the polynomials $f_1,\ldots,f_{s'-1},g$ all have degree strictly smaller than $\delta_{s'}(V)$, as $\deg f_i\leq C\delta_i(V)\leq C\delta_{s'-1}(V)<\delta_{s'}(V)$ for $i<s'$. We showed that each irreducible component of $Z(f_1,\ldots,f_{s'-1},g)$ which contains $V$ has dimension $n-s'$. This implies that $\delta_{s'}(V)\leq\max\{\deg f_1,\ldots,\deg f_{s'-1},\deg g\}$, a contradiction. Thus we have shown that $\deg g\geq\delta_{s'}(V)\gtrsim_{C,n}\delta_s(V)$.

Putting together the upper and lower bounds on $\deg g$, we see that 
\[\max_{0\leq s'\leq s-1}\left(\frac{\deg U}{\prod_{i=1}^{s'}\delta_i(V)}\right)^{\frac{1}{s-s'}}\gtrsim_{C,n}\deg g\gtrsim_{C,n} \delta_s(V).\]
Rearranging, there exists $0\leq s'\leq s-1$ such that
    \[\deg U\gtrsim_{C,n}\delta_s(V)^{s-s'}\prod_{i=1}^{s'}\delta_i(V)\geq \prod_{i=1}^{s}\delta_i(V).\]
By B\'ezout's theorem, we also have $\deg U\leq \deg f_1\cdots\deg f_s\lesssim_{C,n} \delta_1(V)\cdots \delta_s(V)$.

Now all that remains to show that $\delta_i(U)\sim_{C,n} \delta_i(V)$ for all $i\in[s]$.
First, note that $\delta_i(U)\leq \max_{j\in[i]}\deg f_j\leq C\delta_i(V)$. To see this, we need to check that every irreducible component of $Z(f_1,\ldots,f_i)$ which contains $U$ has codimension $i$. By hypothesis, every irreducible component of $Z(f_1,\ldots,f_i)$ which contains $V$ has codimension $i$; as $V\subseteq U$, this suffices.
Now we have the bounds 
    \[\prod_{j=1}^s\delta_j(U)\geq \deg(U)\gtrsim_{C,n}\prod_{j=1}^s \delta_j(V)\gtrsim_{C,n} \delta_i(V)\prod_{j\in[s]\setminus\{i\}}\delta_j(U).\]
    Canceling gives the desired bound, completing the induction.
\end{proof}

\begin{proof}[Proof of \cref{thm:approx-complete-intersection}(1,2)]
For an irreducible codimension-$r$ variety $V\subseteq\F^n$, let $f_1,\ldots, f_r$ be the polynomials produced by \cref{lem:same-partial-degree}. We have that $V$ is an irreducible component of $Z(f_1,\ldots,f_r)$, proving (1). Furthermore, \cref{lem:partial-deg-of-levels-above} with $s=r$ and $U=V$ implies that $\deg V\sim_n \delta_1(V)\cdots\delta_r(V)$, proving (2).
\end{proof}

\subsection{Making the Jacobian full rank}\label{sec:Jacobian}
To prove \cref{thm:approx-complete-intersection}(3), it remains to produce $f_1,\ldots,f_r$ so that the Jacobian of $f_1,\ldots,f_r$ has full rank on a Zariski-dense open subset of $V$. This is trivial if $\codim V=r$ is equal to 0 and easy if $r=1$. We prove the case of $r=2$ in this subsection.

Let $f_1,f_2$ be the polynomials of degree $\delta_1(V),\delta_2(V)$ which are produced by \cref{lem:same-partial-degree}. It may be the case that the Jacobian of $f_1,f_2$ has full rank (i.e., rank 2). However, we do not know how to prove that, nor can we find a counterexample. Instead, we will show how to modify the polynomials so that their Jacobian has full rank without increasing their degrees too much. First we give a sketch of the proof strategy.

If $Z(f_1,f_2)$ has ``multiplicity one'' at $V$ then it turns out that the Jacobian does not vanish identically on $V$ and we are immediately done. (We will define multiplicity and prove this fact in the following pages.)
Thus we have to deal with the case when $Z(f_1,f_2)$ has multiplicity more than one at $V$.
Our goal is to ``differentiate $f_2$'' in some way so that the multiplicity goes down, and eventually the multiplicity at $V$ becomes one.

For any $i\in [n]$ we use $\partial_i$ to denote the derivative operator $\partial/\partial x_i$.
Now suppose that the Jacobian of $f_1,f_2$ is never full rank on $V$. Without loss of generality, assume that $\partial_1f_1$ is not identically zero on $V$. Then for each $x\in V$, the gradient $\vec{\nabla}f_2(x)$ is a multiple of $\vec{\nabla}f_1(x)$. Then replacing $f_2$ by \[g_2=(\partial_1f_1)f_2-(\partial_1f_2)f_1,\] we see that both $g_2$ and $\vec{\nabla}g_2$ are identically zero on $V$. The first claim is immediate since $f_1,f_2$ vanish identically on $V$. To see the second, note that
\[\partial_i g_2=(\partial_i\partial_1f_1)f_2-(\partial_i\partial_1f_2)f_1+(\partial_1f_1)(\partial_i f_2)-(\partial_1 f_2)(\partial_i f_1)=(\partial_i\partial_1f_1)f_2-(\partial_i\partial_1f_2)f_1,\]
where the last equality follows since $\vec\nabla f_1,\vec\nabla f_2$ are multiples of each other.

The strategy is to replace $f_2$ by one of the components of $\vec{\nabla}g_2$ in such a way that $V$ remains an irreducible component of $Z(f_1,f_2)$ but the multiplicity of $V$ in this intersection drops by one.
If this is always possible, we can thus bound the number of replacements of $f_2$ by the multiplicity of $Z(f_1,f_2)$ at $V$. 
In addition, B\'ezout's theorem implies that $m \deg V\leq\deg f_1\deg f_2$ where $m$ is the multiplicity of $Z(f_1,f_2)$ at $V$. By \cref{lem:partial-deg-of-levels-above}, we have $\deg f_1\deg f_2=\delta_1(V)\delta_2(V)\lesssim \deg V$, showing that $m\lesssim 1$. Thus this procedure halts after a bounded number of steps, producing polynomials $f_1,f_2$ whose Jacobian is of full rank on a dense open subset of $V$. Furthermore, the degrees of these polynomials have doubled only a bounded number of times throughout this procedure.

We now give some standard notation, definitions, and results from commutative algebra that will be necessary for the proof. A Noetherian ring $A$ is \emph{local} if it has a unique maximal ideal. For a multiplicatively closed subset $U\subset A$, we use $A[U^{-1}]$ to denote the \emph{localization} of $A$ at $U$. For a prime ideal $\fp\subset A$, we use $A_{\fp}$ to denote the localization of $A$ at $A\setminus\fp$. For any prime ideal $\fp$, the localization $A_{\fp}$ is a local ring.\footnote{The ideals of the localization $A_{\fp}$ are in one-to-one correspondence with the ideals of $A$ which are contained in $\fp$. Thus the unique maximal ideal of $A_{\fp}$ is $\fp A_{\fp}$, which corresponds to $\fp\subset A$.}

For a variety $V\subseteq\CC^n$, its \emph{ideal} $I(V)\subseteq\CC[x_1,\ldots,x_n]$ is the ring of functions which vanish identically on $V$. The quotient $\CC[x_1,\ldots,x_n]/I(V)$ is the \emph{coordinate ring} of $V$ (which can be viewed as the functions $V\to\CC$ which are the restriction of a polynomial). The localization $\CC[x_1,\ldots,x_n]_{I(V)}$ is the ring of rational functions whose denominator does not vanish identically on $V$. Finally $K(V)=\CC[x_1,\ldots,x_n]_{I(V)}/I(V)\CC[x_1,\ldots,x_n]_{I(V)}$ is the \emph{rational function field} of $V$. It is also the field of fractions of the coordinate ring.

The (Krull) \emph{dimension} of a ring is the length of the longest chain of prime ideals. A local ring $A$ is \emph{regular} if its unique maximal ideal can be generated by $\dim A$ elements. A regular local ring of dimension 1 is known as a \emph{discrete valuation ring} (DVR). Let $R$ be a DVR and let $g$ be a generator of its unique maximal ideal. A standard fact in commutative algebra (see, e.g. \cite[Proposition 11.1]{Eis95}) is that every DVR $R$ has a \emph{valuation} $v_R\colon R\setminus\{0\}\to\ZZ_{\geqslant 0}$ which is defined so that for each $x\in R\setminus\{0\}$, there exists a unit $u$ of $R$ such that $x=ug^{v_R(x)}$.

For $V$ an irreducible component of $Z(f_1,f_2)$, we need to measure the multiplicity of $V$ in the ideal $\sang{f_1,f_2}$, which is given by the following definition. The reader less familiar with commutative algebra may want to skip this definition and instead read \cref{lem:valuation-is-mult} where we give a definition that is equivalent in the codimension 2 setting and is easier to work with.

\begin{definition}
\label{defn:multiplicity}
Let $I\subseteq\CC[x_1,\ldots,x_n]$ be an ideal and let $V$ be an irreducible component of $Z(I)$. Then the \emph{multiplicity} of $I$ at $V$, denoted $m_V(I)$, is the positive integer defined by\footnote{The length of an $R$-module $M$, denoted $\len_R(M)$, is the largest $\ell$ such that there exists a chain of $R$-modules $0=M_0\subsetneq M_1\subsetneq\cdots\subsetneq M_\ell=M$. Since $I\subseteq I(V)$, we see that $I\CC[x_1,\ldots,x_n]_{I(V)}\neq \CC[x_1,\ldots,x_n]_{I(V)}$, so this module is not the zero module and thus has positive length.}
\[\len_{\CC[x_1,\ldots,x_n]_{I(V)}}\paren{\CC[x_1,\ldots,x_n]_{I(V)}/I\CC[x_1,\ldots,x_n]_{I(V)}}.\]
If $I=\sang{f_1,\ldots,f_r}$, we write $m_V(f_1,\ldots,f_r)$ for $m_V(I)$.
\end{definition}

The next result, known as the Jacobian criterion, relates the rank of the Jacobian to properties of the ring $\CC[x_1,\ldots,x_n]_{I(V)}/I\CC[x_1,\ldots,x_n]_{I(V)}$ which appears in the definition of multiplicity. This result is essentially standard; we give a formal derivation from results appearing in the literature in \cref{sec:appendix-jacobian}.

\begin{restatable}[{Jacobian criterion \cite[Theorem 16.19]{Eis95}}]{theorem}{jacobian}
\label{thm:jacobian-criterion}
Let $I\subseteq\CC[x_1,\ldots,x_n]$ be an ideal and let $V$ be an irreducible variety lying in $Z(I)$ so that every irreducible component of $Z(I)$ containing $V$ has codimension $r$ in $\CC^n$. Then $R=\CC[x_1,\ldots,x_n]_{I(V)}/I\CC[x_1,\ldots,x_n]_{I(V)}$ is a local ring of dimension $\codim_{\CC^n}(V)-r$. Furthermore, the Jacobian of $I$ has rank $r$ on a Zariski-dense open subset of $V$ if and only if $R$ is regular. 
\end{restatable}

We next give an easy corollary of B\'ezout's theorem with multiplicity (see, e.g., \cite[Chapter I, Theorem 7.7]{Har77}). The usual version of B\'ezout's theorem is an equality when summed over all components of the intersection in projective space. Restricting to a single component gives the following inequality.

\begin{theorem}[B\'ezout]
\label{thm:bezout-with-multiplicity}
Let $V\subset\CC^n$ be a codimension 2 variety which is an irreducible component of $Z(f_1,f_2)$. Then
\[\deg f_1\deg f_2\geq m_V(f_1,f_2)\deg V.\]
\end{theorem}

For our final preparation, we give a simpler definition of multiplicity which applies in the codimension 2 setting. The proof, which is a simple exercise in commutative algebra, is also deferred to \cref{sec:appendix-jacobian}.

\begin{restatable}{lemma}{valuation}
\label{lem:valuation-is-mult}
Let $V\subset \CC^n$ be an irreducible variety of codimension $2$ and let $f_1$ be a nonzero polynomial vanishing identically on $V$.
Suppose that $R=\CC[x_1,\ldots, x_n]_{I(V)}/\langle f_1\rangle\CC[x_1,\ldots, x_n]_{I(V)}$ is a DVR.  For any $f_2\in\CC[x_1,\ldots, x_n]$ vanishing identically on $V$, let $[f_2]$ be its image in $R$. If $[f_2]\neq 0$, then $m_V(f_1,f_2) = v_R([f_2])$.
\end{restatable}

We are now ready to prove the main tool used in the proof, which shows that differentiating decreases the multiplicity.

\begin{proposition}\label{lem:codim-2-key}
Let $V\subset\CC^n$ be an irreducible variety of codimension 2 which is an irreducible component of $Z(f_1,f_2)$. Suppose that $\deg f_1=\delta_1(V)$.
If $\vec{\nabla}f_2$ is identically zero on $V$, then there exists a polynomial $f_2'\in\CC[x_1,\ldots,x_n]$
such that the following hold:
    \begin{enumerate}
        \item $V$ is an irreducible component of $Z(f_1,f_2')$;
        \item $m_V(f_1,f_2')<m_V(f_1,f_2)$; and
        \item $\deg f_2' < \deg f_2$.
    \end{enumerate}
\end{proposition}

\begin{proof}
As $f_1$ is non-constant, $\vec{\nabla}f_1$ has a nonzero component.
As each component of $\vec\nabla f_1$ has strictly smaller degree than $f_1$, by the definition of $\delta_1(V)=\deg f_1$, we see that $\vec{\nabla}f_1$ does not identically vanish on $V$.

Now we apply the Jacobian criterion to $V\subset Z(f_1)$. Since $\vec\nabla f_1$ does not vanish identically on $V$ (in other words, the Jacobian of $\sang{f_1}$ has rank 1 on a Zariski-dense open subset of $V$), \cref{thm:jacobian-criterion} implies that $\CC[x_1,\ldots,x_n]_{I(V)}/\sang{f_1}\CC[x_1,\ldots,x_n]_{I(V)}$ is a regular local ring of dimension 1.

Call this ring $R$. We just showed that $R$ is a DVR. For any $h\in \CC[x_1,\ldots, x_n]_{I(V)}$, let $[h]$ be its congruence class in $R$.
Choose some $g\in \CC[x_1,\ldots, x_n]$ so that $[g]$ generates the maximal ideal in $R$.\footnote{To see that this is possible, note that the maximal ideal of $R$ is generated by $[g']$ for some $g'\in \CC[x_1,\ldots, x_n]_{I(V)}$; then $g'$ can be written as $g/g''$ for some $g\in\CC[x_1,\ldots,x_n]$ and $g''\in\CC[x_1,\ldots,x_n]\setminus I(V)$. Then $[1/g'']$ is a unit in $R$, so $\sang{[g']}=\sang{[g]}$, as desired.}

With this choice of $g$, we claim that the following holds:
\[R/\sang{[g]}=\CC[x_1,\ldots,x_n]_{I(V)}/\sang{f_1,g}\CC[x_1,\ldots,x_n]_{I(V)}=K(V).\]
To see this, note that the quotient map $\psi\colon \CC[x_1,\ldots,x_n]_{I(V)}\to R$ gives a one-to-one correspondence between ideals $I$ of $R$ and ideals $\psi^{-1}(I)$ of $\CC[x_1,\ldots,x_n]_{I(V)}$ which contain $\sang{f_1}\CC[x_1,\ldots,x_n]$. Since $\sang{[g]}$ is a maximal ideal of $R$, then $\psi^{-1}(I)=\sang{f_1,g}\CC[x_1,\ldots,x_n]_{I(V)}$ is a maximal ideal of $\CC[x_1,\ldots,x_n]_{I(V)}$. The latter ring is local and its unique maximal ideal is $I(V)\CC[x_1,\ldots,x_n]_{I(V)}$. Thus the middle term is $\CC[x_1,\ldots,x_n]_{I(V)}/I(V)\CC[x_1,\ldots,x_n]_{I(V)}$, i.e., $K(V)$.

Now $K(V)$, as it is a field, is a regular local ring of dimension 0. Therefore by another application of the Jacobian criterion (\cref{thm:jacobian-criterion}), we see that the Jacobian of $f_1,g$ has rank 2 on a Zariski-dense open subset of $V$. Explicitly, this means that there is some $2\times2$ minor of the Jacobian whose determinant does not vanish identically on $V$. This is the same as saying that $[\vec{\nabla}f_1]$ and $[\vec{\nabla}g]$ are $K(V)$-linearly independent. To see this formally, suppose that there were $\mu,\nu\in \CC[x_1,\ldots,x_n]_{I(V)}$ so that $[\mu][\vec\nabla f_1]+[\nu][\vec\nabla g]=\vec 0$ in $K(V)^n$; this means that there are $h_i\in I(V)\CC[x_1,\ldots,x_n]_{I(V)}$ so that $\mu\partial_i f_1+\nu\partial_i g=h_i$ for each $i\in[n]$. For this to be a non-trivial linear dependency, at least one of $[\mu],[\nu]\neq 0$. Without loss of generality, suppose $[\mu]\neq 0$. This means that $\mu\not\in I(V)\CC[x_1,\ldots,x_n]_{I(V)}$, so we can write $\partial_i f_1=(h_i-\nu\partial_i g)/\mu$, an equality in $\CC[x_1,\ldots,x_n]_{I(V)}$ for each $i\in[n]$. Then the determinant of the $2\times 2$ minor of the Jacobian consisting of rows $i,j$ is
\[(\partial_i f_1)(\partial_j g)-(\partial_j f_1)(\partial_i g)=\frac{h_i-\nu\partial_i g}\mu(\partial_j g)-\frac{h_j-\nu\partial_j g}\mu(\partial_i g)=\frac{h_i\partial_j g-h_j\partial_i g}{\mu}.\]
As $h_i,h_j\in I(V)\CC[x_1,\ldots,x_n]_{I(V)}$, this function vanishes identically on $V$, a contradiction.\footnote{Explicitly, the left-hand side is a polynomial which we showed does not vanish identically on $V$. The right-hand side is a rational function, which agrees with the left-hand side when evaluated on a Zariski-dense open subset of $V$. As the right-hand side evaluates to 0 on the dense open subset of $V$ where the evaluation is well-defined, this is a contradiction.}

Setting $m=m_V(f_1,f_2)\geq 1$, by \cref{lem:valuation-is-mult} we know that $v_R([f_2])=m$. In other words, $[f_2]=[u][g]^m$ where $u\in \CC[x_1,\ldots, x_n]_{I(V)}$ is such that $[u]$ is a unit in $R$. This means that $f_2-ug^m\in\sang{f_1}\CC[x_1,\ldots,x_n]_{I(V)}$, so we can write $f_2=ug^m+f_1h$ for some $h\in \CC[x_1,\ldots, x_n]_{I(V)}$. Note also that $u\not\in I(V)\CC[x_1,\ldots,x_n]_{I(V)}$.\footnote{Since $[u]$ is a unit in $R$, there must exist some $u'$ so that $[u][u']=1$ in $R$. Thus for some $h'$, we have $uu'=1+f_1h'$ in $\CC[x_1,\ldots,x_n]_{I(V)}$. Now if $u\in I(V)\CC[x_1,\ldots,x_n]_{I(V)}$, then the left-hand side would be identically zero on $V$ while the right-hand side would take the value 1 identically on $V$.}
By the chain rule, we have
\[\vec{\nabla}f_2 = g^m\vec{\nabla}u+mug^{m-1}\vec{\nabla}g+f_1\vec{\nabla}h+h\vec{\nabla}f_1.\]

Suppose for the sake of contradiction that all coordinates of $[\vec{\nabla}f_2]$ lie in $\sang{[g]^m}$.
This implies that all coordinates of $m[u][g]^{m-1}[\vec{\nabla}g]+[h][\vec{\nabla}f_1]$ lie in $\sang{[g]^m}$, so all coordinates of $[h][\vec{\nabla}f_1]$ lie in $\sang{[g]^{m-1}}$.
Now observe there is at least one coordinate of $\vec{\nabla}f_1$ that does not vanish identically on $V$, so the corresponding coordinate of $[\vec{\nabla}f_1]$ is a unit. This implies that $[h]\in \sang{[g]^{m-1}}$, say $[h] = [g]^{m-1}[h']$ for some $h'\in \CC[x_1,\ldots, x_n]_{I(V)}$.
However, now we get that all coordinates of
\[[g]^{m-1}\left(m[u][\vec{\nabla}g]+[h'][\vec{\nabla}f_1]\right)\]
lie in $\sang{[g]^m}$, or equivalently that all coordinates of
\[m[u][\vec{\nabla}g]+[h'][\vec{\nabla}f_1]\]
lie in $\sang{[g]}$.
This is a linear dependency between $[\vec\nabla g]$ and $[\vec\nabla f_1]$ over $R/\sang{[g]}$. Furthermore, since $m[u]\not\in\sang{[g]}$ (and $[\vec\nabla f_1]\neq \vec 0$), it is a non-trivial linear dependency. Since $R/\sang{[g]}=K(V)$, we have shown that $\vec{\nabla}g$ and $\vec{\nabla}f_1$ are $K(V)$-linearly dependent, contradicting our earlier deduction.
This contradiction implies that some coordinate of $[\vec{\nabla}f_2]$ does not lie in $\sang{[g]^m}$.

Let $f_2'$ be a coordinate of $\vec{\nabla}f_2$ with $[f_2']\not\in \sang{[g]^m}$.
In particular, $[f_2']\neq 0$, showing that $f_2'$ is not a multiple of $f_1$. By assumption, $\vec\nabla f_2$ is identically zero on $V$, so $V\subseteq Z(f_1,f_2')$. Furthermore, since $f_2'$ is not a multiple of $f_1$ (and $f_1$ is irreducible) we conclude that $V$ is an irreducible component of $Z(f_1,f_2')$.

Moreover, by \cref{lem:valuation-is-mult}, we know 
\[m_V(f_1,f_2') = v_R(f_2') < m = v_R(f_2) = m_V(f_1,f_2).\]
Lastly, by the definition of $f_2'$ we know that $\deg f_2'< \deg f_2$.
\end{proof}

\begin{remark}
This is the first time in this section that we have used crucially the fact that we are working over $\CC$.
The proof here fails in general if the multiplicity $m$ is divisible by the characteristic of the underlying field.
\end{remark}

Now we can prove the codimension $2$ case of \cref{thm:approx-complete-intersection}(3) by iteratively applying the previous lemma.

\begin{proof}[Proof of \cref{thm:approx-complete-intersection} for $r=2$]
Let $f_1,f_2$ be the polynomials produced by applying \cref{lem:same-partial-degree} to $V$. Thus $V$ is an irreducible component of $Z(f_1,f_2)$ and $\deg f_i=\delta_i(V)$ for $i=1,2$. Furthermore, by the definition of $\delta_1(V)$, any nonzero polynomial $f$ vanishing identically on $V$ has degree at least $\deg f_1$. Without loss of generality, suppose that $\partial_1f_1$ is not the zero polynomial. By the minimality of $\deg f_1$, we see that $\partial_1f_1$ is not identically zero on $V$.

By \cref{thm:bezout-with-multiplicity,lem:partial-deg-of-levels-above}, we have
\[m_V(f_1,f_2)\deg V\leq \deg f_1\deg f_2=\delta_1(V)\delta_2(V)\lesssim_n\deg V,\]
so $m_V(f_1,f_2)\lesssim_n 1$.

We construct a sequence of polynomials $f_2=f_{2,1},f_{2,2},\ldots$ such that:
\begin{enumerate}[(i)]
    \item $V$ is an irreducible component of $Z(f_1,f_{2,k})$; and
    \item $\deg f_{2,k}\leq 2 \deg f_{2,k-1}$.
\end{enumerate}
    
Suppose that the Jacobian of $f_1,f_{2,k}$ does not have rank 2 on a Zariski-dense open subset of $V$. Since rank at most $s$ is a Zariski-closed condition (as it is defined by the vanishing of all $s\times s$ minors), we see that this implies that the Jacobian has rank at most 1 on all of $V$. First note that by defining ${g} = (\partial_1f_1)f_{2,k}-(\partial_1f_{2,k})f_1$, we see that $g$ also satisfies conditions (i,ii). To check (i), note that $\sang{f_1,g}=\sang{f_1,(\partial_1 f_1)f_{2,k}}$, so $Z(f_1,g)=Z(f_1,f_{2,k})\cup Z(f_1,\partial_1 f_1)$. Since $\partial_1 f_1$ does not vanish identically on $V$, we see that $\dim Z(f_1,\partial_1f_1)=n-2$, so $V$ remains an irreducible component of $Z(f_1,g)$. To check (ii), note that $\deg f_1\leq\deg f_{2,k}$ by the minimality of $\deg f_1$ and thus $\deg g\leq \deg f_1+\deg f_{2,k}-1<2\deg f_{2,k}$.

We also claim that $m_V(f_1,g)=m_V(f_1,f_{2,k})$. As $\sang{f_1,g}=\sang{f_1,(\partial_1 f_1)f_{2,k}}$, it suffices to check that 
$\sang{f_1,(\partial_1f_1)f_{2,k}}\CC[x_1,\ldots,x_n]_{I(V)}=\sang{f_1,f_{2,k}}\CC[x_1,\ldots,x_n]_{I(V)}$. However, as $\partial_1 f_1$ does not vanish identically on $V$, it is a unit in $\CC[x_1,\ldots,x_n]_{I(V)}$. Thus the two ideals are equal, as are the two multiplicities.

Next, note that 
\[\partial_i{g} = (\partial_i \partial_1f_1)f_{2,k}-(\partial_i\partial_1f_{2,k})f_1+(\partial_1 f_1)(\partial_i f_{2,k})-(\partial_1 f_{2,k})(\partial_i f_1)= (\partial_i \partial_1f_1)f_{2,k}-(\partial_i\partial_1f_{2,k})f_1,\]
where the latter equality follows since the Jacobian of $f_1,f_{2,k}$ has rank at most 1 on $V$, meaning that $\vec\nabla f_1,\vec\nabla f_{2,k}$ are multiples of each other. In particular, $(\partial_1 f_1)\vec\nabla f_{2,k}-(\partial_1 f_{2,k})\vec\nabla f_1=\vec 0$ on $V$.

Thus we can apply \cref{lem:codim-2-key} to $f_1,g$ to produce a polynomial $f_{2,k+1}$. We see that (i) holds as does (ii), since $\deg f_{2,k+1}<\deg g<2\deg f_{2,k}$. Furthermore, $m_V(f_1,f_{2,k+1})<m_V(f_1,g)=m_V(f_1,f_{2,k})$.

We halt this process when the Jacobian of $f_1,f_{2,k}$ has rank 2 on a Zariski-dense open subset of $V$. This occurs for some $k\leq m_V(f_1,f_2)$ since the fact that $V$ is an irreducible component of $Z(f_1,f_{2,k})$ implies that $m_V(f_1,f_{2,k})\geq 1$. We have produced polynomials $f_1,f_{2,k}$ such that $V$ is an irreducible component of $Z(f_1,f_{2,k})$ with $\deg f_1=\delta_1(V)$ and $\deg f_{2,k}\lesssim_n \delta_2(V)$ and such that the Jacobian of $f_1,f_{2,k}$ has rank 2 on a Zariski-dense open subset of $V$.
\end{proof}

\begin{corollary}
\label{cor:small-deg-sing-locus}
For any irreducible variety $V\subseteq \CC^n$, there exists a polynomial $f$ of degree $\lesssim_n \delta(V)$ such that $f$ does not vanish identically on $V$, but $f$ vanishes at all singular points of $V$.
\end{corollary}

As with \cref{thm:approx-complete-intersection}, in this paper we only prove and use this result for varieties of codimension at most 2. The general version of this result follows from the general version of \cref{thm:approx-complete-intersection}, proved in the companion paper \cite{TYZ26b}.

\begin{proof}
Suppose $V$ has codimension $r$ and let $f_1,\ldots, f_r$ be the polynomials produced by \cref{thm:approx-complete-intersection}.
Then $\deg f_1,\ldots, \deg f_r\lesssim_n \delta(V)$.
We can take $f$ to be one of the $r\times r$-minors of the Jacobian that does not vanish identically on $V$.
This $f$ satisfies the condition.
\end{proof}
\section{Cayley--Monge--Salmon for flats}
 \label{sec:csm-for-flats}

In applications of the polynomial method in incidence geometry, one typically needs to study lines that lie in an algebraic variety. An important tool is the following result discovered both by Monge and by Cayley and Salmon. A variety $V$ is \emph{ruled} (by lines) if for a generic point $p\in V$, there exists a line through $p$ that is contained in $V$. For a local weakening of this property, a point $p\in V$ is called a \emph{flecnode} if there exists a line through $p$ that is tangent to $V$ to third order.

\begin{theorem}[Cayley--Monge--Salmon]
\label{thm:cms-for-lines}
Let $V\subset\CC^3$ be an irreducible surface of degree $d$. There exists a polynomial $g$ of degree $O(d)$ such that the following hold:
\begin{enumerate}
    \item every point $p\in V\cap Z(g)$ is a flecnode;
    \item if $g$ vanishes identically on $V$, then $V$ is ruled.
\end{enumerate}
\end{theorem}

Salmon proved that the polynomial $g$, called the \emph{flecnodal polynomial}, has degree at most $11d-24$. One of many useful corollaries of this result is that an unruled surface in $\CC^3$ contains at most $11d^2-24d$ lines. For more exposition on this classical fact see, e.g., \cite{Kat14,Kol15,Tao14blogpost}.

Various generalizations of \cref{thm:cms-for-lines} are known, primarily for varieties ruled by lines in $\CC^n$ \cite{Lan99,SS17,BDSW20,GZ18}. Our goal is to understand varieties ruled by positive 3-flats. In this section we prove an analogue of the Cayley--Monge--Salmon theorem for this problem. We work in a fairly general setting; this does not increase the difficulty of the proofs but abstracts away the particular formulas used to define the collection of positive 3-flats.

Fix positive integers $n,m,k$ with $1\leq k<n$. We consider a family of $k$-flats in $\CC^n$, whose directions form an $m$-dimensional family. Precisely, fix a polynomial map $\varphi\colon\CC^m\to\CC^{(n-k)\times k}$ from $\CC^m$ to the space of $(n-k)\times k$ matrices. In this section we often view $\CC^n=\CC^k\times\CC^{n-k}$ so for any $x\in \CC^n$ we write $x^{(1)}$ for its first $k$ coordinates and $x^{(2)}$ for its last $n-k$ coordinates.
Now for $v\in \CC^n$ and $\theta\in \CC^m$, let $L_{v,\theta}\subset\CC^n$ be the $k$-flat
\begin{equation}
\label{eq:l-v-theta}
L_{v,\theta}=\left\{v+x:x\in \CC^n, x^{(2)} = \varphi(\theta)x^{(1)}\right\}.
\end{equation}
$L_{v,\theta}$ is a flat through $v$ whose direction is given by $\theta$. In particular, for all $v'\in L_{v,\theta}$, we have $L_{v',\theta}=L_{v,\theta}$.

Given a variety $V\subseteq\CC^n$, our goal is to study the set of pairs $(v,\theta)$ such that $L_{v,\theta}\subseteq V$.

\begin{definition}
Given an irreducible variety $V\subseteq\CC^n$, an irreducible variety $\Gamma\subseteq V\times\CC^m$ \emph{rules} $V$ if the projection $\pi\colon \Gamma\to V$ is dominant\footnote{A map $U\to V$ is \emph{dominant} if its image is Zariski-dense in $V$.} and if $L_{v,\theta}\times\{\theta\}\subseteq\Gamma$ for every $(v,\theta)\in\Gamma$.
\end{definition}

As in the classical Cayley--Monge--Salmon theorem, our goal is to detect a local property of $V$ and show that if this local property holds everywhere, then $V$ is genuinely ruled. Unlike the classical case, we make no attempt to give a geometric interpretation for our local property analogous to the definition of a flecnode.

Write $\vec{D}_{\theta}$ for the differential operator $\vec{\nabla}^{(1)}+\varphi(\theta)^{\intercal}\vec{\nabla}^{(2)}$.
Here $\vec{\nabla}^{(1)}$ is the gradient with respect to $x^{(1)}$ and $\vec{\nabla}^{(2)}$ is the gradient with respect to $x^{(2)}$.
We show that the vanishing of the operators $\vec D_\theta$ detects whether a variety $\Gamma\subseteq V\times\CC^m$ rules $V$.

\begin{proposition}\label{lem:Frobenius}
Given polynomials $g_1,\ldots,g_t\in\CC[x_1,\ldots, x_n, y_1,\ldots, y_m]$, let $\Gamma\subseteq \CC^n\times \CC^m$ be an irreducible component of $Z(g_1,\ldots,g_t)$.
Set \[\Gamma'={\{(v,\theta)\in\Gamma:\forall i\in[t],\,\vec{D}_{\theta}g_i(v,\theta)=0\}}.\]
If $(v,\theta)\in \CC^n \times\CC^m$ satisfies $L_{v,\theta}\times\{\theta\}\subseteq\Gamma$, then $(v,\theta)\in \Gamma'$.
Moreover, if $\Gamma'=\Gamma$ and if the Jacobian of $g_1,\ldots,g_t$ is of rank $\codim \Gamma$ on a Zariski-dense open subset of $\Gamma$, then $\Gamma$ rules $\overline{\pi(\Gamma)}\subseteq\CC^n$ where $\pi\colon\CC^n\times\CC^m\to\CC^n$ is the projection to the first $n$ coordinates.
\end{proposition}

Intuitively, for a variety $V\subseteq\CC^n$, this proposition allows us to cut the space $V\times \CC^m$ iteratively while still keeping all $(v,\theta)\in V\times \CC^m$ with $L_{v,\theta}\subseteq V$.
Moreover, if we cannot cut some irreducible component any further, then that irreducible component rules the variety.
This will be the main tool in the proof of Cayley--Monge--Salmon for flats.

\begin{proof}[Proof of \cref{lem:Frobenius}]
Given $\theta\in\CC^m$, write $\vec u_i(\theta)=(\hat e_i,\varphi(\theta)\hat e_i,0)\in\CC^{k}\times\CC^{n-k}\times\CC^m$ for each $i\in[k]$. These vectors are defined so that $L_{v,\theta}\times\{\theta\}$ is the $k$-flat in $\CC^n\times\CC^m$ through $v$ in directions spanned by $\vec u_1(\theta),\ldots,\vec u_k(\theta)$. In addition, $\vec D_\theta$ is the vector of directional derivatives in directions $\vec u_1(\theta),\ldots,\vec u_k(\theta)$.

Let $(v,\theta)$ be such that $L_{v,\theta}\times\{\theta\}\subseteq \Gamma$. Now $g_i$ is identically zero on $L_{v,\theta}\times\{\theta\}\subseteq\Gamma\subseteq Z(g_1,\ldots,g_t)$, implying that $\vec{D}_{\theta}g_i(v,\theta)=0$, since $\vec D_\theta$ is a vector of directional derivatives in directions contained in $L_{v,\theta}\times\{\theta\}$. As a consequence, we see that $(v,\theta)\in \Gamma'$, as desired.

Now suppose that $\Gamma'=\Gamma$ and the Jacobian of $I$ is of rank $\codim \Gamma$ on a Zariski-dense open subset of $\Gamma$. We will show that $\Gamma$ rules $\overline{\pi(\Gamma)}$. Trivially, $\pi\colon\Gamma\to\overline{\pi(\Gamma)}$ is dominant so all that remains to show is that $L_{v,\theta}\times\{\theta\}\subseteq \Gamma$ for each $(v,\theta)\in\Gamma$. Fix some $(v,\theta)\in\Gamma$ so that the Jacobian of $I$ has rank $\codim \Gamma$ at $(v,\theta)$. We claim that $L_{v,\theta}\times\{\theta\}\subseteq \Gamma$. This will suffice to prove that $\Gamma$ rules $V$ since the condition $L_{v,\theta}\times\{\theta\}\subseteq \Gamma$ holds for a Zariski-closed subset of $\Gamma$, yet is true for every $(v,\theta)$ in the Zariski-dense open subset where the Jacobian of $I$ has full rank.

For each nonzero $\vec c\in\CC^k$ set $\vec x(\theta)=c_1\vec u_1(\theta)+\cdots+c_k\vec u_k(\theta)$. We will prove that $\Gamma$ contains a short (real) line segment through $(v,\theta)$ in direction $\vec x(\theta)$. Since $\Gamma$ is an algebraic variety, this will imply that it contains the entire (complex) line in this direction. As this is true for all $\vec c\in\CC^k\setminus\{0\}$, the variety $\Gamma$ thus contains all of $L_{v,\theta}\times\{\theta\}$.

First note that since $\Gamma\subseteq Z(g_1,\ldots,g_t)$ and the Jacobian of $g_1,\ldots,g_t$ has rank $\codim\Gamma$ at $(v,\theta)$, then $T_{(v,\theta)}\Gamma\subseteq T_{(v,\theta)} Z(g_1,\ldots,g_t)$ has codimension at least $\codim\Gamma$. But this means that $\codim T_{(v,\theta)}\Gamma=\codim\Gamma$ and so $(v,\theta)$ is a smooth point of $\Gamma$. Furthermore, we see that the tangent space $T_{(v,\theta)}\Gamma$ is equal to the kernel of the Jacobian of $g_1,\ldots,g_t$. Now in a (Euclidean) neighborhood of $(v,\theta)$, the variety $\Gamma$ is a smooth manifold. Consider the smooth vector field on $\Gamma$ that assigns to each point $(v,\theta)$ the vector $\vec x(\theta)$. This is a polynomial map, so clearly it is smooth. To see that this is actually a vector field on $\Gamma$, note that
\[\begin{split}
T_{(v,\theta)}\Gamma=&\left\{\paren{\vec x^{(1)},\vec x^{(2)},\vec x^{(3)}}\in\CC^{k}\times\CC^{n-k}\times\CC^{m}:\forall i\in[t],\right.\\&\qquad\left.\vec x^{(1)}\cdot \vec\nabla^{(1)}g_i(v,\theta)+\vec x^{(2)}\cdot \vec\nabla^{(2)}g_i(v,\theta)+\vec x^{(3)}\cdot \vec\nabla^{(3)}g_i(v,\theta)=0\right\},
\end{split}\]
where $\vec \nabla^{(1)},\vec \nabla^{(2)},\vec \nabla^{(3)}$ refer to the gradient with respect to the first $k$, the next $n-k$, and the last $m$ variables, respectively. Now we have $\vec x(\theta)\in T_{(v,\theta)}\Gamma$ since
\begin{align*}
\vec x(\theta)^{(1)}\cdot \vec\nabla^{(1)}g_i(v,\theta)+&\vec x(\theta)^{(2)}\cdot \vec\nabla^{(2)}g_i(v,\theta)+\vec x(\theta)^{(3)}\cdot \vec\nabla^{(3)}g_i(v,\theta)\\
&=\vec c\cdot \vec\nabla^{(1)}g_i(v,\theta)+(\varphi(\theta)\vec c)\cdot \vec\nabla^{(2)}g_i(v,\theta)\\
&=\vec c\cdot \vec\cD_\theta g_i(v,\theta)=0.
\end{align*}

Viewing $\Gamma$ as a real manifold of twice the dimension and $\vec x(\theta)$ as a smooth vector field on this manifold, consider an integral curve through $(v,\theta)$. This is a map $\gamma\colon(-\epsilon,\epsilon)\to\Gamma$ so that $\gamma(0)=(v,\theta)$ and $\gamma'(s)=\vec x(\gamma(s)^{(3)})$ for all $s\in(-\epsilon,\epsilon)$. (As before, $\gamma(s)^{(3)}$ refers to the last $m$ coordinates of $\gamma(s)\in\CC^{k}\times\CC^{n-k}\times\CC^m$.) We claim that this integral curve must be a line segment in direction $\vec x(\theta)$. To see this, we compute using the chain rule
\[\gamma_i''(s)=\gamma'(s)^{(3)}\cdot \vec\nabla^{(3)} x_i(\gamma(s)^{(3)})=0\cdot \vec\nabla^{(3)} x_i(\gamma(s)^{(3)})=0\]
for all $i\in[n]$. Thus $\gamma'(s)=\gamma'(0)=\vec x(\theta)$, as desired.
\end{proof}

\begin{remark}
Alternatively, the result can be proved via an application of the Frobenius theorem in differential geometry (see, e.g., \cite[Theorem 19.12]{Lee13}). Consider the distribution on $\Gamma$ which assigns to each point $(v,\theta)\in\Gamma$ the $k$-dimensional subspace $\spn\{\vec u_1(\theta),\ldots,\vec u_k(\theta)\}\subseteq T_{(v,\theta)}\Gamma$. This distribution is easily seen to be involutive, so the Frobenius theorem implies that it is integrable. In other words, there exists an integral manifold for this distribution that passes through $(v,\theta)$ and is contained in $\Gamma$. This integral manifold (locally) must agree with $L_{v,\theta}\times\{\theta\}$, giving an alternative proof that $L_{v,\theta}\times\{\theta\}\subseteq\Gamma$.
\end{remark}

We are now ready to state the Cayley--Monge--Salmon theorem for flats. Comparing this result to \cref{thm:cms-for-lines}, when applied to an unruled surface $V$, \cref{thm:general-CSM} produces a polynomial $g$ that -- like the flecnodal polynomial -- detects the lines lying in $V$. When applied to a ruled surface $V$, the polynomial $g$ -- unlike the flecnodal polynomial -- does not vanish identically on $V$, but rather detects the ``exceptional lines'' in $V$: those which are not part of the ruling. Furthermore, it is crucial that we have the bound $\deg g\lesssim \delta(V)$, instead of the weaker $\deg g\lesssim \deg V$. 

\begin{theorem}[{Cayley--Monge--Salmon for flats}]
\label{thm:general-CSM}
Let $V$ be an irreducible variety in $\CC^n$.
There exist $O_{\varphi}(1)$ irreducible varieties $\Gamma_1,\ldots, \Gamma_{\ell}$ in $V\times \CC^m$ and a polynomial $g\in\CC[x_1,\ldots,x_n]$ of degree $ O_{\varphi}(\delta(V))$ satisfying the following:
    \begin{enumerate}
        \item each $\Gamma_i$ rules $V$;
        \item the polynomial $g$ does not vanish identically on $V$; and
        \item if $L_{v,\theta}\subseteq V$ for some $(v,\theta)\in V\times \CC^m$, then either $(v,\theta)\in\Gamma_1\cup\cdots\cup\Gamma_\ell$ or $v\in Z(g)$; and
        \item for each $i$, the fiber of $\Gamma_i\to V$ over $x$ has degree $O_{\varphi}(1)$ for a generic choice of $x\in V$.
    \end{enumerate}
\end{theorem}

This result is proved using \cref{thm:approx-complete-intersection}. In this paper we will only apply \cref{thm:general-CSM} when $\codim V\leq 2$; the full version of \cref{thm:general-CSM} follows from the full version of \cref{thm:approx-complete-intersection}, proved in the companion paper \cite{TYZ26b}.

\begin{remark}
\cref{thm:general-CSM} has a number of immediate corollaries. For example, it implies that there exists a constant $C$ depending only on $\varphi$ so that if an irreducible variety $V$ is $C$-ruled then it is infinitely ruled. This generalizes the classical fact that if a surface in $\CC^3$ is triply ruled then it is infinitely ruled (and thus a plane).
\end{remark}

We first outline the proof strategy. Say that an irreducible variety $\Gamma\subseteq V\times\CC^m$ is \emph{transverse} if the projection $\pi\colon \Gamma\to V$ is dominant. Suppose $\codim V=r$ and let $f_1,\ldots, f_r$ be the polynomials given by \cref{thm:approx-complete-intersection} applied to $V$. In other words, $V$ is an irreducible component $Z(f_1,\ldots,f_r)$ and the Jacobian of $f_1,\ldots,f_r$ has rank $r$ on a Zariski-dense open subset of $V$.

To prove \cref{thm:general-CSM}, we will construct a (not necessarily irreducible) variety $\Gamma\subseteq V\times\CC^m$, a polynomial $g\in \CC[x_1,\ldots,x_n]$, and a list $\cL$ of polynomials $g_1,\ldots,g_t\in\CC[x_1,\ldots,x_n,y_1,\ldots,y_m]$. We will iteratively refine $\Gamma,g, \cL$ maintaining the following properties:
    \begin{enumerate}[(i)]
        \item $\Gamma$ is the union of transverse components of $(V\times \CC^m)\cap Z(\cL)$;
        \item the polynomial $g$ does not vanish identically on $V$;
        \item if $L_{v,\theta}\subseteq V$ for some $(v,\theta)\in V\times \CC^m$, then either $(v,\theta)\in \Gamma$ or $v\in Z(g)$.
    \end{enumerate}
As we update this data, we will ensure that $g$ has degree $O_{\varphi}(\delta(V))$ and that $g_1,\ldots,g_t$ have degree in $x_1,\ldots,x_n$ bounded by $O_\varphi(\delta(V))$ and degree in $y_1,\ldots,y_m$ bounded by $O_\varphi(1)$. We will continue to update $\Gamma,g,\cL$ until each component of $\Gamma$ rules $V$. At this point, we will show that the $\Gamma$ and $g$ we have produced satisfy the desired properties.

The iteration works as follows. Viewing $f_1,\ldots,f_r\in\CC[x_1,\ldots,x_n]\subset\CC[x_1,\ldots,x_n,y_1,\ldots,y_m]$, we see $\Gamma\subseteq Z(f_1,\ldots,f_r,\cL)\subseteq\CC^n\times\CC^m$. We would like to apply \cref{lem:Frobenius} to cut down the components of $\Gamma$ which do not rule $V$. To do this, we first must ensure that the Jacobian of $f_1,\ldots,f_r,\cL$ has full rank on a Zariski-dense open subset of each irreducible component $\Gamma_i$ of $\Gamma$.

It is not hard to show that for any ideal $I$, every component of $Z(I)$ has multiplicity 1 in the radical ideal $\sqrt{I}$ and the Jacobian of $\sqrt{I}$ has full rank. Thus to make sure that the Jacobian has full rank, we could replace $f_1,\ldots,f_r,\cL$ by a new list of polynomials $\tilde\cL$ so that $\sqrt{\langle f_1,\ldots,f_r,\cL\rangle}=\langle\tilde\cL\rangle$. However, this is quantitatively too wasteful: the polynomials in $\tilde\cL$ may have degree double exponential in the degree of the original polynomials, which have degree as large as $\delta(V)$. 

To fix this, note that all the original polynomials have bounded degree in $y_1,\ldots,y_m$. Viewing these polynomials as elements of $K(V)[y_1,\ldots,y_m]$ (polynomials in $y_1,\ldots,y_m$ whose coefficients lie in the field of rational functions on $V$), we can replace $\langle\cL\rangle$ by its radical ideal in a quantitatively efficient way.\footnote{In algebraic geometry, $\langle\cL\rangle\subseteq K(V)[y_1,\ldots,y_m]$ is the ideal of the \emph{generic fiber} of $\pi\colon\Gamma\to V$.} 
A short calculation shows that this radical ideal has multiplicity 1 at each of its transverse components and so the Jacobian has full rank. We will then be in a position to apply \cref{lem:Frobenius} to $f_1,\ldots,f_r,\tilde\cL$, producing a new list of polynomials.

Call this new list $f_1,\ldots,f_r,\cL^{\new}$ and define $\Gamma^{\new}$ to be the union of the transverse components of $(V\times\CC^m)\cap Z(\cL^{\new})$. Any component $\Gamma_i$ of $\Gamma$ that ruled $V$ is still a component of $\Gamma^{\new}$. However a component which did not rule $V$ is cut by $\cL^{\new}$. If it is cut down to a subvariety which is no longer transverse, it no longer lies in $\Gamma^{\new}$. Thus to maintain property (iii), we have to deal with any $L_{v,\theta}\subseteq V$ where $(v,\theta)\in\Gamma_i$ for some $\Gamma_i$ that is cut by $\cL^{\new}$ into a non-transverse subvariety. We do this by setting $g^{\new}=gh$ for some polynomial $h$ which vanishes at such $v$. We show that we can find some $h\in\CC[x_1,\ldots,x_n]$ that does not vanish identically on $V$, yet vanishes identically on all non-transverse components of $(V\times\CC^m)\cap Z(\cL^{\new})$ and has degree $O_{\varphi}(\delta(V))$.

With these definitions of $\Gamma^{\new},g^{\new},\cL^{\new}$ we repeat until the process terminates. At each step, the components of $\Gamma$ which rule $V$ are unchanged, while those which do not are replaced by lower-dimensional subvarieties. Thus the process terminates after a bounded number of steps at which point each component of $\Gamma$ rules $V$.

Two key steps in this proof are producing $\tilde\cL$ so that $\sqrt{\sang{\cL}}=\sang{\tilde\cL}$ in $K(V)[y_1,\ldots,y_m]$ and computing the polynomial $h$ both in quantitatively efficient ways. Using the theory of Gr\"obner bases, there are known algorithms for doing these computations. By analyzing the algorithms, one can easily bound the degree of the polynomials produced. So as to not interrupt the flow of the proof, we defer the preliminaries on Gr\"obner bases as well as the proofs of these two results (\cref{claim:radical-bounded,claim:grobner-saturation} below) to \cref{sec:appendix-grobner}. 

\begin{proof}[Proof of \cref{thm:general-CSM}]
Let $f_1,\ldots, f_r$ be the polynomials given by \cref{thm:approx-complete-intersection} applied to $V$ where $r=\codim V$. Set $\Gamma=V\times \CC^m$ and $g=1$ and $\cL=\emptyset$. It is clear that properties (i-iii) are met.

We do the following procedure. Suppose that we have produced $\Gamma,g,\cL$ satisfying properties (i-iii). We update these quantities in the following way.

Write $S=\paren{\CC[x_1,\ldots,x_n]/I(V)}\setminus\{0\}$ and consider
\begin{align*}
\CC[x_1,\ldots,x_n,y_1,\ldots, y_m]
&\overset{\psi}{\twoheadrightarrow}\CC[x_1,\ldots,x_n,y_1,\ldots, y_m]/I(V)\\
&\hookrightarrow (\CC[x_1,\ldots,x_n,y_1,\ldots, y_m]/I(V))[S^{-1}]=K(V)[y_1,\ldots,y_m].
\end{align*}
In words, first we quotient by $I(V)$, then we localize at $S$.
For an element $h\in \CC[x_1,\ldots,x_n,y_1,\ldots, y_m]$, write $[h]\in K(V)[y_1,\ldots,y_w]$ for its image under this map and similarly $[\cL]$ for a list of polynomials $\cL\subset\CC[x_1,\ldots,x_n,y_1,\ldots,y_m]$.

\begin{restatable}{claim}{grobnerRadical}
\label{claim:radical-bounded}
Fix $d,D\geq 1$ and $\cL\subset\CC[x_1,\ldots,x_n,y_1,\ldots,y_m]$, an arbitrary finite list of polynomials such that every $g\in\cL$ has degree in $x_1,\ldots,x_n$ bounded by $D$ and its degree in $y_1,\ldots,y_m$ bounded by $d$. Then there exists a list of polynomials $\tilde\cL\subset\CC[x_1,\ldots,x_n,y_1,\ldots,y_m]$ such that $\sqrt{\sang{[\cL]}}=\sang{[\tilde\cL]}$ and such that for each $\tilde g\in\tilde\cL$, the degree of $\tilde g$ in $x_1,\ldots,x_n$ is bounded by $O_{d,m}(D)$ and its degree in $y_1,\ldots,y_m$ is bounded by $O_{d,m}(1)$. 
\end{restatable}

The proof of this claim is found in \cref{sec:appendix-grobner}. Let $\tilde\cL$ be the list of polynomials produced by applying \cref{claim:radical-bounded} to $\cL$. Let $f_{r+1},\ldots,f_s\in\CC[x_1,\ldots,x_n]$ be polynomials so that $I(V)=\sang{f_1,\ldots,f_s}$. Recall that $\Gamma$ is the union of the transverse components of $(V\times\CC^m)\cap Z(\cL)=Z(f_1,\ldots,f_s,\cL)$ by (i). We claim that $\Gamma\subseteq Z(f_1,\ldots,f_s,\tilde\cL)$.\footnote{In fact, it turns out that the transverse components of $Z(f_1,\ldots,f_s,\cL)$ and $Z(f_1,\ldots,f_s,\tilde\cL)$ agree.}

For each $f\in \sang{f_1,\ldots,f_s,\tilde\cL}$, we have $[f]\in \sang{[f_1],\ldots,[f_s],[\tilde\cL]}=\sang{[\tilde\cL]}=\sqrt{\sang{[\cL]}}$. Thus $[f]^k\in\sang{[\cL]}$ for some $k\geq 1$. Clearing denominators, this means $\psi(f^kh)\in\psi(\sang{\cL})$ for some $h$ with $\psi(h)\in S$. In other words, $f^kh\in \psi^{-1}(\psi(\sang{\cL}))=\sang{I(V),\cL}$ and $h\in\psi^{-1}(S)=\CC[x_1,\ldots,x_n]\setminus I(V)$. Now since
\[\Gamma\subseteq (V\times\CC^m)\cap Z(\cL)=Z(I(V),\cL),\]
we know that $f^kh$ vanishes identically on $\Gamma$. Furthermore, $h$ is a function only of $x_1,\ldots,x_n$ which does not vanish identically on $V$. As every component of $\Gamma$ is transverse (meaning the projection onto the $x_1,\ldots,x_n$ coordinates is a Zariski-dense subset of $V$), we see that $h$ does not vanish identically on any irreducible component of $\Gamma$. Thus we conclude that $f$ vanishes identically on $\Gamma$, as desired.

\begin{lemma}
\label{lem:radical-jacobian}
For an irreducible variety $V\subseteq\CC^n$ and an ideal $J\subseteq \CC[x_1,\ldots,x_n,y_1,\ldots,y_m]$, write $I=\sang{I(V),J}$. Suppose that $[J]$ is a radical ideal in $K(V)[y_1,\ldots,y_m]$. Then if $X$ is a transverse irreducible component of $Z(I)$, the Jacobian of $I$ has full rank on a Zariski-dense open subset of $X$.
\end{lemma}

\begin{proof}
We have $I\subseteq I(X)$, so
\[I\CC[x_1,\ldots,x_n,y_1,\ldots,y_m]_{I(X)}\subseteq I(X)\CC[x_1,\ldots,x_n,y_1,\ldots,y_m]_{I(X)}.\]
Our goal is to show that this is an equality.

Fix $f\in I(X)$. Let $g$ be a polynomial vanishing identically on $Z(I)\setminus X$ without vanishing identically on $X$. As the union of the irreducible components of $Z(I)$ other than $X$ is a Zariski-closed set which does not contain $X$, such a $g$ exists.

The product $fg$ vanishes identically on $Z(I)$, so by Hilbert's Nullstellensatz we have $fg\in I(Z(I))=\sqrt{I}$. Thus $(fg)^k\in I$ for some $k\geq 1$.

Now we have $[fg]^k\in [I]=[J]$, which is a radical ideal by hypothesis. Thus $[fg]\in [J]$. This implies that there exists some $h\in \CC[x_1,\ldots,x_n]\setminus I(V)$ so that $fgh\in I$. Indeed, every element of $[J]$ can be written as $a/b$ for some $a\in \psi(J)$ and $b\in S$. Thus $[fg]\in [J]$ implies that $\psi(fgh)\in \psi(J)$ for some $h$ with $\psi(h)\in S$. This gives the claimed $fgh\in \psi^{-1}(\psi(J))=\sang{I(V),J}=I$ and $h\in \psi^{-1}(S)=\CC[x_1,\ldots,x_n]\setminus I(V)$. 

By construction, neither $g$ nor $h$ vanishes identically on $X$. Thus $gh\not\in I(X)$, implying that \[\frac f1=\frac {fgh}{gh} \in I\CC[x_1,\ldots,x_n,y_1,\ldots,y_m]_{I(X)},\]
as desired.

Now the equality gives
\[\CC[x_1,\ldots,x_n,y_1,\ldots,y_m]_{I(X)}/I\CC[x_1,\ldots,x_n,y_1,\ldots,y_m]_{I(X)}= K(X),\]
the field of rational functions on $X$.

As every field is a regular local ring of dimension 0, the Jacobian criterion (\cref{thm:jacobian-criterion}) shows that the Jacobian of $I$ has full rank on a Zariski-dense open subset of $X$.\footnote{As every field has length 1, this also implies that $m_X(I)=1$, though we will not use this fact.}
\end{proof}

We apply \cref{lem:radical-jacobian} with $J=\sang{\tilde\cL}$ for each irreducible component $\Gamma_i$ of $\Gamma$, showing that the Jacobian of $\sang{I(V),J}$, i.e., the Jacobian of $f_1,\ldots,f_s,\tilde\cL$, is full rank on a Zariski-dense open subset of $\Gamma$. We next claim that the same is true for the Jacobian of $f_1,\ldots,f_r,\tilde\cL$.

The Jacobian of $f_1,\ldots,f_s,\tilde\cL$  is the following $(n+m)\times(s+|\tilde\cL|)$ block matrix
\[\begin{pmatrix}
\paren{\frac{\partial f_j}{\partial x_i}}_{i\in[n],j\in[s]} & \paren{\frac{\partial \tilde g}{\partial x_i}}_{i\in[n],\tilde g\in\tilde\cL} \\
0 & \paren{\frac{\partial \tilde g}{\partial y_i}}_{i\in[m],\tilde g\in\tilde\cL}
\end{pmatrix}.\]
Note that the upper-left block is the Jacobian of $f_1,\ldots,f_s\in\CC[x_1,\ldots,x_n]$. As $\sang{f_1,\ldots,f_s}=I(V)$, this Jacobian has rank at most $n-\dim V$ on $V$ (in fact this is an equality almost everywhere). As $f_1,\ldots, f_r$ were produced by \cref{thm:approx-complete-intersection}, the Jacobian of $f_1,\ldots,f_r$ has rank equal to $n-\dim V$ for a Zariski-dense open subset of $V$. Thus we see that removing the polynomials $f_{r+1},\ldots,f_s$ does not decrease the rank of the Jacobian matrix on a Zariski-dense open subset of $\Gamma_i$.

We are now in a position to apply \cref{lem:Frobenius} to the list of polynomials $f_1,\ldots,f_r,\tilde\cL$. Let $\Gamma_i$ be an irreducible component of $\Gamma$, in other words, a transverse irreducible component of $(V\times\CC^m)\cap Z(\cL)$. This means that $\Gamma_i\subseteq Z(f_1,\ldots,f_s,\cL)=Z(f_1,\ldots,f_s,\tilde \cL)\subseteq Z(f_1,\ldots,f_r,\tilde\cL)$. It is not immediately apparent that $\Gamma_i$ is an irreducible component of the final variety. However, this follows from the fact that the Jacobian of $f_1,\ldots,f_r,\tilde\cL$ has full rank on a Zariski-dense open subset of $\Gamma_i$.

To see this, note that the rank of the Jacobian of $f_1,\ldots,f_r,\tilde\cL$ at $p$ is at most the codimension of the tangent space $T_pZ(f_1,\ldots,f_r,\tilde\cL)$. However, on a Zariski-dense open subset of $\Gamma_i$ this rank is equal to the codimension of $\Gamma_i$. Thus the tangent spaces of $\Gamma_i$ and $Z(f_1,\ldots,f_r,\tilde\cL)$ agree on a Zariski-dense subset of $\Gamma_i$, meaning that $\Gamma_i$ must actually be an irreducible component of $Z(f_1,\ldots,f_r,\tilde\cL)$.

Thus \cref{lem:Frobenius} applies to each component $\Gamma_i$ of $\Gamma$. This proposition produces a list of polynomials which are the coordinates of $\vec D_y f(x,y)$ for $f\in \{f_1,\ldots,f_r\}\cup\tilde\cL$. Define $\cL^{\new}$ to be the union of $\tilde\cL$ and all these polynomials produced. 

The first part of \cref{lem:Frobenius}, when combined with condition (iii), implies that if $L_{v,\theta}\subset V$, then either $v\in Z(g)$ or $(v,\theta)\in\Gamma\cap Z(\cL^{\new})$. Indeed, as $L_{v,\theta}=L_{v',\theta}$ for all $v'\in L_{v,\theta}$, we have $(v',\theta)\in\Gamma$ for all $v'\in L_{v,\theta}\setminus Z(g)$. In particular, either $L_{v,\theta}\subseteq Z(g)$, implying $v\in Z(g)$ or $L_{v,\theta}\times\{\theta\}\subseteq\Gamma_i$ for some irreducible component $\Gamma_i$ of $\Gamma$. In the latter case, \cref{lem:Frobenius} implies that $(v,\theta)\in Z(\cL^{\new})$.

We define $\Gamma^{\new}$ to be the union of the transverse components of $(V\times\CC^m)\cap Z(\cL^{\new})$ so that property (i) is maintained. The second part of \cref{lem:Frobenius} implies that if a component $\Gamma_i$ of $\Gamma$ lies in $Z(\cL^{\new})$ -- equivalently if $\Gamma_i$ remains a component of $\Gamma^{\new}$ -- then $\Gamma_i$ rules $V$. To maintain property (iii), we need to define $g^{\new}=gh$ in such a way that if $(v,\theta)$ lies in a component of $\Gamma\cap Z(\cL^{\new})$ which is no longer transverse, then $h(v)=0$.

Write $W^{\new}=(V\times\CC^m)\cap Z(\cL^{\new})$. We will compute a polynomial $h\in\CC[x_1,\ldots,x_n]\setminus I(V)$ so that\footnote{For an ideal $I$ and an element $h$ of a ring $R$ (here taken to be $\CC[x_1,\ldots,x_n,y_1,\ldots,y_m]/I(V)$) the \emph{saturation} is defined to be $(I:h^\infty)=\{x\in R:\exists k\, xh^k\in I\}$.}
\begin{equation}
\label{eq:saturation}
\psi(\sang{\cL^{\new}})[S^{-1}]\cap\CC[x_1,\ldots,x_n,y_1,\ldots,y_m]/I(V)\subseteq \paren{\psi(\sang{\cL^{\new}}):\psi(h)^\infty}.
\end{equation}
In other words, consider the ideal $\psi(\sang{\cL^{\new}})$ in $\CC[x_1,\ldots,x_n,y_1,\ldots,y_m]/I(V)$. For each element $f\in \psi(\sang{\cL^{\new}})$, clearly $f/1$ lies in the localized ideal $\psi(\sang{\cL^{\new}})[S^{-1}]$. However, there may be some $f\not\in\psi(\sang{\cL^{\new}})$ such that $f/1$ still lies in the localized ideal; in particular, this occurs when $f/1=fg/g$ for some $g\in S$ such that $fg\in\psi(\sang{\cL^{\new}})$. We want to find a single polynomial $\psi(h)\in S$ whose powers detects when this happens. In other words, whenever $f/1$ lies in the localized ideal, it is because $f\psi(h)^k/\psi(h)^k=f/1$ for some $k\geq 1$ where $f\psi(h)^k\in \psi(\sang{\cL^{\new}})$.

Such an $h$ can be found by a Gr\"obner basis computation on $\sang{\cL^{\new}}$; we prove this claim in \cref{sec:appendix-grobner}.

\begin{restatable}{claim}{grobnerSaturation}
\label{claim:grobner-saturation}
Fix $d,D\geq 1$ and $\cL^{\new}\subset\CC[x_1,\ldots,x_n,y_1,\ldots,y_m]$, a finite list of polynomials such that every $g\in\cL^{\new}$ has degree in $x_1,\ldots,x_n$ bounded by $D$ and degree in $y_1,\ldots,y_m$ bounded by $d$. Then there exists a polynomial $h\in \CC[x_1,\ldots,x_n]\setminus I(V)$ satisfying \cref{eq:saturation} such that the degree of $h$ in $x_1,\ldots,x_n$ is bounded by $O_{d,m}(D)$ and its degree in $y_1,\ldots,y_m$ is bounded by $O_{d,m}(1)$.
\end{restatable}

Let $h$ be the polynomial produced by this claim which satisfies \cref{eq:saturation}.\footnote{In fact, \cref{eq:saturation} will be an equality, since the reverse containment is true for every $h$ with $\psi(h)\in S$.} We now show that this $h$ has the desired property.

Suppose that $X$ is an irreducible non-transverse component of $W^{\new}$. Thus there exists a polynomial $h_X\in\CC[x_1,\ldots,x_n]\setminus I(V)$ that vanishes identically on $X$. Let $g_X\in\CC[x_1,\ldots,x_n,y_1,\ldots,y_m]$ be a polynomial that vanishes identically on $W^{\new}\setminus X$ without vanishing identically on $X$. (As in the proof of \cref{lem:radical-jacobian}, such a $g_X$ exists since the union of the irreducible components of $W^{\new}$ other than $X$ is a Zariski-closed subset which does not contain $X$.) Now $h_Xg_X$ vanishes identically on $W^{\new}$, implying that $(h_Xg_X)^k\in \sang{I(V),\cL^{\new}}$ for some $k$. In particular, $(\psi(h_X)\psi(g_X))^k\in \psi(\sang{\cL^{\new}})$. Note that $\psi(h_X)\in S$, implying that 
\[\psi(g_X)^k=\frac{(\psi(h_X)\psi(g_X))^k}{\psi(h_X)^k}\in\psi(\sang{\cL^{\new}})[S^{-1}]\cap\CC[x_1,\ldots,x_n,y_1,\ldots,y_m]/I(V).\]
By \cref{eq:saturation}, we have $\psi(g_X)^k\in\paren{\psi(\sang{\cL^{\new}}):\psi(h)^\infty}$. In other words, there exists $\ell$ so that $\psi(g_X)^k\psi(h)^\ell\in\psi(\sang{\cL^{\new}})$. This implies that $g_X^kh^\ell\in\psi^{-1}(\psi(\sang{\cL^{\new}}))=\sang{I(V),\cL^{\new}}$, implying that $g_X^kh^\ell$ vanishes identically on $Z(I(V),\cL^{\new})=W^{\new}$. As $g_X$ does not vanish identically on $X$, we conclude that $h$ vanishes identically on $X$, as desired.

Thus with this definition of $\Gamma^{\new},g^{\new},\cL^{\new}$ we see that properties (i-iii) still hold. We defined $\Gamma^{\new}$ so that (i) holds. We picked $h$ that does not vanish identically on $V$, so defining $g^{\new}=gh$ maintains property (ii). Finally, we just verified that property (iii) holds.

We repeat this process $m+1$ times. We claim that each component of the resulting $\Gamma$ rules $V$. We deduced earlier from \cref{lem:Frobenius} that if a component $\Gamma_i\subseteq \Gamma$ does not rule $V$, then it is cut by this process. Thus after $k$ iterations, the maximum dimension of a component of $\Gamma$ that does not rule $V$ is $\dim V+m-k$. Therefore after $m+1$ iterations, every remaining transverse component rules $V$.

At the end we produce a variety $\Gamma\subseteq V\times\CC^m$ each of whose irreducible components rules $V$. We also have a polynomial $g\in\CC[x_1,\ldots,x_n]$ which, by property (ii), does not vanish identically on $V$. By property (iii), we see that conclusion (3) follows as well. 

Since $\Gamma$ is defined to be the union of the transverse components of $W=(V\times\CC^m)\cap Z(\cL)$, we see that for a generic choice of $x\in V$, the fiber of the map $\pi\colon W\to V$ above $x$ is equal to the fiber of the map $\pi\colon\Gamma\to V$. This fiber is the common zero set of the polynomials $f(x,\cdot)$ for $f\in\cL$. We claim that these polynomials, viewed as elements of $\CC[y_1,\ldots,y_m]$, have degree $O_{\varphi}(1)$. To see this, note that the list $f_1,\ldots,f_r,\cL\subset\CC[x_1,\ldots,x_n,y_1,\ldots,y_m]$ begins with $\cL=\emptyset$ and $f_1,\ldots,f_r$ having degree 0 in $y_1,\ldots,y_m$. During this procedure, we replace $\cL$ by $\tilde\cL$ which, by \cref{claim:radical-bounded}, only increases the $(y_1,\ldots,y_m)$-degree in a bounded way. In addition we add to this list the components of $\vec D_y f(x,y)$ for $f\in\{f_1,\ldots,f_r\}\cup\tilde\cL$. By the definition of $\vec D_y$, this involves differentiating with respect to $x_1,\ldots,x_n$ and multiplying by the components of $\varphi(y)$, i.e., polynomials of degree $O_{\varphi}(1)$. Thus this operation also only increases the $(y_1,\ldots,y_m)$-degree in a bounded way. Thus after all $m+1$ iterations, we have only produced polynomials in $\cL$ whose $(y_1,\ldots,y_m)$-degree is $O_{\varphi}(1)$.

Now by B\'ezout's theorem, we see that the sum of the degrees of the irreducible components of the fiber is bounded by $O_{\varphi}(1)$ (see, e.g., \cite[Example 8.4.6]{Ful98} for this version of B\'ezout's theorem). This immediately implies that $\Gamma$ has at most $O_{\varphi}(1)$ irreducible components, and each irreducible component $\Gamma_i$ has fiber of degree $O_{\varphi}(1)$ over a generic point $x\in V$.
\end{proof}
\section{Bounding the algebraic concentration}
\label{sec:concentration}

In this section, we show that no variety $V\subseteq\CC^6\cong E^+_o(\CC)$ contains too many of the 3-flats $F_{pq}^+$. In particular, the (codimension 0) algebraic concentrations of this collection of 3-flats are as small as possible.
 
\begin{theorem}\label{thm:concentration-codim-0}
Let $\cP\subset \RR^3$ be a set of size $N$ such that $\abs{\cP\cap \gamma}\leq N^{2/3}\deg \gamma$ for every curve $\gamma\subset\CC^3$. 
Let $V\subseteq E^+_o(\CC)$ be an irreducible variety of dimension $d$.
Then
\[\abs{\left\{(p,q)\in \cP\times\cP:F_{pq}^+\subseteq V\right\}}\lesssim N^{2(d-3)/3}\deg V.\]
The symmetric result holds in $E^-_o(\CC)$.
\end{theorem}

The proof has three steps. First, we classify the varieties which are ruled by physical 3-flats. Second, we use this classification to prove a version of \cref{thm:concentration-codim-0} in the case when $V$ is ruled: \cref{thm:concentration-ruled-varieties}. Finally, we use the Cayley--Monge--Salmon theorem for flats (\cref{thm:general-CSM}) to bootstrap this result to the full strength of \cref{thm:concentration-codim-0}.

\subsection{Ruled variety theory}
\label{ssec:ruled-variety-theory}
In this section, we will classify the varieties that are ruled by the $3$-flats $F_x^+$ in $E^+_o(\CC)$.
To specialize the general ruled surface theory of \cref{sec:csm-for-flats} to the case of the physical 3-flats $F_x^+$ in positive space, we set $n=6$, $m=k=3$, and 
\[\varphi(\theta) =
\begin{bmatrix}
    0&\theta_3&-\theta_2\\
    -\theta_3&0&\theta_1\\
    \theta_2&-\theta_1&0
\end{bmatrix}.\]
With this definition of $\varphi$, we have $\varphi(\theta)x=-\theta\times x$.
For $\rho=(\vec\rho,\vec\rho_\ast)\in E^+_o(\CC)\cong \CC^6$ and $\theta\in \CC^3$, define $c(\rho,\theta) = (\theta, -\theta\times \vec{\rho}-\vec{\rho}_*)\in E^{\Phy}_o$. Comparing \cref{eq:l-v-theta} against \cref{eq:f-x-positive}, we see that $L_{\rho,\theta} = F_{c(\rho,\theta)}^+$.

Our goal is to understand the pairs $(V,\Gamma)$ where $V\subseteq\CC^6$ is an irreducible variety and $\Gamma\subseteq V\times\CC^3$ rules $V$. We do a case analysis based on $\dim V$ and $\dim \Gamma$.
Note that the case $\dim V \leq  3$ and $\dim V = 6$ are both trivial, so we just need to consider the two cases $\dim V = 4$ and $\dim V = 5$.

\begin{remark}
While in this section we focus on varieties in positive space ruled by physical 3-flats, the same results hold for varieties in negative space that are ruled by physical 3-flats. Indeed, inspecting \cref{eq:f-x-positive,eq:f-x-negative}, we see that defining $i\colon \CC^6\to\CC^6$ by $i(\vec x,\vec x_\ast)=(\vec x_\ast,\vec x)$, then $i\colon E_o^+\to E_o^-$ is a linear isomorphism satisfying $i(F_x^+)=F_{i(x)}^-$. Thus the theory of varieties ruled by physical 3-flats in positive space and negative space are identical.
\end{remark}

\subsubsection{4-dimensional ruled varieties}
\label{sssec:ruled-variety-theory-4-dim}
We start with the easier case where $\dim V = 4$.
Suppose that $\Gamma$ is an irreducible variety in $V\times \CC^3$ that rules $V$.
By definition, the projection $\Gamma\to V$ is dominant, so $\dim \Gamma \geq 4$.

For a generic point $\rho\in V$, the tangent space $T_{\rho}V$ is $4$-dimensional. 
For such $\rho$, if there are $\theta\neq \theta'$ such that $(\rho,\theta)$ and $(\rho,\theta')$ are both in $\Gamma$, then $F_{c(\rho,\theta)}^+$ and $F_{c(\rho,\theta')}^+$ are both $3$-flats in $T_{\rho}V$.
This shows $F_{c(\rho,\theta)}^+\cap F_{c(\rho,\theta')}^+$ contains a $2$-flat and hence $F_{c(\rho,\theta)}^+=F_{c(\rho,\theta')}^+$ by \cref{prop:intersection-from-same-family}, which is a contradiction.
This shows that for a generic point $\rho\in V$, the fiber of $\Gamma\to V$ above $\rho$ is a single point.
By the fiber dimension theorem (see, e.g., \cite[Corollary 11.13]{Harris92}), we conclude $\dim \Gamma = \dim V = 4$. (Another consequence of this argument is that all ruled 4-dimensional varieties are actually singly ruled.)

Define $\Lambda=\overline{c(\Gamma)}$, the Zariski closure of $c(\Gamma)$. By definition, if $(v,\theta)\in\Gamma$, then $L_{v,\theta}\subset V$ and $(v',\theta)\in\Gamma$ for all $v'\in L_{v,\theta}$. For all of these $v'$, we have $L_{v,\theta}=L_{v',\theta}=F^+_{c(v,\theta)}=F^+_{c(v',\theta)}$. In other words, the map $c\colon \Gamma \to \Lambda\subset E_o^{\Phy}$ has 3-dimensional fibers. By the fiber dimension theorem, we again conclude that $\dim \Lambda = \dim\Gamma-3=1$.

\begin{remark}
Ruled 4-dimensional varieties $V$ with a corresponding curve $\Lambda$ are ubiquitous. Given any irreducible curve $\Lambda\subset E^{\Phy}_o$, define $\Gamma=c^{-1}(\Lambda)$ and take $V$ to be the variety swept out by $(F_x^+)_{x\in \Lambda}$, i.e.,
\[V=\overline{\bigcup_{x\in\Lambda}F_x^+}.\]
Then $V$ is a $4$-dimensional variety ruled by $\Gamma$.
\end{remark}

While this remark shows that there are many ruled $4$-dimensional varieties with $\dim \Gamma = 4$, it suffices for our purpose to be able to bound $\deg \Lambda$ in terms of $\deg V$.

\begin{theorem}\label{thm:ruled-4-dim-4-dim}
    Let $V$ be a $4$-dimensional irreducible variety ruled by some irreducible variety $\Gamma$.
    Set $\Lambda=\overline{c(\Gamma)}$.
    Then $\dim\Lambda=1$ and $\deg \Lambda \leq \deg V$.
\end{theorem}

To show this, we use the following lemma that relates the degree of a projective variety to the degree of its image under a well-behaved projection.

\begin{lemma}\label{lem:projection}
    Let $n>m$ be two positive integers. Let $H_0$ be an $(n-m-1)$-plane in $\PP_{\CC}^n$ and let $H_1$ be an $m$-plane in $\PP_{\CC}^n$ which is disjoint from $H_0$.
    Define the projection $\pi\colon \PP^n\setminus H_0\to H_1\cong \PP^m$ as follows: for each $p\in\PP^n\setminus H_0$, let $\pi(p)$ be the unique intersection of $H_1$ and the $(n-m)$-plane spanned by $H_0\cup \{p\}$. 
    Suppose that $V\subseteq \PP_{\CC}^n$ is an irreducible variety satisfying the following:
\begin{enumerate}
    \item $V$ is disjoint from $H_0$; and
    \item $\pi\colon V\to \pi(V)$ has generic fibers consisting $k<\infty$ points.
\end{enumerate}
Then $\pi(V)$ is an irreducible variety of the same dimension as $V$ and of degree $\tfrac1k\deg V$.
\end{lemma}

As $V$ is irreducible, the same is true for $\pi(V)$. Note that if the generic fibers of $\pi\colon V\to\pi(V)$ are finite, then there always is some $k$ (called the degree of $\pi$) so that a generic fiber consists of $k$ points. (See, e.g., \cite[Proposition 7.16]{Harris92}.)

\begin{proof}
We first prove this lemma in the case when $m=n-1$, i.e., $H_0=\{h_0\}$ is a single point and $H_1$ is a hyperplane.

As the generic fibers of $\pi$ have dimension 0, the fiber dimension theorem implies that $\dim \pi(V) = \dim V$.
Let $d$ be the dimension of $V$.
If $d=0$ then the statement clearly holds, so we may assume that $d>0$.

Write $\G(n-d,n)$ for the Grassmannian of $(n-d)$-planes in $\PP_{\CC}^n$. Define $W\subset \G(n-d,n)$ to be the subvariety consisting of those $(n-d)$-planes which contain $h_0$. Then $\dim W = ((n-d+1)-1)(n-(n-d)) = (n-d)d$. This is because an $(n-d)$-plane can be written as the span of $n-d+1$ points chosen in $\PP_{\CC}^n$; for this plane to contain $h_0$, choose $h_0$ as this first point, then space of choices for the remaining $n-d$ points is $(n-d)n$-dimensional. Each such $(n-d)$-plane is chosen many times: the space of choices for a given $(n-d)$-plane is $(n-d)(n-d)$-dimensional.

Recall that for an $(n-d)$-plane $H$, we say that $V,H$ intersect transversally if for every $p\in H\cap V$, the variety $V$ is smooth at $p$ and $T_pV$ and $H$ are transversal. By B\'ezout's theorem, if $V,H$ intersect transversally, then $\deg V=|H\cap V|$ (see, e.g., \cite[Theorem 18.3]{Harris92} for this version of B\'ezout's theorem). Similarly, if $\pi(V),H$ intersect transversally, $\deg \pi(V)=|H\cap \pi(V)|$. To relate $\deg V,\deg \pi(V)$, we wish to show that for a generic $H\in W$, the $(n-d)$-plane $H$ intersects both $V,\pi(V)$ transversally. To do that, we bound the dimension of the space of $H$ which do not intersect each transversally.

For a $d$-dimensional irreducible variety $U$ which does not contain $h_0$, let $\Gamma_U\subseteq W\times U$ be the set of points $(H,u)\in W\times U$ where $H$ and $U$ do not intersect transversally at $u$. (One can see that $\Gamma_U$ is a variety as the condition $u\in H$ is clearly closed while the condition $\dim(H\cap T_uU)\geq 1$ can be written as the vanishing of the minors of a matrix defined in terms of $H$ and $\vec\nabla f(u)$ for $f\in I(U)$.) We will analyze the components of $\Gamma_U$ that have dimension at least $(n-d)d$, with the aim of showing that $\dim \Gamma_U<\dim W$ for $U=V,\pi(V)$.

Let $X$ be an irreducible component of $\Gamma_U$ with dimension at least $(n-d)d$, and set $X'=\pi_2(X)$ where $\pi_2$ is the projection $\pi_2\colon W\times U\to U$.
Note that by assumption $\pi_2(X)\subseteq U$ does not contain $h_0$.

If $\dim X'=d$ (i.e., $X'=U$), then for a generic point $u\in U$, the fiber dimension theorem implies that the fiber $X_u$ above $u$ has dimension $\dim X-\dim U\geq (n-d)d-d=(n-d-1)d$. Furthermore, for a generic point of $u\in U$, we know that $u$ is a smooth point of $U$, i.e., $\dim T_uU=d$. We noted above that $U$ does not contain $h_0$. Finally, we claim that for generic $u\in U$, that $h_0\not\in (u+T_uU)$, i.e., the affine tangent space through to $U$ at $u$ does not contain $h_0$. If this were not the case, then $U$ would be a cone with apex $h_0$.\footnote{To see this, note that in a small (Euclidean) neighborhood of a generic point $u\in U$, we can view $U$ as a smooth real manifold where the vector field pointing radially toward $h_0$ is a smooth real vector field on this manifold. The integral curves of this vector field are lines through $h_0$ and thus $U$ contains a small Euclidean neighborhood of $u$ on the line connecting $h_0$ and $u$. As $U$ is a variety, this implies that it contains the line connecting $h_0$ and $u$ for generic $u\in U$.} However, this would contradict the fact that $U$ is disjoint from $h_0$.

In the case $n=d-1$, this argument suffices. Indeed, the above calculation shows in this case that $\dim X_u\geq 0$ for generic $u\in U$, meaning that the line $H$ through $h_0$ and $u$ does not intersect $U$ transversely at $u$ for generic $u\in U$. This means that this line is contained in the affine tangent space to $U$ at $u$ for generic $u\in U$, contradicting what we showed above.

Now assume $n\leq d-2$. If $H\in W$ is an $(n-d)$-plane through $u$ that does not intersect $U$ transversally at $u$, then $H\cap (u+T_uU)$ has dimension at least 1. Thus for such a $u$, we conclude that \[\dim X_u \leq (d-1)+((n-d+1)-3)(n-(n-d)) = (d-1)+(n-d-2)d < (n-d-1)d.\]
Here the $d-1$ counts the number of dimensions to pick a 1-dimensional subspace through $u$ in $u+T_uU$ and the second term counts the number of dimensions to extend the $2$-plane spanned by $h_0$ and this 1-dimensional subspace to an $(n-d)$-plane. (We can do this since $n-d\geq 2$.) This is a contradiction, since we showed earlier that, in this case, $\dim X_u\geq (n-d-1)d$ for a generic $u\in U$.

Next if $\dim X'<d$, then for a generic $u\in X'$, the fiber $X_u$ above $u$ has dimension more than $(n-d)d-d = (n-d-1)d$. On the other hand, if $(H,u)\in X$, we simply have that $H$ contains the line spanned by $h_0,u$. This already implies that $\dim X_u \leq ((n-d+1)-2)(n-(n-d)) = (n-d-1)d$, which is still a contradiction.
As a consequence, $\dim \Gamma_U<(n-d)d$ for any irreducible $d$-dimensional $U$ which does not contain $h_0$.

Since $\dim \Gamma_V,\dim \Gamma_{\pi(V)}<\dim W$, we conclude that for a generic $H\in W$, i.e., a generic $(n-d)$-flat $H$ that contains $h_0$, the intersections of $V$ and $\pi(V)$ with $H$ are transversal. Thus by B\'ezout's theorem, we see that $|H\cap V|=\deg V$ and $|H\cap \pi(V)|=\deg \pi(V)$ for a generic $H\in W$.

Now we can write
\[|H\cap V|=\sum_{q\in \pi(V)\cap H}\abs{\pi|_V^{-1}(q)}.\]
By assumption, the generic fibers of $\pi\colon V\to \pi(V)$ have size $k$. Thus for a generic choice of $H$ containing $h_0$, each point $q\in \pi(V)\cap H$ will satisfy $\abs{\pi|_V^{-1}(\pi(p))}=k$. Thus we see that $\abs{H\cap V} = k\abs{H\cap \pi(V)}$, giving the desired result.

Finally, we are ready to prove \cref{lem:projection} in full. Let $h_1,h_2,\ldots,h_{n-m}\in H_0$ be points whose span is $H_0$. Let $H^{(i)}$ be the $(m+i)$-plane which is the span of $H_1\cup\{h_1,\ldots,h_i\}$. Then let $\pi_i\colon H^{(i)}\setminus \{h_{i}\}\to H^{(i-1)}$ be the projection defined in the statement when applied to $\{h_{i}\}, H^{(i-1)}\subset H^{(i)}\cong\PP^{m+i}$, a disjoint 0-plane and $(m+i-1)$-plane.

Note that $\pi=\pi_1\circ\cdots\circ\pi_{n-m}$. Furthermore, for $V$ disjoint from $H_0$, write $V_{n-m}=V$ and $V_{i-1}=\pi_i(V_i)$. Then $V_i\subset H^{(i)}$ does not contain $h_i$ (if it did, the plane spanned by $\{h_i,\ldots,h_{n-m}\}$ would intersect $V$, contradicting assumption (1)). By assumption (2), the generic fibers of $\pi\colon V=V_{n-m}\to \pi(V)=V_0$ are finite of size $k$. Thus the generic fibers of $\pi_i\colon V_i\to V_{i-1}$ are finite, say of size $k_i$. We have $k=k_1k_2\cdots k_{n-m}$. Finally, applying \cref{lem:projection} to $\pi_i\colon V_i\to V_{i-1}$, we conclude that $\dim V_{i-1}=\dim V_i$ and $\deg V_{i-1}=\tfrac1{k_i}\deg V_i$. Combining these gives $\dim V_0=\dim V_{n-m}$ and  $\deg V_0=\tfrac 1k\deg V_{n-m}$.
\end{proof}

\cref{lem:projection} allows us to relate $\deg \Lambda$ and $\deg V$ in \cref{thm:ruled-4-dim-4-dim} by exhibiting an appropriate linear projection.

\begin{proof}[Proof of \cref{thm:ruled-4-dim-4-dim}]
We are given a 4-dimensional variety $V\subset\CC^6\cong E_o^+(\CC)$ ruled by some irreducible variety $\Gamma$. By the discussion at the beginning of \cref{sssec:ruled-variety-theory-4-dim}, we know that $\Gamma$ has dimension 4 and $\Lambda=\overline{c(\Gamma)}$ has dimension 1.

Recall that we write $\rho\in E_o^+(\CC)$ as $\rho=(\vec\rho,\vec\rho_\ast)$ where $\vec\rho,\vec\rho_\ast\in\CC^3$. For $\vec s \in\CC^3$, define $V_{\vec{s}}=\{\rho\in V:\vec{\rho}=\vec{s}\}$. Writing $\pi_1\colon\CC^6\to\CC^3$ for the projection onto the first three coordinates, we have $V_{\vec s}=\pi_1|_V^{-1}(\vec s)$.

We claim that $\pi_1(V)=\CC^3$; by the fiber dimension theorem, this would imply that $\dim V_{\vec s}=\dim V-\dim\CC^3=1$ for a generic $\vec s\in\CC^3$. To see the claim, note that $\pi_1(F_x^+)=\CC^3$ for every $x\in\CC^6\cong E_o^{\Phy}$. (This follows directly from the formula for given in \cref{eq:f-x-positive}.) As $V$ is ruled, it clearly contains at least one such $F^+_x$. Thus we conclude that $\dim V_{\vec s}=1$ for a generic $\vec s\in\CC^3$. As $V_{\vec s}$ is the intersection of $V$ with a 3-flat, B\'ezout's theorem implies that $\deg V_{\vec{s}}\leq \deg V$ for all $\vec{s}\in\CC^3$.

Now consider the map $\phi_{\vec{s}}\colon \Lambda \to E^+_o\cong \CC^6$ sending $x=(\vec x,\vec x_\ast)$ to $(\vec{s},- \vec{x}\times \vec{s}-\vec{x}_*)$.
As this latter point is $F_x^+\cap \{\rho\in E^+_o:\vec{\rho}=\vec{s}\}$, we see that $\phi_{\vec{s}}(\Lambda)$ completely lies in $V_{\vec{s}}$. Set $\Lambda_{\vec{s}}=\overline{\phi_{\vec{s}}(\Lambda)}$. We claim that for a generic $\vec s\in\CC^3$, both $\dim \Lambda_{\vec{s}}=1$ and $\deg \Lambda = \deg \Lambda_{\vec{s}}$ hold. This will suffice to complete the proof, as for a generic $\vec s\in\CC^3$ we will have $\Lambda_{\vec{s}}\subseteq V_{\vec s}$ are both 1-dimensional varieties. Thus $\deg\Lambda_{\vec s}\leq\deg V_{\vec s}$ implies
\[\deg\Lambda=\deg\Lambda_{\vec s}\leq\deg V_{\vec s}\leq \deg V,\]
as desired.

Both parts of the claim follow by applying \cref{lem:projection} to $\phi_{\vec s}\colon \Lambda\subset E_o^{\Phy}\to \Lambda_{\vec s}\subset E_o^+$. To apply this result, we need to extend $E^{\Phy}_o\cong \CC^6$ and $E^+_o\cong\CC^6$ to $\PP_{\CC}^6$. We add appropriate signs in this extension so that $\phi_{\vec s}$ extends to a map which is a linear projection. In particular, embed $x=(\vec x,\vec x_\ast)\in E^{\Phy}_o\cong\CC^6$ to $[1:\vec{x}:-\vec{x}_*]\in\PP^6$ and embed $\rho=(\vec \rho,\vec \rho_\ast)\in E^{+}_o\cong\CC^6$ to $[1:\vec{\rho}:\vec{\rho}_*]\in\PP^6$.
Then $\phi_{\vec{s}}$ extends to the map $\PP^6\to\PP^6$ that sends $[x_0:\cdots:x_6]$ to $[x_0:s_1x_0:s_2x_0:s_3x_0:-s_3x_2+s_2x_3+x_4:-s_1x_3+s_3x_1+x_5:-s_2x_1+s_1x_2+x_6]$. One can easily verify that setting $n=6$ and $m=3$ and
    \[H_0 = \{[x_0:\cdots :x_6]\in \PP_{}^6: x_0=-s_3x_2+s_2x_3+x_4=-s_1x_3+s_3x_1+x_5=-s_2x_1+s_1x_2+x_6=0\}\]
    and
    \[H_1 = \{[x_0:\cdots :x_6]\in \PP_{}^6:x_1-s_1x_0=x_2-s_2x_0=x_3-s_3x_0=0\},\]
then the map $\pi\colon \PP^6\setminus H_0\to H_1\subset\PP^6$ defined in the statement of \cref{lem:projection} is exactly $\phi_{\vec{s}}$. 
Let $\Tilde{\Lambda}$ be the Zariski closure of $\Lambda$ in $\PP^6$.

We claim that for a generic $\vec s\in\CC^3$, the corresponding 2-plane $H_0$ does not intersect $\Tilde{\Lambda}$. In other words, we claim that condition (1) in the statement of \cref{lem:projection} holds.

Consider $X = Z(x_0,-s_3x_2+s_2x_3+x_4,-s_1x_3+s_3x_1+x_5,-s_2x_1+s_1x_2+x_6)\cap (\Tilde{\Lambda}\times \CC^3)$ in $\PP^6\times \CC^3$.
Write $\pi_1\colon X\to \PP^6$ and $\pi_2\colon X\to\CC^3$ for the projection maps onto the first and second coordinates.
For a fixed $x\in Z(x_0)\cap\Tilde\Lambda\subset \PP^6$, we see that $(x,\vec{s})$ is in the fiber of $\pi_1$ above $x$ if and only if $(x_1,x_2,x_3)\times \vec{s}  = (x_4,x_5,x_6)$.
Note that if $(x_1,x_2,x_3)\neq 0$, then the set of $\vec s$ satisfying this equation is at most 1-dimensional. Furthermore, if $(x_1,x_2,x_3)=0$, then $(x_4,x_5,x_6)\neq 0$, showing that there are no solutions.

Thus the fibers of $\pi_1\colon X\to\PP^6$ have dimension at most 1. As the image of this projection lies in $\Tilde{\Lambda}\cap Z(x_0)$, which satisfies $\dim(\tilde\Lambda\cap Z(x_0))=\dim\Lambda-1=0$, we see that $\dim X \leq 1$. In particular, $\dim X<3$, meaning that a generic $\vec s\in\CC^3$ does not lie in $\pi_2(X)$; in other words, for a generic $\vec{s}\in \CC^3$, the corresponding $H_0$ has no intersection with $\Tilde{\Lambda}$.

Finally, we claim that for generic $\vec s$, the projection $\phi_{\vec{s}}\colon \Lambda \to \Lambda_{\vec{s}}$ has generic fibers of size 1. In particular, in the statement of \cref{lem:projection}, condition (2) holds with $k=1$.

To see this, note that if $x,y\in\Lambda$ satisfy $\phi_{\vec s}(x)=\phi_{\vec s}(y)=p\in V$, then by definition $p\in F_x^+,F_y^+\subset V$. As we showed at the beginning of \cref{sssec:ruled-variety-theory-4-dim}, this implies that $\dim T_pV\geq 5$, so $p$ is a singular point of $V$. Now it cannot be the case that a generic point of $\Lambda_{\vec s}\subseteq V_{\vec s}$ is singular for a generic $\vec s$, as otherwise the singular locus of $V$ would be 4-dimensional. Thus we see that for a generic $\vec s\in\CC^3$, the projection $\phi_{\vec s}$ has generic fibers of size 1, implying that we can apply \cref{lem:projection} with $k=1$. This shows that for generic $\vec s\in\CC^3$, both $\dim \Lambda_{\vec s}=\dim V_{\vec s}$ and $\deg \Lambda_{\vec s}\leq\deg V_{\vec s}$ hold, as desired.
\end{proof}

\subsubsection{5-dimensional ruled varieties}
\label{sssec:ruled-variety-theory-5-dim}
Now we study the case $\dim V = 5$.
Again let $\Gamma$ be an irreducible variety in $V\times \CC^3$ that rules $V$.
We once again have $\dim \Gamma \geq 5$.

We first deal with the case $\dim \Gamma>5$.

\begin{lemma}
\label{lem:ruled-5-dim-6-dim-intersecting}
Let $V\subset E^+_o$ be an irreducible hypersurface that is ruled by an irreducible variety $\Gamma\subseteq V\times \CC^3$ of dimension at least $6$.
Set $\Lambda=\overline{c(\Gamma)}$. Then any line passing through two points of $\Lambda$ is isotropic.
\end{lemma}

\begin{proof}
As in \cref{sssec:ruled-variety-theory-4-dim}, the map $c\colon \Gamma\to \Lambda$ has 3-dimensional fibers, so $\dim \Lambda = \dim \Gamma-3$.
Now by assumption $\dim \Gamma>\dim V$ and the projection $\pi_1\colon \Gamma\to V$ is dominant, so for a generic point $p\in V$, the fiber of $\pi_1\colon \Gamma\to V$ above $p$ has dimension $\dim\Gamma-\dim V>0$. Thus for a generic point $x\in \Lambda$, the fiber of $\pi_1\colon \Gamma\to V$ above a generic point in $p\in F_x^+\subset V$ has dimension $\dim \Gamma-\dim V>0$. Fix such an $x\in\Lambda$.

Now $\pi_1^{-1}(F_x^+)\subset\Gamma$ has dimension $\dim F_x^+ + (\dim \Gamma-\dim V)=\dim\Gamma-2$. For any $x'\in \Lambda\setminus\{x\}$, by \cref{prop:intersection-from-same-family} we have $F^+_{x'}\cap F^+_x$ is either one-dimensional or empty. Thus we see that $c\colon \pi_1^{-1}(F_x^+)\to\Lambda$ has 1-dimensional fibers (except over $x$). We conclude that $\dim \overline{c(\pi_1^{-1}(F_x^+))}=\dim\pi_1^{-1}(F_x^+)-1=\dim\Gamma-3=\dim\Lambda$. As $\Gamma$ is irreducible, so is $\Lambda$, implying that $\overline{c(\pi_1^{-1}(F_x^+))}=\Lambda$. Now by Chevalley's theorem, $c(\pi_1^{-1}(F_x^+))$ is a constructible set. As its Zariski-closure is equal to $\Lambda$, we conclude that a generic point $x'\in\Lambda$ lies in $c(\pi_1^{-1}(F_x^+))$. Tracing through the definitions, this implies that $F_{x'}\cap F_x\neq\emptyset$ for a generic point $x'\in \Lambda$.
By \cref{prop:intersection-from-same-family}, we can write this as
\[0=\sang{x,x'}=(\vec x-\vec x')\cdot(\vec x_\ast-\vec x'_\ast)\]
for two generic points $x,x'\in \Lambda$.
Since the vanishing of a polynomial is a Zariski-closed condition, we conclude that the polynomial vanishes for any two $x,x'\in \Lambda$; by another application of \cref{prop:intersection-from-same-family}, this gives the desired result.
\end{proof}

This leads us to classify the subvarieties $X$ of $E^{\Phy}_o$ where any line passing through two points is isotropic.
An obvious example is a variety that lie completely in $F_{\rho}^{\Phy}$ for some $\rho\in E^+\sqcup E^-$.
We will show that this is the only possibility.

\begin{lemma}
\label{lem:q-vanishing-extend-to-3-flat}
Let $X\subseteq E^{\Phy}_o$ be a variety so that any line passing through two points is isotropic.
Then there exists some $\rho\in E^+\sqcup E^-$ so that $X\subseteq F_{\rho}^{\Phy}$.
\end{lemma}

\begin{proof}
View the set $X$ as lying in $X\subseteq E^{\Phy}_o\subset E^{\Phy}\subset \PP(\OO_s)$. Let $\tilde X\subset \OO_s$ be the set of $x\in\OO_s$ such that $[x]\in X\subset\PP(\OO_s)$. By assumption, for any $[x]\neq[x']\in X$, the line through $[x],[x']$ is isotropic. As we showed in \cref{sec:graph-and-compactification}, this implies that $\langle x,x'\rangle=0$. By definition, $\langle x,x\rangle = 2N(x)=0$ whenever $[x]\in E^{\Phy}$. Thus we conclude that $\langle x,x'\rangle=0$ for all $x,x'\in \tilde X$. Since this condition is bilinear, defining $F$ to be the linear span of $X\subset\OO_s$, we conclude that $\langle x,x'\rangle=0$ for all $x,x'\in F$. In particular, $N(x)=\tfrac12 \langle x,x\rangle=0$ for all $x\in F$.

This means that $F$ is a totally isotropic subspace of $\OO_s$. Let $F'$ be a maximal totally isotropic subspace of $\OO_s$ which contains $F$. By \cref{thm:4-dim-isotropic-characterization,prop:all-3-flats}, the subspace $F'$ is 4-dimensional and $\PP(F')\subset\PP(\OO_s)$ is a positive or negative 3-plane in physical space. In other words, $X\subseteq \PP(F')=F_{\rho}^{\Phy}$ for some $\rho\in E^+\sqcup E^-$.
\end{proof}

Combining \cref{lem:ruled-5-dim-6-dim-intersecting,lem:q-vanishing-extend-to-3-flat} immediately lets us classify all the possible $\Lambda$ in this case.

\begin{corollary}
\label{thm:ruled-5-dim-6-dim-Lambda}
Let $V$ be an irreducible hypersurface ruled by some irreducible $\Gamma \subseteq V\times \CC^3$ with $\dim \Gamma\geq 6$.
Set $\Lambda=\overline{c(\Gamma)}$.
Then $\Lambda=F_{\rho}^{\Phy}$ for some $\rho\in E^+\sqcup E^-$.
\end{corollary}

\begin{proof}
Recall that $c\colon\CC^6\times\CC^3\to\CC^6$ has fibers which are 3-flats. Thus $c|_{\Gamma}\colon \Gamma\to\Lambda$ has fibers of dimension at most 3. By the fiber dimension theorem, $\dim\Lambda+3\geq\dim\Gamma\geq 6$, so $\dim\Lambda\geq 3$. By \cref{lem:ruled-5-dim-6-dim-intersecting,lem:q-vanishing-extend-to-3-flat}, we conclude that $\Lambda$ is contained in a positive or negative 3-flat in physical space. Combining the two, we see that $\Lambda$ must in fact be equal to this 3-flat.
\end{proof}

We will not need this fact, but we point out that taking $\Lambda=F_{\rho}^{\Phy}$ for $\rho\in E^+$ does correspond to a bounded-degree ruled hypersurface, but taking $\rho\in E^-$ does not.

\begin{remark}
For $\rho_0\in E^+$, defining
\[V=\overline{\bigcup_{x\in F_{\rho_0}^{\Phy}}F_x^+},\]
one can show that $V=\{\rho\in E^+_o:\langle\rho,\rho_0\rangle=0\}$ is a quadratic hypersurface that is ruled by the 3-flats $\{F_x^+\}_{x\in F_{\rho_0}^{\Phy}}$. (Geometrically, $V$ is the set of positive rigid motions which agree with $\rho_0$ on at least one point: a 5-dimensional family.) Each point of $V$ is contained in a one-dimensional family of these 3-flats. (The set of $\rho\in\PP(\OO_s)$ such that $\langle\rho,\rho_0\rangle=0$ is a hyperplane. However, under the identification of $E^+_o(\CC)=\CC^6$ with an affine open subset of $E^+(\CC)\subset\PP(\OO_s)$ by $(\vec\rho,\vec\rho_\ast)\mapsto[\vec\rho\cdot\vec\rho_\ast:\vec\rho:1:\vec\rho_\ast]$, this hyperplane maps to a degree 2 hypersurface.) 

For $\tau_0\in E^-$, defining
\[V=\overline{\bigcup_{x\in F_{\tau_0}^{\Phy}}F_x^+},\]
one instead finds $V=\CC^6$. (Geometrically, $V$ is the set of positive rigid motions which agree with $\tau_0$ on at least one point: a 6-dimensional family.) This means that there is no ruled hypersurface $V$ with corresponding $\Lambda=F_{\tau_0}^{\Phy}$.
\end{remark}

This concludes the discussion of the case $\dim \Gamma > 5$.
Similar to before, the case $\dim \Gamma = 5$ is hard to classify as precisely, but we are able to relate $\deg \Lambda$ and $\deg V$.
This time the projection $\Gamma \to V$ does not necessarily have generic fibers of size $1$.
Fortunately since we only work with $\Gamma$ that are produced by \cref{thm:general-CSM}, we can assume safely that the generic fibers of $\Gamma \to V$ have size $O(1)$.

\begin{theorem}\label{thm:ruled-5-dim-5-dim}
Let $V$ be an irreducible hypersurface ruled by an irreducible 5-dimensional $\Gamma\subset V\times\CC^3$. Suppose that a generic fiber of $\Gamma\to V$ has size at most $C$. Set $\Lambda=\overline{c(\Gamma)}$. Then $\dim\Lambda=2$ and $\deg \Lambda \leq C \deg V$.
\end{theorem}

\begin{proof}
The proof proceeds the same as the one for \cref{thm:ruled-4-dim-4-dim}. As before, we know that the map $c\colon\Gamma\to\Lambda$ has 3-dimensional fibers, so $\dim\Lambda=\dim\Gamma-3=2$. Define $V_{\vec s}$ as in the proof of \cref{thm:ruled-4-dim-4-dim}. By the same argument, we have $\dim V_{\vec s}=\dim V-3=2$ for a generic $\vec s\in\CC^3$ and $\deg V_{\vec s}\leq\deg V$ for all $\vec s\in\CC^3$.

Define the projection $\phi_{\vec{s}}$ as before. Set $\Lambda_{\vec s}=\phi_{\vec s}(\Lambda)$ and let $\tilde\Lambda$ be the Zariski closure of $\Lambda$ in $\PP^6$ (again embedding $\Lambda\subset\CC^6$ into $\PP^6$ with appropriate signs). We claim that for a generic $\vec s\in\CC^3$ we can apply \cref{lem:projection} to $\tilde \Lambda$ with $k\leq C$ to conclude that $\Lambda_{\vec s}\subseteq V_{\vec s}$ are 2-dimensional varieties satisfying
\[\frac 1C\deg \Lambda =\frac 1C\deg\tilde\Lambda\leq \deg \Lambda_{\vec s} \leq \deg V_{\vec{s}}\leq \deg V.\]

To show that \cref{lem:projection} applies, define $X$ as before. We still have that $\pi_1\colon X\to\PP^6$ has fibers of dimension at most 1 and the image is contained in $\tilde\Lambda\cap Z(x_0)$ which has dimension $\dim \Lambda-1=1$, so $\dim X\leq 2$. In particular, $\dim X<3$, so we again conclude that for generic $\vec s\in\CC^3$, condition (1) of \cref{lem:projection} holds. Finally, by hypothesis, we see that a generic fiber of $\phi_{\vec s}$ has size at most $C$, so condition (2) holds with some $k\leq C$.
\end{proof}

\subsection{Non-concentration in ruled varieties}
\label{ssec:concentration-ruled}

Given $\cP\subset\CC^3$, we consider the collection of physical 3-flats $F^+_x$ for $x\in\phi(\cP\times\cP)$. In the previous subsection, we related a ruled variety $V\subseteq \CC^6\cong E^+_o(\CC)$ to the variety $\Lambda\subseteq\CC^6\cong E^{\Phy}_o(\CC)$ which parametrizes the $x$ so that $F^+_x\subseteq V$. Now we study how many points of $\phi(\cP\times\cP)$ can lie in $\Lambda$. 

We begin by showing that the non-concentration assumption on $\cP$ -- the fact that no curve contains too many points of $\cP$ -- implies a non-concentration property of $\cP\times\cP$. With the results of the previous subsection, this suffices to show non-concentration of the ruled varieties in positive space. (Recall that $\phi$ is just the linear transformation $\phi(p,q)=((p+q)/2,(p-q)/2)$, so its presence is not important to this discussion.)

\begin{lemma}
\label{thm:grid-surface-intersection}
Given parameters $A\geq N^{1/2}$, let $\cP\subset\CC^3$ be a set of size $N$ such that $|\cP\cap\gamma|\leq A\deg \gamma$ for every curve $\gamma\subset\CC^3$. Then \[\abs{X\cap (\cP\times \cP)}\lesssim A^d\deg X\]
for every irreducible variety $X\subset\CC^6$ with $\dim X=d<3$.
\end{lemma}

\begin{proof}
For appropriate constants $1=K_0<K_1<K_2$, we will show that the bound $\abs{X\cap (\cP\times \cP)}\leq K_dA^d\deg X$ holds for all irreducible $X$ of dimension $d$. We prove this by induction on $d$. The case of $d=0$ is trivial.

We write $\pi_1\colon\CC^6\to\CC^3$ for the projection onto the first three coordinates. For an irreducible variety $X\subset\CC^6$ of dimension $d$, define $X_1=\overline{\pi_1(X)}$ and $D_1=\deg X_1$ and $d_1=\dim X_1$. Let $\cP_1=X_1\cap\cP$. It is clear that $X\cap(\cP\times\cP)=X\cap(\cP_1\times\cP)$.

First we handle the $d_1=0$ case. Since $X$ is irreducible, the same is true of $X_1$. Thus if $d_1=0$, then $X_1$ is a single point $p$ and $X=\{p\}\times Y$ for some variety $Y$ of the same degree and dimension as $X$. If $\dim Y=1$, then we have $|X\cap(\cP\times\cP)|\leq |Y\cap\cP|\leq A\deg Y$ by assumption. If $\dim Y=2$, then we have the bound $|X\cap(\cP\times\cP)|\leq |\{p\}\times \cP|= N$, which is strong enough, since $A^2\geq N$.

From now on, we assume that $1\leq d_1\leq d<3$.

\begin{claim}
\label{thm:main-bezout-claim}
Suppose that there exists a polynomial $f\in\CC[x_1,x_2,x_3]$ that vanishes on $\cP_1$ but does not vanish identically on $X_1$. Then $\abs{X\cap (\cP\times \cP)}\leq K_{d-1}A^{d-1}\deg(X)\deg(f)$.
\end{claim}

\begin{proof}
Write $\pi_1^*f$ for the polynomial $f$, viewed as an element of $\CC[x_1,x_2,x_3,y_1,y_2,y_3]$. Since $f$ does not vanish identically on $X_1$, we know that $\pi_1^*f$ does not vanish identically on $X$. By Krull's Hauptidealsatz, $X\cap Z(\pi_1^*f)$ is a variety of pure dimension $d-1$; let $Y_1,\ldots,Y_m$ be its irreducible components. By B\'ezout's theorem, \[\sum_i \deg Y_i\leq \deg(X)\deg(f).\] By the inductive hypothesis, 
\[\abs{Y_i\cap (\cP\times \cP)}\leq K_{d-1} A^{d-1} \deg(Y_i) \] for all $i$. Thus
\[\abs{X\cap (\cP\times \cP)}=\abs{X\cap(\cP_1\times\cP)}\leq \sum_i \abs{Y_i\cap (\cP\times \cP)}\leq K_{d-1}A^{d-1}\deg(X)\deg(f).\qedhere\]
\end{proof}

Next we handle the case $d_1=2$. By parameter counting, \cref{lem:Walsh-param-counting}, there is a polynomial $f$ that vanishes identically on $\cP_1$ without vanishing identically on $X_1$ which satisfies
\[\deg f\lesssim \paren{\frac{|\cP_1|}{\deg X_1}}^{1/2}+\paren{\frac{|\cP_1|\delta_1(X_1)}{\deg X_1}}^{1/3}\lesssim N^{1/2}\leq A.\]
The second inequality follows from the easy bounds $\delta_1(X_1)=\deg(X_1)\geq 1$ and $|\cP_1|\leq|\cP|\leq N$.
By \cref{thm:main-bezout-claim}, this bound suffices to complete the induction, taking $K_d$ sufficiently large in terms of $K_{d-1}$.

Finally we handle the case $d_1=1$. By \cref{lem:Walsh-param-counting}, there is a polynomial $f$ that vanishes identically on $\cP_1$ without vanishing identically on $X_1$ which satisfies
\[\deg f\lesssim \frac{|\cP_1|}{\deg X_1}+\paren{\frac{|\cP_1|\delta_2(X_1)}{\deg X_1}}^{1/2}+\paren{\frac{|\cP_1|\delta_1(X_1)\delta_2(X_1)}{\deg X_1}}^{1/3}.\]
By hypothesis, $|\cP_1|=|\cP\cap X_1|\leq A\deg X_1$. Thus to prove the desired $\deg f\lesssim A$ bound, it suffices to show that $\delta_1(X_1)\leq \delta_2(X_1)\lesssim A$.

Let $f_1\in\CC[x_1,x_2,x_3]$ be a minimal degree nonzero polynomial which vanishes identically on $\cP_1$. By parameter counting, $\deg f_1\lesssim |\cP_1|^{1/3}\leq N^{1/3}\leq A$. If $f_1$ does not vanish identically on $X_1$, we can complete the induction by \cref{thm:main-bezout-claim}. In the remaining case, $f_1$ vanishes identically on $X_1$. Thus $\delta_1(X_1)\leq \deg f_1\lesssim |\cP_1|^{1/3}\leq N^{1/3}$. 

The minimality of $f_1$ implies that it is irreducible. Indeed, if $f_1$ were not irreducible, it would have an irreducible factor which vanished identically on the irreducible variety $X_1\supset \cP_1$, contradicting minimality. Now let $f_2\in\CC[x_1,x_2,x_3]$ be a minimal degree polynomial which vanishes identically on $\cP_1$ without vanishing identically on $Z(f_1)$. Again by \cref{lem:Walsh-param-counting},
\[\deg f_2\lesssim \paren{\frac{|\cP_1|}{\deg g_1}}^{1/2}+|\cP_1|^{1/3}\lesssim |\cP_1|^{1/2}.\]
As before, if $f_2$ does not vanish identically on $X_1$, we can complete the induction by \cref{thm:main-bezout-claim}. Otherwise, $X_1$ is an irreducible component of $Z(f_1,f_2)$ which has pure dimension 1 and so $\delta_2(X_1)\leq\max\{\deg f_1,\deg f_2\}\lesssim |\cP_1|^{1/2}\leq N^{1/2}\leq A$.

Thus we have shown that the polynomial $f$ which vanishes identically on $\cP_1$ without vanishing identically on $X_1$ satisfies
\begin{align*}
\deg f
&\lesssim \frac{|\cP_1|}{\deg X_1}+\paren{\frac{|\cP_1|\delta_2(X_1)}{\deg X_1}}^{1/2}+\paren{\frac{|\cP_1|\delta_1(X_1)\delta_2(X_1)}{\deg X_1}}^{1/3}\\
&\leq A+\paren{A\delta_2(X_1)}^{1/2}+\paren{A\delta_1(X_1)\delta_2(X_1)}^{1/3}\lesssim A.\qedhere
\end{align*}
\end{proof}

\begin{theorem}
\label{thm:concentration-ruled-varieties}
Let $\cP\subset \RR^3$ be a set of size $N$ such that $\abs{\cP\cap \gamma}\leq N^{2/3}\deg \gamma$ for every curve $\gamma\subset\CC^3$. Let $V\subseteq E^+_o(\CC)$ be an irreducible variety of dimension $d$. Suppose that $V$ is ruled by $\Gamma\subseteq V\times\CC^3$ where the generic fibers of the projection $\Gamma\to V$ has degree at most $C$. Setting $\Lambda=\overline{c(\Gamma)}$, then
\[|\phi^{-1}(\Lambda)\cap(\cP\times\cP)|\lesssim C N^{2(d-3)/3}\deg V.\]
\end{theorem}

Recall that $\phi\colon \CC^3\times \CC^3\to \CC^6$ is the map $\phi(p,q)=((p+q)/2, (p-q)/2)$.

\begin{proof}
The result is only non-trivial for $d=4,5$.

When $d=4$, by \cref{thm:ruled-4-dim-4-dim}, we know that $\Lambda$ is a curve of degree at most $\deg V$. As $\phi$ is an invertible linear map, the same is true of $\phi^{-1}(\Lambda)$. By \cref{thm:grid-surface-intersection}, applied with $A=N^{2/3}$, we conclude that that $\abs{\phi^{-1}(\Lambda) \cap (\cP\times\cP)}\lesssim N^{2/3}\deg V$.

When $d=5$, by \cref{thm:ruled-5-dim-6-dim-Lambda,thm:ruled-5-dim-5-dim}, either $\Lambda=F_{\rho}^{\Phy}$ for some $\rho\in E^+\sqcup E^-$ or $\Lambda$ is a $2$-dimensional surface of degree at most $C\deg V$. In the latter case, by \cref{thm:grid-surface-intersection} with $A=N^{2/3}$, we conclude $\abs{\phi^{-1}(\Lambda) \cap (\cP\times\cP)}\lesssim N^{4/3}\deg \Lambda\leq C N^{4/3}\deg V$.

In the former case, we claim that for each $p\in\CC^3$, there is at most one $q\in\R^3$ so that $(p,q)\in\phi^{-1}(F_{\rho}^{\Phy})$. As $\cP\subset\R^3$, this will imply that $|\phi^{-1}(F_{\rho}^{\Phy})\cap (\cP\times\cP)|\leq|\cP|=N$. This claim is obvious if $\rho$ is a genuine rigid motion, as in this case $\phi^{-1}(F_{\rho}^{\Phy})=\{(p,q):\rho(p)=q\}$. However, it remains true even if $\rho$ is not genuine. Indeed, suppose $\phi(p,q),\phi(p,q')\in F_{\rho}^{\Phy}$ for some distinct $q,q'\in\R^3$. Since it lies in $F_{\rho}^{\Phy}$, the line $\aff\{\phi(p,q),\phi(p,q')\}$ is isotropic, meaning that
\[0=((p+q)-(p+q'))\cdot((p-q)-(p-q'))=-(q-q')\cdot(q-q').\]
As $q,q'\in\R^3$, this implies $q=q'$, as desired.
\end{proof}

\subsection{Proof of \texorpdfstring{\cref{thm:concentration-codim-0}}{Theorem 6.1}}
\label{ssec:concentration-main}
We now use the Cayley--Monge--Salmon theorem for flats to bootstrap \cref{thm:concentration-ruled-varieties} to prove \cref{thm:concentration-codim-0} in full.

\begin{proof}[Proof of \cref{thm:concentration-codim-0}]
We prove the result by induction on $d$. The result is only non-trivial for $d=4,5$.
    
Now suppose that $d=4$ or $5$, and the statement holds for varieties of dimension $d-1$.
Let $C$ be a sufficiently large constant to be determined.
Let $i$ be the largest index such that $\delta_i(V)\leq C^iN^{2/3}$. (As $\delta_0(V)=0$ and $\delta_{7-d}(V)=\infty$, such an $i$ exists.) 
Let $f_1,\ldots, f_i$ be the polynomials produced by \cref{lem:same-partial-degree}: this means $\deg f_j = \delta_j(V)$ for each $j\in[i]$ and every irreducible component $W$ of $Z(f_1,\ldots, f_i)$ containing $V$ has codimension $i$.
    
Applying \cref{thm:general-CSM} to each such irreducible component $W$ of $Z(f_1,\ldots,f_i)$ containing $V$ produces a polynomial $g_W$.
We claim that one of these polynomials does not vanish identically on $V$.
Otherwise, observe that $\deg g_W\lesssim \delta_i(V)$ and $g_W$ does not vanish identically on $W$ for each $W$. Thus by \cref{lem:lin-alg}, a generic linear combination $g$ of the $g_W$ does not vanish identically on any of the $W$. In particular $g$ vanishes identically on $V$ but does not vanish identically on any component of $Z(f_1,\ldots, f_i)$ containing $V$. This shows that $\delta_{i+1}(V)\leq \max\{\delta_i(V),\deg g\} \lesssim \delta_i(V)$.
 By taking $C$ sufficiently large, this contradicts the choice of $i$.

Fix a choice $W$ of an irreducible component of $Z(f_1,\ldots, f_i)$ containing $V$ so that $g_W$ does not vanish identically on $V$.
Let $\Gamma_1,\ldots, \Gamma_{\ell}$ be the corresponding irreducible varieties in $W\times \CC^m$ that rule $W$, and let $\Lambda_j=\overline{c(\Gamma_j)}$ for each $j\in[\ell]$.
Then we know that for every $(p,q)\in \cP\times \cP$, if $F_{pq}^+\subset V$, then either $\phi(p,q)\in \bigcup_{j\in[\ell]}\Lambda_j$ or $g_W$ vanishes identically on $F_{pq}^+$.
Applying the inductive hypothesis to $V\cap Z(g_W)$, a variety of pure dimension $d-1$ and degree at most $\deg g_W\deg V$, we see that there are at most $O(N^{2(d-4)/3}\deg g_W\deg V)$ flats $F_{pq}^+$ of the latter type.
    
To bound the number of flats of the former type, we apply \cref{thm:concentration-ruled-varieties} to $\Gamma_j$ ruling $W$. Note that $\dim W=6-i$ and \cref{thm:general-CSM} guarantees that the generic fiber of $\Gamma\to V$ has degree $O(1)$. Thus there are at most $O(N^{2(3-i)/3}\deg W)$ flats $F_{pq}^+$ with $\phi(p,q)\in\Lambda_j$. Since there are $\ell=O(1)$ choices for $j$, the same bound holds once we sum over all $\Lambda_j$.
    
Therefore we conclude that the total number of $F_{pq}^+$ contained in $V$ is bounded by
\[\lesssim N^{2(3-i)/3}\deg W+N^{2(d-4)/3}\deg g_W\deg V\lesssim N^{2(d-3)/3}\deg V.\]
The last inequality holds since $\deg g_W\lesssim \delta_i(V)\lesssim N^{2/3}$ and $\deg V \gtrsim \delta_{i+1}(V)\cdots \delta_{6-d}(V)\deg W\geq N^{2(6-d-i)/3}\deg W$.
\end{proof}
\section{Incidences of 3-flats}
\label{sec:codim-0-incidences}

In this section, we prove two incidence bounds, the first for line--3-flat incidences and the second for 3-flat--3-flat incidences. The first main result, stated below, follows by combining a general line--3-flat incidence bound (\cref{thm:general-codim-0-incidence-bound}) with the codimension 0 concentration bound on $\cF$, proved in \cref{sec:concentration}.

For a set $\cF$ of 3-flats in $\CC^6$, recall that we write $\cL_2(\cF)$ for the set of lines contained in at least two distinct $F,F'\in\cF$. For a set of surfaces $\cS$ where each $S\in\cS$ is contained in a unique $F\in\cF$, we write $\cL_2(\cS;\cF)$ for the set of lines $\ell$ contained in at least two distinct $F,F'\in\cF$ such that there exist $S,S'\in\cS$ with $\ell\subset S\subset F$ and $\ell\subset S'\subset F'$. Also recall that $\deg(\cS)$ refers to the sum of the degrees of the surfaces in $\cS$.

\begin{theorem}
\label{thm:codim-0-incidence-bound}
Let $\cP\subset\R^3$ be a set of size $N$ such that $|\cP\cap\gamma| \leq N^{2/3}\deg \gamma$ for every curve $\gamma\subset\CC^3$. For $H\subseteq\cP^2$ define $\cF=\{F_{pq}^+:(p,q)\in H\}$, a set of physical 3-flats in positive space, $E^+_o(\CC)$. Then there exist $r\ge 1$ and $\cS_1,\ldots,\cS_r$, sets of irreducible surfaces contained in $\cF$, such that for each $i\in[r]$ and each $F\in\cF$, the total degree of the surfaces $S\in \cS_i$ that lie in $F$ is at most $\lesssim N^{2/3}$. Furthermore,
\[\sum_{i=1}^r\deg(\cS_i)\lesssim N^{2/3}|H|\log N\]
and
\[\abs{\cL_2(\cF)\setminus\bigcup_{i=1}^r\cL_2(\cS_i;\cF)}\lesssim N^{4/3}|H|\log^2N.\]
The analogous result holds in $E^-_o(\CC)$.
\end{theorem}

\subsection{Preliminaries and proof strategy}

Let $\F$ be an algebraically closed field. We will apply the incidence bounds in this section with $\F=\CC$, but the same proofs go through for any algebraically closed $\F$.

\begin{definition}
A set $\cF$ of $t$-flats is \emph{$s$-transverse} if $\dim(F\cap F')\leq s$ for all distinct $F, F'\in\cF$.
\end{definition}

\cref{thm:codim-0-incidence-bound} is about 1-transverse 3-flats, though many of the results in this section are more natural to state in the more general setting of $s$-transverse $t$-flats.

\begin{definition}
For a set of flats $\cF$ and a polynomial $f$, define the partition $\cF=\cF_{\ss f}\sqcup\cF_{\ns f}$ where $\cF_{\ss f}$ and $\cF_{\ns f}$ are the sets of flats contained and not contained in $Z(f)$, respectively.
\end{definition}

We employ a surprisingly simple proof strategy to bound $|\cL_2(\cF)|$. We prove the desired bound by induction on $|\cF|$. For the inductive step, the goal is to find a low-degree polynomial $f$ so that 
\[0<|\cF_{\ns f}|\leq \paren{1- c}|\cF|\qquad\qquad\text{and}\qquad\qquad c|\cF|\leq|\cF_{\ss f}|<|\cF|.\]
In other words, a bounded fraction of the problem lies outside of $Z(f)$, while the rest of the problem lies inside $Z(f)$. Now consider a line $\ell\in\cL_2(\cF)$. If $\ell\not\subseteq Z(f)$, then $\ell$ is the intersection of two flats $F,F'\in\cF_{\ns f}$. However, if $\ell\subseteq Z(f)$, then the two flats $F,F'$ can lie in either $\cF_{\ss f}$ or $\cF_{\ns f}$. In the case where, e.g., $F\in\cF_{\ns f}$, then $\ell$ lies in $F\cap Z(f)$, a collection of surfaces. Thus we can write
\[|\cL_2(\cF)|\leq|\cL_2(\cF_{\ns f})|+|\cL_2(\cF_{\ss f}\cup\cS)|,\]
where $\cS$ is the set of irreducible components of $F\cap Z(f)$ for $F\in\cF_{\ns f}$. We apply the inductive hypothesis to bound the first term while we iterate this argument on the second term.\footnote{To formally run this argument, we induct on a stronger statement which involves terms such as $\cL_2(\cF_{\ss f}\cup\cS;\cF)$ that we have not defined yet.} This iteration incurs a loss in that it produces many surfaces $\cS$, but it creates a gain in that the inductive hypothesis is only ever applied to collections of flats of size at most $(1-c)|\cF|$. Balancing these two, the induction closes with only polylogarithmic loss.

The main difficulty in carrying out this proof strategy is finding the polynomial $f$. Parameter counting (\cref{cor:Walsh-param-counting}) easily produces a low-degree polynomial $f$ so that $|\cF_{\ss f}|\geq c|\cF|$. However, it is surprisingly hard to ensure that $f$ does not vanish on all of $\cF$. To do this, we use the following approach, inspired by a similar argument of Walsh \cite{Wal23}. Suppose for contradiction that every low-degree polynomial vanishing on a $c$-fraction of $\cF$ actually vanishes on all of $\cF$. This implies that $\cF$ lies inside some low-degree hypersurface $Z(f_1)$; furthermore, one can also show that $f_1$ is irreducible. Indeed, if not, one of the components of $Z(f_1)$ would contain at least half, but not all of $\cF$.

Continuing this strategy, one can find a codimension 2 variety $Z(f_1,f_2)$ which contains $\cF$. In this case, we will not be able to prove that $Z(f_1,f_2)$ is irreducible, but it will have a property which is close enough for our purposes which we call \emph{combinatorially irreducibility}. Proceeding for $(n-t)$ steps, we find a $t$-dimensional variety contained in $Z(f_1,\ldots,f_{n-t})$ which contains $\cF$ and is combinatorially irreducible. As $\cF$ is a collection of $t$-flats, each flat in $\cF$ must be an irreducible component of the variety. On the other hand, an easy consequence of combinatorial irreducibility will be that the variety contains $\lesssim_n 1$ components. This contradiction proves the desired result: that there exists a low-degree polynomial $f$, vanishing a $c$-fraction of $\cF$ without vanishing on all of $\cF$.

\begin{definition}\label{def:beta-comb-irreducible}
A (possibly reducible) variety $W\subseteq \F^n$ of pure dimension $d$ is \emph{$\beta$-combinatorially irreducible} if there exist polynomials $f_1,\ldots,f_{n-d}$ so that every irreducible component $W_i$ of $W$ is an irreducible component of $Z(f_1,\ldots,f_{n-d})$ and satisfies $\deg(W_i)\geq\beta\deg(f_1)\deg(f_2)\cdots\deg(f_{n-d})$.
\end{definition}

By B\'ezout's theorem, the sum of the degrees of the irreducible components of $Z(f_1,\ldots,f_{n-d})$ is at most $\deg(f_1)\deg(f_2)\cdots\deg(f_{n-d})$ (see, e.g., \cite[Example 8.4.6]{Ful98} for this version of B\'ezout's theorem). This immediately implies that a $\beta$-combinatorially irreducible variety has at most $\beta^{-1}$ irreducible components.

\begin{lemma}
\label{thm:parameter-counting-cor-simple}
Suppose $W$ is a $\beta$-combinatorially irreducible variety of dimension $d$ in $\F^n$, defined by polynomials $f_1,\ldots, f_{n-d}$ of degree at most $D$. Let $X$ be a variety of pure dimension $\ell<d$. If $\deg(X)\leq\rho\deg(W)D^{d-\ell}$ for some $\rho\in(0,1]$, then there exists a polynomial of degree at most $\lesssim_{n} \rho^{1/(n-\ell)}\beta^{-1}D$ that vanishes identically on $X$ without vanishing identically on any irreducible component of $W$.
\end{lemma}
\begin{proof}
We do not need the full strength of combinatorial irreducibility. We just use the fact that each component $W_i$ of $W$ satisfies $\deg(W_i)\geq\beta\deg(W)$ and $\delta(W_i)\leq D$, both of which are immediate from the definition.

We apply \cref{cor:Walsh-param-counting} with $d_1=d_2=\cdots=d_{n-d}=D$. This produces the desired polynomial $g$ which vanishes identically on $X$ without vanishing identically on any irreducible component of $W$ and satisfies
\begin{align*}
\deg(g)
&\lesssim_n\max_{0\leq s \leq n-d}\paren{\frac{\deg(X)D^{n-d-s}}{\beta\deg(W)}}^{\frac1{n-s-\ell}}\\
&\leq\max_{0\leq s \leq n-d}\paren{\frac{\rho\deg(W)D^{n-s-\ell}}{\beta\deg(W)}}^{\frac1{n-s-\ell}}\\
&\leq \rho^{1/(n-\ell)}\beta^{-1}D.\qedhere
\end{align*}
\end{proof}

\subsection{The incidence bound}
\label{ssec:proof-of-codim-0-incidence-bound}

First we show that any set of flats has a good partition: a low-degree polynomial that contains a positive fraction, but not all of the flats. We will use this result for $t=3$ and $n=6$.

\begin{theorem}
\label{thm:good-partitions-exist}
Given $n\geq 2$, there exists a constant $K\geq 1$ such that the following holds for all $t < n$ and $D\ge 1$. Let $\cF$ be a set of $t$-flats in $\FF^n$ satisfying $\cD_k(\cF)\leq D^{k-t}$ for all $t<k\leq n$. If $|\cF|\geq K$, there exists a polynomial $f$ with $\deg f\leq DK$ such that
\[0<|\cF_{\ns f}|\leq \paren{1-\tfrac1{2n}}|\cF|.\]
\end{theorem}

This proof essentially follows from techniques developed by Walsh \cite[Proposition 4.4]{Wal23}.

\begin{proof}
Suppose for contradiction that no such polynomial $f$ exists. We will construct the following data: sets $\cF=\cF_0\supseteq \cF_{1}\supseteq\cdots\supseteq \cF_{n-t}$, polynomials $g_1,\ldots,g_{n-t}$ and $h_1,\ldots,h_{n-t}$, and varieties $\F^n=W_0\supset W_{1}\supset \cdots\supset W_{n-t}$. We will also define some constants
\begin{align*}
1& = \rho_0=B_0^{-1}=\sigma_0=\beta_0\\
&\gg \rho_{1}\gg B_{1}^{-1}\gg\sigma_{1}\gg\beta_{1}\\
&\gg \rho_{2}\gg B_{2}^{-1}\gg\sigma_{2}\gg\beta_{2}\\
&\gg\cdots\\
&\gg \rho_{n-t}\gg B_{n-t}^{-1}\gg\sigma_{n-t}\gg\beta_{n-t}
\\&\gg K^{-1}>0
\end{align*}
which only depend on $n$. Here the notation $\gg$ means that each subsequent parameter is upper bounded by an appropriate function of the preceding parameter (in fact, the dependence will be polynomial). 
The data will have the following properties for all $1\leq i\leq n-t$:
\begin{enumerate}
    \item $g_{i}$ is a minimal degree polynomial vanishing identically on $\cF_{i-1}$ without vanishing identically on any component of $W_{i-1}$;
    \item $\deg(g_{i})\leq B_{i}D$ and $B_{i}\deg(g_{i})\geq\max_{1\leq j\leq i}\deg(g_j)$;
    \item $\cF_{i}=(\cF_{i-1})_{\ns h_i}$;
    \item $|\cF_{i-1}\setminus\cF_{i}|\leq\tfrac1{2n}|\cF|$;
    \item $W_{i}$ is the set of irreducible components of $Z(g_1,\ldots,g_{i})$ that contain an element of $\cF_{i}$;
    \item $W_{i}$ is a variety of pure dimension $n-i$ and is $\beta_{i}$-combinatorially irreducible.
\end{enumerate}
For all $1\leq i\leq n-t$, property (1) implies that $\cF_{i-1}$ lies in $Z(g_1,\ldots,g_i)$; with property (5) this implies that $\cF_i$ lies in $W_i$. From property (4), we immediately have $|\cF_i|\geq|\cF|/2$.

Note that this suffices to produce a contradiction: $W_{n-t}$ is $t$-dimensional and has at most $\beta_{n-t}^{-1}$ irreducible components, yet contains $\cF_{n-t}$ which consists of $|\cF_{n-t}|\geq|\cF|/2\geq K/2$ $t$-flats. Taking $K>2\beta_{n-t}^{-1}$ gives the desired contradiction.

Suppose that we have defined $\cF_i,g_i,h_i,W_i$ for all $i<k$ so that (1-6) hold for all $i<k$. We now construct $\cF_k,g_k,h_k,W_k$ so that (1-6) hold for $i=k$. 

Define $g_{k}$ as in (1). To upper bound $\deg(g_{k})$, we apply \cref{thm:parameter-counting-cor-simple}. By hypothesis, $\deg(\cF_{k-1})\leq \deg(W_{k-1})\cD_{n-k+1}(\cF)\leq \deg(W_{k-1})D^{n-k+1-t}$. Furthermore, by (2), we know \[\max_{1\leq j\leq k-1}\deg(g_j)\leq B_{k-1}\deg(g_{k-1})\leq B_{k-1}^2D.\] We apply \cref{thm:parameter-counting-cor-simple} to $W_{k-1},\cF_{k-1}$ with parameters $(\beta,\rho,D)=(\beta_{k-1},1,B_{k-1}^2D)$, giving the bound
\[\deg(g_{k})\lesssim_n \beta_{k-1}^{-1}B_{k-1}^2D.\]
This gives the first part of (2), taking $B_k\gg_n \beta_{k-1}^{-1} B_{k-1}^2$.

Define $X=W_{k-1}\cap Z(g_{k})$, a variety of pure dimension $n-k$. (As $W_{k-1}$ has pure dimension $n-k+1$ and $g_k$ does not vanish identically on any irreducible component of $W_{k-1}$, this follows from Krull's Hauptidealsatz.) For $k=1$, the second part of (2) is true for any $B\geq 1$. For $k>1$, we will show the bound
\begin{equation}    
\label{eq:deg-X-lower}
\deg(X)\geq\rho_{k}\deg(W_{k-1})\deg(g_{k-1}).
\end{equation}
B\'ezout's theorem gives the bound $\deg(W_{k-1})\leq\deg(W_{k-2})\deg(g_{k-1})$ as (5) implies that every component of $W_{k-1}$ is a component of $Z(g_1,\ldots,g_{k-2})\cap Z(g_{k-1})$ which contains an element of $\cF_{k-1}\subseteq\cF_{k-2}$ and thus is a component of $W_{k-2}\cap Z(g_{k-1})$. For the sake of contradiction, suppose that \cref{eq:deg-X-lower} fails. Combined with the previous bound, we have $\deg(X)<\rho_{k}\deg(W_{k-2})\deg(g_{k-1})^2$.

Now we apply \cref{thm:parameter-counting-cor-simple} to $W_{k-2},X$ with parameters $(\beta,\rho,D)=(\beta_{k-2},\rho_k,B_{k-1}\deg(g_{k-1}))$. Since (2) gives $\max_{1\leq j\leq k-2}\deg(g_j)\leq B_{k-1}\deg(g_{k-1})$, \cref{thm:parameter-counting-cor-simple} applies to produce a polynomial $f$ that vanishes on $X$ without vanishing identically on any irreducible component of $W_{k-2}$ and which satisfies
\[\deg(f)\lesssim_n \rho_k^{1/k}\beta_{k-2}^{-1}B_{k-1}\deg(g_{k-1}).\]
Choosing $\rho_k\ll_n \beta_{k-2}^k B_{k-1}^{-k}$, we see $\deg(f)<\deg(g_{k-1})\leq B_{k-1}D\leq DK$. Since $f$ vanishes on $X$, it vanishes on all of $\cF_{k-1}$. By hypothesis (and the fact that $|\cF_{k-1}|\geq|\cF|/2\geq|\cF|/2n$) we conclude that $f$ vanishes on all of $\cF\supseteq\cF_{k-2}$. However, this contradicts the minimality of $g_{k-1}$. This contradiction proves \cref{eq:deg-X-lower}.

B\'ezout's theorem and \cref{eq:deg-X-lower} give $\rho_{k}\deg(W_{k-1})\deg(g_{k-1})\leq \deg(X)\leq\deg(W_{k-1})\deg(g_k)$. Together with property (2) for $i=k-1$, this implies that
\[\max_{1\leq j\leq k-1}\deg(g_j)\leq B_{k-1}\deg(g_{k-1})\leq B_{k-1}\rho_k^{-1}\deg(g_k),\] proving the second part of property (2) for $i=k$, taking $B_{k}\geq B_{k-1}\rho_{k}^{-1}$.

To continue, we use the following lemma. (The proof appears next in this subsection.)
\begin{lemma}
\label{lem:going-down-a-level}
    Let $\cF$ be a set of $t$-flats contained in a $\beta$-combinatorially irreducible variety $W$ of pure dimension $d$, and let $g$ be a polynomial with minimal degree so that $g$ vanishes identically on $\cF$ but not identically on any irreducible component of $W$.
    Let $B\geq 1$ be such that $W$ is defined (in the sense of \cref{def:beta-comb-irreducible}) by polynomials $f_1,\ldots, f_{n-d}$ of degree at most $B\deg(g)$.
    Then there exist $\sigma\in(0,1]$ depending only on $\beta,B,n$ and a polynomial $h$ only depending on $W,g$ so that the following hold for any $T \le \abs{\cF}/2$:
    \begin{itemize}
        \item $\deg(h)<3\sigma^{-1}\deg(g)$;
        \item any irreducible component of $W\cap Z(g)$ which is not contained in $Z(h)$ has degree at least $\sigma\deg(W)\deg(g)$; and
        \item either $\abs{\cF_{\subseteq h}}<T$ or there exists $h'$ dividing $h$ such that $T\leq \abs{\cF_{\subseteq h'}}<\abs{\cF}$.
    \end{itemize}
\end{lemma}

We apply this lemma to the set of $t$-flats $\cF_{k-1}$, the $\beta_{k-1}$-combinatorially irreducible variety $W_{k-1}$, and the polynomial $g_{k}$ with the parameters $B=B_k\rho_k^{-1}$ and $T=\abs{\cF}/2n$.
This produces a parameter $\sigma_{k}$ and a polynomial $h_k$.
Since $\deg(h_k) <3\sigma_{k}^{-1}\deg(g_{k})\leq DK$ as long as $K$ is sufficiently large with respect to all other parameters, by hypothesis we see that $T\leq \abs{(\cF_{k-1})_{\subseteq h'}}<\abs{\cF_{k-1}}\leq\abs{\cF}$ is not possible for any $h'$ with $\deg(h')\leq\deg(h_k)$. Thus we conclude that $\abs{(\cF_{k-1})_{\subseteq h_{k}}} < \abs{\cF}/2n$.
Defining $\cF_{k}$ as in (3), we see that $\abs{\cF_{k-1}\setminus \cF_{k}}=\abs{(\cF_{k-1})_{\subseteq h_k}}<\abs{\cF}/2n$, proving (4).

Define $W_{k}$ as in (5). We now check (6). By property (5) for $i=k$, each irreducible component $U$ of $W_{k}$ is an irreducible component of $Z(g_1,\ldots,g_{k-1})\cap Z(g_k)$ which contains a flat $F\in\cF_k\subseteq\cF_{k-1}$. Thus by property (5) for $i=k-1$, we see that $U$ is contained in some irreducible component $U'$ of $Z(g_1,\ldots,g_{k-1})$ which also contains $F$, i.e., $U'$ is a component of $W_{k-1}$. As $g_k$ does not vanish identically on any irreducible component of $W_{k-1}$ and $W_{k-1}$ has pure dimension $n-k+1$, we conclude that $\dim U=n-k$. Furthermore, as $F\in\cF_k$ lies in $U$, we know that $h_k$ does not vanish identically on $F$ and thus also does not vanish identically on $U$. Thus, \cref{lem:going-down-a-level} implies that
\[\deg(U)\geq \sigma_{k}\deg(W_{k-1})\deg(g_{k})\geq\sigma_{k}\beta_{k-1}\deg(g_1)\deg(g_2)\cdots\deg(g_{k}),\]
where the final inequality is the definition of $\beta_{k-1}$-combinatorially irreducible. Thus we see that $W_{k}$ has pure dimension $n-k$ and is $\beta_{k}$-combinatorially irreducible for $\beta_{k}\leq\sigma_{k}\beta_{k-1}$.

Thus we have constructed $\cF_k,g_k,h_k,W_k$ satisfying properties (1-6) for $i=k$. Therefore we can reach the desired contradiction.
\end{proof}

\begin{proof}[Proof of \cref{lem:going-down-a-level}]
Let $X = W\cap Z(g)$, and let $X_1,\ldots,X_s$ be the irreducible components of $X$, ordered in decreasing order of degree. 
Also let $\sigma$ be some sufficiently small parameter to be chosen later.
Then pick $r$ such that $\deg(X_i)\geq\sigma\deg(W)\deg(g)$ if and only if $i\leq r$. We write $X=X^{(1)}\cup X^{(2)}$ where $X^{(1)}=X_1\cup\cdots\cup X_r$ and $X^{(2)}=X_{r+1}\cup\cdots\cup X_s$.

Now we group the irreducible components of $X^{(2)}$ as $X^{(2)}=X_1^{(2)}\cup\cdots\cup X_m^{(2)}$ where we define $X_i^{(2)}=X_{j_{i-1}+1}\cup X_{j_{i-1}+2}\cup\cdots\cup X_{j_{i}}$ for some choice of $r=j_0<j_1<j_2<\cdots<j_{m}=s$. We do this so that $\deg(X_i^{(2)})\leq\sigma\deg(W)\deg(g)$ and $m\leq 3\sigma^{-1}$. To see why this is possible, we define $j_0,\ldots,j_m$ greedily. Then each $X_i^{(2)}$, except possibly the last, has degree at least $\sigma\deg(W)\deg(g)/2$. Together with the fact that $\deg(X^{(2)})\leq\deg(X)\leq\deg(W)\deg(g)$, this implies that $m\leq2\sigma^{-1}+1\leq3\sigma^{-1}$.

Let $h_{i}$ be a minimal degree polynomial that vanishes identically on $X_i^{(2)}$ without vanishing identically on any irreducible component of $W$. 
Now we apply \cref{thm:parameter-counting-cor-simple} to $W, X_i^{(2)}$ with parameters $(\beta,\rho,D)=(\beta,\sigma,B\deg(g))$ we have
\[\deg(h_{i})\lesssim_n \sigma^{\frac{1}{n-d+1}}\beta^{-1}B \deg (g)\]
and so $\deg(h_i)<\deg(g)$ as long as $\sigma$ is chosen sufficiently small in terms of $\beta, B,n$.

Define $h=\prod_{i=1}^m h_{i}$. Then we have $\deg(h) < 3\sigma^{-1}\deg(g)$.
Moreover, by construction, it is clear that each irreducible component of $W\cap Z(g)$ not in $Z(h)$ has degree at least $\sigma\deg(W)\deg(g)$.
It remains to prove the last item.

By the minimality of $g$, we see that $h_{i}$ cannot vanish identically on all of $\cF$; however it does vanish on all of $\cF_{\ss h_{i}}$. If $|\cF_{\ss h_{i}}|\geq T$ for some $i$, then we are done.
Therefore we will assume that $|\cF_{\ss h_{i}}|<T$ for every $i\in[m]$. 

Suppose there is some $q\in[m]$ with
\[\abs{\cF_{\ss\prod_{i=1}^q h_{i}}}=\abs{\bigcup_{i=1}^q\cF_{\ss h_{i}}}\geq T, \]
and consider the smallest such $q$.
The minimality of $q$ and the assumption that $|\cF_{\ss h_q}|<T$ implies that $\abs{\cF_{\subseteq h'}}\leq 2T<\abs{\cF}$ where $h' = \prod_{i=1}^{q}h_i$. This gives the desired conclusion in this case that some $q$ exists.
If no such $q$ exists, then $\abs{\cF_{\subseteq h}}<T$, which also suffices to prove the last item.
\end{proof}

Using \cref{thm:good-partitions-exist} iteratively, we can turn an incidence problem involving 3-flats into one involving surfaces. We start with a set of  1-transverse 3-flats $\cFb$. Throughout the process we remove some flats, maintaining a current set $\cF\subseteq\cFb$ of active flats and a set $\cS$ of surfaces, each produced by cutting a flat $F\in\cFb\setminus\cF$ by a polynomial of degree at most $DK$.

\begin{definition}
Let $\cFb$ be a set of 1-transverse 3-flats. Fix parameters $D,K$. For $\cF\subseteq \cFb$, a set of surfaces $\cS$ is \emph{$\cF$-good} if $\cS$ is a set of irreducible surfaces, each contained in a (unique) $F\in \cFb\setminus\cF$ such that each 3-flat $F\in\cFb\setminus\cF$ contains elements of $\cS$ with total degree at most $DK$. For $\cF=\emptyset$, we simply say that $\cS$ is \emph{good}.

Define $\cL_2(\cS\cup\cF;\cFb)$ to be the set of lines $\ell$ such that $\ell=F\cap F'$ for distinct $F,F'\in\cFb$ such that there exist $T,T'\in\cS\cup \cF$ with $\ell\subset T\subseteq F$ and $\ell\subset T'\subseteq F'$. 
\end{definition}

\begin{theorem}
\label{thm:general-codim-0-incidence-bound}
Given $n\geq 5$, there exist constants $C,K\geq 1$ such that the following holds. Let $\cFb$ be a set of 1-transverse 3-flats in $\F^n$ satisfying $\cD_k(\cFb)\leq D^{k-3}$ for all $3<k\leq n$. For $\cF\subseteq\cFb$ and a $\cF$-good set of surfaces $\cS$ there exist good sets of surfaces $\cS_1,\ldots,\cS_r$ with the following properties:
\[\sum_{i=1}^{r}\deg(\cS_i)\leq \deg(\cS)+(CD\log\abs{\cF})\abs{\cF} \]
and
\[\abs{\cL_2(\cS\cup \cF; \cFb)\setminus \bigcup_{i=1}^{r} \cL_2(\cS_i;\cFb)} \leq (CD\log \abs{\cF})\deg(\cS)+(CD\log \abs{\cF})^2|\cF|.\]
\end{theorem}

\begin{proof}
Take $K=K(n)$ as in \cref{thm:good-partitions-exist}.  We prove this result by induction on $\abs{\cF}$. Our base case is $\abs{\cF}<2K$.
For the base case, we have
\[\abs{\cL_2(\cS\cup \cF;\cFb)\setminus \cL_2(\cS;\cFb)}\leq \abs{\cF}\left(\deg(\cS)+\abs{\cF}\right)\leq 2K\deg(\cS)+2K\abs{\cF},\]
and so it suffices to take $\cS_1 = \cS$ and $C$ sufficiently large in terms of $K=K(n)$. (As $\cS$ is $\cF$-good, it certainly is good.)

Now suppose that $\abs{\cF}\geq 2K$ and the statement is proven for all smaller $\abs{\cF}$. 
We will construct sets $\cF=\cF_0\supset \cF_1\supset\cF_2\supset\cdots\supset\cF_m$ and polynomials $f_0,f_1,\ldots,f_{m-1}$ as follows. Suppose we have constructed $\cF_k$. If $|\cF_k|<K$, halt the process and set $m=k$. Otherwise, let $f_k$ be the polynomial of degree at most $DK$ produced by applying \cref{thm:good-partitions-exist} to $\cF_k$. Note that since $\cF_k\subseteq\cFb$, we know $\cD_r(\cF_k)\leq\cD_r(\cFb)\leq D^{r-3}$. Therefore the theorem applies, producing $f_k$ with $\deg(f_k)\leq DK$ such that $0<|(\cF_k)_{\ns f_k}|\leq (1-1/2n)|\cF_k|\leq (1-1/2n)|\cF|$. Then set $\cF_{k+1}=(\cF_k)_{\ss f_k}$. Since $0<|(\cF_k)_{\ns f_k}|$, the $\cF_k$ are strictly decreasing in size and the process eventually terminates with some $\cF_m$ with $|\cF_m|<K$.

We define sets $\cS_1,\ldots,\cS_m$ so that $\cS_i$ is $\cF_i$-good as follows.
Set $\cS_0=\cS$ and for $k=0,1,\ldots, m-1$, set $\cS_{k+1}=(\cS_k)_{\ss f_k}\cup(\cF_k)_{\cap f_k}$ where $(\cF_k)_{\cap f_k}$ is the set of irreducible surfaces which are a component of $F\cap Z(f_k)$ for some $F\in(\cF_k)_{\ns f_k}$. Note that since $\cS_k$ is $\cF_k$-good, we know that $\cS_{k+1}$ is $\cF_{k+1}$-good since by B\'ezout's theorem we only add $\deg(f_k)\leq DK$ surfaces in each $F\in(\cF_k)_{\ns f_k}=\cF_k\setminus\cF_{k+1}$. 

Define $\cF'_k=\cF_k\setminus\cF_{k+1}=(\cF_k)_{\ns f_k}$ and $\cS'_k=\cS_k\setminus\cS_{k+1}=(\cS_k)_{\ns f_k}$ for each $0\leq k<m$. Write $\cF'_m=\cF_m$ and $\cS'_m=\cS_m$. We claim that the sets we have produced satisfy 
\[|\cF'_k|\leq\paren{1-\tfrac1{2n}}|\cF|\]
for all $0\leq k\leq m$, and
\[\sum_{k=0}^m\deg(\cS'_k)\leq \deg(\cS)+DK|\cF|,\]
and
\[\abs{\cL_2 (\cS\cup\cF;\cFb)\setminus \bigcup_{k=0}^m \cL_2 (\cS'_k\cup \cF'_k;\cFb)}\leq DK\deg(\cS)+(DK)^2|\cF|.\]

The first is immediate for $k<m$ from the property of $f_k$ guaranteed by \cref{thm:good-partitions-exist}, and for $k=m$ it follows from $(1-\frac{1}{2n})\abs{\cF}\geq K>\abs{\cF_m}$. The second follows since we start with varieties of total degree $\deg(\cS)$ and add varieties of total degree at most $DK$ for each $F\in\cF$. For the third, consider some $\ell\in\cL_2(\cS\cup\cF;\cFb)$. This line remains an element of $\cL_2(\cS_k\cup\cF_k;\cFb)$ until one of the following occurs: either $\ell\not\subset Z(f_k)$, meaning that $\ell$ moves to $\cL_2(\cS'_k\cup\cF'_k;\cFb)$, or $\ell\subset Z(f_k)$ but it is not an element of $\cL_2(\cS_{k+1}\cup\cF_{k+1};\cFb)$. The latter means that there is some surface $S\in\cS_k\setminus\cS_{k+1}$ that contains $\ell$ (because for any $F\in\cF_k$ which contains $\ell$, either $F\in\cF_{k+1}$ or some component of $F\cap Z(f_k)$ which contains $\ell$ is in $\cS_{k+1}$). Thus $\ell$ is an irreducible component of $S\cap Z(f_k)$ for some $S\in\cS_k\setminus\cS_{k+1}=\cS'_k$. By B\'ezout's theorem, each such $s$ contributes at most $DK\deg(S)$ lines. This shows that the number of 2-rich lines that we lose is at most $DK\sum_{k=0}^m\deg(\cS'_k)$, giving the third bound.

Now we may apply induction to each pair $(\cF_k',\cS_k')$ for $k$ from $0$ to $m$.
Let $c = -\log\paren{1-\frac{1}{2n}}>0$.
For each $k$, we get a collection of good sets $\cS_{k,1}'',\ldots, \cS_{k,r_k}''$ with 
\[\sum_{i=1}^{r_k}\deg(\cS_{k,i}'')\leq \deg(\cS_k')+(CD\log\abs{\cF_k'})\abs{\cF_k'}\leq \deg(\cS_k')+CD(\log \abs{\cF}-c)\abs{\cF'_k} \]
and
\begin{align*}
    \abs{\cL_2(\cS_k'\cup \cF_k'; \cFb)\setminus \bigcup_{i=1}^{r_k} \cL_2(\cS_{k,i}'';\cFb)} \leq& (CD\log \abs{\cF_k'})\deg(\cS_k')+(CD\log \abs{\cF_k'})^2|\cF_k'|\\
    \leq &CD(\log \abs{\cF}-c)\deg(\cS_k')+C^2D^2(\log \abs{\cF}-c)^2\abs{\cF_k'}.
\end{align*}
Now by picking $C\geq K/c$, we see that
\[\sum_{k=0}^{m}\sum_{i=1}^{r_k}\deg(\cS_{k,i}'')\leq \deg(\cS)+DK\abs{\cF}+CD(\log\abs{\cF}-c)\abs{\cF}\leq \deg(\cS)+(CD\log\abs{\cF})\abs{\cF}.\]
Moreover,
\begin{align*}
    \abs{\cL_2(\cS\cup \cF;\cFb)\setminus \bigcup_{k,i}\cL_2(\cS_{k,i}'';\cFb)}
    &\leq 
    \abs{\cL_2 (\cS\cup\cF;\cFb)\setminus \bigcup_{k=0}^m \cL_2 (\cS'_k\cup \cF'_k;\cFb)} \\
    &\qquad\qquad+ \sum_k \abs{\cL_2(\cS_k'\cup \cF_k'; \cFb)\setminus \bigcup_{i=1}^{r_k} \cL_2(\cS_{k,i}'';\cFb)}
    \\ 
    & \leq  DK\deg(\cS)+(DK)^2\abs{\cF}+CD(\log\abs{\cF}-c)(\deg(\cS)+DK\abs{\cF})\\&\qquad\qquad+C^2D^2(\log\abs{\cF}-c)^2\abs{\cF}\\
    &\leq (CD\log \abs{\cF})\deg(\cS)+\left(DK+CD(\log\abs{\cF}-c)\right)^2\abs{\cF}\\
    &\leq (CD\log \abs{\cF})\deg(\cS)+(CD\log\abs{\cF})^2\abs{\cF}.
\end{align*}
Therefore we can take $\{\cS_{k,i}''\}_{0\leq k\leq m, 1\leq i\leq r_k}$ as our good sets, closing the induction.
\end{proof}

\begin{proof}[Proof of \cref{thm:codim-0-incidence-bound}]
By \cref{prop:intersection-from-same-family}, the set $\cF$ is a 1-transverse. Define $\cF_0=\{F^+_{pq}:(p,q)\in\cP^2\}$. By \cref{thm:concentration-codim-0}, there exists a constant $C_0$ so that
\[\cD_4(\cF)\leq\cD_4(\cF_0)\leq C_0N^{2/3}\]
and
\[\cD_5(\cF)\leq\cD_5(\cF_0)\leq C_0N^{4/3}.\]
Furthermore,
\[\cD_6(\cF)=|\cF|\leq N^2.\]

Now we apply \cref{thm:general-codim-0-incidence-bound} with $\cF=\cFb$ and $\cS=\emptyset$ with parameter $D=C_0N^{2/3}$. This gives the desired result.
\end{proof}

\subsection{Weak bounds on very rich partial symmetries}
\label{ssec:3-flat-3-flat-bounds}

In this subsection we will prove the following bound on very rich rigid motions. As discussed in the introduction, this bound is too weak to handle the contribution of $k$-rich orientation reversing rigid motions for $k\gg N^{2/3}$. Nevertheless, it is the starting point for the arguments in \cref{part:iii}.

\begin{theorem}
\label{thm:very-rich-bound-main}
There exists a constant $C>0$ such that the following holds.
Let $\cP\subset\R^3$ be a set of size $N$ such that $|\cP\cap\gamma| \leq N^{2/3}\deg \gamma$ for every curve $\gamma\subset\CC^3$. For $H\subseteq\cP^2$ and $k\geq CN^{2/3}\log N$, let $\cT$ be the set of orientation-reversing rigid motions such that for each $\tau\in \cT$, there are at least $k$ pairs $(p,q)\in H$ with $\tau(p)=q$. Then
\[|\cT|\lesssim \frac{N\abs{H}\log^{3/2}N}{k^{3/2}}.\]
The same holds for orientation-preserving rigid motions.
\end{theorem}

Set $\cF=\{F^+_{pq}:(p,q)\in H\}$ and $\cG=\{F^+_\tau:\tau\in \cT\}$, sets of physical and negative 3-flats in positive space, respectively. Note that each set is 1-transverse, while the pairwise intersections $F\cap G$ for $F\in\cF$ and $G\in\cG$ can have dimension 0 or 2 by \cref{prop:intersection-from-diff-family}. Combined with the codimension 0 concentration bounds on $\cF$, we can use this information to bound the number of incidences between $\cF,\cG$: that is, the number of pairs which intersect at a 2-flat. (Note however that we do not believe that this bound is tight.) In the setting of \cref{thm:very-rich-bound-main}, each element of $\cG$ is $k$-rich for $\cF$, so this suffices to bound $|\cG|=|\cT|$. 

\begin{definition}
Let $\cF,\cG$ be $(t-2)$-transverse sets of $t$-flats in $\F^n$. We write $\cI(\cF,\cG)$ for the set of $(t-1)$-flats $S=F\cap G$ for some $F\in\cF$ and $G\in\cG$. Note that $F,G$ are uniquely determined from $S$ since $\cF,\cG$ are 1-transverse. Write $I(\cF, \cG) = \abs{\cI(\cF, \cG)}$.
\end{definition}

\begin{theorem}
\label{thm:incidence-bound-codim-1}
Let $\cF,\cG$ be $(t-2)$-transverse sets of $t$-flats in $\F^n$. Suppose that $\cD_k(\cF)\leq D^{k-t}$ for all $t<k\leq n$. Then
\[I(\cF,\cG)\lesssim_n \log(|\cF||\cG|)D\paren{|\cG|+|\cG|^{\tfrac1{n-t}}|\cF|^{1-\tfrac1{n-t}}}.\]
\end{theorem}

\begin{remark}
Morally, this result is an incidence bound between $(t-1)$-flats and $t$-flats. As this is a codimension 1 incidence problem -- unlike the codimension 2 problem studied in \cref{ssec:proof-of-codim-0-incidence-bound} -- in theory it could be deduced, without the logarithmic loss, using the strategy that Walsh uses to prove \cite[Theorem 1.1]{Wal23}. As such a bound does not follow directly from any of the results stated in \cite{Wal23} (and as the logarithmic loss is small compared to the losses we will incur in \cref{part:iii}), we instead choose to prove \cref{thm:incidence-bound-codim-1} using the same strategy that we used in \cref{ssec:proof-of-codim-0-incidence-bound}.
\end{remark}

We do not believe that \cref{thm:incidence-bound-codim-1} is tight, even without the logarithmic loss. In the situation that is relevant to the distinct distances problem, we conjecture the following bound

\begin{conjecture}
\label{conj:2-flat--3-flat}
Let $\cF,\cG$ be sets of physical and negative 3-flats in positive space $E^+_o(\R)$, respectively. If $\cD_4(\cF)\leq |\cF|^{1/3}$ and $\cD_5(\cF)\leq |\cF|^{2/3}$, then
\[I(\cF,\cG)\lesssim |\cF|^{1/3}|\cG|+|\cF|^{5/6}|\cG|^{1/2}+|\cF|.\]
\end{conjecture}

\begin{remark}
This conjecture implies that for $k\gg N^{2/3}$ and any set $\cP\subset\R^3$ of size $N$ with $|\cP\cap\gamma| \leq N^{2/3}\deg \gamma$ for every curve $\gamma\subset\CC^3$, the number of $k$-rich partial symmetries of $\cP$ is bounded by $\lesssim N^{10/3}/k^2$. Such a bound is strong enough on its own to replace \cref{part:iii}.
\end{remark}

We first deduce the main result from this incidence bound.

\begin{proof}[Proof of \cref{thm:very-rich-bound-main}]
We apply \cref{thm:incidence-bound-codim-1} with $\cG = \{F^+_\tau: \tau \in \cT\}$ and $\cF = \{F^+_{pq}:(p, q) \in H\}$, both sets of 3-flats in positive space, $E^+_o(\CC)$. By \cref{prop:intersection-from-same-family}, each set is 1-transverse. By \cref{thm:concentration-codim-0}, we have 
    \[
    \cD_4(\cF) \lesssim N^{2/3}\qquad\cD_5(\cF)\lesssim N^{4/3},
    \]
and we trivially have $\cD_6(\cF) = |\cF| = |H|$. Let $\cG'\subseteq\cG$ be a subset of size $\min\{|\cG|,N^2\}$. Now $\log (\abs{\cF}\abs{\cG'}) = O(\log N)$. 
Plugging this into Theorem \ref{thm:incidence-bound-codim-1} gives
    \begin{align*}
        I(\cF,\cG') \lesssim \log N \cdot N^{2/3}\left(|\cG'| + |\cG'|^{1/3}|H|^{2/3}\right).
    \end{align*}
Since $I(\cF,\cG') \ge k |\cG'|$ and $k\geq CN^{2/3}\log N$, the first term cannot dominate for $C$ chosen sufficiently large. Thus we conclude that
    \[|\cG'|\lesssim \frac{N\abs{H}\log^{3/2}N}{k^{3/2}}\leq\frac{|H|}{C^{3/2}}.\]
In particular, for $C$ sufficiently large, we see that $|\cG'|<|H|\leq N^2$, implying that $|\cG'|=|\cG|=|\cT|\lesssim N|H|\log^{3/2}N/k^{3/2}$.
\end{proof}

To prove \cref{thm:incidence-bound-codim-1}, we use an analogous but more complicated version of \cref{thm:good-partitions-exist}.

\begin{definition}
Let $\cF,\cG$ be $(t-2)$-transverse sets of $t$-flats in $\F^n$. A flat $F\in\cF\cup\cG$ is \emph{$r$-rich for $(\cF,\cG)$} if there are at least $r$ flats $S\in\cI(\cF,\cG)$ which are contained in $F$.
\end{definition}

\begin{theorem}
\label{thm:good-partitions-exist-incidence}
Given $n>t\geq1$, there exists a constant $K\geq 1$ such that the following holds. Let $\cF,\cG$ be nonempty $(t-2)$-transverse sets of $t$-flats in $\F^n$. Suppose that $\cD_k(\cF)\leq D_0^{k-t}$ for each $t<k\leq n$. Define
\[D=D_0\cdot \min\set{1,\paren{\frac{|\cG|}{|\cF|}}^{1/(n-t)}}.\]
If each $F\in\cF\cup\cG$ is $DK$-rich for $(\cF,\cG)$, then there exists a polynomial $f$ with $\deg f\leq DK$ such that either
\[0<|\cF_{\ns f}|\leq \paren{1-\tfrac1{2n}}|\cF|\qquad\text{or}\qquad0<|\cG_{\ns f}|\leq \paren{1-\tfrac1{2n}}|\cG|.\]
\end{theorem}

We next show how to prove the incidence bound from this result.

\begin{proof}[Proof of \cref{thm:incidence-bound-codim-1}]
Let $K$ be the constant from \cref{thm:good-partitions-exist-incidence}. Fix a parameter $D_0$. We will prove that
\[I(\cF,\cG)\leq 4Kn\log(|\cF||\cG|)D_0\paren{|\cG|+|\cG|^{\tfrac1{n-t}}|\cF|^{1-\tfrac1{n-t}}}\]
holds for all $\cF,\cG$ which are $(t-2)$-transverse sets of $t$-flats in $\F^n$ that satisfy $\cD_k(\cF)\leq D_0^{k-t}$ for each $t<k\leq n$. The proof will proceed by induction on $|\cF|+|\cG|$.

Write $\cF=\cF_0$ and $\cG=\cG_0$. Set 
\[D=D_0\cdot \min\set{1,\paren{\frac{|\cG|}{|\cF|}}^{1/(n-t)}}.\]

Given $\cF_i,\cG_i$, form $\cF'_i\subseteq\cF_i$ and $\cG'_i\subseteq\cG_i$ by repeatedly removing any flat $F\in\cF'_i\cup\cG'_i$ which is not $2DK$-rich for $(\cF'_i,\cG'_i)$. Thus 
\[I(\cF_i,\cG_i)\leq I(\cF'_i,\cG'_i)+2DK(|\cF_i\setminus\cF'_i|+|\cG_i\setminus\cG'_i|).\]
If $|\cF'_i|>(1-1/2n)|\cF|$, then let $f_i$ be the polynomial produced by applying \cref{thm:good-partitions-exist-incidence} to $\cF'_i,\cG'_i$. Note that the theorem applies since $\cD_k(\cF_i')\leq\cD_k(\cF)\leq D_0^{k-t}$ for all $t<k\leq n$ and 
\[D_0\cdot \min\set{1,\paren{\frac{|\cG_i'|}{|\cF_i'|}}^{1/(n-t)}}\leq D_0\cdot \min\set{1,\paren{\frac{|\cG|}{(1-1/2n)|\cF|}}^{1/(n-t)}}\leq 2D.\]
Since each $F\in\cF_i'\cup\cG_i'$ is $2DK$-rich for $(\cF_i',\cG_i')$ \cref{thm:good-partitions-exist-incidence} applies, producing a polynomial $f_i$ with  $\deg(f_i)\leq 2DK$. (Note that if $\cG_i'$ were empty, since every $F\in\cF_i'$ is $2DK$-rich, then $\cF_i'$ would also be empty. However, we assumed that $|\cF_i'|>(1-1/2n)|\cF|$.)

Define $\cF_{i+1}=(\cF'_i)_{\ss f_i}$ and $\cG_{i+1}=(\cG'_i)_{\ss f_i}$. Then
\[I(\cF'_i,\cG'_i)\leq I(\cF_{i+1},\cG_{i+1})+I((\cF'_i)_{\ns f_i},(\cG'_i)_{\ns f_i})+2DK(|\cF'_i\setminus \cF_{i+1}|+|\cG'_i\setminus\cG_{i+1}|),\]
where the last term comes from B\'ezout's theorem. 
We know that either $0<|\cF'_i\setminus\cF_{i+1}|$ or $0<|\cG'_i\setminus\cG_{i+1}|$ holds, so $|\cF_{i+1}|+|\cG_{i+1}|<|\cF_i|+|\cG_i|$ is strictly decreasing. Thus the process eventually terminates with some $\cF'_m$ satisfying $|\cF'_m|\leq(1-1/2n)|\cF|$. Define $\tilde \cF_i=(\cF'_i)_{\ns f_i}$ for $0\leq i<m$ and $\tilde\cF_m=\cF'_m$. Define $\tilde \cG_i=(\cG_i')_{\ns f_i}$ for $0\leq i<m$ and $\tilde\cG_m=\cG_m'$. Putting everything together, we have shown that
\[I(\cF,\cG)\leq 2DK(|\cF|+|\cG|)+\sum_{i=0}^{m} I(\tilde\cF_i,\tilde\cG_i),\]
where $\tilde\cF_0,\ldots,\tilde\cF_m$ are disjoint subsets of $\cF$ and $\tilde\cG_0,\ldots\tilde\cG_m$ are disjoint subsets of $\cG$. Furthermore, for each $i$, we either have $|\tilde\cF_i|\leq(1-1/2n)|\cF|$ or $|\tilde\cG_i|\leq(1-1/2n)|\cG|$.

From the definition of $D$, the first term above may be bounded as
\[2DK(|\cF|+|\cG|)\leq 2D_0K\paren{|\cG|^{\tfrac1{n-t}}|\cF|^{1-\tfrac1{n-t}}+|\cG|}.\]

Applying the inductive hypothesis to each summand in the second term above, H\"older's inequality then gives
\begin{align*}
\sum_{i=0}^{m} I(\tilde\cF_i,\tilde\cG_i)
&\leq 4Kn\log((1-1/2n)|\cF||\cG|)D_0\sum_{i=0}^m\paren{|\tilde\cG_i|+|\tilde \cG_i|^{\tfrac1{n-t}}|\tilde\cF_i|^{1-\tfrac1{n-t}}}\\
&\leq 4Kn\log((1-1/2n)|\cF||\cG|)D_0\paren{|\cG|+|\cG|^{\tfrac1{n-t}}|\cF|^{1-\tfrac1{n-t}}}\\
&\leq (4Kn\log(|\cF||\cG|)-2K)D_0\paren{|\cG|+|\cG|^{\tfrac1{n-t}}|\cF|^{1-\tfrac1{n-t}}}.
\end{align*}
Combining these two bounds, we see that the induction closes.
\end{proof}

\begin{proof}[Proof of \cref{thm:good-partitions-exist-incidence}]
Suppose for contradiction that no such polynomial $f$ exists. We will construct the following data: sets $\cF=\cF_0\supseteq \cF_{1}\supseteq\cdots\supseteq \cF_{n-t}$ and $\cG=\cG_0\supseteq \cG_{1}\supseteq\cdots\supseteq \cG_{n-t}$, polynomials $g_1,\ldots,g_{n-t}$ and $h_1,\ldots,h_{n-t}$, and varieties $\F^n=W_0\supset W_{1}\supset \cdots\supset W_{n-t}$. We will also define some constants
\begin{align*}
1& = \rho_0=B_0^{-1}=\sigma_0=\beta_0\\
&\gg \rho_{1}\gg B_{1}^{-1}\gg\sigma_{1}\gg\beta_{1}\\
&\gg \rho_{2}\gg B_{2}^{-1}\gg\sigma_{2}\gg\beta_{2}\\
&\gg\cdots\\
&\gg \rho_{n-t}\gg B_{n-t}^{-1}\gg\sigma_{n-t}\gg\beta_{n-t}
\\&\gg K^{-1}>0
\end{align*}
which only depend on $n$. The data will have the following properties for all $1\leq i\leq n-t$:
\begin{enumerate}
    \item $g_{i}$ is a minimal degree polynomial which vanishes identically on $\cF_{i-1}$ or vanishes identically on $\cG_{i-1}$ and in either case does not vanish identically on any component of $W_{i-1}$;
    \item $\deg(g_i)\leq B_iD$ and $B_i\deg(g_i)\geq\max_{1\leq j\leq i}\deg(g_j)$;
    \item $h_i$ satisfies $\deg(h_i)\leq \tfrac1{2n}DK$ and $\cF_{i}=(\cF_{i-1})_{\ns h_i}$ and $\cG_{i}=(\cG_{i-1})_{\ns h_i}$;
    \item $|\cF_{i-1}\setminus\cF_{i}|\leq\tfrac1{2n}|\cF|$ and $|\cG_{i-1}\setminus\cG_{i}|\leq\tfrac1{2n}|\cG|$;
    \item $W_{i}$ is the set of irreducible components of $Z(g_1,\ldots,g_i)$ that contain an element of $\cF_{i}\cup \cG_{i}$;
    \item $W_{i}$ is a variety of pure dimension $n-i$ and is $\beta_{i}$-combinatorially irreducible;
\end{enumerate}
By property (1), we see that $\cF_{i-1}\cup\cG_{i-1}$ lies in $Z(g_1,\ldots,g_i)$, while property (5) implies that $\cF_i\cup\cG_i$ lies in $W_i$. By property (4), we have $|\cF_i|\geq|\cF|/2$ and $|\cG_i|\geq|\cG|/2$. Furthermore property (3) implies that each $F\in\cF_i\cup\cG_i$ is $DK/2$-rich for $(\cF_i,\cG_i)$. Indeed, for $F\in\cF_i$, there are at least $DK$ flats $S\in\cI(\cF,\cG)$ which lie in $F$. Each such $S$ is of the form $F\cap G$ for some $G\in\cG$. Now by (3), if $S\not\subseteq Z(h_1)\cup\cdots\cup Z(h_i)$, then $G\in\cG_i$. However, $F\cap(Z(h_1)\cup\cdots \cup Z(h_i))$ contains at most $\deg(h_1)+\cdots+\deg(h_i)\leq (i/2n)DK\leq DK/2$ of these $(t-1)$-flats, implying that $F$ remains $DK/2$-rich for $(\cF_i,\cG_i)$.

Note that this suffices to produce a contradiction since $W_{n-t}$ is $t$-dimensional, has at most $\beta_{n-t}^{-1}$ irreducible components, yet contains $\cF_{n-t}\cup \cG_{n-t}$. However $\cF_{n-t}$ is nonempty and each $F\in\cF_{n-t}$ is $DK/2$-rich, implying that $|\cG_{n-t}|\geq DK/2$ Taking $DK/2\geq K/2>\beta_{n-t}^{-1}$ gives the desired contradiction.

Suppose that we have defined $\cF_i,\cG_i,g_i,h_i,W_i$ for all $i<k$ so that (1-6) hold for all $i<k$. We now construct $\cF_k,\cG_k,g_k,h_k,W_k$ so that (1-6) hold for $i=k$. 

Define $g_k$ as in (1). To upper bound $\deg(g_k)$, we use the following claim.

\begin{claim}
\label{thm:parameter-counting-more-complicated}
Let $\cX$ be a set of $t$-flats with $|\cX|\leq|\cF|$. There exists a polynomial $f$ with
\[\deg(f)\lesssim_{n}\beta_{k-1}^{-1}B_{k-1}D_0\paren{\frac{|\cX|}{|\cF|}}^{\tfrac1{n-t}}\]
that vanishes identically on $\cX$ without vanishing identically on any irreducible component of $W_{k-1}$.
\end{claim}

\begin{proof}
Let $W'$ be an irreducible component of $W_{k-1}$. For $1\leq r\leq k-1$, by property (5) every irreducible component of $Z(g_1,\ldots,g_r)$ which contains $W'$ is an irreducible component of $W_r$. By property (6), we know that $W_r$ has pure dimension $n-r$. Thus the definition of partial degree implies that 
\[\delta_r(W')\leq \max_{1\leq i\leq r}\deg(g_i)\leq B_{k-1}\deg(g_r)\]
for all $1\leq r\leq k-1$.

Since $\cF_r$ is contained in $W_r$, a variety of pure dimension $n-r$, we see
\[|\cF|/2\leq|\cF_{r}|\leq \deg(W_{r})\cD_{n-r}(\cF)\leq\deg(g_1)\cdots\deg(g_r)D_0^{n-r-t}\]
by our hypothesis on the algebraic concentration of $\cF$.

Now by \cref{cor:Walsh-param-counting}, there exists a polynomial $f$ that vanishes identically on $\cX$ without vanishing identically on any irreducible component of $W_{k-1}$ with
\begin{align*}
\deg(f)&\lesssim_n \max_{0\leq r\leq k-1}\paren{\frac{|\cX|\prod_{i=r+1}^{k-1} B_{k-1}\deg(g_i)}{\beta_{k-1}\deg(g_1)\cdots\deg(g_{k-1})}}^{\tfrac1{n-r-t}}\\
&\leq\beta_{k-1}^{-1} B_{k-1}\max_{0\leq r\leq k-1}\paren{\frac{|\cX|}{\prod_{i=1}^r\deg(g_i)}}^{\tfrac1{n-r-t}}\\
&\leq \beta_{k-1}^{-1} B_{k-1}\max_{0\leq r\leq k-1}\paren{\frac{2D_0^{n-r-t}|\cX|}{|\cF|}}^{\tfrac1{n-r-t}}\\
&\lesssim\beta_{k-1}^{-1} B_{k-1}D_0\paren{\frac{|\cX|}{|\cF|}}^{\tfrac1{n-t}}.\qedhere
\end{align*}
\end{proof}

Applying this claim to the smaller of $\cF_{k-1},\cG_{k-1}$, we find a polynomial $f$ vanishing identically on one of $\cF_{k-1},\cG_{k-1}$. Furthermore, taking $B_k$ sufficiently large in terms of $\beta_{k-1},B_{k-1}$, we see that $\deg(f)\leq B_kD$. By the minimality of $g_k$, this gives the first part of (2).

We claim that $g_k$ actually vanishes on $\cF_{k-1}\cup\cG_{k-1}$. We know that it vanishes on one of these sets, say $\cF_{k-1}$. Then if there exists $G\in\cG_{k-1}$ that is not contained in $Z(g_{k})$, since $\cF_{k-1}$ is contained in this hypersurface, we know that $G$ is at most $\deg(g_k)$-rich for $(\cF_{k-1},\cG_{k-1})$. As $\deg(g_k)\leq B_kD<DK/2$, this contradicts the fact that $G$ is $DK/2$-rich. The same argument applies in the other case.

Now we claim that the second part of (2) holds for $\rho_k, B_k$ chosen appropriately. Indeed, this follows from the identical argument used in the proof of \cref{thm:good-partitions-exist} since $g_k$ vanishes identically on $\cF_k$ while no polynomial of degree strictly smaller than $g_{k-1}$ can vanish on $\cF_{k-2}$ without vanishing identically on any component of $W_{k-2}$.

Next, we construct $h_k$ by applying \cref{lem:going-down-a-level}. Note that $\cF_{k-1}\cup\cG_{k-1}$ lies in $W_{k-1}\cap Z(g_k)$ and that no polynomial of lower degree vanishes identically on either $\cF_{k-1}$ or on $\cG_{k-1}$. Thus we can apply \cref{lem:going-down-a-level} to $\cF_{k-1},W_{k-1},g_k$ with parameters $(\beta,B,T)=(\beta_{k-1},B_k,|\cF|/2n)$. We can also apply \cref{lem:going-down-a-level} to $\cG_{k-1},W_{k-1},g_k$ with parameters $(\beta,B,T)=(\beta_{k-1},B_k,|\cG|/2n)$. This produces a parameter $\sigma_k$ which only depends on $\beta_{k-1},B_k,n$ and a polynomial $h_k$ which only depends on $W_{k-1},g_k$. In particular, since neither depends on $\cF_{k-1},\cG_{k-1}$, the two applications of \cref{lem:going-down-a-level} produce the same data. The $\sigma_k,h_k$ thus have the following properties:
\begin{itemize}
    \item $\deg(h_k)<3\sigma_k^{-1}\deg(g_k)$;
    \item any irreducible component of $W_{k-1}\cap Z(g_k)$ which is not contained in $Z(h_k)$ has degree at least $\sigma_k\deg(W_{k-1})\deg(g_k)$; 
    \item either $|(\cF_{k-1})_{\ss h_k}|<\tfrac1{2n}|\cF|$ or there exists $h'$ dividing $h_k$ so that $\tfrac1{2n}|\cF|\leq|(\cF_{k-1})_{\ss h'}|<|\cF_{k-1}|$; and
    \item either $|(\cG_{k-1})_{\ss h_k}|<\tfrac1{2n}|\cG|$ or there exists $h'$ dividing $h_k$ so that $\tfrac1{2n}|\cG|\leq|(\cG_{k-1})_{\ss h'}|<|\cG_{k-1}|$.
\end{itemize}

The first item implies $\deg(h_k)<3\sigma_k^{-1}B_kD\leq DK/2n$, taking $K$ sufficiently large in terms of all the other parameters. Defining $\cF_k,\cG_k$ as in (3), we see that property (3) holds for $i=k$. By assumption, every polynomial of degree at most $DK$ which vanishes on a $1/2n$-fraction of either $\cF$ or $\cG$ vanishes on all of $\cF$ or $\cG$. Thus we conclude that we are in the first case in the third and fourth items above. This implies property (4) holds for $i=k$.

Finally define $W_k$ as in (5). We now check (6). Any irreducible component $U$ of $W_k$ is an irreducible component of $Z(g_1,\ldots,g_{k-1})\cap Z(g_k)$ which contains some flat $F\in\cF_k\cup\cG_k$. This means that $U$ is contained in some irreducible component $U'$ of $Z(g_1,\ldots,g_{k-1})$ which also contains $F$. Thus $U'$ is an irreducible component of $W_{k-1}$ and as $g_k$ does not vanish identically on any irreducible component of $W_{k-1}$, we conclude that $\dim U=n-k$. Furthermore, since $F\in\cF_{k-1}\cup\cG_{k-1}$, we know that $h_k$ does not vanish identically on $F$, so $U$ is not contained in $Z(h_k)$. Thus the second item above implies \[\deg(U)\geq \sigma_k\deg(W_{k-1})\deg(g_k)\geq\sigma_k\beta_{k-1}\deg(g_1)\deg(g_2)\cdots\deg(g_k),\]
Thus we see that $W_k$ is $\beta_k$-combinatorially irreducible for $\beta_k\leq\sigma_k\beta_{k-1}$.

We have constructed $\cF_k,\cG_k,g_k,h_k,W_k$ satisfying properties (1-6) for $i=k$. Therefore this iteration produces the desired contradiction.
\end{proof}

\newpage

\part{Line-surface incidences}
\label{part:ii}
The results of \cref{part:i} reduce the problem of bounding incidences between lines and 3-flats to one of bounding incidences between lines and surfaces. Starting with a collection of 3-flats $\{F^+_{pq}:(p,q)\in H\}$ for each $(p,q)$, we have produced a collection of irreducible surfaces $\cS_{pq}$ all contained in the 3-flat $F^+_{pq}$. Our task is now essentially to bound the number of lines that lie in at least two of these surfaces. The main result of \cref{part:ii} is \cref{thm:ruled-ruled}, which accomplishes this task in the case that the surfaces are ruled.

Classical results in ruled surface theory describe the structure of lines lying in an algebraic surface. At various points in the paper we will draw on these tools; we collect them together in \cref{sec:ruled-surface-classical}. In short, a surface is either unruled, singly ruled, doubly ruled (a regulus), or infinitely ruled (a plane). The Cayley--Monge--Salmon theorem allows us to easily handle the incidences involving an unruled surface; in \cref{part:iii} we will study the incidences involving a plane. The remaining case is handled in this part by \cref{thm:ruled-ruled} where we study incidences involving two surfaces, both of which are singly or doubly ruled.

As discussed in \cref{ssec:outline-ii}, we use a novel argument to study the line--ruled-surface incidence problem. We transform the original line--ruled-surface incidence problem in positive space to a point-cone incidence problem in physical space which we then project down to a point-cone incidence problem in $\CC^3$. Using the polynomial method to cluster our point set on a low-degree collection of curves, we view our problem as an incidence problem between curves and cones in $\CC^3$. In \cref{sec:multiplicity} we introduce some tools from algebraic geometry that we need to carry out this argument, which is then presented in \cref{sec:ruled-ruled}.

\section{Classical ruled surface theory}
\label{sec:ruled-surface-classical}

The results in the section are all essentially standard. We collect them together here for the reader's convenience. Other good sources for this material from a discrete geometry perspective include \cite{Kol15,SS18}.

We write $\G(1,3)$ for the Grassmanian of lines in $\PP_{\CC}^3$. That is, there is a one-to-one correspondence between points $y\in\G(1,3)$ and lines $\ell_y\subset\PP_{\CC}^3$. Recall from \cref{ssec:plucker} that we view $\G(1,3)$ as a subset of $\PP^5$ using the Pl\"ucker embedding.

\begin{definition}
An irreducible algebraic surface $V\subset \PP_\CC^3$ is \emph{$k$-ruled} if for a generic point $x\in V$, there are exactly $k$ lines through $x$ which are contained in $V$. For a variety $V$, its \emph{Fano variety} $\Lambda_V\subset \G(1,3)$ is the collection of lines contained in $V$.
\end{definition}

In ruled surface theory, there are three classes of surfaces that play a special role: planes, cones, and reguli. A \emph{cone} is an irreducible surface $V\subset\PP^3$ which is not a plane and that has a point $p\in V$ called the \emph{cone point} so that for every $x\in V\setminus\{p\}$, the surface $V$ contains the line through $p$ and $x$. (Since we insist that planes are not cones, every cone has a unique cone point.) A \emph{regulus} is a smooth surface of degree 2. Note that in $\PP_{\CC}^3$, every surface of degree at most 2 is either a plane, a cone, or a regulus. Planes are infinitely ruled, cones are singly ruled, and reguli are doubly ruled. (The first two are obvious, the third is an easy computation that we will do shortly.)

\begin{proposition}
\label{thm:classical-ruled-surface}
For irreducible surface $V\subset\PP_\CC^3$, its Fano variety $\Lambda_V\subset\PP_{\CC}^5$ is an algebraic variety. Furthermore, exactly one of the following holds:
\begin{enumerate}
    \item $V$ is unruled (i.e., 0-ruled). In this case $\Lambda_V$ consists of at most $11(\deg V)^2$ isolated points;
    \item $V$ is singly ruled (i.e., 1-ruled). In this case $\Lambda_V$ consists of a curve of degree at most $\deg(V)$ and at most 2 isolated points;
    \item $V$ is doubly ruled (i.e., 2-ruled). In this case $V$ is a regulus and $\Lambda_V$ consists of two irreducible curves of degree 2 and no isolated points;
    \item $V$ is infinitely ruled. In this case $V$ is a plane and $\Lambda_V$ is a plane.
\end{enumerate}
Furthermore, in cases (2) and (3), $\Lambda_V$ does not contain any line.
\end{proposition}

Before we prove this result, we collect some easy facts about the Pl\"ucker embedding of $\G(1,3)\subset\PP^5$.

\begin{lemma}
\label{clm:plucker-facts}
\begin{enumerate}
    \item For a line $\ell_y\subset\PP^3$, the set of lines incident to $\ell_y$ is parametrized by the hyperplane section $H_y=\{x\in\G(1,3)\subset\PP^5:\vec x\cdot \vec y_{\ast}+\vec x_{\ast}\cdot\vec y=0\}$.
    \item The set of lines contained in a plane $H\subset\PP^3$ is parametrized by a 2-plane contained in $\G(1,3)\subset\PP^5$.
    \item No 3-plane in $\PP^5$ lies entirely in $\G(1,3)$.
    \item For any line $\ell\subset\G(1,3)$, there exist a point $p$ and a plane $H$ in $\PP^3$ so that $\{\ell_y:y\in\ell\}$ is the collection of lines through $p$ which are contained in $H$.
\end{enumerate}
\end{lemma}

\begin{proof}
For (1), suppose $\ell_x$ intersects $\ell_y$ at the point $w$. Let $u,v\neq w$ be other points on $\ell_x,\ell_y$, respectively. Then $x=[x_{01}:x_{02}:x_{03}:x_{23}:x_{31}:x_{12}]$ is a vector of $2\times 2$ minors of the matrix with rows $u,w$ and similarly $y$ is a vector of $2\times 2$ minors of the matrix with rows $v,w$. Now expanding the following determinant in terms of its $2\times 2$ minors, we conclude that
\[0=\det\begin{bmatrix} -\,u\,-\\-\,w\,-\\-\,v\,-\\-\,w\,-\end{bmatrix}=x_{01}y_{23}+x_{02}y_{31}+x_{03}y_{12}+x_{23}y_{01}+x_{31}y_{02}+x_{12}y_{03}=\vec x\cdot \vec y_{\ast}+\vec x_{\ast}\cdot\vec y,\]
as desired.

For (2) see, e.g., \cite[Exercise 6.5]{Harris92}. For completeness, one way to see this is to consider three points $u,v,w\in H$ which span $H$. The lines in $H$ are exactly those which intersect each of the lines $\ell_x,\ell_y,\ell_z$, defined to be the lines spanned by $uv,uw,vw$, respectively. By part (1), the Fano variety $\Lambda_H$ is exactly $\G(1,3)\cap H_{x}\cap H_{y}\cap H_{z}$, where the latter three terms are hyperplanes in $\PP^5$. Since $u,v,w$ are linearly independent, it easily follows that the same is true of $x,y,z$. Thus $H_{x}\cap H_{y}\cap H_{z}$ is a 2-plane in $\PP^5$. Since there is a 2-dimensional family of lines in $H$, we must have that $\G(1,3)\cap H_{x}\cap H_{y}\cap H_{z}$ is 2-dimensional so it must be equal to the plane $H_{x}\cap H_{y}\cap H_{z}$.

For (3), consider $\G(1,3)\subset\PP^5=\P(\CC^6)$. It suffices to show that no linear 4-dimensional subspace of $\CC^6$ lies in the zero set of the quadratic $\vec x\cdot \vec x_{\ast}=0$. Consider the nondegenerate bilinear form $\sang{x,y}=\vec x\cdot\vec y_{\ast}+\vec y\cdot \vec x_{\ast}$. Now if $W$ is a totally isotropic subspace of $\CC^6$ (i.e., if $\sang{x,x}=0$ for all $x\in W$) we conclude that $W\subseteq W^{\perp}$, so $\dim W\leq\dim W^{\perp}=6-\dim W$. Thus $\dim W\leq 3$ showing that $\CC^6$ contains no totally isotropic subspaces of dimension 4, so $\G(1,3)$ contains no 3-planes.

For (4) see, e.g., \cite[Exercise 6.4]{Harris92}. For completeness, one way to see this is to consider two points $x,y\in\G(1,3)$. If the line $\lambda x+\mu y$ lies in $\G(1,3)$ we have $(\lambda \vec x+\mu\vec y)\cdot(\lambda \vec x_{\ast}+\mu \vec y_{\ast})=0$, implying that $\vec x\cdot \vec y_{\ast}+\vec x_{\ast}\cdot \vec y=0$ for all $\lambda,\mu\in\CC$. By part (1), this implies that $\ell_x,\ell_y$ are incident, say at some point $w$. Let $u,v\neq w$ be other points on $\ell_x,\ell_y$, respectively. Now since the Pl\"ucker embedding is bilinear, we see that $\ell_{\lambda x+\mu y}$ is the line through $w$ and $\lambda u+\mu v$. Thus the line through $x,y$ in $\G(1,3)$ is the set of lines passing through $w$ and some point of the form $\lambda u+\mu v$, i.e., the set of lines through $\ell_x\cap\ell_y$ that are contained in the plane spanned by $\ell_x,\ell_y$. 
\end{proof}

\begin{proof}[Proof of \cref{thm:classical-ruled-surface}]
It is well-known that $\Lambda_V$ is an algebraic variety, see, e.g., \cite[Example 6.19]{Harris92}.

Parts of the following argument are easier to carry out in the affine setting. To that end, let $H\subset\PP^3$ be a 2-plane (not equal to $V$). Then consider the affine variety $V_o=V\setminus H\subset\PP^3\setminus H\cong \CC^3$. Say $V_o=Z(f)$ for some irreducible polynomial $f\in\CC[x,y,z]$ with $\deg f=\deg V$.

We break into two cases: $\deg V\geq 3$ and $\deg V\leq 2$. In the former case we will show that $V$ either falls into (1) or (2); in the latter case $V$ is either a cone, regulus, or plane and we will show that it falls into (2), (3), or (4), respectively.

If $\deg V\geq 3$, there is a polynomial $\Flec(f)$ with the properties that $\deg(\Flec(f))\leq 11\deg f-24$ and $\Flec(f)$ vanishes identically on every line $V_o$. Furthermore, if $\Flec(f)$ vanishes identically on $V_o$, then for a generic point $x\in V_o$, there is a line through $x$ which is contained in $V$. (See \cref{thm:cms-for-lines} or \cite[Theorem 2.1]{SS18}).

Now if $\Flec(f)$ does not vanish identically on $V$, then $Z(f,\Flec(f))$ is a 1-dimensional variety of degree at most $11(\deg V)^2-24\deg(V)$. Thus this subvariety, which contains all lines in $V_o$, contains at most $11(\deg V)^2-24\deg V$ lines. By B\'ezout's theorem, there are at most $\deg V$ lines in $V\setminus V_o=V\cap H$.\footnote{By being slightly more careful, one can choose $H$ so that it contains no lines of $V$.} This implies that such $V$ fall into case (1).

Next suppose $\deg V\geq 3$ but $\Flec(f)$ vanishes identically on $V$, implying that $V$ is ruled. Now \cite[Theorem 2.7(a)]{SS18} implies that $\Lambda_V$ consists of an irreducible curve and finitely many isolated points.\footnote{Note that this result is stated in $\CC^3$, but, as argued above, $V\setminus V_o$ only contains finitely many lines, so the claimed result still holds.} Furthermore, \cite[Theorem 2.7(b)]{SS18} implies that $V$ is singly ruled. Finally \cite[Theorem 2.8]{SS18} implies that $\Lambda_V$ has at most 2 isolated points.\footnote{Again this result is stated in $\CC^3$, but since $V$ contains only finitely many isolated points, we can choose an affine open $\PP^3\setminus H$ that does not contain any of the finitely many lines which correspond to these points.}

Other than bounding the degree of the pure dimension 1 part of $\Lambda_V$, this shows that such $V$ fall into case (2).

Next we handle the case $\deg V\leq 2$. This degree assumption implies that $V$ is either a plane, cone, or regulus. If $V$ is a plane, clearly it is infinitely ruled. In this case $\Lambda_V$ is also a plane by \cref{clm:plucker-facts}(2), showing that $V$ falls into case (4). If $V$ is a cone, it is ruled by definition, so by the same argument as above, we see that $V$ falls into case (2). 

Now suppose that $V$ is a regulus, i.e., a smooth degree 2 surface. Up to a linear change of coordinates, we can write $V=Z(x_0x_1-x_2x_3)$. Note that $V$ contains the lines $\ell^1_{a}=\{[x_0:x_1:x_2:x_3]:a_0x_0=a_1x_2,\,a_1x_1=a_0x_3\}$ for each $a=[a_0:a_1]\in\PP^1$ as well as the lines $\ell^2_{b}=\{[x_0:x_1:x_2:x_3]:b_0x_0=b_1x_3,\,b_1x_1=b_0x_2\}$ for each $b\in\PP^1$. Clearly each point in $V$ is passed through by two distinct lines, one from each of these families. This shows that $\Lambda_V$ contains two irreducible curves. We now claim that $\Lambda_V$ has no other points.

If not, there would be some line $\ell$ contained in $V$ which does not come from either ruling. Take $p\in\ell\subset V$. Then $p$ is passed through by three lines, $\ell^1_a,\ell^2_b,\ell$ for some $a,b$. Since $V$ is smooth, the tangent space to $V$ at $p$ is 2-dimensional, yet $T_pV\cap V$ contains the three lines $\ell^1_a,\ell^2_b,\ell$. By B\'ezout's theorem this is impossible unless $V$ is equal to $T_pV$. This contradicts the fact that $V$ is irreducible of degree 2.

Next we bound the degree of the pure dimension 1 part of $\Lambda_V$ in cases (2) and (3).

\begin{claim}
If $V$ is $k$-ruled for some finite $k>0$, then $\deg\Lambda_V=k\deg(V)$. 
\end{claim}

\begin{proof}
Since $V$ is $k$-ruled for some finite $k>0$, we know that $\Lambda_V$ is 1-dimensional. We compute its degree (i.e., the degree of its pure dimension 1 part) by counting intersections between $\Lambda_V$ and an appropriately chosen hyperplane.

Let $y=[\vec y:\vec y_\ast]$ be a generic point in $\G(1,3)$, so that $\ell_y$ is a generic line in $\PP^3_{\CC}$. Then $\abs{\ell_y\cap V}=\deg V$. Now consider the hyperplane $H_y=\{x=[\vec x:\vec x_\ast]\in\PP^5_{\CC}:\vec x\cdot \vec y_\ast+\vec x_\ast\cdot\vec y=0\}$. We claim that for $y\in\G(1,3)$ generic, $\abs{H_y\cap \Lambda_V}=\deg\Lambda_V$. By B\'ezout's theorem, this holds if we can show that for generic $y\in\G(1,3)$, each intersection $p\in H_y\cap \Lambda_V$ is transversal. In other words, that $p$ is a smooth point of $\Lambda_V$ and $H_y$ intersects the tangent space $T_p\Lambda_V$ transversally.

Since $\Lambda_V$ is 1-dimensional, it has a finite number of singular points. For each singular point $p=[\vec p:\vec p_\ast]\in \Lambda_V$, the space of $y\in\G(1,3)$ so that $H_y$ passes through $p$ has codimension 1: it is the intersection of $\G(1,3)\subset\PP^5$ with the hyperplane $H_p$ defined above. (Since $\G(1,3)$ is an (irreducible) quadric hypersurface in $\PP^5$, this intersection must be of codimension 1 inside $\G(1,3)$.) Thus for a generic $y\in\G(1,3)$, the hyperplane $H_y$ avoids all the finitely many singular points of $\Lambda_V$. 

Now for $y\in\G(1,3)$ and a smooth point $p\in\Lambda_V$, we want to understand when the intersection of $H_y$ and $\Lambda_V$ is non-transversal at $p$. This is equivalent to $H_y$ containing the tangent line to $\Lambda_V$ at $p$, which, similar to before, is a codimension 2 linear condition on $y$. Say $P_p\subset\G(1,3)$ is the 3-plane so that $y\in P_p$ if and only if $H_y$ contains the tangent line to $\Lambda_V$ at $p$. Now by \cref{clm:plucker-facts}(3), we know that $\G(1,3)$ does not contain any 3-planes, and so we see that $P_p\cap\G(1,3)$ is 2-dimensional. Thus $\bigcup_p P_p$, where the union is taken over the nonsingular points of $\Lambda_V$ is 3-dimensional, so for a generic $y\in\G(1,3)$, we conclude that $y\not\in\bigcup_p P_p$, meaning that all the intersections of $H_y$ and $\Lambda_V$ are transversal and thus $\abs{H_y\cap \Lambda_V}=\deg \Lambda_V$.

By \cref{clm:plucker-facts}(1), we know that $x\in H_y\cap\G(1,3)$ if and only if the lines $\ell_x$ and $\ell_y$ are incident. Now for a generic point $p\in V$, by hypothesis, there are exactly $k$ lines through $x$ which are contained in $V$. Since $\ell_y$ is a generic line in $\PP^3$, we conclude that every point of $\ell_y\cap V$ is of this type. Thus there are exactly $k\abs{\ell_y\cap V}$ lines contained in $V$ which are incident to $\ell_y$. (These are clearly distinct lines else they are equal to $\ell_y$. However, since $y$ is generic, $\ell_y\not\subset V$.) Thus we conclude that $k\abs{\ell_y\cap V}=\abs{H_y\cap \Lambda_V}$, implying the desired result.
\end{proof}

To prove the final claim, note that for any line $\ell\subset\G(1,3)$, the union of the lines it corresponds to, $\bigcup_{y\in\ell}\ell_y$ is a plane in $\PP^3$ by \cref{clm:plucker-facts}(4). Thus if $\ell\subset\Lambda_V$, it must be the case that $V$ contains (and thus is equal to) a plane.
\end{proof}

\begin{definition}
\label{defn:in-the-ruling}
For an irreducible variety $V\subseteq\PP^3_{\CC}$, we say that a line $\ell\subset V$ is \emph{in the ruling} of $V$ if it corresponds to a non-isolated point of $\Lambda_V$.
\end{definition}

In particular if $V$ is unruled, then no line is in the ruling; if $V$ is singly ruled, then \cref{thm:classical-ruled-surface} implies that at most two lines are not in the ruling; and if $V$ is doubly ruled or infinitely ruled, then every line is in the ruling.

\begin{corollary}
\label{cor:lines-not-in-ruling}
Let $V\subset\PP^3_{\CC}$ be an irreducible variety, then the number of lines in $V$ not in the ruling is bounded by $O((\deg V)^2)$.
\end{corollary}

We remark that this result can also be deduced from the general ruled surface theory that we developed in \cref{sec:csm-for-flats}.

\begin{proposition}
\label{thm:classical-ruled-surface-incidences}
For an irreducible surface $V\subset\PP_\CC^3$, 
\begin{enumerate}
    \item if $V$ is unruled, the number of pairwise intersecting lines in $V$ is at most $66(\deg V)^3$;
    \item if $V$ is singly ruled and not a cone, then for all but at most two lines $\ell$ in $V$, the line $\ell$ intersects at most $\deg V$ other lines in $V$;
    \item if $V$ is doubly ruled, say $\lambda_1, \lambda_2$ are the two irreducible components of $\Lambda_V$. Then for every line corresponding to a point of $\lambda_1$ and every line corresponding to a point of $\lambda_2$, the two lines intersect. Furthermore, no two lines corresponding to points of the same component $\lambda_i$ intersect; and
    \item if $V$ is a cone or a plane, then every pair of lines in $V$ intersects.
\end{enumerate}
\end{proposition}

\begin{proof}
(1) is \cite[Corollary 21(2)]{Kol15}.

For (2), if $V$ is singly ruled and not a cone, then by \cite[Proposition 55(4)]{Kol15}, all but at most two lines in $V$ have a property known as being \emph{non-special}. Then by \cite[Proposition 55(5)]{Kol15}, each non-special line intersects at most $\deg V-2$ non-special lines, and thus at most $\deg V$ lines in total.

For (3), suppose $V$ is a regulus. As in the proof of \cref{thm:classical-ruled-surface}, change coordinates so $V=Z(x_0x_1-x_2x_3)$. Now a line in $\lambda_1$ is of the form $\ell^1_{a}=\{[x_0:x_1:x_2:x_3]:a_0x_0=a_1x_2,\,a_1x_1=a_0x_3\}$ for some $a=[a_0:a_1]\in\PP^1$ and a line in $\lambda_2$ is of the form $\ell^2_{b}=\{[x_0:x_1:x_2:x_3]:b_0x_0=b_1x_3,\,b_1x_1=b_0x_2\}$ for some $b\in\PP^1$. These two lines intersect at the point $[a_1b_1:a_0b_0:a_0b_1:a_1b_0]$, as desired. Now if two lines from the same ruling, say $\ell^1_a,\ell^1_{a'}$ intersected at some point $p$, that point would also be passed through by a line of the form $\ell^2_b$. As in the proof of \cref{thm:classical-ruled-surface}, this is impossible, since $T_pV\cap V$ cannot contain all three lines $\ell^1_a,\ell^1_{a'},\ell^2_b$. This proves the second part of (3).

(4) is clearly true if $V$ is a plane. Suppose $V$ is an irreducible cone with cone point $p\in V$. If $V$ contains a line $\ell$ not through $p$, then $V$ must be the plane spanned by $p$ and $\ell$ which is forbidden by definition. Otherwise, every line in $V$ passes through $p$, again giving the desired result.
\end{proof}
\section{Multiplicities and fibers}
\label{sec:multiplicity}
In this section, we collect some tools from algebraic geometry that we will need to carry out the arguments in \cref{sec:ruled-ruled}. We will often be in the following situation: we have curves $\gamma_1,\gamma_2\subset\CC^n$ and a map $\psi\colon\gamma_1\to\gamma_2$ that is dominant, meaning that its image is Zariski-dense in $\gamma_2$. Suppose that this map is generically $d$-to-1, meaning that for a generic point $p\in\gamma_2$, the fiber of $\psi$ above $p$ has size $d$, i.e., $\abs{\psi^{-1}(p)}=d$. (It is a standard fact in algebraic geometry that any such map is generically $d$-to-1 for some finite $d$. Indeed, since $\dim \gamma_1=\dim\gamma_2=1$, the fiber dimension theorem gives that the size of a generic fiber of $\psi$ is finite. Then by \cite[Proposition 7.16]{Harris92}, there is a well-defined number $d$ so that a generic fiber of $\psi_1$ has size $d$.) The situation we find ourself in is that for a specific point $p\in\gamma_2$, we wish to bound the size of the fiber above $p$. It turns out that if $p$ is non-singular, then this fiber will still have size at most $d$. However for singular points of $\gamma_2$, this can fail.

For example, consider the map $\psi\colon \CC^1_t\to Z(y^2-x^2(x+1))\subset \CC^2_{x,y}$ defined by $\psi(t)=(t^2-1, t(t^2-1))$. (This is a normalization of the nodal singularity of a cubic curve.) Now $\psi$ is generically $1$-to-$1$, but the fiber above the singularity $(0,0)$ consists of the two points $t=\pm 1$. It turns out that this is essentially the worst that can happen: for any point $p\in\gamma_2$, the size of the fiber above $p$ is at most $d$ times the \emph{multiplicity} of $p$ in $\gamma_2$: this is a measure of how singular $\gamma_2$ is at $p$.

In this section we collect several standard tools for working with multiplicity of points on a curve, first in the general case $\gamma\subset\CC^n$ and then in the simpler planar case $\gamma\subset\CC^2$. These results are all essentially standard, though occasionally we need to do a little work to translate them into our setting. The reader who is less familiar with algebraic geometry is encouraged to skip this section and only refer to the results as needed when reading \cref{sec:ruled-ruled}.

\subsection{Multiplicity of a point on a curve}

For an irreducible curve $\gamma\subset\CC^n$ and a polynomial $f\in\CC[x_1,\ldots,x_n]$ which does not vanish identically on $\gamma$, consider $p\in \gamma\cap Z(f)$. Recall that the intersection multiplicity of $\gamma$ and $f$ at $p$ (see \cref{defn:multiplicity}) is defined by
\[m_p(\gamma,f):=m_p(I(\gamma)+\sang{f})=\len_{\CC[x_1,\ldots,x_n]_{\fm}}\paren{\CC[x_1,\ldots,x_n]_{\fm}/(I(\gamma)+\sang{f})\CC[x_1,\ldots,x_n]_{\fm}},\]
where $\fm=I(p)\subset\CC[x_1,\ldots,x_n]$ is the maximal ideal consisting of polynomials vanishing at $p$.

As $p$ is a point, it is a standard fact that this definition simplifies to 
\[m_p(\gamma,f)=\dim_{\CC}\CC[x_1,\ldots,x_n]_{\fm}/(I(\gamma)+\sang{f})\CC[x_1,\ldots,x_n]_{\fm},\]
the dimension of this ring as a $\CC$-vector space.\footnote{To see this, suppose that $0=M_0\subsetneq M_1\subsetneq\cdots\subsetneq M_\ell=M$ is a maximal chain of $\CC[x_1,\ldots,x_n]_{\fm}$-modules. By maximality, each $M_{i+1}/M_i$ is a simple $\CC[x_1,\ldots,x_n]_{\fm}$-module which must be $\CC[x_1,\ldots,x_n]_{\fm}/I$ for some maximal ideal $I$ of $\CC[x_1,\ldots,x_n]_{\fm}$  (see, e.g., \cite[Section 2.4]{Eis95}). Now as $\CC[x_1,\ldots,x_n]_{\fm}$ has a unique maximal ideal $\fm\CC[x_1,\ldots,x_n]_{\fm}$, all of these quotients must be $(\CC[x_1,\ldots,x_n]/\fm\CC[x_1,\ldots,x_n]_{\fm})_{\fm}=K(p)=\CC$, i.e., the field of rational functions on $p$. Thus in a maximal chain, each $M_i$ is a finite-dimensional $\CC$-vector space satisfying $\dim_{\CC} M_{i+1}=\dim _{\CC}M_i+1$.}

Now we define the multiplicity of a point on a curve. We give two definitions, one for curves in $\CC^n$ and in \cref{ssec:planar-multiplicity} we give an equivalent definition that is easier to work with in $\CC^2$.

\begin{definition}
\label{defn:mult-of-point-on-curve}
For an irreducible curve $\gamma\subset\CC^n$ and a point $p\in\gamma$, the \emph{multiplicity} of $p$ in $\gamma$ is the positive integer $\mu_p(\gamma)$ such that 
\[\dim_{\CC}\cO_{\gamma,p}/\fm^N = \mu_p(\gamma)N+c\]
for some integer $c\in\ZZ$ and all sufficiently large $N$.

In the above equation, $\mathcal{O}_{\gamma,p}=\CC[x_1,\ldots,x_n]_{\fm}/I(\gamma)\CC[x_1,\ldots,x_n]_{\fm}$ is the \emph{local ring} of $\gamma$ at $p$ and $\mathfrak m$ is the maximal ideal consisting of polynomials vanishing at $p$. As $\cO_{\gamma,p}$ is a 1-dimensional local ring, it is a standard fact that $\dim_{\CC}\cO_{\gamma,p}/\fm^N=\len_{\CC[x_1,\ldots,x_n]_\fm}(\cO_{\gamma,p}/\fm^N)$ is equal to a linear polynomial for all sufficiently large $N$ (see, e.g, \cite[Theorem 12.4]{Eis95}), so $\mu_p(\gamma)$ is well-defined.
\end{definition}

The following two results relate this definition to intersection multiplicity. 

\begin{lemma}\label{lem:mult-of-all-surface}
For an irreducible curve $\gamma\subset \CC^n$, a point $p\in\gamma$, and a polynomial $f$ vanishing at $p$ but not vanishing identically on $\gamma$, we have
\[m_p(\gamma,f)\geq\mu_p(\gamma).\]
\end{lemma}
\begin{proof}
Consider the map $\psi_f\colon \cO_{\gamma,p}/\fm^{N-1}\to \cO_{\gamma,p}/\fm^N$ defined by multiplication by $f$ (this is well-defined as $f\in\fm$).
The cokernel of $\psi_f$ is $\cO_{\gamma,p}/(\langle f\rangle+\fm^N)$, showing that
\[\dim_{\CC}\cO_{\gamma,p}/(\langle f\rangle+\fm^N)=\dim_{\CC}\cO_{\gamma,p}/\fm^{N}-\dim_{\CC}(\mathrm{im}(\psi_f))\geq  \dim_{\CC}\cO_{\gamma,p}/\fm^{N}-\dim_{\CC}\cO_{\gamma,p}/\fm^{N-1}.\]

The intersection multiplicity satisfies 
$m_p(\gamma,f)=\dim_{\CC}\cO_{\gamma,p}/\langle f\rangle\geq \dim_{\CC}\cO_{\gamma,p}/(\langle f\rangle+\fm^N)$ for all $N$, while the above inequality implies 
\[\dim_{\CC}\cO_{\gamma,p}/(\langle f\rangle+\fm^N)\geq (\mu_p(\gamma)N+c)-(\mu_p(\gamma)(N-1)+c) = \mu_p(\gamma),\]
for $N$ sufficiently large.
\end{proof}

Furthermore, the above result is an equality if $f$ is a generic linear function vanishing at $p$.

\begin{proposition}\label{prop:mult-alt-def}
For an irreducible curve $\gamma\subset \CC^n$ and a point $p\in\gamma$, let $f$ be a generic linear function vanishing at $p$. Then
\[m_p(\gamma,f)=\mu_p(\gamma).\]
\end{proposition}

This is a standard property of multiplicity. For completeness, we show how one can deduce this from standard results in the textbook \cite{Mat89}.

\begin{proof}
First, \cite[Formula 14.1]{Mat89} defines the notation $e(\fm, \cO_{\gamma,p})$ so that it is equal to $\mu_p(\gamma)$. Note that $\fm=\sang{x_1-p_1,\ldots, x_n-p_n}$ and let $f$ be a generic $\CC$-linear combination of $x_1-p_1,\ldots,x_n-p_n$. By \cite[Theorem 14.14]{Mat89}, setting $\mathfrak b=\sang{f}$, we see that $\mathfrak{b}$ is a \emph{reduction} of $\fm$ and $f$ is a \emph{system of parameters} of $\cO_{\gamma,p}$.
By \cite[Theorem 14.13]{Mat89}, since $\mathfrak b$ is a reduction of $\fm$, we have $\mu_p(\gamma) = e(\fm, \cO_{\gamma,p}) = e(\mathfrak{b},\cO_{\gamma,p})$.

To complete the proof, we apply \cite[Theorem 14.11]{Mat89} with $A=M=\cO_{\gamma,p}$ and $x_1=f$ and $\mathfrak{q} = \mathfrak{b}$. This result applies as $f$ is a system of parameters of $\cO_{\gamma,p}$, so we conclude that 
$e(\mathfrak{b},\cO_{\gamma,p}) = e(\sang{0}, \cO_{\gamma,p}/\langle f\rangle)$ as long as $f$ is not a zero-divisor of $\cO_{\gamma,p}$. The latter holds as  $f$ is a nonzero element of $\cO_{\gamma,p}$ which is an integral domain.

As $\cO_{\gamma,p}/\langle f\rangle$ is zero-dimensional, by \cite[Formula 14.1]{Mat89} we get that 
\[e(\sang{0}, \cO_{\gamma,p}/\langle f\rangle)=\mathrm{length}_{\CC[x_1,\ldots, x_n]_{\fm}}(\cO_{\gamma,p}/\langle f\rangle) = m_p(\gamma,f),\]
as desired.
\end{proof}

One of our main uses for multiplicity is to upper bound the size of a fiber of a map between curves.

\begin{theorem}\label{thm:upper-bound-fiber}
For irreducible curves $\gamma_1\subset \CC^n$ and $\gamma_2\subset \CC^m$, let $\psi\colon \gamma_1\to \gamma_2$ be a morphism so that the fiber above a generic point has size $d>0$. Then for every $p\in \gamma_2$, \[\abs{\psi^{-1}(p)}\leq d\cdot \mu_p(\gamma_2).\]
\end{theorem}

\begin{proof}
This follows from fairly standard tools in intersection theory (see, e.g., \cite{Ful98} for a general reference). To apply these tools, we first modify $\psi$ so that it is a so-called projective morphism. A morphism $f\colon X\to Y$ between complex varieties $X, Y$ is \emph{projective} if it factors as a closed embedding $X\hookrightarrow \PP^d_{\CC}\times Y$ composed with the projection map $\pi_2\colon \PP^d_{\CC}\times Y\twoheadrightarrow Y$ (see the definition in \cite[pg. 103]{Har77}).

To do this, first let $\Tilde{\gamma_1}$ be the closure of $\gamma_1\subset \PP_{\CC}^n$, and let $\Gamma_\psi\subset \gamma_1\times \gamma_2\subset\Tilde{\gamma_1}\times \gamma_2\subset\PP^n_{\CC}\times\gamma_2$ be the graph of $\psi$. Let $\gamma_1'$ be the closure of $\Gamma_{\psi}$ in $\Tilde{\gamma_1}\times \gamma_2$.

Finally, write $\psi'\colon\gamma_1'\to \gamma_2$ for the restriction of the projection map $\pi_2\colon\PP^n_{\CC}\times\gamma_2\twoheadrightarrow\gamma_2$ to $\gamma_1'$. Then $\psi'$ is the composition of a closed embedding (the identity map) with the projection map $\pi_2$, meaning that it is a projective morphism. Since $\psi'$ is projective it is also proper \cite[Theorem 4.9]{Har77}, in particular, the image of $\psi'$ is closed. This image contains the image of $\psi$ which is a dense subset of $\gamma_2$ by hypothesis. Thus $\psi'\colon\gamma_1'\to\gamma_2$ is a surjective proper morphism.

We now apply standard intersection theory tools. Apply \cite[Example 4.3.6]{Ful98} to the proper surjective morphism $\psi'\colon \gamma_1'\to \gamma_2$ with $X=V=p$ for an arbitrary point $p\in \gamma_2$. By hypothesis $\psi'^{-1}(p)$ is pure dimension 0, so \cite[Example 4.3.6]{Ful98} applies to give
    \[\deg(\gamma_1'/\gamma_2)(e_p\gamma_2)_p = \sum_{q\in \psi'^{-1}(p)}\deg(q/p)(e_{\psi'^{-1}(p)}\gamma_1')_q.\]
Here the notation $\deg(\gamma_1'/\gamma_2)$ and $\deg(q/p)$ refers to the size of the generic fibers (see \cite[Section 1.4]{Ful98}), so these are $d$ and 1 respectively. The notation $(e_p\gamma_2)_p$ and $(e_{\psi'^{-1}(p)}\gamma_1')_q$ refers to the multiplicity of the point on the curve (see \cite[Example 4.3.1]{Ful98}), so these are $\mu_p(\gamma_2)$ and at least one, respectively. Therefore, the above equation gives
\[d\cdot\mu_p(\gamma_2)=\sum_{q\in\psi'^{-1}(p)}(e_{\psi'^{-1}(p)}\gamma_1')_q\geq\abs{\psi'^{-1}(p)}\geq\abs{\psi^{-1}(p)}.\]
To see the last inequality, note that $\psi=\psi'\circ\iota$ where $\iota\colon\gamma_1\hookrightarrow\gamma_1'$ is the injection $(\Id_{\gamma_1},\psi)\colon \gamma_1\to \Gamma_{\psi}\subset\gamma_1'$, so the fiber of $\psi$ above $p$ injects into the fiber of $\psi'$. 
\end{proof}

\subsection{Multiplicity for planar curves}
\label{ssec:planar-multiplicity}

For a nonzero polynomial $f\in\CC[x_1,x_2]$ and a point $p\in Z(f)$, we give an alternate definition of the multiplicity of $p$ in $Z(f)$. This agrees with \cref{defn:mult-of-point-on-curve} in the case that $f$ is irreducible.\footnote{In fact it always agrees with \cref{defn:mult-of-point-on-curve} if that definition is interpreted scheme-theoretically.}

\begin{definition}
\label{defn:mult-of-point-on-curve-planar}
For a nonzero polynomial $f\in\CC[x_1,x_2]$ and a point $p\in Z(f)$, define the \emph{multiplicity} of $p$ in $f$, denoted $\mu_p(f)$, to be the degree of the lowest-degree term in the Taylor expansion of $f$ at $p$.
\end{definition}

\begin{proposition}\label{prop:planar-mult-alt-def}
For an irreducible polynomial $f\in \CC[x_1,x_2]$ and $p\in Z(f)$, \cref{defn:mult-of-point-on-curve,defn:mult-of-point-on-curve-planar} agree.
\end{proposition}

\begin{proof}
Without loss of generality, assume that $p=(0,0)$. Suppose that the multiplicity of $p$ in $f$, in the sense of \cref{defn:mult-of-point-on-curve-planar}, is $k$. By a linear change of coordinate, we may assume that the coefficient of $x_2^k$ in $f$ is equal to 1. Then write $f(x_1,x_2)=x_2^k-g(x_1,x_2)$ where all terms in $g$ have degree at least $k$ and $g$ does not have a $x_2^k$ term. 

We claim
\[\mathcal{O}_{Z(f),p}/\mathfrak m^N\cong\CC[x_1,x_2]/(\langle f\rangle+\langle x_1,x_2\rangle^N).\]
To see this, note that the left-hand side is $(\CC[x_1,x_2]/(\sang{f}+\fm^N))_{\fm}$ where $\fm=\sang{x_1,x_2}$. Thus it suffices to show that localizing at $\fm$ does not change the ring. In other words, we need to check that for each $h_1\in\CC[x_1,x_2]\setminus\fm$, there exists $h_2\in\CC[x_1,x_2]$ so that $h_1h_2\in 1+\sang{x_1,x_2}^N$. Scaling $h_1$ if necessary, assume that $h_1(0,0)=1$. Then defining
\[h_2=\sum_{i=0}^{N-1}(1-h_1)^i,\]
we see that $h_1h_2=1-(1-h_1)^N\in1+\sang{x_1,x_2}^N$, as desired. This shows that every element of $\CC[x_1,x_2]\setminus\fm$ was already a unit in $\CC[x_1,x_2]/(\langle f\rangle+\langle x_1,x_2\rangle^N)$, so the localization has no effect.

Now we claim the right-hand side has a basis (as a $\CC$-vector space)
$B_1:=\{x_1^ix_2^j:(i,j)\in \ZZ_{\geqslant 0}^2,\, j<k,\, i+j<N\}$. To see this, note that clearly $\CC[x_1,x_2]/\langle x_1,x_2\rangle^N$ has a basis consisting of $B_0:=\{x_1^ix_2^j:(i,j)\in \ZZ_{\geqslant 0}^2,\, i+j<N\}$. Furthermore, by our choice of representation $f=x_2^k-g$, we see that any nonzero multiple of $f$, when expressed in terms of the basis $B_0$ must include a monomial $x_1^ix_2^j$ with $j\geq k$ with nonzero coefficient. This shows that $B_1$ is linear independent in $\CC[x_1,x_2]/(\langle f\rangle+\langle x_1,x_2\rangle^N)$

Now one can quotient by $\sang{f}$ by repeatedly replacing any $x_1^ix_2^{j}$ term for $j\geq k$ with $x_1^ix_2^{j-k}g(x_1,x_2)$ (and removing any terms with degree at least $N$). This procedure replaces monomials with strictly larger monomials under the order where monomial are ordered lexicographically by degree in $x_1$ then degree in $x_2$. Thus this procedure terminates with a linear combination of $B_1$, showing that they also span.

When $N\geq k$, this basis has size
\[\binom{N+1}{2}-\binom{N-k+1}{2} = kN-\binom{k}{2}.\]
This shows that the multiplicity of $p$ in $f$, in the sense of \cref{defn:mult-of-point-on-curve}, is also equal to $k$.
\end{proof}

We relate this definition with intersection multiplicity as follows.

\begin{lemma}\label{lem:mult-of-two-planar-curves}
For relatively prime polynomials $f_1,f_2\in\CC[x_1,x_2]$ and a point $p\in Z(f_1)\cap Z(f_2)$, the intersection multiplicity satisfies
\[m_p(f_1,f_2)\geq \mu_p(f_1)\mu_p(f_2).\]
\end{lemma}

\begin{proof}
Let $\cO_{\CC^2,p}$ be the local ring of $\CC^2$ at $p$, let $\fm$ be the maximal ideal, and let $m_i = \mu_p(f_i)$ for $i=1,2$.
By definition, $f_i\in \fm^{m_i}$ for $i=1,2$.
Now consider the map
\[\psi_f\colon \cO_{\CC^2,p}/\fm^{m_1}\times \cO_{\CC^2,p}/\fm^{m_2}\to \cO_{\CC^2,p}/\fm^{m_1+m_2}\]
sending $(g_1,g_2)$ to $f_2g_1+f_1g_2$.
This is well defined as $f_2\in \fm^{m_2}$ and $f_1\in \fm^{m_1}$.
The cokernel of the map is $\cO_{\CC^2,p}/(\langle f_1,f_2\rangle+\fm^{m_1+m_2})$, so
\begin{align*}
\dim_{\CC}\cO_{\CC^2,p}/(\langle f_1,f_2\rangle+\fm^{m_1+m_2})
&=\dim_{\CC}\cO_{\CC^2,p}/\fm^{m_1+m_2}-\dim_{\CC}(\mathrm{im}(f_\psi))\\
&\geq \dim_{\CC}\cO_{\CC^2,p}/\fm^{m_1+m_2}-\dim_{\CC}\cO_{\CC^2,p}/\fm^{m_1}-\dim_{\CC}\cO_{\CC^2,p}/\fm^{m_2}.
\end{align*}
As $\cO_{\CC^2,p}/\fm^{m_1+m_2}=\CC[x_1,x_2]/\sang{x_1-p_1,x_2-p_2}^{m_1+m_2}$ and similarly for the other two terms, we see that the right-hand side evaluates to
\[\binom{m_1+m_2+1}{2}-\binom{m_1+1}{2}-\binom{m_2+1}{2}=m_1m_2.\]
To conclude, $m_p(f_1,f_2)=\dim_{\CC}\cO_{\CC^2,p}/\langle f_1,f_2\rangle\geq\dim_{\CC}\cO_{\CC^2,p}/(\langle f_1,f_2\rangle+\fm^{m_1+m_2})=m_1m_2$.
\end{proof}

Finally, we can bound the total multiplicity on a curve using B\'ezout's theorem.

\begin{theorem}[B\'ezout's theorem]
\label{thm:bezout-with-multiplicity-projective}
Let $\gamma\subset\PP^n_{\CC}$ be an irreducible curve and $f$ be a homogeneous polynomial that does not vanish identically on $\gamma$. Then
\[\sum_{p\in\gamma\cap Z(f)}m_p(\gamma,f)=\deg\gamma\deg f.\]
\end{theorem}
This is the standard form of B\'ezout's theorem (with multiplicity in projective space over an algebraically closed field), see, e.g., \cite[Chapter I, Theorem 7.7]{Har77}. (Note that to apply the definition of intersection multiplicity in projective space, simply restrict to an affine open set containing $p$.)

\begin{proposition}\label{prop:planar-easy-multiplicity}
Let $\gamma\subset \CC^2$ be an irreducible curve.
Then
\[\sum_{p\in \gamma}\mu_p(\gamma)(\mu_p(\gamma)-1)\leq (\deg \gamma)(\deg \gamma -1).\]
\end{proposition}
\begin{proof}
Let $\gamma = Z(f)$ where $f\in \CC[x_1,x_2]$ is an irreducible polynomial.
Define $f'=\alpha \partial_1 f+\beta \partial_2 f$ for generic $\alpha,\beta\in\CC$. As at least one of $\partial_1f,\partial_2f$ is not the zero polynomial, we conclude the same for $f'$. Since $f$ is irreducible, we see that $f,f'$ are relatively prime.

By definition, for $p\in\gamma$ with $\mu_p(\gamma)>1$, the lowest-degree term in the Taylor expansion of $f$ at $p$ has degree $\mu_p(\gamma)>1$.
Therefore the lowest-degree term in the Taylor expansion of $f'$ at $p$ has degree $\mu_p(\gamma)-1$. (In particular, this means that if $\mu_p(\gamma)>1$, then $p\in Z(f')$.)
By \cref{lem:mult-of-two-planar-curves}, $m_p(f,f')\geq\mu_p(f)\mu_p(f')=\mu_p(\gamma)(\mu_p(\gamma)-1)$.

Summing over all $p\in\gamma$, B\'ezout's theorem (\cref{thm:bezout-with-multiplicity-projective}) gives
\[\sum_{p\in \gamma}\mu_p(\gamma)(\mu_p(\gamma)-1)\leq \sum_{p\in Z(f)\cap Z(f')}m_p(f,f')\leq \deg f \deg f' = (\deg \gamma)(\deg \gamma -1).\qedhere\]
\end{proof}

\section{Incidences between ruled surfaces}
\label{sec:ruled-ruled}

In this section, we prove the following theorem on line--ruled-surface incidences.

\begin{theorem}\label{thm:ruled-ruled}
    Let $\cP\subset \RR^3$ be a set of size $N$ such that $\abs{\cP\cap \gamma}\leq N^{2/3}\deg \gamma$ for every curve $\gamma\subset\CC^3$, and such that every plane contains at most $N^{2/3}$ points of $\cP$.
    Let $C\geq 2$ be a constant and let $H\subseteq \cP\times \cP$ be any subset.
    For every $(p,q)\in H$, let $\cS_{pq}$ be a set of singly or doubly ruled surfaces in $F_{pq}^+$ of total degree $O(N^{2/3})$.
    Let $\cE$ be the set of $(p,q,p',q')\in H^2$ satisfying the following:
    \begin{enumerate}
        \item $(p,q)\neq (p',q')$ and the line $\ell^{\Phy}$ connecting $\phi(p,q),\phi(p',q')$ is an isotropic line satisfying $|\ell^{\Phy} \cap \phi(H)| \le C$; and
        \item $\ell^+$ is in the ruling of $S$ and $S'$ for some $S\in \cS_{pq}$ and $S'\in \cS_{p'q'}$. 
    \end{enumerate}
    Then
    \[\abs{\cE}\lesssim CN^{2/3}\log N\sum_{(p,q)\in H} \deg \cS_{pq}.\]
\end{theorem}

Recall that for a line $\ell^{\Phy}\subset E^{\Phy}(\CC)$ in physical space, we use $\ell^+\subset E^+(\CC)$ to denote the corresponding line in positive space defined in \cref{prop:all-lines}. Recall also that $\ell^+$ is in the ruling of $S$ if it corresponds to a non-isolated point in the Fano variety $\Lambda_S$, as defined in \cref{defn:in-the-ruling}.

\subsection{Moving from positive space to \texorpdfstring{$\CC^3$}{C3}}

We first move the entire problem from positive space to physical space and then project down to $\CC^3$.
For this subsection, fix a pair $(p,q)\in H$.

By \cref{prop:all-lines}, the isotropic lines in $F_{pq}^+$ are in one-to-one correspondence with the isotropic lines in physical space through $x=\phi(p,q)$.
We parameterize the lines through $x$ via a copy of projective space that we denote $\PP^5_{(p,q)}$. Here $[d]\in\PP^5_{(p,q)}$ corresponds to the line $x+d\CC$.
Recall from \cref{defn:isotropic} that for $d=(\vec d,\vec d_\ast)$, the line corresponding to $[d]\in\PP^5_{(p,q)}$ is isotropic if and only if $\vec{d}\cdot \vec{d}_*=0$.
For each $[d]\in \PP^5_{(p,q)}$ with $\vec{d}\cdot \vec{d}_*=0$, denote by $\ell^{\Phy}_{[d]}$ the isotropic line through $x$ in direction $d$, and let $\ell^+_{[d]}$ be the corresponding isotropic line in $E^+_o$.

\begin{definition}
Assume the same setup as \cref{thm:ruled-ruled}.
Fix $(p,q)\in H$. For each $S\in \cS_{pq}$, write $C_S\subset \PP^5_{(p,q)}$ for the set of $[d]\in\PP^5_{(p,q)}$ so that $\ell_{[d]}^+$ lies in the ruling of $S$.
Define $\cC_{pq}$ to be the collection of all irreducible components of $C_S$ as $S$ ranges through $\cS_{pq}$. 
\end{definition}

As $\cS_{pq}$ is a set of singly or doubly ruled surfaces, and $C_S$ only consists of lines in the ruling of $S$, we see that $\cC_{pq}$ is a set of irreducible curves.

From \cref{thm:classical-ruled-surface}, we immediately deduce the following bound on $\deg\cC_{pq}$.

\begin{lemma}
\label{prop:pos-to-phys}
We have $\deg C_S\leq 2\deg S$ for each $S\in \cS_{pq}$. Thus $\deg \cC_{pq}\leq 2\deg \cS_{pq}$.
\end{lemma}

\begin{proof}
By \cref{thm:plucker-embeding}, there is a linear isomorphism $i\colon F^+_{pq}\to\PP^3_\CC$ so that the parametrization of lines in $\PP^3_{\CC}$ using Pl\"ucker coordinates agrees with the parametrization of lines in $F^+_{pq}$ defined above by geometric triality. Identifying $\PP^5$ with $\PP^5_{(p,q)}$ directly, we see that $C_S$ is the pure dimension 1 part of $\Lambda_{i(S)}$. Then the desired result follows from \cref{thm:classical-ruled-surface}.
\end{proof}

\begin{definition}
Define $\psi_{(p,q)}\colon \CC^6\setminus \{(p,q)\}\to \PP^5_{(p,q)}$ by
$\psi_{(p,q)}(p',q')=[\phi(p'-p,q'-q)]$.
\end{definition}

Note that $\psi_{(p,q)}$ is the composition of $\phi$ with the radial projection centered at $\phi(p,q)$.

With these definitions, $\cE$ is the set of quadruples $(p,q,p',q')\in H^2$ such that the following hold:
\begin{itemize}
    \item $(p,q)\neq (p',q')$ and $|\aff\{\phi(p,q),\phi(p',q')\}\cap \phi(H)|\leq C$;
    \item in $\PP_{(p,q)}^5$, we have $\psi_{(p,q)}(p',q')\in \gamma$ for some $\gamma\in \cC_{pq}$; and
    \item in $\PP_{(p',q')}^5$, we have $\psi_{(p',q')}(p,q)\in \gamma'$ for some $\gamma'\in \cC_{p'q'}$.
\end{itemize}

\begin{lemma}\label{prop:phys-easy-fact}
For any $(p,q)\in H$, let $\cE_{pq}\subseteq H\setminus\{(p,q)\}$ be a set so that $\{(p,q)\}\times\cE_{pq}\subseteq\cE$, then $\abs{\cE_{pq}}\leq C\abs{\psi_{(p,q)}(\cE_{pq})}$ and $\psi_{(p,q)}(\cE_{pq})\subset \bigcup_{\gamma\in \cC_{pq}}\gamma$.
\end{lemma}
\begin{proof}
By the first property above, we have that the map $\psi_{(p,q)}\colon \cE_{pq}\to \PP_{(p,q)}^5$ is at most $C$-to-1. Furthermore, by the second property, its image is contained in $\bigcup_{\gamma\in\cC_{pq}}\gamma$.
\end{proof}

So far we have shown that we can move the problem for positive space to physical space, losing only a factor of $C$. Next we move the problem from physical space to $\CC^3$ by projecting from onto the first or second coordinate. More precisely, let $\PP_p^2$ be the copy of projective space that parameterizes the lines in $\CC^3$ through $p$ where $[d_1]\in\PP^2_p$ corresponds to $p+d_1\CC$. Define $\PP^2_q$ analogously. We view $\PP^2_p,\PP^2_q$ as subspaces of $\PP^5_{(p,q)}$ via the embeddings $[d_1]\mapsto [\phi(d_1,0)]$ and $[d_2]\mapsto[\phi(0,d_2)]$, respectively.

\begin{definition}
Define the projections $\pi_1\colon\PP^5_{(p,q)}\setminus\PP^2_q\to\PP^2_p$ and $\pi_2\colon\PP^5_{(p,q)}\setminus\PP^2_p\to\PP^2_q$ by $\pi_1([d_1:d_2])=[d_1+d_2]$ and $\pi_2([d_1:d_2])=[d_1-d_2]$, i.e., we apply $\phi^{-1}$ and then project onto the first three or last three coordinates.

Define $\psi_p\colon\CC^3\setminus\{p\}\to\PP_p^2$ to be the radial projection $\psi_p(p')=[p'-p]$. Define $\psi_q\colon\CC^3\setminus\{q\}\to\PP_q^2$ analogously.
\end{definition}

We make these definitions so that $\pi_1(\psi_{(p,q)}(p',q'))=\psi_p(p')$ and $\pi_2(\psi_{(p,q)}(p',q'))=\psi_q(q')$ hold.

Instead of studying the incidences $\psi_{(p,q)}(p',q')\in\gamma$ for some point $(p',q')\in\cE_{pq}$ and some curve $\gamma\in\cC_{pq}$, we wish to apply either $\pi_1$ or $\pi_2$ to this containment, to conclude that $\pi_i(\psi_{p,q}(p',q'))\in\overline{\pi_i(\gamma)}$ for $i=1,2$. As the left-hand side is equal to $\psi_p(p')$ or $\psi_q(q')$, this will be a simpler problem. There are two things that could go wrong: potentially $\psi_{(p,q)}(p',q')$ lies in $\PP^2_p$ or $\PP^2_q$ so the left-hand side is undefined and potentially $\overline{\pi_i(\gamma)}$ is a point instead of a curve. We show that the former cannot happen and -- while the latter can occur for one value of $i$ -- it cannot happen for both $i=1,2$.

\begin{lemma}\label{prop:phys-to-C3-edge-case}
For any point $(p',q')\in \RR^6\setminus \{(p,q)\}$ such that $\aff\{\phi(p,q),\phi(p',q')\}$ is an isotropic line, we have $\psi_{(p,q)}(p',q')\not\in \PP^2_p\cup \PP^2_q$.
\end{lemma}
\begin{proof}
Without loss of generality, assume that $\psi_{(p,q)}(p',q')\in \PP^2_p$. By definition, this implies that $q'=q$. Then since $\aff\{\phi(p,q),\phi(p',q')\}$ is assumed to be isotropic $\sum_{i=1}^3(p_i-p'_i)^2-(q_i-q'_i)^2=\sum_{i=1}^{3}(p_i-p_i')^2=0$. Since $p,p'$ are assumed to be real points, this implies that $(p,q)=(p',q')$, which is a contradiction.
\end{proof}

\begin{definition}
Define $\tilde\cC_{pq}\subseteq\cC_{pq}$ to be the set of curves $\gamma\in\cC_{pq}$ so that $\gamma\not\subset\PP_p^2$ and $\gamma\not\subset\PP_q^2$. 
\end{definition}

\cref{prop:phys-to-C3-edge-case} shows that there are no incidences of the form $\psi_{(p,q)}(p',q')\in\gamma$ for $\gamma\in\cC_{pq}\setminus\tilde\cC_{pq}$, so we can safely pass to $\tilde\cC_{pq}$.

\begin{definition}
For $\gamma\in\tilde\cC_{pq}$, write $\gamma_1=\overline{\pi_1(\gamma\setminus\PP_q^2)}$ and $\gamma_2=\overline{\pi_2(\gamma\setminus\PP_p^2)}$.
\end{definition}

We defined $\tilde\cC_{pq}$ so that both $\gamma_1,\gamma_2$ are non-empty. Each is either an irreducible curve or a single point.

\begin{definition}
Define the set of ``vertical'' and ``horizontal'' curves, $\cC^v_{pq},\cC^h_{pq}\subseteq\tilde\cC_{pq}$ as follows. A curve $\gamma\in\tilde\cC_{pq}$ lies in $\cC^v_{pq}$ if $\gamma_1$ is a single point and in $\cC^h_{pq}$ if $\gamma_2$ is a single point.
\end{definition}

\begin{lemma}\label{prop:disjoint-exceptional}
The sets $\cC^v_{pq}$ and $\cC^h_{pq}$ are disjoint.
\end{lemma}
\begin{proof}
For $\gamma\in\tilde\cC_{pq}$, suppose $\gamma_1=\{[d_1]\}$ and $\gamma_2=\{[d_2]\}$. Then $\gamma\setminus\PP_q^2\subseteq\pi_1^{-1}([d_1])$ and also $\gamma\setminus\PP_p^2\subseteq\pi_2^{-1}([d_2])$, meaning that $\gamma\setminus (\PP^2_p\cup \PP^2_q)\subseteq\pi_1^{-1}([d_1])\cap\pi_2^{-1}([d_2])$. The latter set is the line $\{[sd_1+td_2:sd_1-td_2]:s,t\in\CC\setminus\{0\}\}$. (Technically this is a 1-plane in $\PP^5_{pq}$ with the two points where it intersects $\PP_p^2$ and $\PP_q^2$ removed.) Now since $\gamma$ is an irreducible curve that is neither contained in $\PP_p^2$ nor $\PP_q^2$, it must be the line $\{[sd_1+td_2:sd_1-td_2]:s,t\in\CC\}$.

However, by \cref{thm:classical-ruled-surface}, no singly or doubly ruled surface contains a line in its Fano variety, implying that no $\gamma\in\cC_{pq}$ is a line. 
\end{proof}

\begin{definition}
For $\gamma\in\tilde\cC_{pq}\setminus\cC^v_{pq}$, define $m_1(\gamma)$ to be the size of the fiber of $\pi_1\colon\gamma\setminus\PP_q^2\to\gamma_1$ above a generic point. Similarly, for $\gamma\in\tilde\cC_{pq}\setminus\cC^h_{pq}$, define $m_2(\gamma)$ to be the size of the fiber of $\pi_2\colon\gamma\setminus\PP_p^2\to\gamma_2$ above a generic point.
\end{definition}

As explained at the start of \cref{sec:multiplicity}, both $m_1(\gamma),m_2(\gamma)$ are well-defined finite numbers.
Next, we show that we can control $\deg \gamma_i$ and $m_i(\gamma)$ in terms of $\deg\gamma$.

\begin{lemma}\label{prop:individual-curve-deg-bound}
For $\gamma\in \tilde\cC_{pq}\setminus\cC^v_{pq}$, we have $m_1(\gamma)\deg \gamma_1\leq \deg \gamma$. Similarly, for $\gamma\in\tilde\cC_{pq}\setminus\cC^h_{pq}$, we have $m_2(\gamma)\deg \gamma_2\leq \deg \gamma$.
\end{lemma}

\begin{proof}
We will prove the first result; the second is analogous.

Take $\gamma\in\tilde\cC_{pq}\setminus\cC^v_{pq}$. By definition, the fiber of $\pi_1\colon\gamma\setminus\PP_q^2\to\gamma_1$ above $x\in\gamma_1$ has size $m_1(\gamma)$ for $x$ in a Zariski-dense subset of $\gamma_1$. In other words, there is a finite collection of $x$ for which this property fails. Now let $\ell$ be a generic line in $\PP^2_p$, so $\abs{\ell\cap \gamma_1} = \deg \gamma_1$. Furthermore, since $\ell$ is generic, it avoids the finite collection of points which do not have fiber of the correct size. This implies that $\abs{\pi_1^{-1}(\ell)\cap\gamma}=m_1(\gamma)\deg\gamma_1$.

Now $\overline{\pi_1^{-1}(\ell)}$ is a $4$-plane in $\PP^5_{(p,q)}$. This 4-plane does not contain $\gamma$, as $\ell$ does not contain $\gamma_1$. Furthermore, the 4-plane contains at least $m_1(\gamma)\deg \gamma_1$ points of $\gamma$. By B\'ezout's theorem, this implies that $m_1(\gamma)\deg\gamma_1\leq \deg \gamma$, as desired.
\end{proof}

\begin{definition}
Let $\cC^1_{pq}$ be the set of curves $\gamma_1$ (i.e.,$\overline{\pi_1(\gamma\setminus\PP_q^2)}$) for $\gamma\in\tilde\cC_{pq}\setminus\cC^v_{pq}$. Analogously, let $\cC^2_{pq}$ be the set of curves $\gamma_2$ for $\gamma\in\tilde\cC_{pq}\setminus\cC^h_{pq}$.
For each $\gamma_1\in \cC^1_{pq}$, define
\[m_{pq}(\gamma_1) = \sum_{\gamma\in \tilde\cC_{pq}\setminus\cC^v_{pq}:\overline{\pi_1(\gamma\setminus\PP_q^2)} = \gamma_1}m_1(\gamma).\]
For $\gamma_2\in\cC^2_{pq}$, define $m_{pq}(\gamma_2)$ analogously.
    
\end{definition}

Combining \cref{prop:pos-to-phys,prop:individual-curve-deg-bound}, we get the following.

\begin{lemma}\label{prop:total-curve-bound}
    For each $(p,q)\in H$ and $i=1,2$, we have
    \[\sum_{\gamma_i\in \cC^i_{pq}}m_{pq}(\gamma_i)\deg \gamma_i \leq 2\deg \cS_{pq}.\]
\end{lemma}

Moreover, ignoring the incidences coming from $\cC^v_{pq}$, the following proposition says that we can project the problem via $\pi_1$ down to $\CC^3$. (Symmetrically, ignoring the incidences coming from $\cC^h_{pq}$, we can project the problem via $\pi_2$ down to $\CC^3$.)

\begin{proposition}\label{prop:pos-to-C3}
Assume the same setup as \cref{thm:ruled-ruled}. Fixing $p,q\in \cP$ and $\cP_2\subseteq\cP$, let $\cE_{pq}\subseteq \cP_2\times\cP$ be a set so that $\{(p,q)\}\times\cE_{pq}\subseteq\cE$ and so that for each $(p',q')\in\cE_{pq}$ we have $\psi_{(p,q)}(p',q')\in\gamma$ for some $\gamma\in\tilde\cC_{pq}\setminus\cC^v_{pq}$. Then
\[\abs{\cE_{pq}}\lesssim CN^{2/3}\deg\cS_{pq}+C\sum_{\gamma_1\in \cC^1_{pq}}m_{pq}(\gamma_1)\abs{\psi_p(\cP_2\setminus\{p\})\cap \gamma_1}.\]
\end{proposition}

\begin{proof}
By \cref{prop:phys-easy-fact}, we know that $\abs{\cE_{pq}}\leq C\abs{\psi_{(p,q)}(\cE_{pq})}$. By hypothesis $\psi_{(p,q)}(\cE_{pq})\subset\bigcup_{\gamma\in\tilde\cC_{pq}\setminus\cC^v_{pq}}\gamma$. Thus it suffices to show
\[\sum_{\gamma\in \tilde\cC_{pq}\setminus \cC^v_{pq}}\abs{\psi_{(p,q)}(\cE_{pq})\cap \gamma}\lesssim N^{2/3}\deg\cS_{pq}+\sum_{\gamma_1\in \cC^1_{pq}}m_{pq}(\gamma_1)\abs{\psi_p(\cP_2\setminus\{p\})\cap \gamma_1}.\]

For each point $y\in \psi_{(p,q)}(\cE_{pq})\cap \gamma$, by \cref{prop:phys-to-C3-edge-case} we know $y\not\in\PP^2_q$, so $\pi_1(y)$ is well-defined. Since $\gamma\in\tilde\cC_{pq}\setminus\cC^v_{pq}$, we know $\pi_1(y)\in\gamma_1$ for some $\gamma_1\in\cC^1_{pq}$, say $\pi_1(y)=z$. As we computed above, if $y=\psi_{(p,q)}(p',q')$, then $z=\pi_1(y)=\psi_p(p')$. Since $z$ is well-defined, we must have $p'\neq p$. In other words, $z\in\psi_p(\cP_2\setminus\{p\})$. Thus we have shown $z=\pi_1(y)\in\psi_p(\cP_2\setminus\{p\})\cap\gamma_1$.

Now we bound how many $y\in \psi_{(p,q)}(\cE_{pq})\cap \gamma$ can correspond to the same $z$. By \cref{thm:upper-bound-fiber}, for each such $\gamma,\gamma_1,z$, the fiber of $\pi_1\colon\gamma\setminus\PP_q^2\to\gamma_1$ above $z$ has size at most $m_1(\gamma)\mu_z(\gamma_1)$. From the definition of $m_{pq}$ we thus have
\[\sum_{\gamma\in \tilde\cC_{pq}\setminus \cC^v_{pq}}\abs{\psi_{(p,q)}(\cE_{pq})\cap \gamma}\leq \sum_{\gamma_1\in \cC^1_{pq}}m_{pq}(\gamma_1)\sum_{z\in \psi_p(\cP_2\setminus\{p\})\cap \gamma_1}\mu_z(\gamma_1).\]
    
Now we have
\[\sum_{z\in \psi_p(\cP_2\setminus\{p\})\cap \gamma_1}\mu_z(\gamma_1)\leq\abs{\psi_p(\cP_2\setminus\{p\})\cap\gamma_1}+\sum_{z\in\gamma_1:\mu_z(\gamma_1)>1}\mu_z(\gamma_1).\]
Using \cref{prop:planar-easy-multiplicity}, we can bound the latter sum as
\[\sum_{z\in\gamma_1:\mu_z(\gamma_1)>1}\mu_z(\gamma_1)\leq\sum_{z\in\gamma_1}\mu_z(\gamma_1)(\mu_z(\gamma_1)-1)\leq\deg\gamma_1(\deg\gamma_1-1)\lesssim N^{2/3}\deg\gamma_1.\]
The final inequality follows from the very rough bound that since $\gamma_1=\overline{\pi_1(\gamma\setminus\PP_q^2)}$ for some $\gamma\in\cC_{pq}$, we have $\deg\gamma_1\leq\deg\gamma\leq\deg\cC_{pq}\leq 2\deg\cS_{pq}\lesssim N^{2/3}$.

Thus we have shown
\[\sum_{\gamma\in \tilde\cC_{pq}\setminus \cC^v_{pq}}\abs{\psi_{(p,q)}(\cE_{pq})\cap \gamma}\lesssim N^{2/3}\sum_{\gamma_1\in\cC^1_{pq}}m_{pq}(\gamma_1)\deg\gamma_1+\sum_{\gamma_1\in \cC^1_{pq}}m_{pq}(\gamma_1)\abs{\psi_p(\cP_2\setminus\{p\})\cap \gamma_1}.\]
Applying \cref{prop:total-curve-bound} to bound the first term thus completes the proof.
\end{proof}

\subsection{Curves and cones in \texorpdfstring{$\CC^3$}{C3}}
By \cref{prop:pos-to-C3}, the bulk of the problem is to bound terms of the form $|\psi_p(\cP_2\setminus\{p\})\cap\gamma_1|$. Taking the preimage under $\psi_p$, this quantity is upper bounded by $|\cP_2\cap\psi_p^{-1}(\gamma_1)|$. 

\begin{definition}
For an irreducible curve $\gamma_1\subset\PP^2_p$, write $C(p,\gamma_1)$ for the cone $\psi_p^{-1}(\gamma_1)\subset\CC^3$. 
\end{definition}

To control the incidences between $\cP_2\subseteq\cP$ and the collection of cones, we first cluster $\cP$ into a low-degree collection of curves.

\begin{proposition}\label{prop:curves-containing-points}
There exists a collection $\cZ$ of irreducible curves in $\CC^3$ such that $\deg \cZ\lesssim N^{2/3}$ and each $p\in\cP$ is contained in some $\gamma\in\cZ$.
\end{proposition}

\begin{proof}
By parameter counting, \cref{lem:Walsh-param-counting}, there exists a nonzero polynomial $f\in\CC[x,y,z]$ that vanishes on $\cP$ with $\deg f\lesssim N^{1/3}$. Write $V_1,\ldots,V_t$ for the irreducible components of $Z(f)$ and define a partition $\cP=\bigsqcup_{i=1}^t \cP_i$ so that $\cP_i\subset V_i$. (If a point lies in multiple surfaces, assign it to one of the corresponding parts arbitrarily.)

By a second application of \cref{lem:Walsh-param-counting}, for each $i$, there is a polynomial $g_i\in\CC[x,y,z]$ that vanishes on $\cP_i$ without vanishing identically on $V_i$ such that $\deg g_i\lesssim\max\{\abs{P_i}^{1/3},(\abs{P_i}/\deg V_i)^{1/2}\}$. Setting $\cZ=\bigcup_{i=1}^t V_i\cap Z(g_i)$, we clearly have that this collection of curves contains $\cP$. Furthermore, by B\'ezout's theorem,
\begin{align*}
\deg\cZ
&\leq\sum_{i=1}^t \deg V_i\deg g_i
\lesssim \sum_{i=1}^t\abs{P_i}^{1/3}\deg V_i+\sum_{i=1}^t\abs{P_i}^{1/2}(\deg V_i)^{1/2}\\
&\leq N^{1/3}\sum_{i=1}^t\deg V_i+\paren{\sum_{i=1}^t \abs{P_i}}^{1/2}\paren{\sum_{i=1}^t\deg V_i}^{1/2}\\
&\lesssim N^{1/3}\cdot N^{1/3}+\paren{N}^{1/2}\paren{N^{1/3}}^{1/2}
\lesssim N^{2/3}.\qedhere
\end{align*}
\end{proof}

For each incidence between a point $p'\in\cP$ and a cone $C(p,\gamma_1)$, there is a curve $\zeta\in\cZ$ containing $p'$. B\'ezout's theorem gives a good bound on the number of incidences where $\zeta\not\subset C(p,\gamma_1)$. The main term is to bound the contribution from incidences where the curve $\zeta$ is contained in the cone $C(p,\gamma_1)$.

In the situation where $\zeta\subset C(p,\gamma_1)$, we can relate $\zeta$ to $\gamma_1$ since usually $\overline{\psi_p(\zeta\setminus\{p\})}=\gamma_1$. (The only time this does not occur is when $\zeta$ is a line through $p$; this case is easy to deal with.) This allows us to relate $\deg\zeta$ to $\deg\gamma_1$ in all circumstances except when $p$ both lies on $\zeta$ and is a very singular point of this curve. We next explain this relation and show that the very singular case is rare. 

\begin{definition}
For an irreducible curve $\zeta\subset \CC^3$, write $a(p,\zeta)$ for the size of a generic fiber of $\psi_p\colon \zeta\setminus \{p\}\to \overline{\psi_p(\zeta\setminus \{p\})}$.
\end{definition}

Note that if $a(p,\zeta)=\infty$, then $\zeta$ is a line through $p$.

\begin{proposition}\label{prop:proj-formula-with-singular-origin}
For $p\in\CC^3$, let $\zeta\subset \CC^3$ be an irreducible curve that is not a line through $p$. Defining  $\overline{\psi_p(\zeta\setminus \{p\})}=\gamma_1$, then
 $\deg \zeta = a(p,\zeta)\deg\gamma_1+\mu_p(\zeta)$.
\end{proposition}

This proof is very similar to \cref{lem:projection,prop:individual-curve-deg-bound}.

\begin{proof}
Let $\ell$ be a generic line in $\PP_p^2$. Then $\ell\cap\gamma_1$ consists of $\deg\gamma_1$ generic points of $\gamma_1$. In particular, $\abs{\psi_p^{-1}(\ell)\cap\zeta}=a(p,\gamma)\deg\gamma_1$.

Now $F:=\overline{\psi^{-1}_p(\ell)}$ is a generic plane through $p$ in $\CC^3$. The plane $F$ intersects $\zeta$ potentially at $p$ as well as at the $a(p,\zeta)\deg\gamma_1$ points we found previously. We first claim that the intersections between $F$ and $\zeta$ not at $p$ are all transversal. As in the proof of \cref{lem:projection}, note that $\zeta$ has a finite number of singular points and a finite number of non-singular points $q\in\zeta$ for which $T_q\zeta$ (the line through $q$ tangent to $\zeta$) also passes through $p$. Indeed, the condition $p\in T_q\zeta$ is a Zariski-closed condition. If it were satisfied for all $q\in\zeta$, then $\zeta$ would be a line through $p$ which we assumed was not the case. Since $F$ is a generic plane through $p$, it avoids this finite set of bad points on $\zeta$.

Thus we have shown that away from $p$, the intersection $F\cap\zeta$ consists of $a(p,\zeta)\deg\gamma_1$ intersections $q$ where $q$ is a non-singular point of $\zeta$ and $F$ intersects $T_q\zeta$ transversally. This means that these $a(p,\zeta)\deg\gamma_1$ intersections occur with multiplicity 1 (see, e.g., \cite[Lecture 18]{Harris92}). 
Next note that the intersection multiplicity of $F$ and $\zeta$ at $p$ is $\mu_p(\zeta)$ (which is defined to be 0 if $p$ does not lie on $\zeta$) by \cref{prop:mult-alt-def}. Finally, we count the number of intersections at infinity.

Defining $\tilde{\zeta}$ and $\tilde{F}$ to be the Zariski-closures of $\zeta$ and $F$ in $\PP^3$, we claim that $\tilde{\zeta}$ and $\tilde{F}$ do not intersect at infinity. We know $\tilde{\zeta}\setminus\zeta$ is finite, so again -- since $\tilde{F}$ is a generic plane through $p$ in $\PP^3$ -- we see that $\tilde{F}$ avoids this finite collection of bad points.

Now by B\'ezout's theorem (\cref{thm:bezout-with-multiplicity-projective})
\[\deg\zeta=\sum_{q\in \tilde{\zeta}\cap\tilde{F}}m_q(\tilde{\zeta},\tilde{F})=\mu_p(\zeta)+\sum_{q\in (\zeta\cap F)\setminus\{p\}}m_q(\zeta, F)=a(p,\zeta)\deg\gamma_1+\mu_p(\zeta).\qedhere\]
\end{proof}

\begin{definition}
For an irreducible curve $\zeta\subset\CC^3$, say that a point $p$ is a \emph{very singular} point of $\zeta$ if $\mu_p(\zeta)>\tfrac12\deg\zeta$.
\end{definition}

The above result implies that $\deg\zeta\leq 2a(p,\zeta) \deg \gamma_1$ as long as $p$ is not a very singular point of $\zeta$. We next show that each curve has few very singular points.

\begin{proposition}\label{prop:unique-very-singular}
Let $\zeta\subset \CC^3$ be an irreducible curve that is not a line.
Then there is at most one very singular point of $\zeta$.
\end{proposition}

\begin{proof}
Suppose for the sake of contradiction that there exist points $p_1,p_2$ with $\mu_{p_i}(\zeta)>\tfrac{1}{2}\deg \zeta$ for $i=1,2$.
Then let $f$ be a linear function that vanishes at $p_1,p_2$ but does not vanish identically on $\zeta$.
By \cref{lem:mult-of-all-surface}, we can bound the intersection multiplicity as 
\[m_{p_i}(\zeta,f)\geq \mu_{p_i}(\zeta)>\frac{1}{2}\deg\zeta\]
for $i=1,2$.
However this contradicts B\'ezout's theorem (\cref{thm:bezout-with-multiplicity-projective}) since
\[\deg\zeta=\deg\zeta\deg f\geq\sum_{p\in \zeta\cap Z(f)} m_p(\zeta,f)\geq \mu_{p_1}(\zeta)+\mu_{p_2}(\zeta)>\deg\zeta.\qedhere\]
\end{proof}

Putting together these tools, we can bound the main contribution in \cref{thm:ruled-ruled} by reducing it to a curve-cone incidence problem.

\begin{theorem}
\label{lem:ruled-ruled-dyadic-general}
Assume the same setup as \cref{thm:ruled-ruled}. Given $\cP_1,\cP_2\subseteq\cP$, let $\cE'\subseteq\cE$ be the set of quadruples $(p,q,p',q')\in\cE$ so that $p\in\cP_1$ and $p'\in\cP_2$ with $\psi_{(p,q)}(p',q')\in\gamma$ for some $\gamma\in\tilde\cC_{pq}\setminus\cC^v_{pq}$. Then
\[|\cE'|\lesssim CN^{2/3}\sum_{(p,q)\in H\cap(\cP_1\times\cP)}\deg \cS_{pq}+C|\cP_1|^{2/3}|\cP_2|\max_{p\in\cP_1}\sum_{q:(p,q)\in H}\deg\cS_{pq}.\]
\end{theorem}

Summing the contribution to $\cE$ over each pair $(p,q)\in H\cap(\cP_1\times\cP)$, this theorem immediately follows from \cref{prop:pos-to-C3} and the next result, a curve cone incidence bound.

\begin{proposition}
\label{prop:ruled-ruled-curve-cone-incidences}
Assume the same setup as \cref{thm:ruled-ruled}. Given $\cP_1,\cP_2\subseteq\cP$, write $H_1=H\cap(\cP_1\times\cP)$. Then
\[\sum_{(p,q)\in H_1}\sum_{\gamma_1\in\cC^1_{pq}}m_{pq}(\gamma_1)\abs{\psi_p(\cP_2\setminus\{p\})\cap\gamma_1}\lesssim N^{2/3}\sum_{(p,q)\in H_1}\deg \cS_{pq}+|\cP_1|^{2/3}|\cP_2|\max_{p\in\cP_1}\sum_{q:(p,q)\in H}\deg\cS_{pq}.\]
\end{proposition}

\begin{proof}
Applying \cref{prop:curves-containing-points} to $\cP_2$, we produce a collection $\cZ$ of irreducible curves in $\CC^3$ so that $\deg\cZ\lesssim N^{2/3}$ such that each $p\in\cP_2$ is contained in some $\zeta\in\cZ$. Define a partition $\cP_2=\bigsqcup_{\zeta\in\cZ}\cP_\zeta$ where $\cP_{\zeta}\subset\zeta$ (if a point lies in multiple curves, assigned it to one of the sets arbitrarily). The term we wish to bound is thus
\[\sum_{(p,q)\in H_1}\sum_{\gamma_1\in\cC^1_{pq}}\sum_{\zeta\in\cZ}m_{pq}(\gamma_1)\abs{\psi_p(\cP_\zeta\setminus\{p\})\cap\gamma_1}.\]

We break this sum into several pieces depending on properties of the curve $\zeta$ and the cone $C(p,\gamma_1)$:
\begin{enumerate}
    \item terms with $\zeta\not\subset C(p,\gamma_1)$: these are handled by B\'ezout's theorem;
    \item terms with $\zeta\subset C(p,\gamma_1)$ where $C(p,\gamma_1)$ is a plane: these are easy since each plane contains few points of $\cP$;
    \item terms with $\zeta\subset C(p,\gamma_1)$ where $C(p,\gamma_1)$ is not a plane, but $\zeta$ is a line: these terms are easy;
    \item terms with $\zeta\subset C(p,\gamma_1)$ where $\zeta$ is not a line but is very singular at $p$ (meaning that $\mu_p(\zeta)>\tfrac12\deg\zeta$): we can handle these terms since each non-line curve contains at most one very singular point;
    \item terms with $\zeta\subset C(p,\gamma_1)$ where $\zeta$ is not a line and is not very singular at $p$ (meaning that $\mu_p(\zeta)\leq \tfrac12\deg\zeta$): this is the main term.
\end{enumerate}

For (3), since $C(p,\gamma_1)$ is not a plane but contains the line $\zeta$, note that $\zeta$ must be a line through $p$. Then $\psi_p(\cP_\zeta\setminus\{p\})\subseteq\psi_p(\zeta\setminus\{p\})$, which is a single point. Thus the contribution of each of these terms is bounded as
\[\sum_{p,q,\gamma_1,\zeta:(3)}m_{pq}(\gamma_1)\abs{\psi_p(\cP_\zeta\setminus\{p\})\cap\gamma_1}\leq \sum_{p,q,\gamma_1,\zeta:(3)}m_{pq}(\gamma_1).\]
The number of lines in $\cZ$ is bounded by $\deg\cZ\lesssim N^{2/3}$. Using the trivial bound $1\leq \deg(\gamma_1)$ with \cref{prop:total-curve-bound}, we can bound the contribution of (3) by
\[\lesssim N^{2/3}\sum_{(p,q)\in H_1}\sum_{\gamma_1\in\cC^1_{pq}}m_{pq}(\gamma_1)\deg\gamma_1\lesssim N^{2/3}\sum_{(p,q)\in H_1}\deg\cS_{pq}.\]

In all the other cases, we use the bound 
\[|\psi_p(\cP_\zeta\setminus\{p\})\cap\gamma_1|\leq |\cP_\zeta\cap C(p,\gamma_1)|.\]

For (1), B\'ezout's theorem gives
\begin{align*}
\sum_{p,q,\gamma_1,\zeta:(1)}m_{pq}(\gamma_1)\abs{\psi_p(\cP_\zeta\setminus\{p\})\cap\gamma_1}
&\leq \sum_{p,q,\gamma_1,\zeta:(1)}m_{pq}(\gamma_1)\abs{\cP_\zeta\cap C(p,\gamma_1)}\\
&\leq \sum_{p,q,\gamma_1,\zeta:(1)}m_{pq}(\gamma_1)\deg\zeta\deg\gamma_1.
\end{align*}
Since $\sum_{\zeta\in\cZ}\deg\zeta =\deg\cZ\lesssim N^{2/3}$, \cref{prop:total-curve-bound} again bounds the contribution of this case by
\[\lesssim N^{2/3}\sum_{(p,q)\in H_1}\deg\cS_{pq}.\]

For (2), since each plane contains at most $N^{2/3}$ points of $\cP$ (and since the $\cP_\zeta$ are disjoint), we have
\[\sum_{p,q,\gamma_1,\zeta:(2)}m_{pq}(\gamma_1)\abs{\cP_\zeta\cap C(p,\gamma_1)}\leq N^{2/3}\sum_{p,q,\gamma_1}m_{pq}(\gamma_1)\lesssim N^{2/3}\sum_{(p,q)\in H_1}\deg \cS_{pq}.\]
The second inequality again follows from the trivial bound $1\leq \deg(\gamma_1)$ and \cref{prop:total-curve-bound}.

For (4), \cref{prop:unique-very-singular} shows that each $\zeta\in\cZ$ which is not a line contains at most one very singular point $p$. Thus 
\begin{align*}
\sum_{p,q,\gamma_1,\zeta:(4)}m_{pq}(\gamma_1)\abs{\cP_\zeta\cap C(p,\gamma_1)}
&=\sum_{p,q,\gamma_1,\zeta:(4)}m_{pq}(\gamma_1)\abs{\cP_\zeta}\\
&\leq\sum_{\zeta\in\cZ}\abs{\cP_\zeta}\max_{p\in\cP_1}\sum_{q:(p,q)\in H}\sum_{\gamma_1\in\cC^1_{pq}}m_{pq}(\gamma_1)\\
&=|\cP_2|\max_{p\in\cP_1}\sum_{q:(p,q)\in H}\sum_{\gamma_1\in\cC^1_{pq}}m_{pq}(\gamma_1)\\
&\lesssim |\cP_2|\max_{p\in\cP_1}\sum_{q:(p,q)\in H}\deg\cS_{pq}.
\end{align*}
The final inequality again follows from the trivial bound $1\leq \deg(\gamma_1)$ and \cref{prop:total-curve-bound}.

Finally we isolate the main term, case (5).

\begin{definition}
\label{defn:zeta-p-gamma}
Fix a set $\cZ$ of irreducible curves in $\CC^3$. For a point $p\in\CC^3$ and an irreducible curve $\gamma_1\subset\PP_p^2$ that is not a line, define $\cZ_{p,\gamma_1}\subseteq\cZ$ to be the set of curves $\zeta\in\cZ$ where $\zeta$ is not a line, $\zeta\subset C(p,\gamma_1)$, and $p$ is not a very singular point of $\zeta$. If $\gamma_1$ is a line, define $\cZ_{p,\gamma_1}=\emptyset$.
\end{definition}

Now the desired bound on case (5) follows from the next claim which we prove in \cref{ssec:trisecant}.

\begin{claim}
\label{claim:curve-cone-incidences}
\[\sum_{(p,q)\in H_1}\sum_{\gamma_1\in\cC^1_{pq}}\sum_{\zeta\in\cZ_{p,\gamma_1}}m_{pq}(\gamma_1)|\cP_\zeta|\lesssim N^{2/3}\sum_{(p,q)\in H_1}\deg\cS_{pq}+|\cP_1|^{2/3}|\cP_2|\max_{p\in\cP_1}\sum_{q:(p,q)\in H}\deg\cS_{pq}.\]
\end{claim}

Combining the bounds in each of the cases gives the desired bound.
\end{proof}

\subsection{Curve-cone incidences via trisecants}
\label{ssec:trisecant}

In this section, we prove \cref{claim:curve-cone-incidences}. To bound the number of curve-cone incidences we show that it is unlikely to have three distinct curves $\zeta_1,\zeta_2,\zeta_3$ all contained in the same cone. Morally, this can be thought of as showing that the curve-cone incidence graph is $K_{3,C}$-free for some large $C$. (If the cones all have degree $d$, we can take $C=O(d^2)$.) Indeed, the expression $|\cP_1|^{2/3}|\cP_2|$ matches with the K\H{o}v\'ari--S\'os--Tur\'an bound for $K_{3,C}$-free graphs.

To make this argument work, we also show that it is unlikely to have a single curve $\zeta$ that is ``triply contained'' in a cone $C(p,\gamma_1)$, i.e., with $a(p,\zeta)\geq 3$. Both of these are proved using the following version of the trisecant lemma from algebraic geometry.

Recall that for a point $y\in\G(1,3)$, we write $\ell_y\subset\PP_{\CC}^3$ for the corresponding line.

\begin{theorem}[{\cite[Theorem 1]{Ran91} for $\ell=n=1, r=3$}]
\label{thm:trisecant}
Let $X\subset \PP_{\CC}^3$ be a (not necessarily irreducible) curve, and let $Y\subseteq \G(1,3)$ be an irreducible variety such that for a generic $y\in Y$, the line $\ell_y$ is a trisecant of $X$, i.e., $\abs{\ell_y\cap X}\geq 3$. Then
\[\dim\paren{\overline{\bigcup_{y\in Y}\ell_y}}\leq 2.\]
\end{theorem}

Using ruled surface theory, we immediately deduce the following.

\begin{corollary}\label{cor:trisecant}
Let $X\subset \CC^3$ be a (not necessarily irreducible) curve, and let $Y\subseteq \G(1,3)$ be an irreducible variety such that for a generic $y\in Y$, the line $\ell_y$ is a trisecant of $X$. If $\dim Y\geq 2$, then $\overline{\bigcup_{y\in Y}\ell _y}$ is a plane.
\end{corollary}

Note that although \cref{thm:trisecant} was stated projectively, there is no issue restricting to $\CC^3$. When we say that $\ell_y\subset\PP^3_{\CC}$ is a trisecant of $X$ we simply mean that it contains at least 3 points of $X\subset\CC^3\subset\PP^3_{\CC}$.

\begin{proof}
Write $V=\overline{\bigcup_{y\in Y}\ell _y}$ and define $W\subseteq Y\times V$ to be the set of pairs $(y,p)\in Y\times V$ for which $p\in\ell_y$. Writing $\pi_1\colon W\to Y$ and $\pi_2\colon W\to V$ for the projection maps, note that $V=\overline{\pi_2(\pi_1^{-1}(Y)})$. Since $Y$ and all the fibers of $\pi_1$ are irreducible, we conclude that $W$ and thus $V$ are also irreducible \cite[Theorem 11.14]{Harris92}. By \cref{thm:trisecant}, we know $\dim V\leq 2$. Now clearly $Y$ is contained in $\Lambda_V$, the Fano variety of $V$. However, for an irreducible surface $V$, we have $\dim\Lambda_V\leq 1$ unless $V$ is a plane by \cref{thm:classical-ruled-surface}. (If $\dim V\leq 1$, then clearly $\Lambda_V$ is a finite set of points.) This shows the desired result.
\end{proof}

We apply \cref{cor:trisecant} in the following slightly different form.

\begin{lemma}
\label{lem:trisecant}
Let $\zeta_1,\zeta_2,\zeta_3\subset\CC^3$ be distinct irreducible curves. The following three statements hold:
\begin{enumerate}
        \item if $\zeta_1$ is not contained in a plane, then for a generic point $p\in\zeta_1$, we have $a(p,\zeta_1)=1$;
        \item if $\zeta_1\cup\zeta_2$ is not contained in a plane, then for a generic point $p\in\zeta_1$, we have $a(p,\zeta_2)=1$ and $\overline{\psi_p(\zeta_1\setminus \{p\})}\neq \overline{\psi_p(\zeta_2\setminus \{p\})}$.
        \item If $\zeta_1\cup\zeta_2\cup\zeta_3$ is not contained in a plane, then for a generic point $p\in\zeta_1$, we have $\overline{\psi_p(\zeta_2\setminus \{p\})}\neq \overline{\psi_p(\zeta_3\setminus \{p\})}$.
    \end{enumerate}
\end{lemma}

\begin{proof}
To prove (1), suppose for the sake of contradiction that $a(p,\zeta_1)\geq 2$ for a generic point $p\in \zeta_1$. This means that the generic fibers of $\psi_p\colon\zeta_1\setminus\{p\}\to\overline{\psi_p(\zeta_1\setminus\{p\}})$ have size at least 2. Thus for a generic point $p_1\in\zeta_1$, the line $\aff\{p,p_1\}$ intersects $\zeta_1$ in at least one other point. This produces a 2-dimensional family of trisecants.

In particular, consider the map
$s\colon\CC^3\times\CC^3\setminus\Delta\to\G(1,3)$ which maps a pair of distinct points $p,q\in\CC^3$ to $s(p,q)\in\G(1,3)$ so that $\ell_{s(p,q)}$ is the line through $p,q$. Defining $Y=\overline{s(\zeta_1\times\zeta_1\setminus\Delta)}$, we just showed that $\ell_y$ is a trisecant of $\zeta_1$ for generic $y\in Y$. To apply the trisecant lemma, we just need to check that $Y$ is irreducible and has dimension 2. The first is clear since $\zeta_1$, and thus $\zeta_1\times\zeta_1$, is irreducible. For the second, $\dim(\zeta_1\times\zeta_1)=2$ while the generic fibers of $s|_{\zeta_1\times\zeta_1\setminus\Delta}$ are finite (since the intersection between any line and $\zeta_1$ is finite), implying that $\dim Y=\dim(\zeta_1\times\zeta_1)$.

Now \cref{cor:trisecant} implies $\overline{\bigcup_{y\in Y}\ell_y}$ is a plane, yet this union contains a (Zariski-dense subset of) $\zeta_1$, which contradicts the non-planarity of $\zeta_1$.

A similar argument shows (2) and (3), applying \cref{cor:trisecant} to $\zeta_1\cup\zeta_2$ and $\zeta_1\cup\zeta_2\cup\zeta_3$ respectively. Indeed, if $a(p,\zeta_2)\geq 2$ for a generic point $p\in\zeta_1$, then for a generic point $p_2\in\zeta_2$, the line $\aff\{p,p_2\}$ intersects $\zeta_2$ in at least one other point. On the other hand, if for a generic point $p\in\zeta_1$, we have $\overline{\psi_p(\zeta_1\setminus \{p\})}=\overline{\psi_p(\zeta_2\setminus \{p\})}$, then for a generic point $p_1\in\zeta_1$, the line $\aff\{p,p_1\}$ also intersects $\zeta_2$. Finally, if for a generic point $p\in\zeta_1$, we have $\overline{\psi_p(\zeta_2\setminus \{p\})}= \overline{\psi_p(\zeta_3\setminus \{p\})}$, then for a generic point $p_2\in\zeta_2$, the line $\aff\{p,p_2\}$ also intersects $\zeta_3$.
\end{proof}

To formalize the ideas from the beginning of the subsection, consider three curves $\zeta_1,\zeta_2,\zeta_3$. The goal is to show that there are few cones containing all three curves. Given such a cone $C(p,\gamma_1)$, consider the line $\aff\{p,p_1\}$ for a generic point $p_1\in\zeta_1$. By definition, this line intersects $\zeta_2$ in $a(p,\zeta_2)$ points and $\zeta_3$ in $a(p,\zeta_3)$ points. If $p$ is not a very singular point of $\zeta_2,\zeta_3$, then these quantities are at large: at least $\deg\zeta_2/2\deg\gamma_1$ and $\deg\zeta_3/2\deg\gamma_1$.

Projecting radially from $p_1$, the line $\aff\{p,p_1\}$ corresponds to the point $\psi_{p_1}(p)$ in the intersection of the curves $\overline{\psi_{p_1}(\zeta_2\setminus\{p_1\})}$ and $\overline{\psi_{p_1}(\zeta_3\setminus\{p_1\})}$. By the trisecant lemma, these two curves are distinct and $\psi_{p_1}$ is a generically 1-to-1 projection onto each of the curves. By \cref{thm:upper-bound-fiber}, this implies that $\psi_{p_1}(p)$ must be a very singular point on each curve, so the intersection multiplicity is very high here. This shows that such cones $C(p,\gamma_1)$ are rare.

\begin{lemma}\label{lem:self-intersection-bound-3}
Let $\cZ$ be a set of irreducible curves in $\CC^3$. For $\zeta_1,\zeta_2,\zeta_3\in\cZ$ distinct curves, we have
\[\sum_{p,\gamma_1:\zeta_1,\zeta_2,\zeta_3\in\cZ_{p,\gamma_1}}\left(\deg \gamma_1\right)^{-2}\leq 4.\]
\end{lemma}

Recall from \cref{defn:zeta-p-gamma} that $\zeta\in\cZ_{p,\gamma_1}$ implies that neither $\zeta$ nor $\gamma_1$ is a line.

\begin{proof}
We may assume that no two of $\zeta_1,\zeta_2,\zeta_3$ is simultaneously contained in any plane. Indeed, suppose that $\zeta_1,\zeta_2$ are contained in the plane $F$. If there exist $p,\gamma_1$ with $\zeta_1,\zeta_2,\zeta_3\in C(p,\gamma_1)$, since $\gamma_1$ is not a line, $C(p,\gamma_1)$ is not a plane. Then $F\cap C(p,\gamma_1)$ is a finite union of lines if $p\in F$ or an irreducible curve if $p\not\in F$. Since $\zeta_1,\zeta_2\subseteq F\cap C(p,\gamma_1)$ are distinct curves, neither of which is a line, this is impossible. Thus we conclude that if any two of the curves is simultaneously contained in a plane, then the desired result is vacuously true.

Let $\cQ$ be the set of points $p$ for which there exists $\gamma_1$ with $\zeta_1,\zeta_2,\zeta_3\in\cZ_{p,\gamma_1}$. Note that given $p\in\cQ$, the curve $\gamma_1$ is determined uniquely: it must be $\gamma_1=\overline{\psi_p(\zeta_1\setminus\{p\})}=\overline{\psi_p(\zeta_2\setminus\{p\})}=\overline{\psi_p(\zeta_3\setminus\{p\})}$. Our first goal is to show that $\cQ$ is finite.

By \cref{lem:trisecant}(3), for a generic point $p_1\in\zeta_1$, we have $\zeta_{2,p_1}=\overline{\psi_{p_1}(\zeta_2\setminus\{p_1\})}$ and $\zeta_{3,p_1}=\overline{\psi_{p_1}(\zeta_3\setminus\{p_1\})}$ are distinct irreducible curves. Thus $\zeta_{2,p_1}\cap\zeta_{3,p_1}$ is a finite subset of $\PP_{p_1}^2$. Set $\cL_{p_1}=\psi_{p_1}^{-1}(\zeta_{2,p_1}\cap\zeta_{3,p_1})$, a finite set of lines through $p_1$. 

Now suppose we have some $p,\gamma_1$ with $\zeta_1,\zeta_2,\zeta_3\in\cZ_{p,\gamma_1}$. Since none of $\zeta_1,\zeta_2,\zeta_3$ are lines, a generic line through $p$ in $C(p,\gamma_1)$ is a trisecant, i.e., intersects each of $\zeta_1,\zeta_2,\zeta_3$. Furthermore, for a generic point $p_1\in\zeta_1$, the line through $p,p_1$ is a trisecant. Note that $\cL_{p_1}$ is the set of trisecants that pass through $p_1$. Thus we know that for a generic point $p_1\in\zeta_1$, we have $p\in\cL_{p_1}$.

Pick two generic points $p_1,\tilde p_1\in\zeta_1$. We argued above that $p\in\cL_{p_1}\cap\cL_{\tilde p_1}$, so $\cQ\subseteq\cL_{p_1}\cap\cL_{\tilde p_1}$. As $\zeta_1$ is not a line, $\cL_{p_1}\cap \zeta_1$ is finite, so $\tilde p_1\in \zeta_1\setminus\cL_{p_1}$. Now $\cL_{p_1}$ is a set of lines through $p_1$ and $\cL_{\tilde p_1}$ is a set of lines through $\tilde p_1$. Since $\cL_{p_1}$ does not include the line connecting $p_1,\tilde p_1$, we conclude that the intersection, and thus $\cQ$, is finite.

Next pick a generic point $p_1\in\zeta_1$. Clearly $p_1\not\in\zeta_2,\zeta_3$. By \cref{lem:trisecant}(2), the maps $\psi_{p_1}\colon \zeta_2\to\zeta_{2,p_1}$ and $\psi_{p_1}\colon \zeta_3\to\zeta_{3,p_1}$ are generically 1-to-1. Furthermore, for each $p,\gamma_1$ so that $\zeta_1,\zeta_2,\zeta_3\in\cZ_{p,\gamma_1}$, the line $\aff\{p,p_1\}$ intersects $\zeta_2$ in $a(p,\zeta_2)$ points, intersects $\zeta_3$ in $a(p,\zeta_3)$ points, and does not intersect $\cQ$ in any other point that $p$. (The last is true simply because $\cQ$ is a finite set of points.)

In particular, the map $\psi_{p_1}\colon \zeta_2\to\zeta_{2,p_1}$ has fiber above $\psi_{p_1}(p)$ of size $a(p,\zeta_2)$ and the map $\psi_{p_1}\colon \zeta_3\to\zeta_{3,p_1}$ has fiber above $\psi_{p_1}(p)$ of size $a(p,\zeta_3)$. Since these maps are generically 1-to-1, by \cref{thm:upper-bound-fiber}, we conclude that $\mu_{\psi_{p_1}(p)}(\zeta_{2,p_1})\geq a(p,\zeta_2)$ and $\mu_{\psi_{p_1}(p)}(\zeta_{3,p_1})\geq a(p,\zeta_3)$. Thus by \cref{lem:mult-of-two-planar-curves}, the intersection multiplicity satisfies
\[m_{\psi_{p_1}(p)}(\zeta_{2,p_1},\zeta_{3,p_1})\geq\mu_{\psi_{p_1}(p)}(\zeta_{2,p_1})\mu_{\psi_{p_1}(p)}(\zeta_{3,p_1})\geq a(p,\zeta_2)a(p,\zeta_3).\]

Recall that $\zeta_2,\zeta_3\in\cZ_{p,\gamma_1}$ implies that $\mu_p(\zeta_2)\leq\tfrac12\deg\zeta_2$ and $\mu_p(\zeta_3)\leq\tfrac12\deg\zeta_3$, so \cref{prop:proj-formula-with-singular-origin} implies that $\deg\zeta_2\leq 2a(p,\zeta_2)\deg\gamma_1$ and $\deg\zeta_3\leq 2a(p,\zeta_3)\deg\gamma_1$. Thus we see that the intersection multiplicity is at least $\deg\zeta_2\deg\zeta_3/(4(\deg\gamma_1)^2)$.

However by B\'ezout's theorem (\cref{thm:bezout-with-multiplicity-projective}), we see
\[\sum_{q\in\zeta_{2,p_1}\cap \zeta_{3,p_1}}m_q(\zeta_{2,p_1},\zeta_{3,p_1})\leq\deg\zeta_{2,p_1}\deg\zeta_{3,p_1}\leq\deg\zeta_2\deg\zeta_3.\]
Finally, since distinct points $p\in \cQ$ project to distinct intersection points (and each choice of $p$ determines $\gamma_1$ uniquely), we conclude that
\[\deg\zeta_2\deg\zeta_3\geq\sum_{p,\gamma_1:\zeta_1,\zeta_2,\zeta_3\in\cZ_{p,\gamma_1}}\frac{\deg\zeta_2\deg\zeta_3}{4(\deg\gamma_1)^2}.\qedhere\]
\end{proof}

\begin{lemma}\label{lem:self-intersection-bound-1}
Let $\cZ$ be a collection of irreducible curves in $\CC^3$. 
For $\zeta\in\cZ$,
\[\sum_{p,\gamma_1:\zeta\in\cZ_{p,\gamma_1}\text{ and }a(p,\zeta)\geq 3}\left(\deg \gamma_1\right)^{-2}\leq 18.\]
\end{lemma}
This follows by a similar argument to \cref{lem:self-intersection-bound-3}.

\begin{proof}
We first show that we may assume that $\zeta$ is not contained in any plane. Suppose that $\zeta$ is contained in the plane $F$ and suppose that $p,\gamma_1$ are such that $\zeta\in\cZ_{p,\gamma_1}$ and $a(p,\zeta)\geq 3$. If $p\not\in F$, then every line through $p$ intersects $F$ (and thus $\zeta\subset F$) at most once. This contradicts the fact that $a(p,\zeta)\geq 3$. Thus we may assume that $p\in F$, so $F\cap C(p,\gamma_1)$ is a finite union of lines. However, this intersection contains $\zeta$ which is not a line. Thus we see that the sum is empty if $\zeta$ is contained in a plane.

Let $\cQ$ be the set of points $p$ for which there exists $\gamma_1$ with $\zeta\in\cZ_{p,\gamma_1}$ and $a(p,\zeta)\geq 3$. Given $p\in\cQ$, the curve $\gamma_1$ is determined uniquely: it must be $\gamma_1=\overline{\psi_p(\zeta\setminus\{p\})}$. As before, we first show that $\cQ$ is finite.

For $p_1\in\zeta$, write $\zeta_{p_1}=\overline{\psi_{p_1}(\zeta\setminus\{p_1\})}$. By \cref{lem:trisecant}(1), for a generic point $p_1\in\zeta$, we know that the map $\psi_{p_1}\colon\zeta\setminus\{p_1\}\to\gamma_{p_1}$ is generically 1-to-1. Thus there is a finite collection of lines $\cL_{p_1}$ through $p_1$, the union of the fibers of $\psi_{p_1}$ for which the projection is not 1-to-1. (By \cref{thm:upper-bound-fiber}, these are all fibers over singular points of $\zeta_{p_1}$.)

If $p,\gamma_1$ is such that $\zeta\in\cZ_{p,\gamma_1}$ and $a(p,\zeta)\geq 3$, then a generic line through $p$ in $C(p,\gamma_1)$ is a trisecant, i.e., intersects $\zeta$ in at least three points. This means that for a generic point $p_1\in\zeta$, we have $p\in\cL_{p_1}$. As before, this implies that $\cQ$ is finite, since it is contained in $\cL_{p_1}\cap\cL_{\tilde p_1}$ where $p_1,\tilde p_1\in\zeta$ are two generic points.

Pick a generic point $p_1\in\zeta$. The map $\psi_{p_1}\colon\zeta\setminus\{p_1\}\to\zeta_{p_1}$ is generically 1-to-1 and for each $p,\gamma_1$ so that $\zeta\in \cZ_{p,\gamma_1}$ and $a(p,\zeta)\geq 3$, we have that the line $\aff\{p,p_1\}$ intersects $\zeta$ in $a(p,\zeta)$ points and does not intersect $\cQ$ in any point other than $p$.

In particular, the map $\psi_{p_1}\colon\zeta\setminus\{p_1\}\to\zeta_{p_1}$ has fiber above $\psi_{p_1}(p)$ of size $a(p,\zeta)-1$ (since the intersection at $p_1$ is not counted in the fiber). Then \cref{thm:upper-bound-fiber} gives $\mu_{\psi_{p_1}(p)}(\zeta_{p_1})\geq a(p,\zeta)-1$. By \cref{prop:planar-easy-multiplicity},
\[(\deg\zeta)^2\geq (\deg\zeta_{p_1})^2\geq\sum_{q\in\zeta_{p_1}}\mu_q(\zeta_{p_1})(\mu_q(\zeta_{p_1})-1)\geq \sum_{p\in\cQ}(a(p,\gamma)-1)(a(p,\gamma)-2)\geq\frac29\sum_{p\in\cQ}a(p,\zeta)^2,\]
where the last inequality follows from the assumption that $a(p,\zeta)\geq 3$.

As before, whenever $\zeta\in\cZ_{p,\gamma_1}$, we have $\deg\zeta\leq 2a(p,\zeta)\deg\gamma_1$, so this rearranges to the desired bound.
\end{proof}

We are now ready to prove \cref{claim:curve-cone-incidences}.

\begin{proof}[{Proof of \cref{claim:curve-cone-incidences}}]
We wish to bound
\[\sum_{p,q,\gamma_1}m_{pq}(\gamma_1)\sum_{\zeta\in\cZ_{p,\gamma_1}}\abs{\cP_\zeta}=\sum_{p,\gamma_1}\paren{\sum_{q:(p,q)\in H_1\text{ and }\gamma_1\in\cC^1_{pq}}m_{pq}(\gamma_1)}\paren{\sum_{\zeta\in\cZ_{p,\gamma_1}}\abs{\cP_\zeta}}.\]

We break into three cases. The first two cases are when a single term dominates the final sum. For a given $p,\gamma_1$ the cases are as follows:
\begin{enumerate}[(5a)]
    \item there exists $\zeta_{p,\gamma_1}\in\cZ_{p,\gamma_1}$ with $a(p,\zeta_{p,\gamma_1})\leq 2$ such that
\[\abs{\cP_{\zeta_{p,\gamma_1}}}\geq\frac14\sum_{\zeta\in\cZ_{p,\gamma_1}}|\cP_{\zeta}|;\]
    \item there exists $\zeta_{p,\gamma_1}\in\cZ_{p,\gamma_1}$ with $a(p,\zeta_{p,\gamma_1})\geq 3$ such that
\[\abs{\cP_{\zeta_{p,\gamma_1}}}\geq\frac14\sum_{\zeta\in\cZ_{p,\gamma_1}}|\cP_{\zeta}|;\]
    \item for all $\zeta\in\cZ_{p,\gamma_1}$,\[\abs{\cP_{\zeta}}<\frac14\sum_{\zeta'\in\cZ_{p,\gamma_1}}|\cP_{\zeta'}|.\]
\end{enumerate}

Recall that by definition, for $\zeta\in\cZ_{p,\gamma_1}$, we have $\mu_p(\zeta)\leq\tfrac12\deg\zeta$. Thus \cref{prop:proj-formula-with-singular-origin}, implies that $\deg\zeta\leq 2a(p,\zeta)\deg\gamma_1$.

In case (5a) we have 
\[\sum_{\zeta\in\cZ_{p,\gamma_1}}|\cP_{\zeta}|\leq 4|\cP_{\zeta_{p,\gamma_1}}|\leq 4N^{2/3}\deg\zeta_{p,\gamma_1}\leq 8N^{2/3}a(p,\zeta_{p,\gamma_1})\deg\gamma_1\leq 16N^{2/3}\deg\gamma_1.\]
Thus we can bound the contribution from this case as
\begin{align*}
\sum_{p,\gamma_1:(5a)}\paren{\sum_{q:(p,q)\in H_1\text{ and }\gamma_1\in\cC^1_{pq}}m_{pq}(\gamma_1)}\paren{\sum_{\zeta\in\cZ_{p,\gamma_1}}\abs{\cP_\zeta}}
&\lesssim N^{2/3}\sum_{p,q,\gamma_1}m_{pq}(\gamma_1)\deg\gamma_1\\
&\lesssim N^{2/3}\sum_{(p,q)\in H_1}\deg\cS_{pq},
\end{align*}
by \cref{prop:total-curve-bound}.

In case (5b) we wish to use \cref{lem:self-intersection-bound-1}. We have
\[\sum_{\zeta\in\cZ_{p,\gamma_1}}|\cP_{\zeta}|\leq 4\abs{\cP_{\zeta_{p,\gamma_1}}},\]
so we can bound
\begin{align*}
\sum_{p,\gamma_1:(5b)}&\paren{\sum_{q:(p,q)\in H_1\text{ and }\gamma_1\in\cC^1_{pq}}m_{pq}(\gamma_1)}\paren{\sum_{\zeta\in\cZ_{p,\gamma_1}}\abs{\cP_\zeta}}\\
&\lesssim \sum_{p,\gamma_1:(5b)}\paren{\sum_{q:(p,q)\in H_1\text{ and }\gamma_1\in\cC^1_{pq}}m_{pq}(\gamma_1)}\abs{\cP_{\zeta_{p,\gamma_1}}}\\
&= \sum_{p,\gamma_1:(5b)}\paren{\sum_{q:(p,q)\in H_1\text{ and }\gamma_1\in\cC^1_{pq}}m_{pq}(\gamma_1)\deg\gamma_1}\abs{\cP_{\zeta_{p,\gamma_1}}}(\deg\gamma_1)^{-1}\\
&\leq\paren{\sum_{p,\gamma_1:(5b)}\paren{\sum_{q:(p,q)\in H_1\text{ and }\gamma_1\in\cC^1_{pq}}m_{pq}(\gamma_1)\deg\gamma_1}^{3/2}}^{2/3}\paren{\sum_{p,\gamma_1:(5b)}\abs{\cP_{\zeta_{p,\gamma_1}}}^3(\deg\gamma_1)^{-3}}^{1/3},
\end{align*}
where the last line is H\"older's inequality.

For the first term, we use the easy bound $\sum x_i^{3/2}\leq(\sum x_i)^{3/2}$ together with \cref{prop:total-curve-bound}:
\begin{align*}
\left(\sum_{p,\gamma_1:(5b)}\right.&\left.\left(\sum_{q:(p,q)\in H_1\text{ and }\gamma_1\in\cC^1_{pq}}m_{pq}(\gamma_1)\deg\gamma_1\right)^{3/2}\right)^{2/3}\\
&\leq\paren{\sum_{p\in\cP_1}\paren{\sum_{\gamma_1}\sum_{q:(p,q)\in H\text{ and }\gamma_1\in\cC^1_{pq}}m_{pq}(\gamma_1)\deg\gamma_1}^{3/2}}^{2/3}\\
&\lesssim \paren{\sum_{p\in\cP_1}\paren{\sum_{q:(p,q)\in H}\deg\cS_{pq}}^{3/2}}^{2/3}\\
&\leq |\cP_1|^{2/3}\max_{p\in\cP_1}\sum_{q:(p,q)\in H}\deg\cS_{pq}.
\end{align*}
For the second term, \cref{lem:self-intersection-bound-1} gives
\begin{align*}
\paren{\sum_{p,\gamma_1:(5b)}\abs{\cP_{\zeta_{p,\gamma_1}}}^3(\deg\gamma_1)^{-3}}^{1/3}
&=\paren{\sum_{\zeta\in\cZ}\abs{\cP_\zeta}^3\sum_{p,\gamma_1:\zeta=\zeta_{p,\gamma_1}}(\deg\gamma_1)^{-3}}^{1/3}\\
&\leq\paren{\sum_{\zeta\in\cZ}\abs{\cP_\zeta}^3\sum_{p,\gamma_1:\zeta\in\cZ_{p,\gamma_1}\text{ and }a(p,\zeta)\geq 3}(\deg\gamma_1)^{-3}}^{1/3}\\
&\lesssim\paren{\sum_{\zeta\in\cZ}\abs{\cP_\zeta}^3}^{1/3}\leq\sum_{\zeta\in\cZ}\abs{\cP_\zeta}=|\cP_2|.
\end{align*}
Note we are using  \cref{lem:self-intersection-bound-1} to bound the sum of $(\deg\gamma_1)^{-3}\leq(\deg\gamma_1)^{-2}$. The last line uses the easy bound  $\sum x_i^{3}\leq(\sum x_i)^{3}$.

Finally, in case (5c) we wish to use \cref{lem:self-intersection-bound-3}. First note that by the assumption of this case, we have
\begin{equation}
\label{eq:distinct-triples-dominate}
\sum_{\zeta_1,\zeta_2,\zeta_3\in\cZ_{p,\gamma_1}\text{ distinct}}\abs{\cP_{\zeta_1}}\abs{\cP_{\zeta_2}}\abs{\cP_{\zeta_3}}\gtrsim \sum_{\zeta_1,\zeta_2,\zeta_3\in\cZ_{p,\gamma_1}}\abs{\cP_{\zeta_1}}\abs{\cP_{\zeta_2}}\abs{\cP_{\zeta_3}}=\paren{\sum_{\zeta\in\cZ_{p,\gamma_1}}\abs{\cP_{\zeta}}}^3.
\end{equation}
Indeed, 
\[\sum_{\zeta_1,\zeta_2,\zeta_3\in\cZ_{p,\gamma_1}\text{ distinct}}\abs{\cP_{\zeta_1}}\abs{\cP_{\zeta_2}}\abs{\cP_{\zeta_3}}=\sum_{\zeta_1\in\cZ_{p,\gamma_1}}\paren{\abs{\cP_{\zeta_1}}\sum_{\zeta_2\in\cZ_{p,\gamma_1}\setminus\{\zeta_1\}}\paren{\abs{\cP_{\zeta_2}}\sum_{\zeta_3\in\cZ_{p,\gamma_1}\setminus\{\zeta_1,\zeta_2\}}\abs{\cP_{\zeta_3}}}}.\]
The innermost sum is at least $\tfrac12\sum_{\zeta\in\cZ_{p,\gamma_1}}|\cP_{\zeta}|$ and so the double sum is at least $\tfrac12\cdot\frac34\paren{\sum_{\zeta\in\cZ_{p,\gamma_1}}|\cP_{\zeta}|}^2$.

By a similar application of H\"older's inequality as in case (5b), we have
\begin{align*}
\sum_{p,\gamma_1:(5c)}&\paren{\sum_{q:(p,q)\in H_1\text{ and }\gamma_1\in\cC^1_{pq}}m_{pq}(\gamma_1)}\paren{\sum_{\zeta\in\cZ_{p,\gamma_1}}\abs{\cP_\zeta}}\\
& =\sum_{p,\gamma_1:(5c)}\paren{\sum_{q:(p,q)\in H_1\text{ and }\gamma_1\in\cC^1_{pq}}m_{pq}(\gamma_1)\deg\gamma_1}\paren{\sum_{\zeta\in\cZ_{p,\gamma_1}}\abs{\cP_\zeta}(\deg\gamma_1)^{-1}}\\
&\leq \paren{\sum_{p,\gamma_1:(5c)}\paren{\sum_{q:(p,q)\in H_1\text{ and }\gamma_1\in\cC^1_{pq}}m_{pq}(\gamma_1)\deg\gamma_1}^{3/2}}^{2/3}\\
&\qquad\qquad\cdot \paren{\sum_{p,\gamma_1:(5c)}\paren{\sum_{\zeta\in\cZ_{p,\gamma_1}}\abs{\cP_\zeta}(\deg\gamma_1)^{-1}}^{3}}^{1/3}.
\end{align*}

The first term is bounded as in case (5b). For the second term, \cref{eq:distinct-triples-dominate} combined with \cref{lem:self-intersection-bound-3} gives
\begin{align*}
\left(\sum_{p,\gamma_1:(5c)}\right.&\left.\paren{\sum_{\zeta\in\cZ_{p,\gamma_1}}\abs{\cP_\zeta}(\deg\gamma_1)^{-1}}^{3}\right)^{1/3}\\
&\lesssim \paren{\sum_{p,\gamma_1}(\deg\gamma_1)^{-3}\sum_{\zeta_1,\zeta_2,\zeta_3\in\cZ_{p,\gamma_1}\text{ distinct}}\abs{\cP_{\zeta_1}}\abs{\cP_{\zeta_2}}\abs{\cP_{\zeta_3}}}^{1/3}\\
&= \paren{\sum_{\zeta_1,\zeta_2,\zeta_3\in\cZ\text{ distinct}}\abs{\cP_{\zeta_1}}\abs{\cP_{\zeta_2}}\abs{\cP_{\zeta_3}}\sum_{p,\gamma_1:\zeta_1,\zeta_2,\zeta_3\in\cZ_{p,\gamma_1}}(\deg\gamma_1)^{-3}}^{1/3}\\
&\lesssim \paren{\sum_{\zeta_1,\zeta_2,\zeta_3\in\cZ\text{ distinct}}\abs{\cP_{\zeta_1}}\abs{\cP_{\zeta_2}}\abs{\cP_{\zeta_3}}}^{1/3}\leq \sum_{\zeta\in\cZ}\abs{\cP_\zeta}=|\cP_2|.
\end{align*}

Combining everything gives the bound
\[\sum_{p,q,\gamma_1}m_{pq}(\gamma_1)\sum_{\zeta\in\cZ_{p,\gamma_1}}\abs{\cP_\zeta}\lesssim N^{2/3}\sum_{(p,q)\in H_1}\deg\cS_{pq}+|\cP_1|^{2/3}|\cP_2|\max_{p\in\cP_1}\sum_{q:(p,q)\in H}\deg\cS_{pq}.\qedhere\]
\end{proof}

\subsection{Exceptional quadruples}
To efficiently apply \cref{lem:ruled-ruled-dyadic-general}, we partition $\cP$ into sets where $\sum_{q:(p,q)\in H}\deg\cS_{pq}$ is approximately constant as a function of $p$. Unfortunately the factor of $|\cP_2|$ is too large if we apply this result with $\cP_2=\cP$. Instead we apply this partition both to the $p$ and the $p'$ coordinate of $(p,q,p',q')\in\cE$.  Then \cref{lem:ruled-ruled-dyadic-general} gives a strong bound when $p\in\cP_1$ and $p'\in\cP_2$ with $|\cP_1|\geq|\cP_2|$. Applying the symmetric bound with the roles of $(p,q)$ and $(p',q')$ swapped handles the case $|\cP_2|\geq|\cP_1|$. This accounts for quadruples $(p,q,p',q')\in\cE$ such that there exists $\gamma\in\tilde \cC_{pq}\setminus\cC^v_{pq}$ where $\psi_{(p,q)}(p',q')\in\gamma$ and such that there also exists a curve $\gamma'\in\tilde\cC_{p'q'}\setminus\cC^v_{p'q'}$.

Using the symmetry between $p,q$, we can also deal with the quadruples $(p,q,p',q')\in\cE$ for which there exist curves $\gamma\in\tilde\cC_{pq}\setminus\cC^h_{pq}$ with $\psi_{(p,q)}(p',q')\in\gamma$ and $\gamma'\in\tilde\cC_{p'q'}\setminus\cC^h_{p'q'}$ with $\psi_{(p',q')}(p,q)\in\gamma'$. This handles all of $\cE$ except for the exceptional quadruples: those where $\gamma\in\cC^h_{pq}$ and $\gamma'\in\cC^v_{p'q'}$ or those where $\gamma\in\cC^v_{pq}$ and $\gamma'\in\cC^h_{p'q'}$.

The next lemma deals with the exceptional quadruples. We give a geometric argument that essentially shows that the contribution from the exceptional quadruples is dominated by the terms where $\cP_1=\cP_2$.

\begin{proposition}
\label{lem:ruled-ruled-vertical-horizontal}
Assume the same setup as \cref{thm:ruled-ruled}. Given $\cP_1\subseteq\cP$, let $\cE'\subseteq\cE$ be the set of quadruples $(p,q,p',q')\in\cE$ so that $p\in\cP_1$ and so that $\psi_{(p,q)}(p',q')\in\gamma$ for some $\gamma\in\cC^{h}_{pq}$ and $\psi_{(p',q')}(p,q)\in\gamma'$ for some $\gamma'\in\cC^{v}_{p',q'}$.
Then
\[|\cE'|\lesssim CN^{2/3}\sum_{(p,q)\in H}\deg \cS_{pq}+C|\cP_1|^{5/3}\max_{p\in\cP_1:(p,q)\in H}\deg\cS_{pq}.\]
\end{proposition}

\begin{proof}
Write $H_1=(\cP_1\times\cP)\cap H$.

For $(p',q')\in H$, let $\cC^v_{p'q'} = \cC^{v,\mathrm{rich}}_{p'q'}\sqcup \cC^{v,\mathrm{poor}}_{p'q'}$ where $\gamma'\in \cC^{v,\mathrm{rich}}_{p'q'}$ if and only if the point $\pi_1(\gamma')\in\PP_{p'}^2$ corresponds to a line through $p'$ that passes through at least two points in $\cP_1 \setminus \{p'\}$; that is, this line passes through at least one point of $\cP_1$ other than $p,p'$.

Now for any $(p,q)\in H_1$, let $\cE^{\mathrm{rich}}_{pq}$ be the pairs $(p',q')$ so that $(p,q,p',q')\in\cE'$ and $\psi_{(p',q')}(p,q)\in\gamma'$ for some $\gamma'\in\cC^{v,\mathrm{rich}}_{p'q'}$. Define $\cE^{\textup{poor}}_{pq}$ analogously. We first handle the contribution from the poor curves.

For $(p',q')\in H$ and $\gamma'\in\cC^{v,\mathrm{poor}}_{p'q'}$ and any $(p,q)$ with $(p',q')\in\cE^{\mathrm{poor}}_{pq}$, the definition of $\cC^{v,\mathrm{poor}}_{p'q'}$ implies that $p$ is the unique point of $\cP_1\setminus\{p'\}$ on the line through $p'$ corresponding to the point $\pi_1(\gamma')\in\PP_{p'}^2$. Now $q$ satisfies $\norm{p-p'}^2=\norm{q-q'}^2$ and $\psi_{q'}(q)\in\overline{\pi_2(\gamma'\setminus\PP^2_{p'})}=:\gamma_2'$, implying that $q$ lies in the intersection of the cone $C(q',\gamma_2')$ and the sphere centered at $q'$ with radius $\norm{p-p'}$.\footnote{To be precise, this sphere is the complex surface $Z\paren{\sum_{i\in[3]}(x_i-q_i')^2-(p_i-p_i')^2}$.} This intersection is a curve of degree $\lesssim \deg\gamma_2'\leq \deg\gamma'$. By hypothesis, there are $\lesssim N^{2/3}\deg\gamma'$ points of $\cP$ on this curve. Thus we see that
\[\sum_{(p,q)\in H_1}\abs{\cE_{pq}^{\mathrm{poor}}}\lesssim\sum_{(p',q')\in H}\sum_{\gamma'\in\cC^{v,\mathrm{poor}}_{p'q'}}N^{2/3}\deg\gamma'\lesssim N^{2/3}\sum_{(p',q')\in H}\deg\cS_{p'q'},\]
where the last inequality follows by \cref{prop:pos-to-phys}.

We next study the rich contribution. Consider $\psi_p(\pi_1(p',q'))=\psi_p(p')$ for some $(p',q')\in\cE^{\mathrm{rich}}_{pq}$. A priori, we know that $\psi_p(p')\in\psi_p(\cP\setminus\{p\})$. However, we claim that this point actually lies in $\psi_p(\cP_1\setminus\{p\})$. To see this, note that there exists $\gamma'\in \cC^{v,\textup{rich}}_{p'q'}$ with $\psi_{p',q'}(p,q)\in \gamma'$. By the definition of $\cC^{v,\mathrm{rich}}_{p'q'}$, we know that $\pi_1(\gamma')$ corresponds to a line through $p'$ that contains at least two points of $\cP_1\setminus\{p'\}$. In particular, it contains some point $p''\in \cP_1\setminus \{p\}$, and so $\psi_p(p')=\psi_p(p'')\in \psi_p(\cP_1\setminus \{p\})$.

For each $\gamma\in \cC^h_{pq}$, let $\cE^{\mathrm{rich}}_{pq}(\gamma)\subseteq\cE^{\mathrm{rich}}_{pq}$ be the set of pairs $(p',q')$ with $\psi_{(p,q)}(p',q')\in \gamma$. Clearly we have
\[\abs{\cE^{\mathrm{rich}}_{pq}}\leq\sum_{\gamma\in\cC^h_{pq}}\abs{\cE^{\mathrm{rich}}_{pq}(\gamma)}.\]

Given $(p,q)\in H_1$ and $\gamma\in\cC^h_{pq}$, we claim that there are at most $2C$ choices for $(p',q')\in\cE^{\mathrm{rich}}_{pq}(\gamma)$ which achieve the same value of $\psi_p(p')$. Suppose that $\psi_p(p')=[d_1]$ for some $d_1\in\R^3\setminus\{0\}$. (We have that $d_1$ is real since both $p,p'\in\cP$ are real points.) Without loss of generality, assume that $\norm{d_1}=1$. Since $\psi_{(p,q)}(p',q')\in\gamma$, we know that $\psi_q(q')\in\pi_2(\gamma)$. Since $\gamma\in\cC^h_{pq}$, we have that $\pi_2(\gamma)$ is a single point, say $[d_2]$ for some $d_2\in\R^3\setminus\{0\}$. Without loss of generality, assume that $\norm{d_2}=1$. Now $(p',q')=(p+d_1s,q+d_2t)$ for some $s,t\in\R$. However, since $\norm{p-p'}=\norm{q-q'}$, we conclude that $s=\pm t$. Thus $(p',q')$ lies on one of two (isotropic) lines through $(p,q)$. By hypothesis, each of these lines is at most $C$-rich for $H$, implying that there are at most $2C$ choices for $(p',q')$.

Since $\psi_p(p')\in \psi_p(\cP_1\setminus\{p\})\cap \overline{\pi_1(\gamma\setminus\PP_q^2)}$, we conclude that
\[\sum_{\gamma\in\cC^h_{pq}}\abs{\cE^{\mathrm{rich}}_{pq}(\gamma)}\leq \sum_{\gamma\in\cC^h_{pq}}2C\abs{\psi_p(\cP_1\setminus\{p\})\cap \overline{\pi_1(\gamma\setminus\PP_q^2)}}\leq 2C\sum_{\gamma_1\in\cC^1_{pq}}m_{pq}(\gamma_1)\abs{\psi_p(\cP_1\setminus\{p\})\cap\gamma_1}.\]
Here we use \cref{prop:disjoint-exceptional} again to show that if $\gamma\in\cC^h_{pq}$, then $\gamma\in\tilde\cC_{pq}\setminus\cC^v_{pq}$. (Note that the second inequality is lossy, as we are using $m_{pq}(\gamma_1)$ as an upper bound on the number of curves $\gamma$ which project to the curve $\gamma_1$.)

So far we have shown that
\[\abs{\cE'}\lesssim N^{2/3}\sum_{(p',q')\in H}\deg\cS_{p'q'}+C\sum_{(p,q)\in H_1}\sum_{\gamma_1\in\cC^1_{pq}}m_{pq}(\gamma_1)\abs{\psi_p(\cP_1\setminus\{p\})\cap \gamma_1}.\]
To finish, we bound the second term by \cref{prop:ruled-ruled-curve-cone-incidences} applied with $\cP_1=\cP_2$.
\end{proof}

\subsection{Dyadic decomposition}
\label{ssec:dyadic-decomposition}

\begin{proof}[Proof of \cref{thm:ruled-ruled}]
For each $(p,q,p',q')\in \cE$, by \cref{prop:phys-to-C3-edge-case} we have $\psi_{(p,q)}(p',q')\in \gamma$ for some $\gamma\in \tilde\cC_{pq}$ and $\psi_{(p',q')}(p,q)\in \gamma'$ for some $\gamma'\in \tilde \cC_{p'q'}$. By \cref{prop:disjoint-exceptional}, at most one of the projections of each $\gamma,\gamma'$ can be a point. Thus we can cover $\cE$ as $\cE^{(1,1)}\cup\cE^{(2,2)}\cup\cE^{(1,2)}\cup\cE^{(2,1)}$ where
\begin{itemize}
    \item $(p,q,p',q')\in\cE^{(1,1)}$ if $\gamma\in\tilde\cC_{pq}\setminus\cC^v_{pq}$ and $\gamma'\in\tilde\cC_{p'q'}\setminus\cC^v_{p'q'}$; 
    \item $(p,q,p',q')\in\cE^{(2,2)}$ if $\gamma\in\tilde\cC_{pq}\setminus\cC^h_{pq}$ and $\gamma'\in\tilde\cC_{p'q'}\setminus\cC^h_{p'q'}$; 
    \item $(p,q,p',q')\in\cE^{(1,2)}$ if $\gamma\in\cC^h_{pq}$ and $\gamma'\in\cC^v_{p'q'}$; and
    \item $(p,q,p',q')\in\cE^{(2,1)}$ if $\gamma\in\cC^v_{pq}$ and $\gamma'\in\cC^h_{p'q'}$.
\end{itemize}

For each $p\in \cP$, define
\[t(p) = \sum_{q:(p,q)\in H}\deg \cS_{pq}.\]
By hypothesis, $\deg\cS_{pq}\lesssim N^{2/3}$, so $t(p)\leq C_0N^{5/3}$ for some absolute constant $C_0$. Consider the dyadic partition $\cP=\bigsqcup_{i=1}^T\cP_i$ where $p\in\cP_i$ if $2^{-i}C_0N^{5/3}<t(p)\leq 2^{1-i}C_0N^{5/3}$. Here $T\sim\log N$. Set $H_i=H\cap(\cP_i\times\cP)$ for each $1\leq i\leq T$. 

Then 
\[\abs{\cE^{(1,1)}} = \sum_{1\leq i, j\leq T}\abs{\left\{(p,q,p',q')\in \cE^{(1,1)}: p\in \cP_i, p'\in \cP_j\right\}}.\]
Denote the sets on the right-hand side by $\cE^{(1,1)}_{i,j}$. For each $i,j$, by \cref{lem:ruled-ruled-dyadic-general}, we have
\begin{align*}
\abs{\cE^{(1,1)}_{i,j}}
&\lesssim CN^{2/3}\sum_{(p,q)\in H_i}\deg \cS_{pq}+C\abs{P_i}^{2/3}\abs{P_j}\frac{N^{5/3}}{2^i}.
\end{align*}

Swapping the role of $(p,q)$ with $(p',q')$ and averaging, we can strengthen this bound to
\begin{align*}
&\lesssim CN^{2/3}\left(\sum_{(p,q)\in H_i}\deg \cS_{pq}+\sum_{(p',q')\in H_j}\deg \cS_{p'q'}\right)+C\sqrt{\left(\abs{\cP_i}^{2/3}\abs{\cP_j}\frac{N^{5/3}}{2^i}\right)\cdot \left(\abs{\cP_i}\abs{\cP_j}^{2/3}\frac{N^{5/3}}{2^j}\right)}\\
&\qquad\qquad= CN^{2/3}\left(\sum_{(p,q)\in H_i}\deg \cS_{pq}+\sum_{(p',q')\in H_j}\deg \cS_{p'q'}\right)+C\abs{\cP_i}^{5/6}\abs{\cP_j}^{5/6}\frac{N^{5/3}}{2^{(i+j)/2}}.
\end{align*}
Note that
\[\sum_{1\leq i,j\leq T}CN^{2/3}\left(\sum_{(p,q)\in H_i}\deg \cS_{pq}+\sum_{(p',q')\in H_j}\deg \cS_{p'q'}\right)\lesssim CN^{2/3}\log N\sum_{(p,q)\in H}\deg \cS_{pq}.\]
Next we bound the contribution of the second term. This can be written as
\begin{equation}
\label{eq:e11-second-term}
\sum_{1\leq i,j\leq T}C\abs{\cP_i}^{5/6}\abs{\cP_j}^{5/6}\frac{N^{5/3}}{2^{(i+j)/2}} = CN^{5/3}\left(\sum_{i=1}^{T}\abs{\cP_i}^{5/6}2^{-i/2}\right)^2.
\end{equation}

By two applications of H\"older's inequality, we have
\begin{align*}
\paren{\sum_{i=1}^{T}\abs{\cP_i}^{5/6}2^{-i/2}}^2
&\leq \left(\sum_{i=1}^{T}\abs{\cP_i}2^{-i}\right)\left(\sum_{i=1}^{T}\abs{\cP_i}^{2/3}\right)\\
&\leq \left(\sum_{i=1}^{T}\abs{\cP_i}2^{-i}\right)\paren{\sum_{i=1}^T\abs{\cP_i}}^{2/3}\paren{\sum_{i=1}^{T}1}^{1/3}.
\end{align*}
Now
\begin{equation}
\label{eq:tp-dyadic-sum}
\sum_{(p,q)\in H}\deg\cS_{pq}=\sum_{p\in\cP}t(p)>\sum_{i=1}^{T} 2^{-i}C_0N^{5/3}|\cP_i|,
\end{equation}
implying that the first parenthetical term in the expression on the right-hand side above is $\lesssim N^{-5/3}\sum_{(p,q)\in H}\deg\cS_{pq}$. The second term is $|\cP|=N$, while the third is $T\sim\log N$. Thus the expression in \cref{eq:e11-second-term} can be bounded by $\lesssim CN^{2/3}(\log N)^{1/3}\sum_{(p,q)\in H}\deg\cS_{pq}$.

We have shown the desired bound on $\abs{\cE^{(1,1)}}$; symmetrically the same bound holds for $\abs{\cE^{(2,2)}}$.

For $\abs{\cE^{(1,2)}}$, we use the same dyadic decomposition, just on the $p$ coordinate. Define
\[\cE^{(1,2)}_{i}=\set{(p,q,p',q')\in\cE^{(1,2)}:p\in\cP_i}.\]
Then by \cref{lem:ruled-ruled-vertical-horizontal},
\[\abs{\cE^{(1,2)}_{i}}\lesssim CN^{2/3}\sum_{(p,q)\in H}\deg\cS_{pq}+C\abs{\cP_i}^{5/3} \frac{N^{5/3}}{2^i}.\]
Summing over $i$ gives
\[\abs{\cE^{(1,2)}}\lesssim CN^{2/3}\log N\sum_{(p,q)\in H}\deg\cS_{pq}+CN^{5/3}\sum_{i=1}^T\frac{\abs{\cP_i}^{5/3}}{2^i}.\]
Bounding $\abs{\cP_i}^{5/3}2^{-i}\leq N^{2/3}\abs{\cP_i}2^{-i}$ and using \cref{eq:tp-dyadic-sum} again, we see that the second term is bounded by $CN^{2/3}\sum_{(p,q)\in H}\deg\cS_{pq}$. Thus $\abs{\cE^{(1,2)}}$, and symmetrically $\abs{\cE^{(2,1)}}$, are of the desired size.
\end{proof}

\newpage

\part{Very rich 3-flats}
\label{part:iii}
In this part we complete the proof of \cref{thm:main}. By \cref{thm:beck-trick}, it suffices to prove a bound of $N^{10/3+o(1)}$ on the number of isotropic lines $\ell^{\Phy}\subset E_o^{\Phy}(\R)$ that lie in at least 2 and at most $C$ points of $\phi(\cP\times\cP)$. The results of \cref{part:i,part:ii} imply that such a bound holds for all isotropic lines except those contained in a very rich 3-flat: some $F_\tau^{\Phy}$ with $\tau\in E^-_g$ that is $k$-rich for $k\gg N^{2/3}$, meaning that $\tau(p)=q$ for at least $k$ pairs $(p,q)\in\cP^2$. Applying the symmetric argument in negative space, it remains to study lines that are contained in both such an $F_\tau^{\Phy}$ and an $F_\rho^{\Phy}$ for $\rho\in E^+_g$ and $\tau\in E^-_g$. This deduction is stated as \cref{thm:input}.

The remaining task is to bound the number of these isotropic lines which are contained in a very rich positive and negative 3-flat. This problem has significantly more geometric structure which we heavily exploit.

In \cref{sec:incidence-geometry-r3}, we introduce some tools for incidence geometry in $\R^3$ that will be used in the proof. Some of these are straightforward consequences of results appearing in the literature and some are novel. In \cref{sec:degenerate}, we bound the number of isotropic lines in several cases that are geometrically structured. Finally, in \cref{sec:main-argument}, we put everything together to prove \cref{thm:main}.

\section{Geometric preliminaries}
\label{sec:incidence-geometry-r3}

\subsection{Incidence geometry in \texorpdfstring{$\R^3$}{R3}}

In this section we collect some standard incidence geometry results in $\R^3$ and prove some less-standard results. 
First is the Guth--Katz point-line incidence bound \cite{GK15}.

\begin{theorem}[{\cite[Theorem 12.1]{Guth16}}]
\label{thm:guth-katz-point-line-incidences}
Let $\cP,\cL$ be sets of points and lines in $\R^3$. If at most $B$ lines of $\cL$ lie in any plane, then
\[I(\cP,\cL)\lesssim \abs{\cP}^{1/2}\abs{\cL}^{3/4}+B^{1/3}\abs{\cP}^{2/3}\abs{\cL}^{1/3}+\abs{\cP}+\abs{\cL}.\]
\end{theorem}

\begin{remark}
Note that the theorem as stated in \cite{Guth16} has the additional hypothesis that $B\geq \abs{\cL}^{1/2}$. This is unnecessary, as applying the result with $B'=\max\{B,\abs{\cL}^{1/2}\}$, one recovers the desired bound with an additional term of the form $\abs{\cP}^{2/3}\abs{\cL}^{1/2}$. However this term is dominated by the other terms as it is equal to $(\abs{\cP}^{1/2}\abs{\cL}^{3/4})^{2/3}(\abs{\cP})^{1/3}$.
\end{remark}

\cref{thm:guth-katz-point-line-incidences} implies an estimate on the number of rich lines determined by a set of points that does not concentrate too much on any plane.

\begin{corollary}
\label{thm:rich-lines-r3}
Let $\cP$ be a set of $N$ points in $\R^3$ with at most $A$ in any plane. For $k\geq 2$, let $\cL$ be the set of $k$-rich lines.
Then
\[I(\cP,\cL)\lesssim N^{1/2}\abs{\cL}^{3/4}+A^{2/3}\abs{\cL}^{1/3}N^{2/3}k^{-1}+A^{1/3}\abs{\cL}^{1/3}N^{2/3}k^{-1/3}+\abs{\cL}+N\]
and
\[\abs{\cL}\lesssim \frac{N^2}{k^4}+\frac{NA}{k^3}+\frac{N}{k}.\]
\end{corollary}

\begin{proof}
By the Szemer\'edi--Trotter theorem \cite{ST83}, since each plane contains at most $A$ points of $\cP$, it also contains at most $B:=C(A^2/k^3+A/k)$ lines of $\cL$ for some absolute constant $C$. Now \cref{thm:guth-katz-point-line-incidences} implies
\begin{align*}
k\abs{\cL}&\leq I(\cP,\cL)\lesssim N^{1/2}\abs{\cL}^{3/4}+B^{1/3}\abs{\cL}^{1/3}N^{2/3}+\abs{\cL}+N\\
&\lesssim N^{1/2}\abs{\cL}^{3/4}+A^{2/3}\abs{\cL}^{1/3}N^{2/3}k^{-1}+A^{1/3}\abs{\cL}^{1/3}N^{2/3}k^{-1/3}+\abs{\cL}+N.
\end{align*}
This proves the first inequality, and it remains to prove the second one.

For $k\lesssim 1$, the trivial bound of $\abs{\cL}\leq N^2$ suffices. Otherwise, assume that the hidden constant in the above inequality is at most $k/2$. Thus we have
\[k\abs{\cL}\lesssim N^{1/2}\abs{\cL}^{3/4}+A^{2/3}\abs{\cL}^{1/3}N^{2/3}k^{-1}+A^{1/3}\abs{\cL}^{1/3}N^{2/3}k^{-1/3}+N,\]
which rearranges to
\[\abs{\cL}\lesssim\frac{N^2}{k^4}+\frac{AN}{k^3}+\frac{A^{1/2}N}{k^2}+\frac{N}{k}.\]
The third term is the geometric mean of the second and fourth terms so we can remove it, giving the desired bound.
\end{proof}

The next result is a small strengthening of Beck's theorem that allow one to switch between a set of points and a large subset of it.

\begin{lemma}
\label{thm:beck-variant}
There exist constants $C,c>0$ such that the following holds for all $d,N\ge 2$ and $K\geq 1$ with $K\leq N^3$. Let $\cP,\cQ \subset \R^d$ be point sets with $|\cP| = N$ and $|\cQ| \le KN$ such that $|\cP \cap \ell| \le N/2$ for every line $\ell$. Then there are at least $cN^2$ lines $\ell$ with  $|\ell \cap \cP| \geq 2$ and $|\ell\cap \cQ| \leq C K^{2/3}$.
\end{lemma}

\begin{proof}
By Beck's theorem (\cref{thm:beck}) there are at least $c_BN^2$ lines that contain at least 2 points of $\cP$ for some constant $c_B>0$. On the other hand, the number of lines which are at least $CK^{2/3}$-rich for $\cQ$ is upper bounded by Szemer\'edi--Trotter theorem by 
\[\lesssim \frac{|\cQ|^2}{(CK^{2/3})^3}+ \frac{|\cQ|}{CK^{2/3}}\leq\frac{N^2}{C^3}+\frac{K^{1/3}N}{C}\leq\frac{2N^2}{C}.\]
The last inequality holds for any $C\geq 1$ as $K\leq N^3$. Now taking $C$ large enough, this number of lines will be at most $(c_B/2)N^2$.

Therefore there are at least $c_BN^2-(c_B/2)N^2$ lines which contain at least 2 points of $\cP$ and at most $CK^{2/3}$ points of $\cQ$.
\end{proof}

Next we develop several structural results about line-line incidences and point-plane incidences in $\R^3$. Using the techniques of Guth and Katz, de Zeeuw gave a bipartite incidence bound between two sets of lines in $\R^3$ \cite{deZ16}. He used this bound to give a short proof of a result of Rudnev which bounds point-plane incidences \cite{Rud18}. Using similar techniques, we can get a structural result on configurations of points and planes with many incidences.

\begin{theorem}[{cf. \cite[Theorem 1.1]{deZ16}}]
\label{thm:point-plane-bipartite-graph-decomp}
Let $\cP,\cF$ be sets of points and planes in $\R^3$. Say that $\{|\cP|,|\cF|\}=\{L,M\}$ with $L\leq M$. Then the incidence graph of $\cP,\cF$ can be decomposed into $O(L^{1/2})$ vertex-disjoint complete bipartite graphs and at most $O(L^{1/2}M)$ additional edges.  
\end{theorem}

\begin{proof}
By \cite[Lemma 2.1]{deZ16}, we can map the points $\cP$ to a set of lines $\cL$ in $\R^3$ and the planes $\cF$ to a set of lines $\cM$ in $\R^3$ such that $p\in\pi$ if and only if the line corresponding to $p$ intersects the line corresponding to $\pi$. By inspecting this mapping, we see that no line in $\cL$ is parallel to a line in $\cM$. Indeed, the points map to lines in direction $(u,0,v)$ with $u\neq 0$ and the planes map to lines pointing in direction $(0,s,t)$.

Without loss of generality, assume that $|\cL|=L\leq |\cM|=M$. Let $G=(\cL\sqcup\cM,E)$ be the incidence graph. By parameter counting (see \cref{cor:Walsh-param-counting}), there is a nonzero polynomial $f\in\R[x_1,x_2,x_3]$ which vanishes identically on $\cL$ and satisfies $\deg f\leq 3L^{1/2}$. Let $\tilde f\in\R[x_0,x_1,x_2,x_3]\subset\CC[x_0,x_1,x_2,x_3]$ be the homogenization of $f$. Thus we view $Z(\tilde f)$ as a surface in $\PP_{\CC}^3$. Write $\tilde\cL,\tilde\cM$ for the set of (complex projective) lines in $\PP^3_{\CC}$ produced by taking the Zariski-closure of each line $\ell\in\cL,\cM$, respectively. Clearly $\tilde\cL$ is contained in $Z(\tilde f)$. Furthermore, a line $\ell\in\tilde \cL$ intersects a line $\ell'\in\tilde\cM$ if and only if the corresponding lines in $\cL$ and $\cM$ intersect. This is because no line in $\cL$ is parallel to a line in $\cM$, so the corresponding lines in $\tilde\cL,\tilde\cM$ cannot intersect at infinity. Furthermore, since $\ell,\ell'$ are defined over the reals, their intersection point must also be real.

Now decompose $Z(\tilde f)$ into irreducible components as $V_1\cup\cdots\cup V_t$. For $i\in[t]$, write $\cL_i\subseteq\tilde\cL$ and $\cM_i\subseteq\tilde\cM$ for the lines contained in $V_i$ but not $V_1\cup\cdots\cup V_{i-1}$. Note that $\bigsqcup_{i\in[t]} \cL_i=\tilde\cL$. Defining $\cM_\infty=\tilde\cM\setminus\bigcup_{i\in[r]} \cM_i$, we also have $\bigsqcup_{i\in[t]\cup\{\infty\}} \cM_i=\tilde\cM$.

Define $I$ to be the set of $i\in[t]$ such that $V_i$ is a plane, regulus, or cone. If $V_i$ is a plane or cone, then $G[\cL_i\times\cM_i]$ is complete bipartite since every pair of lines a plane or cone intersects by \cref{thm:classical-ruled-surface-incidences}(4). If $V_i$ is a regulus, write $\cL_i=\cL_{i,1}\sqcup\cL_{i,2}$ and $\cM_i=\cM_{i,1}\sqcup\cM_{i,2}$ where these are the sets of lines contained in the first ruling and the second ruling, respectively. Note that $G[\cL_{i,\epsilon}\times\cM_{i,3-\epsilon}]$ is complete bipartite for $\epsilon\in\{1,2\}$ while $G[\cL_{i,\epsilon}\times\cM_{i,\epsilon}]$ is empty by \cref{thm:classical-ruled-surface-incidences}(3).

Thus we have found at most $\deg f\leq 3L^{1/2}$ vertex-disjoint complete bipartite graphs in $G$.

Now we bound the remaining edges. We claim that are at most $6L^{1/2}M$ incidences not captured in $\bigcup_{i\in[t]} G[\cL_i\times\cM_i]$. To see this, consider an incidence between $\ell\in\cL_i$ and $\ell'\in\cM_j$. If $i<j$, we have $\ell\subset V_i$ while $\ell'\not\subset V_i$. By B\'ezout's theorem, $|\ell'\cap V_i|\leq\deg V_i$. Furthermore, since the lines of $\tilde\cL$ are pairwise skew, $\ell$ is uniquely determined by this intersection points. Thus there are at most $\deg V_i |\cM_j|$ incidences in this case. Similarly, if $i>j$, there are at most $\deg V_j|\cL_i|$ incidences. In total, the number of these incidences is bounded by
\[\sum_{i,j}\deg V_i(|\cL_j|+|\cM_j|)=\deg f(|\tilde\cL|+|\tilde\cM|)\leq 6L^{1/2}M.\]

Next we count incidences in $G[\cL_i\times\cM_i]$ for $i\not\in I$. If $i\not\in I$ and $V_i$ is unruled, then it contributes at most $66(\deg V_i)^3$ incidences by \cref{thm:classical-ruled-surface-incidences}(1). If $i\not\in I$ and $V_i$ is ruled, then $V_i$ must be a singly ruled surface that is not a cone. By \cref{thm:classical-ruled-surface-incidences}(2), all but at most 2 lines in $V_i$ intersect at most $\deg V_i$ other lines in $V_i$. Thus $V_i$ contributes at most $2\max\{|\cL_i|,|\cM_i|\}+(|\cL_i|+|\cM_i|)\deg V_i$ incidences.

Summing all of these, the number of uncaptured incidences is bounded by
\[6L^{1/2}M+\sum_i 66(\deg V_i)^3+\sum_i(L_i+M_i)\deg V_i\lesssim L^{1/2}M.\qedhere\]
\end{proof}

We remark that any complete bipartite graph in the point-plane incidence graph (where both parts have size at least 2) corresponds to a unique line; this line contains the points and is contained in the planes which correspond to the vertices of the complete bipartite graph.

The next result has a similar flavor: it roughly states that in $\R^3$, incidences between points and planes are largely dominated by the contribution from lines. Unlike \cref{thm:point-plane-bipartite-graph-decomp}, the following result uses real polynomial partitioning to obtain better bounds.

\begin{lemma}\label{lem:pair-covering-structure-new}
There exists some an constant $C$ such that the following holds for any $\eta\in(0,1/2)$. Let $\cP \subset \R^3$ be a set of size $N$ such that each plane contains at most $\eta N$ points of $\cP$. Let $\cF$ be the set of $k$-rich planes for $\cP$ where $k\geq C\eta^{-3}N^{1/2}$. For each $F\in\cF$, define $\cP_F\subseteq \cP\cap F$ to be the set of points $p\in\cP\cap F$ for which there does not exist an $N^{1/2}$-rich line $\ell$ that contains $p$ and lies in $F$. Then
\[\sum_{F\in\cF}|\cP_F|^2\lesssim \eta N^2.\]
\end{lemma}

\begin{proof}
Define $\cF=\cF_{\textrm{sm}}\sqcup\cF_{\textrm{lg}}$ where $F\in\cF_{\textrm{sm}}$ if $|\cP_F|\leq\eta |\cP\cap F|$ and $F\in\cF_{\textrm{lg}}$ otherwise. Write $\cL$ for the set of $N^{1/2}$-rich lines. By the Szemer\'edi--Trotter theorem (or double counting), we have $|\cL|\lesssim N^{1/2}$ and $I(\cP,\cL)\lesssim N$.

To bound the contribution from $\cF_{\textrm{sm}}$, note that for $F\in\cF_{\textrm{sm}}$ we have 
\[|\cP_F|^2\leq\frac{\eta}{1-\eta}|\cP_F|(|\cP\cap F|-|\cP_F|)\leq 2\eta |\cP_F|(|\cP\cap F|-|\cP_F|).\]
Thus we can write
\[\sum_{F\in\cF_{\textrm{sm}}}|\cP_F|^2\leq 2\eta \abs{\set{(p,p',F)\in\cP\times\cP\times\cF_{\textrm{sm}}:p\in\cP_F,p'\in(\cP\cap F)\setminus\cP_F}}.\]
By definition, for each $(p,p',F)$ as above, there exists an $N^{1/2}$-rich line $\ell\in\cL$ which passes through $p'$ and lies in $F$. Furthermore, it follows from the definition of $\cP_F$ that $\ell$ does not pass through $p$. Thus $F$ can be recovered from $(p,p',\ell)$, as it is the affine span of $p\cup\ell$. Therefore we can bound the size of the set on the right-hand side by
\[\abs{\set{(p,p',\ell)\in\cP\times\cP\times\cL:p'\in\ell}}\leq N\cdot I(\cP,\cL)\lesssim N^2.\]
We conclude that 
\[\sum_{F\in\cF_{\textrm{sm}}}|\cP_F|^2\lesssim \eta N^2.\]

To bound the contribution from $\cF_{\textrm{lg}}$, we use real polynomial partitioning. Set $D=\eta^{-1}\geq 1$. By real polynomial partitioning \cite[Theorem 4.1]{GK15}, there exists a nonzero polynomial $f$ with $\deg f\leq D$ so that $\R^3\setminus Z(f)$ consists of $\lesssim D^3$ open sets $\Omega_i$, each containing $\lesssim D^{-3}N$ points of $\cP$.

We partition $\cF_{\textrm{lg}}=\cF_{\textrm{cell-rich}}\sqcup\cF_{\textrm{alg-rich}}\sqcup\cF_{\textrm{alg}}$ where $F\in\cF_{\textrm{lg}}$ lies in $\cF_{\textrm{cell-rich}}$ if $|\cP_F\setminus Z(f)|\geq|\cP_F|/2$ and it lies in $\cF_{\textrm{alg-rich}}$ if $F\not\subseteq Z(f)$ yet $|\cP_f\cap Z(f)|>|\cP_F|/2$. Finally $F\in\cF_{\textrm{alg}}$ if $F\subseteq Z(f)$.

For each plane $F\in\cF_{\textrm{cell-rich}}$, let $I_F$ be the set of indices $i$ so that $\abs{\cP_F\cap\Omega_i}\geq 2N^{1/2}$. Since the total number of cells intersecting $F$ is at most $C_WD^2$ for some absolute constant $C_W$ by Warren's theorem \cite{War68}, we have
\[\sum_{i\in I_F}\abs{\cP_F\cap\Omega_i} \geq |\cP_F\setminus Z(f)|-C_WD^2(2N^{1/2})\geq |\cP_F|/2-C_WD^2(2N^{1/2})\geq |\cP_F|/4.\]
The last inequality holds since $|\cP_F|\geq\eta|\cP\cap F|\geq \eta k\geq C\eta^{-2}N^{1/2}\geq 8C_W\eta^{-2}N^{1/2}= 4C_WD^2(2N^{1/2})$ as $D=\eta^{-1}$ and we can take $C$ sufficiently large in terms of $C_W$.

Now combining the above with the Cauchy--Schwarz inequality, we have
\[|\cP_F|^2\lesssim \paren{\sum_{i\in I_F}\abs{\cP_F\cap \Omega_i}}^2\lesssim D^2\sum_{i\in I_F}\abs{\cP_F\cap\Omega_i}^2.\]
Note that for each $i\in I_F$, we have that $\cP_F\cap\Omega_i$ is a set of size at least $2N^{1/2}$ and, by definition of $\cP_F$, does not contain any $N^{1/2}$-rich line. Thus by Beck's theorem (\cref{thm:beck}), the number of 2-rich lines for $\cP_F\cap\Omega_i$ satisfies $\abs{\cL_2(\cP_F\cap\Omega_i)}\geq c\abs{\cP_F\cap\Omega_i}^2$ for some absolute constant $c>0$.

Write $\cL_i=\cL_2(\cP\cap\Omega_i)$. We have shown that for each $F\in\cF_{\textrm{cell-rich}}$, we have $|\cP_F|^2\lesssim D^2\sum_i\abs{\set{\ell\in\cL_i:\ell\subset F}}$, so
\[\sum_{F\in\cF_{\textrm{cell-rich}}}|\cP_F|^2\lesssim D^2\sum_i\abs{\set{(\ell,F)\in\cL_i\times\cF:\ell\subset F}}.\]
We bound the right-hand side by taking the projective dual and applying a point-line incidence bound. Formally, embed $\R^3\subset\P_{\R}^3$. Write $\cF^*$ for the projective dual of $\cF$ (a set of points in $\P_\R^3$) and $\cL_i^*$ for the projective dual of $\cL_i$ (a set of lines in $\P_\R^3$). For each pair $(\ell,F)\in\cL_i\times\cF$ with $\ell\subset F$, by duality we have $F^*\in\ell^*$. Furthermore, for any plane $H\subset\P^3_{\R}$, there is a point $H^*\in\P^3_{\R}$ so that $\ell^*\subset H$ if and only if $H^*\in\ell$.

Now any collection of coplanar lines in $\cL_i^*$ corresponds to a collection of lines in $\cL_i$ which all pass through some point $H^*\in\P^3$. However, each line of $\cL_i$ is 2-rich for $\cP\cap\Omega_i$ and thus passes through some point of this set other than $H^*$. These points must be distinct, showing that any set of coplanar lines has size at most $|\cP\cap\Omega_i|$.

Therefore we have shown that
\[\sum_{F\in\cF_{\textrm{cell-rich}}}|\cP_F|^2\lesssim D^2\sum_iI(\cF^*,\cL_i^*)\]
where $\cL_i^*$ is a set of at most $|\cP\cap\Omega_i|^2$ lines with at most $|\cP\cap\Omega_i|$ lines lying in any plane.
Thus by the Guth--Katz incidence bound (\cref{thm:guth-katz-point-line-incidences}), we have
\begin{align*}
I(\cF^*,\cL_i^*)
&\lesssim \abs{\cF^*}^{1/2}\abs{\cL_i^*}^{3/4}+\abs{\cP\cap\Omega_i}^{1/3}\abs{\cF^*}^{2/3}\abs{\cL_i^*}^{1/3}+\abs{\cF^*}+\abs{\cL_i^*}\\
&\lesssim \abs{\cF_{\textrm{cell-rich}}}^{1/2}\abs{\cP\cap\Omega_i}^{3/2}+\abs{\cF_{\textrm{cell-rich}}}^{2/3}\abs{\cP\cap\Omega_i}+\abs{\cF_{\textrm{cell-rich}}}+\abs{\cP\cap\Omega_i}^{2}.
\end{align*}
Summing over the cells, this gives
\begin{align*}
\sum_i I(\cF^*,\cL_i^*)
&\lesssim D^3\paren{\abs{\cF_{\textrm{cell-rich}}}^{1/2}\paren{\frac N{D^3}}^{3/2}+\abs{\cF_{\textrm{cell-rich}}}^{2/3}\paren{\frac N{D^3}}+\abs{\cF_{\textrm{cell-rich}}}+\paren{\frac N{D^3}}^2}\\
&=\abs{\cF_{\textrm{cell-rich}}}^{1/2}N^{3/2} D^{-3/2}+\abs{\cF_{\textrm{cell-rich}}}^{2/3}N+\abs{\cF_{\textrm{cell-rich}}}D^3+N^2D^{-3}.
\end{align*}
We can ignore the second term as it is dominated by the first and third term. Combining everything gives
\[\sum_{F\in\cF_{\textrm{cell-rich}}}|\cP_F|^2\lesssim \abs{\cF_{\textrm{cell-rich}}}^{1/2}N^{3/2} D^{1/2}+\abs{\cF_{\textrm{cell-rich}}}D^5+N^2D^{-1}.\]
Now the left-hand side is at least $\abs{\cF_{\textrm{cell-rich}}}(\eta k)^2$. Combining these two inequalities lets us bound $\abs{\cF_{\textrm{cell-rich}}}$. Since $(\eta k)^2\geq C^2\eta^{-4}N$ while $D^5=\eta^{-5}$, the second term on the right-hand side is negligible (note that the result we are trying to prove is vacuous for $\eta<N^{-1/6}$); rearranging the other two terms gives
\[\abs{\cF_{\textrm{cell-rich}}}\lesssim \frac{N^3D}{\eta^4k^4}+\frac{N^2}{D\eta^2k^2}.\]
Plugging this bound back into the previous inequality thus gives
\[\sum_{F\in\cF_{\textrm{cell-rich}}}|\cP_F|^2\lesssim \frac{N^3D}{\eta^2k^2}+\frac{N^{5/2}}{\eta k}+\frac{N^3D^6}{\eta^4k^4}+\frac{N^2D^4}{\eta^2k^2}+\frac{N^2}{D}.\]
Using $D=\eta^{-1}$ and $k\geq C\eta^{-3}N^{1/2}$, the right-hand side can be bounded by $\lesssim \eta N^2$.

Next we bound the contribution from $\cF_{\textrm{alg-rich}}$. For each $F\in\cF_{\textrm{alg-rich}}$, let $\gamma_F$ be the union of the non-line components of $F\cap Z(f)$. Formally, we argue over the complex numbers. The intersection $F_{\CC}\cap Z_{\CC}(f)$ is a variety of pure dimension 1 and degree at most $D$. There are at most $D$ components of this variety which are lines and, by definition, each contains at most $N^{1/2}$ points of $\cP_F$. Thus $\abs{\cP_F\cap\gamma_F}\geq \abs{\cP_F\cap Z(f)}-DN^{1/2}\geq \abs{\cP_F}/2-DN^{1/2}\geq \abs{\cP_F}/4$, where the last inequality follows since $\abs{\cP_F}\geq\eta k\geq C\eta^{-2}N^{1/2}\geq 4\eta^{-1}N^{1/2}=4DN^{1/2}$.

Consider the number of incidences between $\cP$ and the set of curves $\Gamma:=\{\gamma_F\}_{F\in\cF_{\textrm{alg-rich}}}$. We just showed that
\[I(\cP,\Gamma)\geq \sum_{F\in\cF_{\textrm{alg-rich}}}\abs{\cP_F\cap\gamma_F}\gtrsim \sum_{F\in\cF_{\textrm{alg-rich}}} \abs{\cP_F}\gtrsim \eta k\abs{\cF_{\textrm{alg-rich}}}.\]
On the other hand, the incidence graph between $\cP,\Gamma$ is $K_{D^2+1,2}$-free. To see this, note that any $\gamma_F,\gamma_{F'}$ do not share any components (as each irreducible component is a non-line curve and thus is contained in a unique plane). Then by B\'ezout's theorem, $\gamma_F\cap\gamma_{F'}$ consists of at most $D^2$ points. Thus by the K\H{o}v\'ari--S\'os--Tur\'an theorem, we have
\[\eta k\abs{\cF_{\textrm{alg-rich}}}\lesssim I(\cP,\Gamma)\lesssim DN^{1/2}\abs{\cF_{\textrm{alg-rich}}}+N.\]
Since $\eta k\geq C\eta^{-2}N^{1/2}$ while $DN^{1/2}=\eta^{-1}N^{1/2}$, the first term on the right-hand side is negligible; rearranging thus gives $|\cF_{\textrm{alg-rich}}|\lesssim N/(\eta k)$.

Now using the hypothesis that each plane is at most $\eta N$-rich, we conclude that
\begin{align*}
\sum_{F\in\cF_\textrm{alg-rich}}|\cP_F|^2&\leq \eta N\sum_{F\in\cF_\textrm{alg-rich}}|\cP_F|\lesssim \eta N\cdot I(\cP,\Gamma)\\&\lesssim \eta DN^{3/2}\abs{\cF_{\textrm{alg-rich}}}+\eta N^2 \lesssim \frac{DN^{5/2}}{k}+\eta N^2.
\end{align*}
Using $D=\eta^{-1}$ and $k\geq C\eta^{-3}N^{1/2}$, the right-hand side can be bounded by $\lesssim \eta N^2$.

Finally each $F\in\cF_{\textrm{alg}}$ lies in $Z(f)$, so there are at most $D$ of these planes. Again since each plane is at most $\eta N$-rich, we have
\[\sum_{F\in\cF_{\textrm{alg}}}\abs{\cP_F}^2\leq D(\eta N)^2=\eta N^2.\]
Combining all the cases gives the desired bound.
\end{proof}

Next we study line-line incidences. Suppose that we have (not necessarily distinct) lines $\ell_1,\ldots, \ell_N$ in $\PP_{\RR}^3$ and consider the number of incident pairs $(i,j)$ with $\ell_i\cap \ell_j\neq \emptyset$. There are several ways for this number to be quadratic in $N$. In the most degenerate setting, all the lines could be equal, creating $N^2$ incident pairs. Another option is to have all $N$ lines distinct but lying in a fixed plane, or all of them distinct but passing through a fixed point. Both constructions still have $N^2$ incident pairs. Lastly, we can fix a regulus and take $N/2$ lines from each of its ruling, creating $N^2/4$ incident pairs. The following lemma says that these four constructions are the only possible ways to get quadratically many crossings.

\begin{lemma}[Structure for line sets with many crossings]
\label{lem:line-structure-new}
Let $\ell_1,\ldots,\ell_N$ be (not necessarily distinct) lines in $\PP_{\RR}^3$ and let $\eta \in (0,1/2)$. There exist a set of points $\fP$, a set of planes $\fF$, a set of reguli $\mathfrak R$, and a set of lines $\mathfrak L$ with the following properties:
\begin{enumerate}[(i)]
    \item $|\fP|, |\fF|, |\fR|,|\fL|\lesssim \eta^{-{8}}$;
    \item there are $\lesssim \eta N^2$ pairs $(i,j)\in[N]^2$ such that $\ell_i,\ell_j$ are incident, but none of the following hold:
    \begin{itemize}
        \item $\ell_i\cap\ell_j=\{p\}$ for some $p\in\fP$;
        \item $\ell_i,\ell_j\subset F$ for some $F\in\fF$;
        \item $\ell_i,\ell_j\subset R$ for some $R\in\fR$;
        \item $\ell_i\in\fL$ or $\ell_j\in\fL$.
    \end{itemize}
\end{enumerate}
\end{lemma}

\begin{proof}
Let $\fL$ be the set of lines which appear in the list $\ell_1,\ldots,\ell_N$ at least $\eta^{8}N$ times. Clearly $|\fL|\leq\eta^{-{8}}$.

We iteratively construct the rest of the required data. Initialize $I=\{i\in[N]:\ell_i\not\in\fL\}$ and $\fP=\fF=\fR=\emptyset$. At all times we will maintain a set $E\subseteq [N]^2$ consisting of the pairs $(i,j)\not\in I^2$ so that $\ell_i,\ell_j$ are incident but none of the following hold:
\begin{itemize}
    \item $p\in \ell_i,\ell_j$ for some $p\in\fP$;
    \item $\ell_i,\ell_j\subset F$ for some $F\in\fF$;
    \item $\ell_i,\ell_j\subset R$ for some $R\in\fR$;
    \item $\ell_i\in\fL$ or $\ell_j\in\fL$.
\end{itemize}
Note that this is almost the same as the number of pairs we are aiming to bound, but the first bullet point is slightly different.

We start with $E=\emptyset$, since if $(i,j)\not\in I^2$, then at least one of $\ell_i,\ell_j$ lies in $\fL$.

Consider the graph $G$ with vertex set $I$ where two vertices are adjacent if the corresponding lines are incident. First suppose that $G$ has at least $\eta N^2$ edges.

We claim that this implies that there exist $a,b,c\in I$ so that $\ell_a,\ell_b,\ell_c$ are distinct lines and the common neighborhood of $a,b,c$ in $G$ has size at least $\eta^3N$. To see this, consider $Y$, the number of tuples $(x,a,b,c)\in I^4$ such that $a,b,c\in N_G(x)$ and $\ell_a,\ell_b,\ell_c$ are distinct. If $\deg_G(x)\geq 2\eta^8N$, then $x$ contributes at least $\deg_G(x)(\deg_G(x)-\eta^{8}N)(\deg_G(x)-2\eta^{8}N)$ quadruples, so each $x$ contributes at least $\max\{0,\deg_G(x)-2\eta^{8}N\}^3$. As this function is convex, we conclude
\[Y\geq \sum_{x\in I}\max\{0,\deg_G(x)-2\eta^{8}N\}^3\geq |I|\paren{\frac{2e(G)}{|I|}-2\eta^{8}N}^3\geq \eta^3 N^4,\]
averaging over all possible triples $a,b,c\in I$, we find one which contributes at least $\eta^3 N$ to $Y$, as desired.

First suppose $\ell_a,\ell_b,\ell_c$ are pairwise skew. In this case, the set of lines incident to all three of $\ell_a,\ell_b,\ell_c$ sweeps out a regulus (see, e.g., \cite[Proposition 8.15]{Guth16}). Call this regulus $R$. In this case we add $R$ to $\fR$ and remove all $i$ with $\ell_i\subset R$ from $I$. By our choice of $a,b,c$, we have removed at least $\eta^3N$ elements from $I$.

We now check how many incidences this adds to $E$. By definition, the only incidences we add involve one line in $R$ and one line that does not lie in $R$. For each $\ell_i$ not lying in $R$, there are at most two points of $R$ which it intersects. Every point of $R$ is contained in exactly two lines, so there are at most four lines contained in $R$ which $\ell_i$ can be incident to. Finally, each line appears with multiplicity at most $\eta^{8}N$, so in total we increase $E$ by at most $4\eta^{8}N^2$.

If $\ell_a,\ell_b,\ell_c$ are not pairwise skew then two of them, say $\ell_a,\ell_b$ intersect. Let $\ell_a,\ell_b$ intersect at the point $p$ and span the plane $F$. Let $I_F\subseteq I$ be the set of $i\in I$ so that $\ell_i\subset F$ and let $\fP_F$ be the set of points $p'\in F$ which are $2\eta^4N$-rich, i.e., there exist at least $2\eta^4N$ indices $i\in I_F$ so that $\ell_i$ passes through $p'$. Dually, let $I_p\subseteq I$ be the set of $i\in I$ so that $\ell_i\ni p$ and let $\fF_p$ be the set of planes $F'\ni p$ for which there exist at least $2\eta^4N$ indices $i\in I_p$ so that $\ell_i$ is contained in $F'$. In this case we add $\{p\}\cup\fP_F$ to $\fP$, add $\{F\}\cup \fF_p$ to $\fF$, and remove $I_F\cup I_p$ from $I$.

We first claim that $|I_F\cup I_p|\geq \eta^3N$. By construction there are at least $\eta^3N$ indices $i\in I$ so that $\ell_i$ intersects all of $\ell_a,\ell_b,\ell_c$. Now every line which intersects both of $\ell_a,\ell_b$ either contains the point $p$ or is contained in the plane $F$, showing that all of these indices lie in $I_F\cup I_p$.

Next we claim that $|\fP_F|,|\fF_p|< 2\eta^{-4}$. To see this, suppose for contradiction that we have distinct points $p_1,p_2,\ldots,p_{2\eta^{-4}}\in\fP_F$. Now $p_i$ is incident to $\ell_j$ for at least $2\eta^4N$ values of $j$. At most $(i-1)\eta^8N$ of these lines can be incident to one of $p_1,\ldots,p_{i-1}$, since the line connecting $p_{i'},p_i$ appears with multiplicity at most $\eta^8N$. Thus we find
\[\sum_{i=1}^{2\eta^{-4}}(2\eta^4N-(i-1)\eta^8N)\geq (2\eta^{-4})(2\eta^4N)-\frac12(2\eta^{-4})^2(\eta^8N)=2N\]
distinct indices $j$ so that $\ell_j$ is incident to one of the $p_i$. This is a contradiction, since there are only $N$ possible indices. The bound on $|\fF_p|$ follows by the same argument.

Finally we bound the number of incidences that we add to $E$. All incidences that we add to $E$ must either involve one element of $I_F$ and one element $I\setminus I_F$ or one element of $I_p$ and one element of $I\setminus I_p$. (Indeed if both indices lie in $I_F$, then both lines are contained in $F$ and the incidence is not counted and if both indices lie in $I_p$, then both lines pass through $p$ and the incidence is not counted.)

Consider a line not in $F$; it intersects $F$ at a unique point. If this point lies in $\fP_F$ then it contributes no incidences to $E$; otherwise it contributes at most $2\eta^4N$ incidences. Thus the incidences between $I_F$ and $I\setminus I_F$ contribute at most $2\eta^4N^2$ to $E$. Similarly, the incidences between $I_p$ and $I\setminus I_p$ contribute at most $2\eta^4N^2$ to $E$.

We continue this process until $G$ has fewer than $\eta N^2$ edges. Note that each step of the process decreases the size of $I$ by at least $\eta^3 N$, so there are at most $\eta^{-3}$ steps. Thus we see that at the end $|\fP|,|\fF|\lesssim \eta^{-7}$ and $|\fR|\lesssim\eta^{-3}$ while $|E|\lesssim \eta N^2$.

Finally we bound the number of incidences counted in $(ii)$. We have $|E|\lesssim \eta N^2$ incidences produced by the iteration and $e(G)\leq \eta N^2$ incidences remaining in $I^2$. The only type of incidence we have not yet counted is those $(i,j)$ where $\ell_i=\ell_j$ contains some point $p\in\fP$. There are at most $\eta^{8}N^2$ of these incidences which we have not yet counted. Indeed, there are $N$ choices for $i$ and at most $\eta^{8}N$ choices for $j$ so that $\ell_i=\ell_j$ and $\ell_j\not\in\fL$. (If $\ell_j\in\fL$, then this incidence was already counted.)
\end{proof}

\begin{remark}
We remark that this result also follows from the regularity lemma for semialgebraic graphs \cite[Theorem 1.3]{TY24}. Parameterizing $\ell_1,\ldots,\ell_N$ using the Pl\"ucker embedding, the property of two lines $\ell_i,\ell_j$ being incident can be written as a single polynomial equation. The semialgebraic regularity lemma decomposes the incidence graph into $\eta^{-O(1)}$ complete bipartite graphs and $\eta N^2$ additional edges. One can check that each of the complete bipartite graphs encodes the structure of either a rich point, plane, regulus, or line. Optimizing this argument gives an alternative proof of \cref{lem:line-structure-new} with the $\eta^{-8}$ replaced by $\eta^{-4}$.
\end{remark}

As in \cref{sec:ruled-ruled}, it will be sometimes helpful to cluster our collection of points onto a small number of low-degree curves. This will be useful when we want to show that most of the incidences in a point-surface incidence problem come from curve-surface incidences. In some sense the following result is a more general version of \cref{prop:curves-containing-points}.

In this result we are working over $\RR^3$ and it will be important to distinguish curves $\gamma\subset\CC^3$ where $\gamma\cap\R^3$ is 1-dimensional from those where $\gamma\cap\R^3$ is 0-dimensional.

\begin{definition}
An irreducible algebraic curve $\gamma\subset\CC^n$ is a \emph{real curve} if $\gamma\cap\R^n$ is infinite. A reducible algebraic curve $\gamma\subset\CC^n$ (i.e., a variety of pure dimension 1) is a real curve if each of its irreducible components is a real curve. We sometimes abuse notation and refer interchangeably to a curve $\gamma\subset\CC^n$ and its set of real points, $\gamma\cap\R^n$. As a real curve $\gamma$ is the Zariski-closure of $\gamma\cap\R^n$, either set can recovered from the other.
\end{definition}

\begin{lemma}\label{lem:rich-variety}
Let $\cP \subset \RR^3$ be a set of size $N$. Let $\mc V$ be a family of irreducible surfaces of degree at most $C_0$ and let $K\geq 2$ be a parameter. Then there exist a constant $C_1$, a collection of irreducible real curves $\Gamma$, and a set $\cP_\gamma\subseteq \cP$ for each $\gamma\in\Gamma$ such that:
    \begin{enumerate}[(i)]
        \item $|\Gamma| \le C_1 N^{2/3}$;
        \item $\cP_{\gamma}\subset \gamma$ for each $\gamma\in \Gamma$ and the sets $\cP_\gamma$ are disjoint;
        \item for all $r \ge C_1 N^{2/3}$, the number of $V \in \mc V$ such that $\abs{(V\cap \cP)\setminus \left(\bigcup_{\gamma\in \Gamma: \gamma\subset V}\cP_{\gamma}\right)}\ge r$ is bounded by $C_1 N/r$;
        \item there are at most $C_1N^{1/3}$ surfaces $V\in \mc V$ that contain more than $C_1N^{1/3}$ curves from $\Gamma$; and
        \item every $\gamma \in \Gamma$ lies in at least $K$ surfaces $V\in\cV$ and has degree at most $C_1$.
    \end{enumerate}
\end{lemma}

\begin{proof}
We will take $C_1$ sufficiently large in terms of $C_0,K$. Let $f$ be a minimal degree nonzero polynomial vanishing on $\cP$. Decompose $Z(f)$ into irreducible components as
    \[
    Z(f) = V_1 \cup \cdots \cup V_t
    \]
and set $d_i = \deg V_i$. Let $\mc V_0$ be the set of components $V_i$ which belong to $\mc V$. By parameter counting (\cref{cor:Walsh-param-counting}), we have $\deg f = \sum d_i \le CN^{1/3}$; in particular $t \le CN^{1/3}$ for some absolute constant $C$.
Later, we will construct $\Gamma$ so that $\gamma\subset Z(f)$ for every $\gamma\in \Gamma$. Note that this implies that for every $V\in \mc V\setminus\cV_0$, we have that $V$ contains at most $\deg f\deg V\leq C_0\deg f\leq C_0CN^{1/3}$ different curves from $\Gamma$. Since $\abs{\cV_0}\leq t\leq CN^{1/3}$, this proves (iv).

Let $\cP=\cP_1\sqcup\cdots\sqcup\cP_t$ be a partition of $\cP$ where $\cP_i\subset V_i$.  Let $g_i$ be a minimal degree polynomial vanishing on $\cP_i$ and not vanishing identically on $V_i$. Again by parameter counting (\cref{cor:Walsh-param-counting}), we have 
    \[
    \deg g_i \le C \max\{(|\cP_i| / d_i)^{1/2}, |\cP_i|^{1/3}\}
    \]
for some absolute constant $C$.

Let $\Gamma_i$ be the set of irreducible components of $V_i \cap Z(g_i)$ and define $(\cP_{i,\gamma})_{\gamma\in \Gamma_i}$ to be some partition of $\cP_i$ so that $\cP_{i,\gamma}\subset \gamma$ for each $\gamma\in \Gamma_i$. Defining $M_i=\deg(\Gamma_i)$, by B\'ezout's theorem we have the bound
\begin{align*}
M &:= \sum_i \deg(\Gamma_i) \le \sum_i d_i\deg g_i \leq C \sum_i \paren{(|\cP_i|d_i)^{1/2} + |\cP_i|^{1/3} d_i} \\
&\le C \paren{\sum_i |\cP_i|}^{1/2}\paren{\sum_i d_i}^{1/2} + C N^{1/3}\sum_i d_i \leq 2C^2 N^{2/3}.
\end{align*}

Let us write $\Gamma'_i$ for the set of real curves $\gamma \in \Gamma_i$ which are contained in at least $K$ surfaces $V\in \mc V$. Let $\Gamma = \bigcup \Gamma_i'$ and for every $\gamma\in \Gamma$, let $\cP_{\gamma} = \bigcup_{i:\gamma\in \Gamma_i}\cP_{i,\gamma}$. We have
\[|\Gamma|\leq\sum_i|\Gamma_i|\leq M\leq 2C^2N^{2/3},\]
proving (i). By construction, $\cP_{i,\gamma}\subset\gamma$ and the $\cP_{i,\gamma}$ partition $\cP$, showing that (ii) holds. Since each $\gamma\in\Gamma$ lies in $K\geq 2$ surfaces $V\in\cV$, B\'ezout's theorem implies that each $\gamma \in \Gamma$ has degree at most $C_0^2$, proving (v).

For any finite $\tilde\cV\subset\cV$, we aim to bound the size of the set
\[\set{(p,V)\in\cP\times\tilde \cV:p\in V, p\not\in\bigcup_{\gamma\in\Gamma:\gamma\subset V}\cP_\gamma}.\]
Note that for each $p\in \cP$, there exists a unique $i$ and $\gamma\in\Gamma_i$ such that $p\in \cP_{i,\gamma}$. Then for each pair $(p,V)$ in the above set, one of the following occurs:
    \begin{enumerate}
        \item $\gamma\not\subset V$;
        \item $\gamma\subset V$ and $\gamma$ is contained in at most $K-1$ surfaces in $\Tilde{\cV}$;
        \item $\gamma\subset V$ and $\gamma$ is contained in at least $K$ surfaces in $\Tilde{\cV}$ but $\gamma\cap\RR^3$ is finite;
        \item $\gamma\subset V$ and $\gamma\in\Gamma$.
    \end{enumerate}
The number of pairs $(p,V)$ of type (1) is bounded by B\'ezout's theorem as $\sum_{V\in\tilde\cV}\deg(V)\deg(\Gamma)\leq C_0M|\tilde\cV|\leq 2C_0C^2N^{2/3}|\tilde\cV|$. The number of pairs $(p,V)$ of type (2) is bounded by $\sum_{i,\gamma\in\Gamma_i}(K-1)|\cP_{i,\gamma}|<KN$.

For the number of pairs $(p,V)$ of type (3), note that $\deg \gamma\leq C_0^2$ and $\gamma\cap\R^3$ is finite, i.e., $\gamma$ is not a real curve. In this case, we claim that $|\gamma\cap\R^3|\leq C_0^4$. 

This is a standard fact in real algebraic geometry. To see this, first suppose that $\gamma\subset\CC^3$ is not defined over the real numbers. Then its complex conjugate $\bar\gamma$ is an irreducible curve, distinct from $\gamma$. However $\gamma\cap\R^3=\gamma\cap\bar\gamma\cap\R^3\subseteq\gamma\cap\bar\gamma$. By B\'ezout's theorem, the last set has size at most $\deg(\gamma)^2\leq C_0^4$. Next suppose that $\gamma$ is defined over the real numbers. By the implicit function theorem, $\gamma\cap\R^3$ is locally a 1-dimensional manifold in the neighborhood of any nonsingular point. Since we assume that $\gamma\cap\R^3$ is finite, every such point must be singular. Again by B\'ezout's theorem, we have $|\gamma\cap\R^3|\leq\deg(\gamma)(\deg(\gamma)-1)<C_0^4$. 

Therefore the number of pairs of type (3) is bounded by $C_0^4\sum_{i}|\Gamma_i||\tilde\cV|\leq C_0^4M|\tilde\cV|\leq 2C_0^4C^2N^{2/3}|\tilde\cV|$. 

Now we claim that there are no pairs of type (4) in the above set. Indeed, $p\in \cP_{i,\gamma}\subseteq \cP_\gamma$, so $\gamma\in\Gamma$ is not possible. Thus we conclude that for each finite $\tilde\cV\subset\cV$, 
\begin{align*}
\abs{\set{(p,V)\in\cP\times\tilde\cV:p\in V, p\not\in\bigcup_{\gamma\in\Gamma:\gamma\subset V}\cP_\gamma}}
&\leq 4C_0^4C^2N^{2/3}|\tilde\cV|+KN.
\end{align*}

Setting $\tilde\cV_r$ to be the set of $V\in\cV$ such that $|(V \setminus \bigcup_{\gamma \in \Gamma: \gamma\subset V} \gamma) \cap \cP| \ge r$, we see that for any $r \ge 8C_0^4C^2N^{2/3}$, we have $|\tilde\cV_r|\leq 2KN/r$. This proves (iii) for $C_1$ sufficiently large in terms of $C,C_0,K$.
\end{proof}

\subsection{Generalizations of Guth--Katz}
\label{sec:gen-GK}

In this section, we give some slight generalizations of the 2-dimensional Guth--Katz bound.

\begin{definition}
For point sets $\cP_1,\cP_2,\cP_3,\cP_4$, denote the set of \emph{distances quadruples}
\[\cE(\cP_1,\cP_2,\cP_3,\cP_4)={\{(p,q,p',q')\in \cP_1\times \cP_2\times \cP_3\times \cP_4:\norm{p-p'}=\norm{q-q'}\}}.\]
The size of this set is the (asymmetric) \emph{distance energy}, denoted
\[E(\cP_1,\cP_2,\cP_3,\cP_4)=\abs{\cE(\cP_1,\cP_2,\cP_3,\cP_4)}.\]
\end{definition}

Note that $E(\cP_1,\cP_2,\cP_3,\cP_4)=E(\cP_1,\cP_4,\cP_3,\cP_2)$.
In the symmetric case, we write
\[\cE(\cP)=\cE(\cP,\cP,\cP,\cP)\qquad\text{ and }\qquad E(\cP)=E(\cP,\cP,\cP,\cP).\]

The original Guth--Katz bound is as follows.

\begin{theorem}[{\cite[Theorem 1.1 and Proposition 2.2]{GK15}}]
\label{thm:guth-katz-distance-energy}
Let $\cP\subset\R^2$ be a set of $N$ points. Then the distance energy satisfies
\[E(\cP)\lesssim N^3\log N.\]
As a consequence, $\cP$ determines $\gtrsim N/\log N$ distinct distances.
\end{theorem}

We will extend their argument to the asymmetric setting on both planes and spheres. It was already observed \cite{Tao11blogpost} that the Guth--Katz argument can also be used to give a lower bound on distinct distances on spheres; here we spell out the details. In the plane, the bipartite distance energy $\cE(\cP_1,\cP_1,\cP_2,\cP_2)=\cE(\cP_1,\cP_2,\cP_2,\cP_1)$ was previously studied by Mathialagan \cite{Mat21}; our argument is essentially equivalent for the fully asymmetric case.

Let $\SS^2$ be the unit sphere in $\RR^3$. Let $0$ be the origin of $\RR^3$, and let $p,q$ be two points in $\SS^2$.
To carry out the Elekes--Sharir reduction in this setting, we parametrize (orientation-preserving) rigid motions of the sphere that send $p$ to $q$. These are the orientation-preserving rigid motions of $\R^3$ that fix $0$ and send $p$ to $q$.
In light of this, we set
\begin{equation}
\label{eq:elekes-sharir-sphere}
\ell_{pq} = F_{pq}^+(\RR)\cap F_{00}^+(\RR)\subset F_{00}^+(\RR)\cong \RR^3.
\end{equation}
By \cref{prop:missing-rigid-motion}, note that $F_{00}^+(\R)\subset E^+_o(\R)$ is the set of all orientation-preserving rigid motions fixing $0$ (and thus acting on $\SS^2$) except for the ones that rotate by angle $\pi$.

To show that those lines determine few rich points, we need to show that they do not concentrate on planes and reguli.
This can be done with a similar argument to \cite[Proposition 2.7]{GK15}.
Given the ruled surface theory that we have developed so far, we will provide the details in a slightly different manner.

\begin{lemma}
\label{lem:line-concentration-bound}
Let $\cP\subset \SS^2\subset \RR^3$ be a set of $N$ points on the sphere. Let $\cL = \{\ell_{pq}\}_{p,q\in \cP}$ be the corresponding set of $N^2$ lines in $F_{00}^+(\RR)\cong \RR^3$. Then any plane or regulus contains $\lesssim N$ lines from $\cL$.
\end{lemma}

\begin{proof}
By \cref{prop:all-2-flats}, any plane in $F_{00}^+$ is of the form $F_{00}^+\cap F_{\tau}^+$ for some $\tau\in E^-(\RR)$.
Since $F_{\tau}^+$ has non-empty intersection with $F_{00}^+\subset E^+_o$, we see that $\tau$ is a genuine orientation-reversing rigid motion by \cref{prop:actual-rigid-motion}. Moreover, by \cref{prop:intersection-from-diff-family}, we know that $\phi(0,0)\in F_{\tau}^{\Phy}$, and so $\tau$ fixes $0$. Now if $\ell_{pq}\subset F_{\tau}^+$, then $F_{pq}^+\cap F_{\tau}^+$ is necessarily a $2$-flat by \cref{prop:intersection-from-diff-family}. By the same argument, we have $\phi(p,q)\in F_{\tau}^{\Phy}$, implying that $\tau(p)=q$. Therefore the number of lines from $\cL$ contained in $F_{00}^+\cap F_{\tau}^+$ can be upper bounded by the number of pairs $(p,q)\in \cP^2$ with $\tau(p)=q$. Since for each $p$ there is at most one such $q$, we see that there are at most $N$ such pairs, as desired.

Next fix a regulus $\fr\subset F_{00}^+(\RR)$. View $\fr\subset F_{00}^+(\RR)\subset F_{00}^+(\CC)\cong\PP^3_{\CC}$ and let $\tilde\fr$ be the Zariski closure of $\fr$ in $\PP^3_{\CC}$. Recall from \cref{ssec:plucker} that we can identify the set of lines $\cL^+(00)$ contained in $F_{00}^+$ with the set of lines $\cL^{\Phy}(00)$ through the origin in $E^{\Phy}$ which is naturally identified with a subset of $\PP^5$. By \cref{thm:plucker-embeding}, there is a linear isomorphism $i\colon F_{00}^+\to\PP^3_{\CC}$ so that the Grassmannian of lines in $F_{00}^+$ agrees with the parametrization defined above.

Now by \cref{thm:classical-ruled-surface}(3), the Fano variety $\Lambda_{\tilde\fr}\subset\G(1,3)\subset\PP^5$ of $\tilde\fr$ consists of two irreducible curves of degree 2. Write $Y\subset\CC^6$ for the cone over $\Lambda_{\tilde\fr}$ with apex the origin, i.e., the set of $x\in\CC^6$ so that $[x]\in\Lambda_{\tilde\fr}$. This is a degree 4 variety of pure dimension 2. Now a line $\ell_{pq}$ lies in $\fr$ only if $\phi(p,q)\in Y$, where $\phi\colon\CC^6\to\CC^6$ is the linear isomorphism defined by $\phi(p,q)=((p+q)/2,(p-q)/2)$. We are investigating the set of pairs $(p,q)\in\cP\times\cP\subset\SS^2\times\SS^2$ which also lie on $\phi^{-1}(Y)$. Such a pair lies in $X\cap(\cP\times\cP)$ where \[X=\phi^{-1}(Y\cap Z((x_1+x_4)^2+(x_2+x_5)^2+(x_3+x_6)^2-1)),\] as the hypersurface encodes the condition that $p\in\SS^2$, so $\norm{p}^2=1$. We claim that $X$ is a variety of pure dimension 1. Indeed, this occurs unless an irreducible component of $Y$ lies in $Z((x_1+x_4)^2+(x_2+x_5)^2+(x_3+x_6)^2-1)$. Since $Y$ is a cone over the origin, each of its irreducible components contains $0$, which does not lie in the latter hypersurface.

Thus $X$ is a variety of pure dimension 1 and, by B\'ezout's theorem, has degree at most 8. 

Finally, we apply \cref{thm:grid-surface-intersection} to bound $|X\cap(\cP\times\cP)|$. We apply this result with $A=N$, so that the hypothesis that $|\cP\cap\gamma|\leq N\deg(\gamma)$ is trivially true. We conclude that the number of lines $\ell_{pq}$ contained in $\fr$ is $\lesssim N$, as desired.
\end{proof}

With this, we may now prove the following.

\begin{theorem}
\label{thm:bipartite-guth-katz}
Let $\cP_1,\cP_2,\cP_3,\cP_4$ be sets of points either all lying in $\R^2$ or all on the unit sphere $\S^2$. Setting $N=\max\{|\cP_1|,|\cP_2|,|\cP_3|,|\cP_4|\}$,
\[E(\cP_1,\cP_2,\cP_3,\cP_4)\lesssim (|\cP_1||\cP_2|+|\cP_3||\cP_4|)N\log N.\]
By the symmetry of $E$, we also have the bound
\[E(\cP_1,\cP_2,\cP_3,\cP_4)\lesssim (|\cP_1||\cP_4|+|\cP_2||\cP_3|)N\log N.\]
\end{theorem}

\begin{proof}
First we perform the Elekes--Sharir reduction. For $p,q\in\S^2$, let $\ell_{pq}\subset\F_{00}^+(\R)\cong\RR^3$ be the line defined in \cref{eq:elekes-sharir-sphere}. For $p,q\in\R^2$, let $\ell_{pq}\subset\R^3$ be the line defined in \cite[Proposition 2.7]{GK15}. In the planar setting, the points of $\R^3$ correspond to the orientation-preserving rigid motions of $\R^2$ other than the translations; $\ell_{pq}$ is the set of those rigid motions which send $p$ to $q$. Note that $\ell_{pq}=\ell_{p'q'}$ means that every rigid motion which sends $p$ to $q$ also sends $p'$ to $q'$. In the plane, this means that $(p,q)=(p',q')$, on the sphere, this means that $(p,q)=\pm(p',q')$. Thus each line corresponds to at most two pairs of points.

Let $\cL'=\{\ell_{pq}:(p,q)\in \cP_1\times \cP_2\}$ and $\cL''=\{\ell_{p'q'}:(p',q')\in \cP_3\times \cP_4\}$. Writing $\cL=\cL'\cup\cL''$, we have $|\cL|\leq |\cP_1||\cP_2|+|\cP_3||\cP_4|$. By \cite[Proposition 2.8]{GK15} or \cref{lem:line-concentration-bound} applied to $\cP_1\cup \cP_2$ and $\cP_3\cup \cP_4$, we see that each plane and regulus contains $\lesssim N$ lines of $\cL$.

Now for each distance quadruple $(p,q,p',q')\in \cP_1\times \cP_2\times \cP_3\times \cP_4$ with $\norm{p-p'}=\norm{q-q'}$, there is an orientation-preserving rigid motion sending $p$ to $q$ and $p'$ to $q'$. This rigid motion corresponds to the intersection point of the lines $\ell_{pq}\in\cL'$ and $\ell_{p'q'}\in\cL''$ unless either $\ell_{pq}=\ell_{p'q'}$ or if the rigid motion is missed by the reduction (in $\R^2$ this means that it is a translation, in $\S^2$ this means that it is a rotation with angle $\pi$). There are at most $2N^2$ distances quadruples of the first type. We claim that there are at most $|\cP_1||\cP_2||\cP_3|+|\cP_1||\cP_3||\cP_4|$ of the second type. Indeed, in $\R^2$, there is a unique translation sending $p$ to $q$. Thus the distance quadruple is uniquely determined by $p,q,p'$. In $\S^2$ there is a unique rotation with angle $\pi$ sending $p$ to $q$ unless $p=-q$. Thus if $p\neq -q$, then the distance quadruple is uniquely determined by $p,q,p'$. If $p=-q$, then the distance quadruple is uniquely determined by $p,p',q'$.

All the other distance quadruples correspond to an intersection point of a line from $\cL'$ and a line from $\cL''$. If some point $\rho$ is passed through by $k'$ lines from $\cL'$ and $k''$ lines from $\cL''$, it corresponds to at most $4k'k''$ distance quadruples. (We gain two factors of 2 since each line can correspond to two pairs of points.) 

Writing $\cP_{=k}(\cL)$ for the set of exactly $k$-rich points of $\cL$, we have the bound
\begin{align*}
E(\cP_1,\cP_2,\cP_3,\cP_4)
&\leq 2N^2+|\cP_1||\cP_2||\cP_3|+|\cP_1||\cP_3||\cP_4|+4\sum_{k=2}^{2N} k^2 \abs{P_{=k}(\cL)}\\
&\lesssim (|\cP_1||\cP_2|+|\cP_3||\cP_4|)N+\sum_{k=2}^{2N}k \abs{\cP_k(\cL)}\\
&\lesssim (|\cP_1||\cP_2|+|\cP_3||\cP_4|)N+\sum_{k=2}^{2N}k \paren{\frac{|\cL|^{3/2}}{k^2}+\frac{|\cL|N}{k^3}+\frac{|\cL|}k}\\
&\lesssim (|\cP_1||\cP_2|+|\cP_3||\cP_4|)N+|\cL|^{3/2}\log N+|\cL|N+|\cL|N.
\end{align*}
In the third line we used the Guth--Katz bound on point-line incidences in $\R^3$, \cite[Theorem 8.3 and Corollary 12.2]{Guth16}, to bound $\abs{\cP_k(\cL)}$. Plugging in the values of $\abs{\cL}, N$ gives the desired bound.
\end{proof}

Applying the standard Cauchy--Schwarz argument (see \cite[Lemma 2.1]{GK15} or \cref{thm:beck-trick}), this immediately implies a distinct distances bound on the sphere.

\begin{corollary}
\label{thm:distinct-distances-sphere}
Let $\cP\subset\S^2$ be a set of $N$ points on the sphere. Then $\cP$ determines $\gtrsim N/\log N$ distinct distances.
\end{corollary}
\section{Structured cases}\label{sec:degenerate}

In this section, we will bound the structured contribution to the distance energy, that is, the contribution coming from a rigid motion $\rho$ which maps a rich line to a rich line, a rich regulus to a rich regulus, or a fixed point to another fixed point.

Recall that we write $E^+_g\subset E^+(\R)$ for the collection of genuine orientation-preserving rigid motions. Though we state all results in this section for orientation-preserving rigid motions, a symmetric argument gives the same results for orientation-reversing rigid motions.

\subsection{Very rich lines}\label{ssec:rich-line}

\begin{proposition}\label{lem:rich-lines}
Let $\cP\subset\R^3$ be a set of size $N$ such that $|\cP\cap\gamma|\leq N^{2/3}\deg \gamma$ for every curve $\gamma\subset\CC^3$, and every plane contains at most $N^{2/3}$ points of $\cP$. Define $\cE$ to be the collection of quadruples $(p,q,p',q')\in \cP^4$ for which there exists a genuine rigid motion $\rho\in E^+_g$ and $N^{1/3}$-rich lines $\ell_p,\ell_q\in \cL_{N^{1/3}}(\cP)$ such that:
    \begin{itemize}
        \item $p\in \ell_p$ and $q\in \ell_q$;
        \item $\rho(p)=q$ and $\rho(p')=q'$ and $\rho(\ell_p)=\ell_q$.
    \end{itemize}
    Then we have the bound
    \[
    \abs{\cE} \lesssim N^{10/3}\log^2 N.
    \]
\end{proposition}

\begin{definition}
For a line $\ell \subset \R^3$ and $r> 0$, define $C_r(\ell) = \{p: \dist(p, \ell) = r\}$ to be the \emph{cylinder} of radius $r$ with axis $\ell$.
\end{definition}

Essentially all of the distance energy $\abs{\cE}$ comes from quadruples $(p,q,p',q')$ with $p,q$ on rich lines $\ell_p,\ell_q$ and $p',q'$ on cylinders $C_r(\ell_p),C_r(\ell_q)$. To bound $\cE$, we study the incidences between $\cP$ and this collection of cylinders. We do this using \cref{lem:rich-variety} which reduces the problem to a curve-cylinder incidence problem. We will prove a geometric result, \cref{clm:only-lines-in-cylinders}, which implies that lines are the only relevant curve for this problem.

To prove this, we will need the following auxiliary lemma.

\begin{lemma}\label{thm:complete-intersection-degrees}
Let $f_1,f_2$ be relatively prime degree $2$ polynomials in $\CC[x_1,x_2,x_3]$ so that $Z(f_1,f_2)$ is a curve of degree $4$.
Then a degree $2$ polynomial vanishes identically on $Z(f_1,f_2)$ if and only if it is a $\CC$-linear combination of $f_1$ and $f_2$.
\end{lemma}

\begin{proof}
Let $H$ be a generic 2-flat in $\CC^3$, so that $Z(f_1|_H, f_2|_H) = H\cap Z(f_1,f_2)$ is a set of 4 points.
We will first show that if $g$ is a degree $2$ polynomial defined on $H$ that vanishes on $Z(f_1|_H,f_2|_H)$, then $g$ is a linear combination of $f_1|_H$ and $f_2|_H$.
    
By B\'ezout's theorem, no three of these four points are collinear. Indeed, if they were $f_1|_H,f_2|_H$ would both vanish identically on the line that they span. 

Next let $W$ be the space of polynomials on $H$ of degree at most $2$, and let $\varphi\colon W\to \CC^4$ be the evaluation map at the four points in $Z(f_1|_H, f_2|_H)$. Now pick any $p\in H\setminus Z(f_1|_H,f_2|_H)$. We claim that 
    \[\dim\left(\ker \varphi \cap \{g\in W: g(p)=0\}\right)\leq 1.\]
To see this, suppose that $g_1,g_2$ both lie in this set and are not multiples of each other. If $g_1,g_2$ are relatively prime, then by B\'ezout's theorem $Z(g_1,g_2)$ consists of at most 4 points, yet contains the five points of $\{p\}\cup Z(f_1|_H,f_2|_H)$. Otherwise, $g_1,g_2$ share a linear factor $\ell$. Now $Z(\ell)$ contains at most two points of $Z(f_1|_H,f_2|_H)$, so $g_1/\ell$ and $g_2/\ell$ are polynomials of degree at most 1 which both vanish on the two other points of $Z(f_1|_H,f_2|_H)$. This implies that $g_1/\ell$ and $g_2/\ell$ are multiples of each other, a contradiction.

Therefore we conclude that
\[\dim \ker\varphi\leq 1+\dim\left(\ker \varphi \cap \{g\in W: g(p)=0\}\right)\leq 2.\]
Since $f_1|_H,f_2|_H$ both lie in $\ker\varphi$ and are not multiples of each other, we conclude that they generate $\ker\varphi$ as a $\CC$-vector space.

Now let $f\in\CC[x_1,x_2,x_3]$ be a degree $2$ polynomial vanishing on $Z(f_1,f_2)$.
Then $f|_H$ is a degree $2$ polynomial on $H$ vanishing on $Z(f_1|_H,f_2|_H)$. In other words, $f|_H\in\ker\varphi$, implying that 
\[f|_H = a_1f_1|_H+a_2f_2|_H.\]
for some $a_1,a_2\in\CC$.

Writing $h$ for the degree $1$ polynomial defining $H$, there exists $q\in\CC[x_1,x_2,x_3]$ so that
\[f = a_1f_1+a_2f_2+hq.\]
It is clear that $\deg q\leq 1$.
Moreover, as $f,f_1,f_2$ all vanish identically on $Z(f_1,f_2)$, we have that $q$ vanishes identically on $Z(f_1,f_2)$ as well.
If $q$ is nonzero, then $q$ is irreducible and thus is relatively prime to $f_i$ for some $i\in[2]$.
By B\'ezout's theorem, $Z(f_i,q)$ is a curve of degree at most $2$, which is a contradiction as it contains $Z(f_1,f_2)$.

Thus we conclude that $q=0$, so $f=a_1f_1+a_2f_2$, as desired.
\end{proof}

\begin{remark}
This same proof can be generalized to the case when $\deg f_1=\deg f_2=d$ and $Z(f_1,f_2)$ consists of $d^2$ points. In fact the result is true for $f_1,f_2$ of arbitrary degree: if $Z(f_1,f_2)$ is a curve of degree $\deg f_1\deg f_2$, then any polynomial $f$ vanishing identically on $Z(f_1,f_2)$ can be written as $f=a_1f_1+a_2f_2$ with $\deg a_i\leq \deg f-\deg f_i$. This is a corollary of Max Noether's $AF+BG$ theorem \cite{Noether73}.
\end{remark}

\begin{proposition}
\label{clm:only-lines-in-cylinders}
Let $\gamma_\CC\subset\CC^3$ be an irreducible real curve such that $\gamma_\CC\cap\R^3$ lies in three distinct cylinders. Then $\gamma_\CC$ is a line.
\end{proposition}

\begin{proof}
Let $C_1,C_2\subset\R^3$ be distinct cylinders and let $\gamma_\CC\subset\CC^3$ be an irreducible real curve so that $\gamma:=\gamma_\CC\cap\R^3$ lies in $C_1\cap C_2$ and is not a line. We will show that $\gamma$ does not lie in any other cylinder.

Let $\tilde C_1, \tilde C_2$ be the complex projective surfaces defined to be the Zariski-closures of $C_1,C_2\subset\RR^3\subset \PP^3_{\CC}$, respectively. We compute the structure of $\tilde C_1,\tilde C_2$: up to a rigid motion and a dilation, any cylinder can be written as $Z(x^2+y^2-1)$. The complex projective closure of this surface is $[x^2+y^2:x:y:z]$. This is a cone with apex $[0:0:0:1]$ at infinity. Furthermore, the part of this closure lying at infinity is the pair of complex lines $[0:x:\pm ix:z]$; each line has exactly one real point, namely the apex of the cone.

Now we claim that $\tilde C_1\cap \tilde C_2$ has no component which is a line. Indeed, since $\tilde C_1, \tilde C_2$ are cones with no planar component, any line in $\tilde C_1$ passes through its apex $p_1$ and any line in $\tilde C_2$ passes through its apex $p_2$. Thus if $\tilde C_1\cap \tilde C_2$ had a component which was a line, then either $p_1=p_2$ or this would be the line connecting $p_1,p_2$. The latter is not possible since this line lies entirely in $\PP^3_{\CC}\setminus\CC^3$ and has two distinct real points $p_1,p_2$; this contradicts our earlier computation. On the other hand, if $p_1=p_2$, then $\tilde C_1,\tilde C_2$ are cones with the same apex, meaning that their intersection is a union of lines through this point. This contradicts the assumption that $\gamma$ is not a line.

Therefore we see that no component of $\tilde C_1\cap \tilde C_2$ is a line. Now by B\'ezout's theorem with multiplicity \cite[Chapter 1, Theorem 7.7]{Har77}, we have
\[\sum_{\gamma_i\subset \tilde C_1\cap \tilde C_2} m_i\deg\gamma_i=\deg \tilde C_1\deg \tilde C_2=4,\]
where the sum is over the irreducible components of $\tilde C_1\cap \tilde C_2$ and $m_i\geq 1$ is the multiplicity of the component in the intersection. Since we showed that $\deg\gamma_i>1$ for all $i$, we conclude that each irreducible component of $\tilde C_1\cap \tilde C_2$ has degree 2 or 4.

Suppose that $\deg\gamma_\CC=4$. Let $f_1,f_2\in\CC[x_1,x_2,x_3]$ be the polynomials defining $C_1,C_2$. By \cref{thm:complete-intersection-degrees}, we conclude that any degree 2 polynomial which vanishes identically on $\gamma_\CC$ must be a $\CC$-linear combination of $f_1,f_2$.

Assume for contradiction that there is a third cylinder $C_3$ containing $\gamma$. Let $\ell_1,\ell_2,\ell_3$ be the axes and $r_1,r_2,r_3$ be the radii of $C_1,C_2,C_3$ respectively. Now $C_i=Z(f_i)$ where $f_i(x)=\dist(x,\ell_i)^2-r_i^2$. We assumed that $f_1,f_2,f_3$ vanish identically on $\gamma=\gamma_\CC\cap\R^3$. However, as this is an infinite subset of $\gamma_\CC$, we conclude that $f_1,f_2,f_3$ vanish identically on $\gamma_\CC$. Thus in this case, \cref{thm:complete-intersection-degrees} implies that $f_3$ is a $\CC$-linear combination of  $f_1,f_2$.

We claim that this implies that there are nonzero real coefficients $a_1,a_2,a_3$ so that $a_1f_1+a_2f_2+a_3f_3$ is the zero polynomial. As $f_1,f_2,f_3$ are defined over the reals, taking either the real or imaginary parts of the $a_i$ produces a tuple $(a_1,a_2,a_3)\neq(0,0,0)$ of real coefficients so that this identity holds. Next if exactly one of the coefficients were zero, then two of the $C_i$ would be equal and if two of the coefficients were zero, then one of the $f_i$ would be the zero polynomial, neither of which can occur. 

Thus there are nonzero reals $a_1,a_2,a_3$ and a constant $a_0$ so that $a_0+\sum_{i=1}^3 a_i\dist(x,\ell_i)^2$ is identically zero. Note first that $\ell_1,\ell_2,\ell_3$ point in pairwise distinct directions. Indeed, if two were parallel, then two of the cylinders would have an intersection which only consists of lines, contradicting the assumption that $\gamma$ is not a line. Without loss of generality, assume that $a_1,a_2>0$. For any $R$, we can take $x\in\ell_3$ so that $\dist(x,\ell_1),\dist(x,\ell_2)>R$. Then we have
\[0=a_0+\sum_{i=1}^3 a_i\dist(x,\ell_i)^2=a_0+a_1\dist(x,\ell_1)^2+a_2\dist(x,\ell_2)^2>R^2(a_1+a_2)+a_0.\]
Taking $R$ sufficiently large, the right-hand side is positive, contradiction.

Finally, suppose that $\deg\gamma_\CC=2$. First note that any degree 2 curve is planar. Indeed any three points on $\gamma_\CC$ span a plane $F$; by B\'ezout's theorem, if $\gamma_\CC\not\subset F$, then $|\gamma_\CC\cap F|\leq \deg\gamma_\CC\deg F=2$. Since $|\gamma_\CC\cap F|\geq 3$, we conclude that $\gamma_\CC\subset F$. Now let $C$ be a cylinder with axis $\ell$ and radius $r$ and $F$ be a plane. If $\ell$ is parallel to $F$, then $C\cap F$ only has linear components. Otherwise, suppose that $\ell$ intersects $F$ at point $p$ and makes an angle of $\theta$ with $F$. If $\theta\neq \pi/2$, elementary geometry implies that $C\cap F$ is an ellipse with center $p$, semimajor axis $\proj_F \ell$ of length $r/\sin\theta$ and semiminor axis of length $r$. If $\theta=\pi/2$, then $C\cap F$ is a circle with center $p$ and radius $r$.

Thus if $\deg\gamma_\CC=2$, then $\gamma$ is an ellipse. From $\gamma$, we can uniquely determine $F$ (the plane it lies in), $p$ (its center), $r$, and $\theta$. Now if $\theta=\pi/2$, i.e., $\gamma$ is a circle, then there is a unique cylinder $C$ so that $\gamma=F\cap C$: this is the cylinder of radius $r$ whose axis $\ell$ is the line through $p$ perpendicular to $F$. If $\theta\neq\pi/2$, we let $\ell_\gamma\subset F$ be the semimajor axis of $\gamma$. We know that if $\gamma=F\cap C$, then $C$ is a cylinder of radius $r$ whose axis $\ell$ passes through $p$, satisfies $\proj_F\ell=\ell_\gamma$, and makes an angle of $\theta$ with $F$. This shows that there are only two possibilities for $C$. Thus in this case $\gamma$ also cannot lie in a third cylinder, as desired.
\end{proof}

\begin{proof}[Proof of \cref{lem:rich-lines}]
Write $\cL=\cL_{N^{1/3}}(\cP)$.
Since each plane contains at most $N^{2/3}$ points of $\cP$, we can use \cref{thm:rich-lines-r3} to bound
\[|\cL|\lesssim \frac{N^2}{N^{4/3}}+\frac{N^{5/3}}{N}+\frac{N}{N^{1/3}}\lesssim N^{2/3},\]
and
\[I(\cP,\cL)\lesssim N^{1/2}N^{1/2}+N^{4/9}N^{2/9}N^{2/3}N^{-1/3}+N^{2/9}N^{2/9}N^{2/3}N^{-1/9}+N^{2/3}+N\lesssim N.\]

Now if $\rho(p')=q'$ and $\rho(\ell_p)=\ell_q$ for some genuine rigid motion $\rho$, then $\dist(p',\ell_p)=\dist(q',\ell_q)$, implying that $p',q'$ lie on the cylinders $C_r(\ell_p),C_r(\ell_q)$ for some $r\geq 0$. (Here we abuse notation and write $C_0(\ell_p)$ for the line $\ell_p$.) Thus we write
\[\cE\subseteq\bigcup_{\ell_p,\ell_q\in\cL}\bigcup_{r\geq 0}\ell_p\times\ell_q\times C_r(\ell_p)\times C_r(\ell_q).\]

We first handle the terms with $r=0$. For lines $\ell_p,\ell_q$ and any quadruple $(p,q,p',q')\in\ell_p\times\ell_q\times\ell_p\times\ell_q$ with $\norm{p-p'}=\norm{q-q'}$, note that given $\ell_p,\ell_q$ and three of the points, there are only two possibilities for the fourth. Thus we can write
\[\abs{\cE\cap(\ell_p\times\ell_q\times\ell_p\times\ell_q)}\leq 2\abs{\ell_p\cap\cP}^2\abs{\ell_q\cap\cP}\lesssim N^{2/3}\abs{\ell_p\cap\cP}\abs{\ell_q\cap\cP}.\]
Therefore
\[\abs{\cE\cap \bigcup_{\ell_p,\ell_q\in\cL}\paren{\ell_p\times\ell_q\times\ell_p\times\ell_q}}\lesssim N^{2/3}\sum_{\ell_p,\ell_q\in\cL}\abs{\ell_p\cap\cP}\abs{\ell_q\cap\cP} = N^{2/3} I(\cP,\cL)^2\lesssim N^{8/3}.\]

From now on we focus on the terms with $r>0$. For lines $\ell_p,\ell_q\in\cL$, a radius $r>0$, and a quadruple $(p,q,p',q')\in\ell_p\times\ell_q\times C_r(\ell_p)\times C_r(\ell_q)$ with $\norm{p-p'}=\norm{q-q'}$, again note that given $\ell_p,\ell_q$ and $(p,p',q')$, there are only two possibilities for $q$. Similarly, given $\ell_p,\ell_q$ and $(q,p',q')$, there are only two possibilities for $p$. Thus the remaining contribution to $|\cE|$ can be bounded by
\begin{equation}\label{eq:boundQ}
\lesssim \sum_{\ell_p,\ell_q\in\cL}\sum_{r>0}\abs{\cP\cap C_r(\ell_p)}\abs{\cP\cap C_r(\ell_q)}\min\set{\abs{\ell_p\cap\cP},\abs{\ell_q\cap\cP}}.
\end{equation}

Let $\cV$ be the family of cylinders in $\R^3$. We apply \cref{lem:rich-variety} to $\cP,\cV$ with $C_0=2$ and $K=3$. This produces a constant $C_1$, a collection of irreducible real curves $\Gamma$, and a set $\cP_\gamma\subseteq\cP$ for each $\gamma\in\Gamma$. 

We split \cref{eq:boundQ} into the contribution coming from curves $\gamma\in\Gamma$ contained in cylinders and curves intersecting cylinders. For a cylinder $C$, define $\cP_C$ to be the points of $\cP\cap C$ which are not captured by any curve; that is, $\cP_C$ satisfies
\[\cP\cap C=\cP_C\sqcup\bigsqcup_{\gamma\in\Gamma:\gamma\subset C}\cP_\gamma.\]
The quantity in \cref{eq:boundQ} involves a pair $(p',q')\in (\cP\cap C_r(\ell_p))\times(\cP\cap C_r(\ell_q))$. We break the sum into three parts: where both $p',q'$ are captured by $\Gamma$, where $p'$ is not captured by $\Gamma$, and where $q'$ is not captured by $\Gamma$. By symmetry, the last two parts are equal. Thus defining
\[E_{\mathrm{curve}}=\sum_{\ell_p,\ell_q\in\cL}\sum_{r>0}\sum_{\substack{\gamma_p,\gamma_q\in\Gamma\\\gamma_p\subset C_r(\ell_p), \gamma_q\subset C_r(\ell_q)}}\abs{\cP_{\gamma_p}}\abs{\cP_{\gamma_q}}\min\set{|\ell_p\cap\cP|,|\ell_q\cap\cP|}\]
and
\[E_{\mathrm{cross}}=\sum_{\ell_p,\ell_q\in\cL}\sum_{r>0}\abs{\cP_{C_r(\ell_p)}}\abs{\cP\cap C_r(\ell_q)}\min\set{|\ell_p\cap\cP|,|\ell_q\cap\cP|},\]
we have that the quantity in \cref{eq:boundQ} is bounded by $E_{\mathrm{curve}}+2E_{\mathrm{cross}}$.

Now note that each $\gamma\in\Gamma$ lies in at least three distinct cylinders. By \cref{clm:only-lines-in-cylinders}, we conclude that $\Gamma$ is a collection of lines. 

To bound $E_{\mathrm{curve}}$, we take a dyadic decomposition. Write $\cL=\bigsqcup_{i\geq 0}\cL_i$ where $\cL_i$ is the set of lines $\ell\in\cL$ so that $\abs{\ell\cap\cP}\in[N^{1/3}2^i,N^{1/3}2^{i+1})$. Similarly, write $\Gamma=\bigsqcup_{j\geq 0}\Gamma_j$ where for $j>0$ we define $\Gamma_j$ to be the set of lines $\gamma\in\Gamma$ so that $\abs{\cP_\gamma}\in[N^{1/3}2^j,N^{1/3}2^{j+1})$. We let $\Gamma_0$ be the remaining lines: those $\gamma\in\Gamma$ so that $\abs{\cP_\gamma}<2N^{1/3}$.

Since $I(\cP,\cL_i)\leq I(\cP,\cL)\lesssim N$, we have $\abs{\cL_i}\lesssim N^{2/3}/2^i$. As the sets $(\cP_{\gamma})_{\gamma\in \Gamma}$ are disjoint, we also know that $\abs{\Gamma_j}\leq N^{2/3}/2^j$ for all $j>0$. For $j=0$ we have $\abs{\Gamma_0}\leq\abs{\Gamma}\lesssim N^{2/3}$.

Now define $\cS\subseteq \cL^2\times \Gamma^2$ be the set of $(\ell_p,\ell_q,\gamma_p,\gamma_q)$ for which there exists $r>0$ with $\gamma_p\subset C_r(\ell_p)$ and $\gamma_q\subset C_r(\ell_q)$. Let $\cS_{ii'jj'}\subseteq \cS$ be the set of tuples for which $(\ell_p,\ell_q,\gamma_p,\gamma_q)\in\cL_i\times\cL_{i'}\times\Gamma_j\times\Gamma_{j'}$. We can write
\[E_{\mathrm{curve}} = \sum_{i,i',j,j'\geq 0}\sum_{(\ell_p,\ell_q,\gamma_p,\gamma_q)\in \cS_{ii'jj'}}|\cP_{\gamma_p}| |\cP_{\gamma_q}| \min\{ |\ell_p \cap \cP|, |\ell_q \cap \cP| \}.\]
Now for $(\ell_p,\ell_q,\gamma_p,\gamma_q)\in \cS_{ii'jj'}$, we have
\[|\cP_{\gamma_p}| |\cP_{\gamma_q}| \min\{ |\ell_p \cap \cP|, |\ell_q \cap \cP| \}\lesssim N2^{j+j'}\min\{2^i,2^{i'}\}\lesssim N^{4/3}2^{\min\{i,i'\}+\min\{j,j'\}}.\]
Thus we can bound
\[E_{\mathrm{curve}}\lesssim N^{4/3}\sum_{i,i',j,j'\geq 0}2^{\min\{i,i'\}+\min\{j,j'\}}|\cS_{ii'jj'}|.\]

To write this sum more efficiently, define \[\cS_{ij}=\bigcup_{\substack{i_1,i_2\geq i\\j_1,j_2\geq j}}\cS_{i_1i_2j_1j_2},\qquad \cL_{\geq i}=\bigcup_{i'\geq i}\cL_{i'},\qquad\text{ and }\qquad\Gamma_{\geq j}=\bigcup_{j'\geq j}\Gamma_{j'}.\]
Then
\[E_{\mathrm{curve}}\lesssim N^{4/3}\sum_{i,j\geq 0}2^{i+j}|\cS_{ij}|.\]
Note that $|\cS_{ij}|$ is the number of quadruples $(\ell_p,\ell_q,\gamma_p,\gamma_q)\in\cL_{\geq i}^2\times\Gamma_{\geq j}^2$ where $\ell_p,\gamma_p$ are parallel and $\ell_q,\gamma_q$ are parallel, and $\dist(\ell_p,\gamma_p)=\dist(\ell_q,\gamma_q)$. We still have the bounds $\abs{\cL_{\geq i}}\lesssim N^{2/3}/2^i$ and $\abs{\Gamma_{\geq j}}\lesssim N^{2/3}/2^j$.

From this definition, note that we can view $|\cS_{i,j}|$ as a distance energy which can be bounded by \cref{thm:bipartite-guth-katz}.

\begin{lemma}
\label{thm:parallel-line-energy}
Let $\cL,\Gamma$ be sets of lines in $\R^3$ with $\abs{\cL},\abs{\Gamma}\leq M$. Define $\cE$ to be the set of quadruples $(\ell_p,\ell_q,\gamma_p,\gamma_q)\in\cL^2\times\Gamma^2$ so that $\ell_p,\gamma_p$ are parallel and $\ell_q,\gamma_q$ are parallel, and $\dist(\ell_p,\gamma_p)=\dist(\ell_q,\gamma_q)$. Then
\[\abs{\cE}\lesssim\abs{\cL}\abs{\Gamma}M\log M.\]
\end{lemma}

\begin{proof}
Given a set of lines $\cM$ in $\R^3$ and a unit vector $\hat v\in\S^2$, write $F_{\hat v}$ for the plane through the origin with normal vector $\hat v$. Then let $\cM(\hat v)\subset F_{\hat v}$ be the set of points $\ell\cap F_{\hat v}$ for each line $\ell\in\cM$ that is parallel to $\hat v$. Finally, pick a generic isometry $F_{\hat v}\cong \R^2$ so we can view $\cM(\hat v)\subset \R^2$. With these definitions, note that if $(\ell_p,\ell_q,\gamma_p,\gamma_q)\in\cE$, then there exist $\hat u,\hat v\in\S^2$ so that the quadruple corresponds to a quadruple $(p,q,p',q')\in\cL(\hat u)\times\cL(\hat v)\times\Gamma(\hat u)\times\Gamma(\hat v)$ with $\norm{p-p'}=\norm{q-q'}$. Thus we can bound
\begin{align*}
\abs{\cE}
&\leq \sum_{\hat u,\hat v\in\S^2}E(\cL(\hat u),\cL(\hat v),\Gamma(\hat u),\Gamma(\hat v))\\
&\leq E\paren{\bigsqcup_{\hat u\in\S^2}\cL(\hat u),\bigsqcup_{\hat v\in\S^2}\cL(\hat v),\bigsqcup_{\hat u\in\S^2}\Gamma(\hat u),\bigsqcup_{\hat v\in\S^2}\Gamma(\hat v)}\\
&\lesssim\abs{\cL}\abs{\Gamma} M\log M.
\end{align*}
In the second line above, since we choose a generic isometry $F_{\hat v}\cong\R^2$ for each $\hat v\in\S^2$, we may assume that the sets $\cL(\hat u)$ are disjoint for distinct $\hat u$ and similarly with the sets $\Gamma(\hat u)$.
The last line follows from \cref{thm:bipartite-guth-katz}.
\end{proof}

Applying \cref{thm:parallel-line-energy} to bound
\[\abs{\cS_{ij}}\lesssim |\cL_{\geq i}||\Gamma_{\geq j}|\max\{|\cL_{\geq i}|,|\Gamma_{\geq j}|\}\log N\lesssim\frac{N^2\log N}{2^{i+j}}\paren{2^{-i}+2^{-j}},\]
we conclude that
\[E_{\mathrm{curve}}\lesssim N^{4/3}\sum_{i,j\geq 0}2^{i+j}|\cS_{ij}|\lesssim N^{10/3}\log N\sum_{i,j\geq 0}\paren{2^{-i}+2^{-j}}\lesssim N^{10/3}\log^2N.\]

It remains to bound $E_{\mathrm{cross}}$. Let $\cV_{k}$ be the set of pairs $(r,\ell)$ so that $\abs{\cP_{C_r(\ell)}}\geq k$. By \cref{lem:rich-variety}, we know that $\abs{\cV_{k}} \lesssim N/k$ for all $k\geq C_1N^{2/3}$. First, for the terms in $E_{\mathrm{cross}}$ where $\abs{\cP_{C_r(\ell_p)}}<C_1N^{2/3}$, we have the bound
\begin{align*}
\sum_{\ell_p,\ell_q\in\cL}\sum_{r>0}C_1N^{2/3}\abs{\cP\cap C_r(\ell_q)}\min\set{|\ell_p\cap\cP|,|\ell_q\cap\cP|}
&\lesssim N^{2/3}\sum_{\ell_p\in\cL}|\ell_p\cap\cP|\sum_{\ell_q,r}\abs{\cP\cap C_r(\ell_q)}\\
&\leq N^{2/3}I(\cP,\cL)\sum_{\ell_q\in\cL}N\\
&\lesssim N^{2/3}I(\cP,\cL)|\cL|N\lesssim N^{10/3}.
\end{align*}

For the remaining terms, we take a dyadic decomposition based on the value of $\abs{\cP_{C_r(\ell_p)}}\in [C_1N^{2/3},N]$. This gives
\begin{align*}
\sum_{i:2^i\in [C_1N^{2/3},N]} \sum_{\ell_p,\ell_q\in\cL}&\sum_{\substack{r>0:\\(r,\ell_p)\in\cV_{2^i}}} 2^i\abs{\cP\cap C_r(\ell_q)}\min\set{|\ell_p\cap\cP|,|\ell_q\cap\cP|}\\
&\leq N\sum_{i:2^i\in [C_1N^{2/3},N]} \sum_{\ell_q\in\cL}\sum_{\substack{(r,\ell_p)\in\cV_{2^i}}} 2^i|\ell_q\cap\cP|\\
&\lesssim N\sum_{i:2^i\in [C_1N^{2/3},N]} 2^i\abs{\cV_{2^i}}\sum_{\ell_q\in\cL}|\ell_q\cap\cP|\\
&\lesssim N\sum_{i:2^i\in [C_1N^{2/3},N]} N\cdot I(\cP,\cL)\\
&\lesssim N^3\log N.
\end{align*}

Thus $E_{\mathrm{cross}}\lesssim N^{10/3}$, completing the proof.
\end{proof}

\subsection{Very rich reguli}
\label{ssec:rich-regulus}
In this subsection we study the distance quadruples which are captured by rich reguli. Here the term regulus refers to a \emph{real regulus} $\fr\subset\R^3$. In $\PP^3_{\CC}$, every smooth degree 2 surface is doubly ruled by lines and thus a regulus. However, in $\R^3$, some of these surfaces are no longer ruled by real lines. The non-degenerate real quadrics (the surfaces which are the affine real points on a smooth degree 2 surface) are ellipsoids, one-sheeted hyperboloids, two-sheeted hyperboloids, elliptic paraboloids, and hyperbolic paraboloids. Up to a rigid motion, these have equations
\begin{align*}
\frac{x^2}{a^2}+\frac{y^2}{b^2}+\frac{z^2}{c^2}&=1,\\
\frac{x^2}{a^2}+\frac{y^2}{b^2}-\frac{z^2}{c^2}&=1,\\
\frac{x^2}{a^2}-\frac{y^2}{b^2}-\frac{z^2}{c^2}&=1,\\
\frac{x^2}{a^2}+\frac{y^2}{b^2}-z&=0,\\
\frac{x^2}{a^2}-\frac{y^2}{b^2}-z&=0,
\end{align*}
respectively. It is easy to see that the ellipsoids, two-sheeted hyperboloids, and elliptic paraboloids contain no real lines. On the other hand, one can check that the one-sheeted hyperboloids and the hyperbolic paraboloids are doubly ruled (by real lines) and in fact contain exactly two lines through each point. (See, e.g., \cite[Section 6.4]{Som34}.)

\subsubsection{Geometry of reguli}
To control the distance quadruples that are captured by rich reguli, we will need some geometric claims.

\begin{claim}
\label{claim:regulus-symmetries}
Let $\fr_1,\fr_2\subset\R^3$ be two real reguli. Given $p\in\fr_1$ and $q\in \fr_2$, there are at most 8 genuine orientation-preserving rigid motions $\rho\in E^+_g$ so that $\rho(\fr_1)=\fr_2$ and $\rho(p)=q$.
\end{claim}

\begin{proof}
Given $p\in\fr_1$, there are exactly two lines $\ell^{(1)}_p,\ell^{(2)}_p$ contained in $\fr_1$ that pass through $p$.
There are also two lines $\ell^{(1)}_q,\ell^{(2)}_q$ through $q$ in $\fr_2$.
Now if $\rho(\fr_1)=\fr_2$ and $\rho(p)=q$, then $\rho$ maps $\{\ell^{(1)}_p,\ell^{(2)}_p\}$ to $\{\ell^{(1)}_q,\ell^{(2)}_q\}$. There are two possibilities for how these lines are paired up. Once this is chosen, there are two possibilities for each of $\rho|_{\ell^{(i)}_p}$. Therefore there are at most 8 options for $\rho|_{\ell^{(1)}_p\cup\ell^{(2)}_p}$. 
This determines $\rho$ completely since any rigid motion is determined by its value on three non-collinear points.
\end{proof}

\begin{claim}
\label{claim:point-reg-reg}
Given two real reguli $\fr_1,\fr_2$ and a point $p\in \fr_1$, the set of $q\in \fr_2$ such that there exists some $\rho\in E^+_g$ with $\rho(\fr_1)=\fr_2$ and $\rho(p)=q$ lies on a curve with degree $4$.
\end{claim}

\begin{proof}
Up to a rigid motion, the regulus $\fr_2$ is either a one-sheeted hyperboloid with equation $x^2/a^2+y^2/b^2-z^2/c^2=1$ or a hyperbolic paraboloid with equation $x^2/a^2-y^2/b^2=z$.

One can compute the Gaussian curvature of the former to be
\[\frac{-1}{a^2b^2c^2\paren{\frac{x^2}{a^4}+\frac{y^2}{b^4}+\frac{z^2}{c^4}}^2}\]
at the point $(x,y,z)$. (See, e.g., \cite[Chapter 3, Section III]{Spi75}.)
The level sets of the Gaussian curvature are the curves $\gamma_k=\fr_2\cap\mathfrak c_k$ where $\mathfrak c_k$ is the vertical cylinder defined by $x^2/a^4+y^2/b^4+(x^2/a^2+y^2/b^2-1)/c^2=k$. The base of this cylinder is an ellipse (or a point or empty). Thus the $\gamma_k$ are curves of degree at most 4 (or a set of at most 2 points). Now the image of $\rho(p)$ is confined to one of these level sets, which is the desired result. 

One can compute the Gaussian curvature of the latter to be
\[\frac{-1}{4a^2b^2\paren{\frac{x^2}{a^4}+\frac{y^2}{b^4}+\frac14}^2}\]
at the point $(x,y,z)$. The level sets of the Gaussian curvature are the curves $\gamma_k=\fr_2\cap\mathfrak c_k$ where $\mathfrak c_k$ is the vertical cylinder defined by $x^2/a^4+y^2/b^4=k$. These are curves of degree at most 4 (or a point or empty). Since the image of $\rho(p)$ is confined to one of these level sets we are done.
\end{proof}

\begin{claim}
\label{claim:point-curve-point-curve}
Given two real irreducible algebraic curves $\gamma_1,\gamma_2\subset\R^3$ and points $p\in\gamma_1$ and $q\in\gamma_2$, let $\rho\in E_g^+$ be a genuine rigid motion such that $\rho(p)=q$ and $\rho(\gamma_1)=\gamma_2$. If $\gamma_1$ is not a line, there are at most $4(\deg\gamma_2)^2$ possibilities for $\rho$.
\end{claim}

\begin{proof}
Pick points $p',p''\in\gamma_1$ so that $p,p',p''$ are not collinear. Now $\norm{q-\rho(p')}=\norm{p-p'}$, so $\rho(p')$ must lie on the intersection of $\gamma_2$ and $\S$, the sphere of radius $\norm{p-p'}$ centered at $q$. Note that $\gamma_2,\S$ are irreducible algebraic varieties and $\gamma_2\not\subset \S$ (since $q\in\gamma_2$ but $q\not\in\S$). Thus by B\'ezout's theorem, $|\gamma_2\cap\S|\leq\deg\gamma_2\deg\S=2\deg\gamma_2$. Thus there are at most $2\deg\gamma_2$ possibilities for $\rho(p')$. Similarly, there are at most $2\deg\gamma_2$ possibilities for $\rho(p'')$.

Since $p,p',p''$ are not collinear, $\rho$ is determined by its image on these three points, so there are at most $4(\deg\gamma_2)^2$ possibilities for $\rho$, as desired.
\end{proof}

\begin{claim}
\label{claim:curve-regulus-1dim-family}
Given a real regulus $\fr\subset\R^3$, a line $\ell\subset\fr$, and a real irreducible algebraic curve $\gamma\subset\R^3$ with $\gamma\neq\ell$, there are at most 4 rigid motions $\rho$ which are rotations about $\ell$ such that $\rho(\gamma)\subset\fr$.
\end{claim}

\begin{proof}
Pick $p\in\gamma\setminus\ell$. Consider the circle $\mathfrak c$ which is the image of $p$ under the rotations about $\ell$. First note that $\mathfrak c\not\subset\fr$. To see this, suppose for contradiction that $\mathfrak c\subset\fr$ and let $F$ be the plane containing $\mathfrak c$. By B\'ezout's theorem, $\fr\cap F$ is a curve of degree at most 2, yet it contains $\mathfrak c$ and the additional point $F\cap\ell$, a contradiction.

Since $\mathfrak c\not\subset\fr$, by B\'ezout's theorem, $\mathfrak c\cap\fr$ consists of at most $\deg\mathfrak c\deg\fr=4$ points. As $\rho$ is determined by $\rho(p)$ which must lie in $\mathfrak c\cap \fr$, we see that there are at most 4 possibilities for $\rho$.
\end{proof}

\begin{claim}
\label{claim:curve-reg}
Let $\gamma\subset\RR^3$ be a (not necessarily irreducible) real algebraic curve of degree $4$, and let $\fr\subset \RR^3$ be a real regulus. Then there are at most $O(1)$ reguli which are congruent to $\fr$ and contain $\gamma$.
\end{claim}

\begin{proof}
If $\gamma$ only lies in at most one regulus which is congruent to $\fr$, then the statement immediately holds.
Therefore we may suppose that $\gamma$ lies in two distinct reguli $\fr_1$ and $\fr_2$ which are both congruent to $\fr$. Write $\fr_1 = Z(f_1)$ and $\fr_2 = Z(f_2)$ where $f_1,f_2\in \RR[x,y,z]$ are polynomials of degree $2$. By \cref{thm:complete-intersection-degrees}, we see that any degree 2 polynomial $f\in\RR[x,y,z]$ vanishing on $\gamma$ is an $\R$-linear combination of $f_1$ and $f_2$. (To be precise, let $\gamma_{\CC}$ be the Zariski-closure of $\gamma$ in $\CC^3$. Then $f_1,f_2,f$ all vanish identically on $\gamma_\CC$, so \cref{thm:complete-intersection-degrees} implies that $f$ is a $\CC$-linear combination of $f_1,f_2$. Since $f_1,f_2,f$ all have real coefficients, this must actually be an $\R$-linear combination.)

Consider the 1-parameter family of surfaces $Z(f_1+tf_2)$ for $t\in\R$. Define $T\subset\R$ to be the set of $t\in\R$ so that $Z(f_1+tf_2)$ is congruent to $\fr_1$. There exists an absolute constant $C$ so that either $|T|\leq C$ or $T$ is infinite.

This follows from standard results in real algebraic geometry. We say that a set $X\subset\R^d$ is \emph{semialgebraic} of complexity $D$, if there exist $D$ polynomials $g_1,\ldots,g_D\in\R[x_1,\ldots,x_d]$ of degree at most $D$ so that $X$ can be expressed as a union of intersections of the sets $\{x\in\R^d:g_i(x)>0\}$ and $\{x\in\R^d:g_i(x)=0\}$. We will show that $T\subset\R$ is semialgebraic of complexity $O(1)$, implying the desired result.

To prove this, we define a certain set $X$ of triples $(t,c,\tilde\rho)\in\R\times\R\times\R^8$. This will be a semialgebraic set of complexity $O(1)$ with the property that $\rho\in E_g^+$ maps $\fr_1$ to $Z(f_1+tf_2)$ if and only if there exist $c\in\R$ and $\tilde\rho\in\R^8$ so that $\rho=[\tilde\rho]$ and $(t,c,\tilde\rho)\in X$. Note that when we write $\rho=[\tilde\rho]$ we are viewing $E_g^+\subset E^+(\R)\subset\P^7_\R$. First note that the set of $\tilde\rho\in\R^8$ so that $[\tilde\rho]\in E_g^+$ is a semialgebraic set: we have $E^+(\R)=\P(Z(N))$ for a degree 2 polynomial $N$ while $E^+(\R)\setminus E_g^+$ is cut out by a bounded number of polynomials of bounded degree (see \cref{prop:actual-rigid-motion}). Then $[\tilde\rho]$ maps $\fr_1=Z(f_1)$ to $Z(f_1+tf_2)$ if and only if $f_1\circ[\tilde\rho]^{-1}$ and $f_1+tf_2$ are multiples of each other. Now define $X$ to be the set of triples $(t,c,\tilde\rho)$ so that $c\cdot f_1\circ[\tilde\rho]^{-1}=f_1+tf_2$. To see that this condition is semialgebraic, note that it can be rewritten as $c\cdot f_1(p)=(f_1+tf_2)(q)$ for all $(p,q)\in\R^6$ with $\phi(p,q)\in F_{[\tilde\rho]}^{\Phy}$. The latter is cut out by a bounded number of polynomials of bounded degree.

To complete the proof, note that $T$ is the projection of $X$ to the first coordinate. Thus $T$ is a semialgebraic set as the Tarski--Seidenberg theorem \cite[Theorem 2.76]{BPR06} states that the projection of a semialgebraic set is semialgebraic. This result follows from (and is equivalent to) the fact that quantifier elimination is possible over $\R$. The existence of efficient algorithms for quantifier elimination imply that the complexity of $T$ is bounded in terms of the complexity of $X$. (For example, \cite[Theorem 14.13]{BPR06} shows that the complexity of $T$ is at most polynomial in the complexity of $X$.)

Now for the sake of contradiction, assume that $T$ is infinite. Letting $\vec{v} = (x,y,z)^{\intercal}$, we may write
\[f_i = \vec{v}^{\intercal} A_i\vec{v}+B_i^{\intercal}\vec{v}+C_i\]
where $A_i$ is $3$-by-$3$ and symmetric, and $B_i$ is $3$-by-$1$.
As before, if $\rho$ maps $\fr_1$ to $Z(f_1+tf_2)$ then $f_1\circ\rho^{-1}$ and $f_1+tf_2$ are multiples of each other; in particular, the homogeneous degree 2 part of these polynomials are multiples of each other. The latter has homogeneous degree 2 part defined by the quadratic form $A_1+tA_2$ and the former is defined by $OA_1O^{-1}$ for an orthogonal matrix $O$; in particular, writing $\rho$ as the composition of a rotation and a translation, $O$ is the matrix of the rotation. Thus for $t\in T$, there exist $c_t\neq 0$ and an orthogonal matrix $O_t$ so that
\[c_t\cdot O_tA_1O_t^{-1}=A_1+tA_2.\]

Raising both sides to the $k$th power and taking the trace, we conclude that
\[\tr((A_1+tA_2)^k)=c_t^k(\lambda_1^k+\lambda_2^k+\lambda_3^k)\]
for all $k\geq 1$ and all $t\in T$, where $\lambda_1,\lambda_2,\lambda_3\in\R$ are the eigenvalues of $A_1$.

Let $m_k = \lambda_1^k+\lambda_2^k+\lambda_3^k$. Suppose that $m_1$ or $m_3\neq 0$. We will show that in this case $A_1,A_2$ are multiples of each other. First assume $m_1\neq 0$. Then define the polynomial
\[P(t)=m_1^2\tr((A_1+tA_2)^2)=m_2\tr(A_1+tA_2)^2.\]
To be precise, for each $t\in T$ both expressions above are equal to $c_t^2m_1^2m_2$; since they are polynomials and $T$ is infinite, they must be equal for all $t\in\R$.

If $\tr A_2\neq 0$, there exist some $t\in\R$ so that $\tr(A_1+tA_2)=0$. For this value of $t$, we have $P(t)=0$; since $m_1\neq 0$, we conclude that $\tr((A_1+tA_2)^2)=0$. However, since $A_1+tA_2$ is a real symmetric matrix, this implies that $A_1+tA_2$ is the zero matrix, as desired. Instead, if $\tr A_2=0$, then $P(t)=m_2\tr(A_1)^2$, a constant. However the coefficient of $t^2$ in the other expression for $P(t)$ is $m_1^2\tr(A_2^2)$, implying that $\tr(A_2^2)=0$. Since $A_2$ is real symmetric, this implies that $A_2$ is the zero matrix, which is a contradiction to the fact that $\deg f_2=2$.

Next assume $m_3\neq 0$. Similar to before, define the polynomial
\[P(t)=m_3^2\tr((A_1+tA_2)^2)^3=m_2^3\tr((A_1+tA_2)^3)^2.\]
As before, both expressions above are equal to $c_t^6m_2^3m_3^2$ for all $t\in T$, meaning that the polynomials must agree. Factoring $P(t)$ over $\CC$ into linear factors, the first expression implies that every root appears with multiplicity divisible by 3 and the second implies the multiplicity is divisible by 2. Thus $P(t)=kg(t)^6$ for $k\in\R$ and $g\in\R[t]$ with $\deg g\leq 1$. First note that if $\deg g=0$, then $P(t)$ is a constant. Now the coefficient of $t^6$ in the first expression for $P(t)$ is $m_3^2\tr(A_2^2)^3$, implying that $\tr(A_2^2)=0$. As before, this implies that $A_2$ is the zero matrix, contradiction. 

Otherwise, $\deg g=1$, so there exists some $t\in\R$ so that $g(t)=0$. For this value of $t$ we have $0=P(t)=m_3^2\tr((A_1+tA_2)^2)^3$, implying that $\tr((A_1+tA_2)^2)=0$. As before, this implies that $A_1+tA_2$ is the zero matrix, as desired.

Thus we have shown that $A_1+t_0A_2$ is the zero matrix for some $t_0\in\R$. Now $\gamma=Z(f_1,f_2)\subseteq Z(f_1,f_1+t_0f_2)$. Since $A_1+t_0A_2$ is the zero matrix, $f_1+t_0f_2$ has degree at most 1. It cannot be the zero polynomial, since we assumed that $Z(f_1)\neq Z(f_2)$. Thus by B\'ezout's theorem, we see $\deg \gamma\leq \deg f_1\deg(f_1+t_0f_2)\leq 2$, a contradiction.

It remains to handle the case $m_1=m_3=0$. Solving $\lambda_1+\lambda_2+\lambda_3=\lambda_1^3+\lambda_2^3+\lambda_3^3=0$ implies one of the eigenvalues is 0, say $\lambda_3=0$ and $\lambda_1+\lambda_2=0$. 
Therefore up to rescaling and an orthogonal change of basis, we may assume that 
        \[A_1 = \begin{bmatrix}
            1&&\\&-1&\\&&0
        \end{bmatrix}.\]
        Note that in this case, the regulus must be a hyperbolic paraboloid, and by translation we may now assume that $f_1(x,y,z) = x^2-y^2-kz$ for some $k\in\RR\setminus \{0\}$. (By translation we can shift the equation $x^2-y^2+b_1x+b_2y+b_3z+c=0$ to the equation $x^2-y^2-kz=0$ or $x^2-y^2-k=0$. The latter is degenerate.)
        Suppose that
        \[A_2=\begin{bmatrix}
            a_{11}&a_{12}&a_{13}\\a_{12}&a_{22}&a_{23}\\a_{13}&a_{23}&a_{33}
        \end{bmatrix}.\]
Since $m_1=0$, we know $\tr(A_1+tA_2)=0$ for all $t\in T$ (and thus all $t\in\R$). This implies $\tr(A_2)=0$, so $a_{11}+a_{22}+a_{33}=0$.
Similarly, since $m_3=0$, we know $\tr((A_1+tA_2)^3)=0$. Taking the coefficient of $t$, we see $3\tr(A_1^2A_2)=0$. Expanding out the definition of $A_1,A_2$ gives $a_{22}=-a_{11}$. Combining the two implies $a_{33}=0$.

Write $A_2=a_{11}A_1+A_2'$. Defining $g_2=f_2-a_{11}f_1$, we have $Z(f_1+tf_2)=Z((1+ta_{11})f_1+tg_2)$. In particular, there are infinitely many values of $s\in\R$ so that $Z(f_1+sg_2)$ is congruent to $\fr_1$ (take any $s=t/(1+ta_{11})$ for $t\in T\setminus\{-1/a_{11}\}$). Replacing $f_2$ by $g_2$, we may assume that $a_{11}=0$. As $a_{22}=-a_{11}$, this implies $a_{22}=0$ as well.
        
We still know that $\tr((A_1+tA_2)^3)$ is the zero polynomial. Taking the coefficient of $t^2$, we see $3\tr(A_1A_2^2)=0$.
        Expanding, we get
        \[\tr(A_1A_2^2) = (a_{12}^2+a_{13}^2)-(a_{12}^2+a_{23}^2) = a_{13}^2-a_{23}^2,\]
        showing that $a_{23}=\pm a_{13}$.
        Finally, taking the coefficient of $t^3$, we see $\tr(A_2^3)=0$, which shows that
        \[6a_{12}a_{23}a_{13}=0.\]
We can rescale $A_2$ so that it is one of the following matrices
        \[A_2=\begin{bmatrix}
            &1&\\1&&\\&&0
        \end{bmatrix}\qquad\text{or}\qquad\begin{bmatrix}
            &&1\\&&1\\1&1&0
        \end{bmatrix}\qquad\text{or}\qquad\begin{bmatrix}
            &&1\\&&-1\\1&-1&0
        \end{bmatrix}.\]

We know that there are infinitely many values of $t$ so that $Z(f_1+tf_2)$ is congruent to $Z(f_1)=Z(x^2-y^2-kz)$. In the first case, $Z(f_1+tf_2)$ has quadratic part $x^2+2txy-y^2$, in particular, this does not depend on $z$. This quadratic form has eigenvalues $\pm\sqrt{1+t^2}$, so we can perform a rotation about the $z$-axis so that $Z(f_1+tf_2)$ maps to
\[Z(\sqrt{1+t^2}x^2-\sqrt{1+t^2}y^2+a_xx+a_yy+a_zz+a)\]
for some coefficients $a_x,a_y,a_z,a\in\R$, where crucially $a_z$ is the coefficient of $z$ in $f_1+tf_2$. Note that $a_z\neq 0$ since otherwise this surface would be a cylinder instead of a regulus. Thus we can perform a translation so that this surface maps to
\[Z(\sqrt{1+t^2}x^2-\sqrt{1+t^2}y^2+a_zz)=Z\paren{x^2-y^2+\frac{a_z}{\sqrt{1+t^2}}z}.\]
We claim that this surface is congruent to $Z(f_1)=Z(x^2-y^2-kz)$ if and only if $k=\pm a_z/\sqrt{1+t^2}$. Indeed, $Z(x^2-y^2-kz)$ has Gaussian curvature
\[\frac{-4}{k^2\paren{1+\frac{4x^2}{k^2}+\frac{4y^2}{k^2}}^2}.\]
This is largest in magnitude at the point $(0,0,0)$ which has curvature $-4/k^2$. Thus if $Z(x^2-y^2-kz)$ and $Z(x^2-y^2-k'z)$ are congruent, we must have $k=\pm k'$, as desired.

Letting $\ell$ be the coefficient of $z$ in $f_2$, note that $a_z$, the coefficient of $z$ in $f_1+tf_2$, is $-k+\ell t$. Thus, in this case, $Z(f_1+tf_2)$ is congruent to $Z(f_1)$ only if
\[k=\pm\frac{k-\ell t}{\sqrt{1+t^2}}.\]
This has infinitely many solutions only when the polynomial identity $k^2(1+t^2)=(k-\ell t)^2$ holds in $\R[t]$. This implies $k=\ell=0$, which means that $Z(f_1)=Z(x^2-y^2-kz)$ is the union of two planes, contradiction.

In the second case, the homogeneous degree 2 parts of $f_1$ and $f_2$ are $x^2-y^2$ and $2xz+2yz$, respectively. Homogenize $f_1,f_2$ to the homogeneous degree 2 polynomials $\tilde f_1,\tilde f_2\in\R[w,x,y,z]$ so that $f_1(x,y,z)=\tilde f_1(1,x,y,z)$ and similarly for $f_2$. Now we see that $\tilde f_1(0,x,-x,z)=\tilde f_2(0,x,-x,z)=0$ for all $x,z$. This implies that $Z(\tilde f_1,\tilde f_2)\subset\P^3_{\R}$ contains both the degree 4 curve $\gamma\subset\R^3\subset\P^3_{\R}$ and the line $\{[0:x:-x:z]:x,z\in\R\}\subset\P^3_\R\setminus\R^3$. This contradicts B\'ezout's theorem.

An analogous argument handles the third case, as the homogeneous degree 2 parts of $f_1$ and $f_2$ are $x^2-y^2$ and $2xz-2yz$, respectively. Thus $\tilde f_1(0,x,x,z)=\tilde f_2(0,x,x,z)=0$ in this case. This gives a contradiction for the same reason.
\end{proof}

\begin{remark}
While the proof of \cref{claim:curve-reg} seems significantly more involved than one would expect given its innocuous statement, we point out that the result is nearly false.

Consider the reguli $\fr_1=Z(f_1)$ and $\fr_2=Z(f_2)$ where
\[f_1=x^2-y^2-z\qquad\text{and}\qquad f_2=xz+yz-y.\]
The intersection $\fr_1\cap\fr_2$ consists of a single degree 3 curve $\gamma\subset\R^3$. In $\PP^3_{\RR}$, the intersection also contains the line $\ell=\{[0:x:-x:z]:x,z\in\RR\}$ at infinity.

Now the infinite family of reguli $Z(f_1+tf_2)$ for $t\in \R$ all contain the degree 3 curve $\gamma$ as well as the line $\ell$ at infinity. Furthermore, one can show that all of these reguli are congruent.
This example shows that any proof of \cref{claim:curve-reg} must use the fact that $\gamma$ actually lies in $\R^3$ and not just in $\PP^3_{\RR}$.

One might attempt to turn this into a counterexample for \cref{claim:curve-reg} by applying a projective transformation to $\PP^3_{\RR}$ to move $\gamma\cup\ell$ into $\RR^3$. However, such a projective transformation will destroy the congruence of this class of reguli.
\end{remark}

\subsubsection{Distances quadruples on rich reguli}

\begin{proposition}\label{lem:rich-reguli}
There exists a constant $C$ such that the following holds. Let $\cP\subset\RR^3$ be a set of size $N$ such that $|\cP\cap \gamma|\leq N^{2/3}\deg\gamma$ for every curve $\gamma\subset\CC^3$, and such that every plane contains at most $N^{2/3}$ points of $\cP$. Fix parameters $k\geq CN^{2/3}\log N$ and $\eta \leq 1$. Let $\cR\subset E^+_g$ be a set of genuine positive rigid motions that are $k$-rich for $\cP$. For each $\rho\in \cR$, let $H_{\rho}$ be a subset of $\{(p,q)\in \cP\times \cP: \rho(p)=q\}$ of size at most $2k$ and let $\fR_{\rho}$ be a collection of at most $\eta^{-1}$ reguli, each contained in the 3-flat $\phi^{-1}(F_{\rho}^{\Phy})\subset\RR^6$ such that each is $\eta k$-rich for $H_{\rho}$.

Let $\cE$ be the set of tuples $(p,q,p',q')\in \cP^4$ for which there exists a rigid motion $\rho\in\cR$ and a regulus $\fr\in\fR_\rho$ so that $(p,q),(p',q')\in\fr\cap H_\rho$. Then
\[\abs{\cE}\lesssim \eta^{-1}N^{10/3}\log^2N.\]
\end{proposition}

The hypothesis that $(p,q),(p',q'),\fr$ all lie in $\phi^{-1}(F_\rho^{\Phy})$ implies that $\rho(p)=q$ and $\rho(p')=q'$. Furthermore, writing $\pi_1,\pi_2$ for the projections onto the first and last three coordinates, $\rho$ also maps the regulus $\pi_1(\fr)$ to the (congruent) regulus $\pi_2(\fr)$.

\begin{proof}
Let $\fR$ be the collection of real reguli in $\R^3$. We apply \cref{lem:rich-variety} to $\cP,\fR$ with parameters $C_0=K=2$. This produces a set of curves $\Gamma$ and disjoint subsets $\cP_{\gamma}\subseteq\cP$. For each regulus $\fr$, define $\Gamma_{\fr}\subseteq\Gamma$ to be the set of curves contained in $\fr$. Also write $\fr\cap\cP=\cP_{\fr\cap\Gamma}\sqcup\cP_{\fr\setminus\Gamma}$ where 
\[\cP_{\fr\cap\Gamma}=\bigcup_{\gamma\in \Gamma_{\fr}}P_{\gamma}\qquad\text{ and }\qquad\cP_{\fr\setminus\Gamma}=(\fr\cap\cP)\setminus\cP_{\fr\cap\Gamma}\]
are the points of $\fr\cap\cP$ captured and not captured by $\Gamma$ respectively.

Let $\fR^{\mathrm{pt}}$ be the set of reguli $\fr$ with $|\cP_{\fr\setminus\Gamma}|>C_1N^{2/3}$ and let $\fR^{\mathrm{curve}}$ be the set of reguli $\fr$ with $\deg \Gamma_{\fr} > C_1^2 N^{1/3}$. By \cref{lem:rich-variety}(iii)(iv), we have $|\fR^{\mathrm{pt}}|,|\fR^{\mathrm{curve}}|\lesssim N^{1/3}$.

Let $\tilde\cE$ be the set of tuples $(p,q,p',q',\rho,\fr_1,\fr_2)\in\cP^4\times\cR\times\fR^2$ so that the regulus $\fr\subset\phi^{-1}(F_\rho^{\Phy})$ defined by $\pi_1(\fr)=\fr_1$ and $\pi_2(\fr)=\fr_2$ (such a regulus exists uniquely if and only if $\rho(\fr_1)=\fr_2$) satisfies $\fr\in\fR_\rho$ and $(p,q),(p',q')\in\fr\cap H_\rho$. Clearly $\cE$ is the projection of $\tilde\cE$ onto the first four coordinates.

We start with two simple estimates.

\begin{claim}\label{claim:rho-given-reg-reg}
Given reguli $\fr_1$ and $\fr_2$, we have
\[\abs{\set{\rho\in \cR: \exists \fr\in \fR_{\rho}\text{ s.t. }\pi_1(\fr)=\fr_1,\,\pi_2(\fr)=\fr_2}}\lesssim \frac{N^{5/3}}{\eta k}.\]
\end{claim}

\begin{proof}
Note that $\fr$ is uniquely determined by $\rho,\fr_1,\fr_2$. Define $\cR'\subseteq \cR$ to be the set of $\rho\in\cR$ that send $\fr_1$ to $\fr_2$ so that the corresponding $\fr$ is $\eta k$-rich for $H_\rho$. Our goal is to show that $\abs{\cR'}\lesssim N^{5/3}/\eta k$.
    
By \cref{claim:regulus-symmetries}, we know that each $(p,q)\in\cP^2$ lies in at most 8 sets of the form $\{(p,q)\in\fr_1\times \fr_2:\rho(p)=q\}$ for $\rho\in\cR'$. Furthermore, by definition, each set has size at least $\eta k$. On the other hand, for each $p\in\cP\cap\fr_1$, there are few possibilities for $q\in\cP\cap\fr_2$ for which there exists some $\rho\in\cR'$ with $\rho(p)=q$. In particular, by \cref{claim:point-reg-reg}, all such $q$ lie on a curve of degree $4$. By hypothesis, this curve contains at most $4N^{2/3}$ points of $\cP$. Thus we see that 
\[\eta k|\cR'|\leq\sum_{\rho\in\cR'}\abs{\{(p,q)\in\fr_1\times \fr_2:\rho(p)=q\}}\lesssim N^{5/3}.\]
Rearranging gives the desired bound on $|\cR'|$.
\end{proof}

\begin{claim}
\label{claim:pt-pt-rho-reg-reg}
We have
\[\abs{\set{(p,q,\rho,\fr_1,\fr_2)\in\cP^2\times\cR\times\fR^2:\exists (p',q')\text{ s.t. }(p,q,p',q',\rho,\fr_1,\fr_2)\in \tilde\cE}}\lesssim \eta^{-1}N^{8/3}\log N.\]
The symmetric statement with $(p,q)$ and $(p',q')$ swapped also holds.
\end{claim}

\begin{proof}
By \cref{thm:very-rich-bound-main}, we have $\abs{\cR} = k^{-3/2}N^{3}\log^{3/2}N$ as long as $C$ is chosen larger than the constant in \cref{thm:very-rich-bound-main}. Given $\rho$, there are at most $\eta^{-1}$ choices for $\fr\in\fR_\rho$, which determines $\fr_1,\fr_2$. Finally, there are at most $2k$ choices for $(p,q)\in H_{\rho}$.
This gives a bound of
\[k^{-3/2}N^3\log^{3/2}N \cdot \eta^{-1} \cdot 2k\lesssim \eta^{-1}N^{8/3}\log N,\]
since $k\gtrsim N^{2/3}\log N$.
\end{proof}

With the preparations out of the way, we bound $|\cE|$. We start with two easy cases; the first is when both $\fr_1,\fr_2$ lie in $\fR^{\mathrm{pt}}\cup\fR^{\mathrm{curve}}$.

\begin{claim}
\label{claim:both-reguli-rich}
We have
\[\abs{\set{(p,q,p',q',\rho,\fr_1,\fr_2)\in \tilde\cE:\fr_1,\fr_2\in \fR^{\mathrm{pt}}\cup \fR^{\mathrm{curve}}}}\lesssim \eta^{-1}N^{10/3}.\]
\end{claim}

\begin{proof}
Since $|\fR^{\mathrm{pt}}\cup \fR^{\mathrm{curve}}|\lesssim N^{1/3}$, there are $\lesssim N^{2/3}$ possibilities for the pair $(\fr_1,\fr_2)$. Each pair gives $\lesssim N^{5/3}/\eta k$ possibilities for $\rho$ by \cref{claim:rho-given-reg-reg}. For each choice of $\fr_1,\fr_2,\rho$, there are at most $(2k)^2$ possibilities for $(p,q,p',q')\in H_{\rho}^2$. As $k\leq N$, this gives a bound of
\[\lesssim N^{2/3}\cdot N^{5/3}\eta^{-1}k^{-1}\cdot k^2 \lesssim \eta^{-1}N^{10/3}.\qedhere\]
\end{proof}

Next we handle the case when one of the reguli $\fr_i$ does not lie in $\fR^{\mathrm{pt}}$ and one of the points on that regulus does not lie in the set of points captured by the curves $\Gamma_{\fr_i}$.

\begin{claim}
\label{claim:poor-regulus-point-off-curve}
We have
\[\abs{\set{(p,q,p',q',\rho,\fr_1,\fr_2)\in\tilde\cE:\fr_1\not\in\fR^{\mathrm{pt}},\text{ $p$ or $p'$ is in }\cP_{\fr_1\setminus\Gamma}}}\lesssim \eta^{-1}N^{10/3}\log N.\]
The symmetric statement with $(p,p',\fr_1)$ swapped with $(q,q',\fr_2)$ also holds.
\end{claim}

\begin{proof}
Without loss of generality, suppose $p'\in\cP_{\fr_1\setminus\Gamma}$; the other cases can be handled symmetrically.
In this case, the definition of $\fR^{\mathrm{pt}}$ implies that for each fixed choice of $(p,q,\rho,\fr_1,\fr_2)$ there are only $C_1N^{2/3}$ possibilities for $p'$. With $p'$ chosen, $q'=\rho(p')$ is uniquely determined.
By \cref{claim:pt-pt-rho-reg-reg}, this gives a bound of
\[\lesssim \eta^{-1}N^{8/3}\log N\cdot N^{2/3} = \eta^{-1}N^{10/3}\log N.\qedhere\]
\end{proof}

Now we handle the case when one of the reguli is special (lies in $\fR^{\mathrm{pt}}\cup\fR^{\mathrm{curve}}$) and the two points on the other regulus are captured by the curves. 

\begin{claim}
\label{claim:one-regulus-rich-pts-on-curves}
We have
\[\abs{\set{(p,q,p',q',\rho,\fr_1,\fr_2)\in \tilde\cE:\fr_1\in \fR^{\mathrm{pt}}\cup \fR^{\mathrm{curve}}, q,q'\in \cP_{\fr_2\cap\Gamma}}}\lesssim \eta^{-1}N^{10/3}\log N.\]
The symmetric statement with $(p,p',\fr_1)$ swapped with $(q,q',\fr_2)$ also holds.
\end{claim}

\begin{proof}
For each pair $(q,q')\in \cP\times\cP$, there is at most one pair of curves $(\gamma,\gamma')\in \Gamma\times \Gamma$ such that $q\in \cP_{\gamma}$ and $q'\in\cP_{\gamma'}$. As we saw earlier, there are at most $\lesssim N^{1/3}$ possibilities for $\fr_1\in \fR^{\mathrm{pt}}\cup \fR^{\mathrm{curve}}$.

We first deal with the case where one of $\rho^{-1}(\gamma), \rho^{-1}(\gamma')$ is not contained in $\fr_1$.
Without loss of generality, assume that $\rho^{-1}(\gamma)\not\subseteq \fr_1$. By \cref{claim:pt-pt-rho-reg-reg}, there are $\lesssim \eta^{-1}N^{8/3}\log N$ choices for $(p',q',\rho,\fr_1,\fr_2)$. There are $|\Gamma|\lesssim N^{2/3}$ choices for $\gamma$ (by \cref{lem:rich-variety}(i)). Then $p\in\rho^{-1}(\gamma)\cap\fr_1$; by B\'ezout's theorem, this is a set of size at most $2\deg(\gamma)\lesssim 1$ (by \cref{lem:rich-variety}(v)). Finally $q=\rho(p)$ is uniquely determined. Thus in this case we get a bound of 
\[\lesssim \eta^{-1}N^{8/3}\log N\cdot N^{2/3} = \eta^{-1}N^{10/3}\log N.\]
From now on, we may assume that $\rho^{-1}(\gamma), \rho^{-1}(\gamma')\subset \fr_1$.

Next consider the case $\gamma=\gamma'$. Given $(p,q,\rho,\fr_1,\fr_2)$, the curve $\gamma$ with $q\in\cP_\gamma$ is uniquely determined. By hypothesis $\gamma'=\gamma$, which contains $\lesssim N^{2/3}\deg\gamma'\lesssim N^{2/3}$ points of $\cP$. Thus there are $\lesssim N^{2/3}$ possibilities for $q'\in\cP_{\gamma'}\subset \gamma'\cap \cP$. This uniquely determines $p'=\rho^{-1}(q')$. Again using \cref{claim:pt-pt-rho-reg-reg}, we see that in this case we get a bound of
\[\lesssim \eta^{-1}N^{8/3}\log N\cdot N^{2/3} = \eta^{-1}N^{10/3}\log N.\]
From now on, we will assume that $\gamma\neq \gamma'$.
    
The next case is that when $\deg\gamma+\deg \gamma'\geq 4$. First choose $(q,q',\fr_1)$. There are $\lesssim N\cdot N\cdot N^{1/3}$ possibilities as $|\fR^{\mathrm{pt}}\cup \fR^{\mathrm{curve}}|\lesssim N^{1/3}$. Recall that $(q,q')$ determines $\gamma,\gamma'$. 
If $\deg \gamma+\deg \gamma'>4$, then there is at most one possibility for $\fr_2\supset\gamma\cup\gamma'$ by B\'ezout's theorem. 
If $\deg \gamma+\deg \gamma'=4$, then by \cref{claim:curve-reg} there are $O(1)$ possibilities for $\fr_2$ once $(q,q',\fr_1)$ are determined. Given $(q,q',\fr_1,\fr_2)$, by \cref{claim:rho-given-reg-reg}, there are at most $N^{5/3}/\eta k$ possibilities for $\rho$. Then given $(q,q',\rho,\fr_1,\fr_2)$, the pair $(p,p')$ is uniquely determined as $\rho(p)=q$ and $\rho(p')=q'$.
In total, in this case we get a bound of 
\[\lesssim N\cdot N\cdot N^{1/3}\cdot \frac{N^{5/3}}{\eta k}\leq \eta^{-1}N^{10/3}.\]
From now on we will assume that $\deg\gamma+\deg\gamma'\leq 3$, implying that one of $\gamma,\gamma'$ is a line.

Without loss of generality, suppose that $\gamma$ is a line (and $\gamma'\neq \gamma$). We claim that given $(p,q,q',\fr_1)$, there are only 16 possibilities for $\rho$. First note that there are 4 possibilities for $\rho^{-1}|_{\gamma}$, since $\rho^{-1}$ must map $q$ to $p$ and must map the line $\gamma$ through $q$ to one of the two lines in $\fr_1$ through $p$. After choosing which of these lines is the image of $\gamma$, there are 2 possibilities for $\rho^{-1}|_\gamma$ sending $q$ to $p$.

Once $\rho^{-1}|_{\gamma}$ is fixed, the possibilities for $\rho^{-1}$ all differ by rotations around a fixed axis. Thus we are in a position to apply \cref{claim:curve-regulus-1dim-family} which states that there are at most 4 possibilities for such a $\rho^{-1}$ with $\rho^{-1}(\gamma')\subset\fr_1$.

Now given $(p,q,q',\fr_1)$, we have shown that there are 16 possibilities for $\rho$ which determines $p'=\rho^{-1}(q')$ and $\fr_2=\rho(\fr_1)$ uniquely. Thus we get a bound of $\lesssim N\cdot N \cdot N \cdot N^{1/3}=N^{10/3}$ in this case as $|\fR^{\mathrm{pt}}\cup \fR^{\mathrm{curve}}|\lesssim N^{1/3}$.
\end{proof}

For our final case, we assume that neither regulus lies in $\fR^{\mathrm{curve}}$ and all four points are captured by the curves. This is the only case where we bound the contribution to $\cE$ instead of to $\tilde\cE$.

\begin{claim}
\label{claim:few-curves-on-reguli-all-pts-on-curves}
    We have
    \begin{multline*}\left|\left\{(p,q,p',q')\in\cE:\exists (\rho,\fr_1,\fr_2) \text{ s.t. }(p,q,p',q',\rho,\fr_1,\fr_2)\in \tilde\cE,\right.\right.\\\left.\left. p,p'\in \cP_{\fr_1\cap\Gamma},\,q,q'\in \cP_{\fr_2\cap\Gamma},\,\fr_1,\fr_2\not\in \fR^{\mathrm{curve}}\right\}\right|\lesssim \eta^{-1}N^{10/3}\log^{2}N.\end{multline*}
\end{claim}
\begin{proof}
Consider a tuple $(p,q,p',q',\rho,\fr_1,\fr_2)\in  \tilde\cE$ so that $\fr_1,\fr_2\not\in\fR^{\mathrm{curve}}$, $p,p'\in \cP_{\fr_1\cap\Gamma}$ and $q,q'\in \cP_{\fr_2\cap\Gamma}$. Then there exist  curves $\gamma_{p},\gamma_{p'}\in \Gamma_{\fr_1}$ and $\gamma_q,\gamma_{q'}\in \Gamma_{\fr_2}$ so that $p\in\cP_{\gamma_p}$ and similarly for the other curves.

First let us consider the case $\rho(\gamma_{p'})\neq \gamma_{q'}$. By \cref{claim:pt-pt-rho-reg-reg}, there are $\lesssim \eta^{-1}N^{8/3}\log N$ possibilities for $(p,q,\rho,\fr_1,\fr_2)$. As $\fr_1,\fr_2$ are both not in $\fR^{\mathrm{curve}}$, there are $\lesssim (N^{1/3})^2$ possibilities for $(\gamma_{p'},\gamma_{q'})$. Now $q'\in\rho(\gamma_{p'})\cap\gamma_{q'}$ which is a collection of $O(1)$ points. Thus in this case we get a bound of
\[\lesssim \eta^{-1}N^{8/3}\log N\cdot N^{2/3}=\eta^{-1}N^{10/3}\log N.\]
From now on we may assume that $\rho(\gamma_{p'})=\gamma_{q'}$ and symmetrically, we may assume $\rho(\gamma_{p})=\gamma_q$.

Next consider the case that $\gamma_p$ is not a line. There are $N^3$ choices for $(p,q,p')$. From $p,q$, the curves $\gamma_p,\gamma_q$ are determined. Since $\gamma_p$ is not a line, \cref{claim:point-curve-point-curve} implies that there are at most $4(\deg\gamma_p)^2\lesssim 1$ possibilities for $\rho$ satisfying $\rho(p)=q$ and $\rho(\gamma_p)=\gamma_q$. Finally $q'=\rho(p')$ determines $q'$, giving a bound of $\lesssim N^3$ possible $(p,q,p',q')$ in this case.
    
We may therefore assume that $\gamma_p$ is a line. As $\gamma_q=\rho(\gamma_p)$ we may also assume that $\gamma_q$ is a line. Symmetrically we may assume that $\gamma_{p'},\gamma_{q'}$ are both lines.

If $\gamma_p$ and $\gamma_q$ are both $N^{1/3}$-rich for $\cP$, then the number of quadruples $(p,q,p',q')$ in this case is $\lesssim N^{10/3}\log^2N$ by \cref{lem:rich-lines}. Therefore we may assume that at least one of $\gamma_p,\gamma_q$ is less than $N^{1/3}$-rich. Symmetrically we can make this assumption for $\gamma_{p'},\gamma_{q'}$.

We next study the case that $\gamma_p,\gamma_{p'}$ are not parallel. There are $|\Gamma|^4\lesssim N^{8/3}$ choices for $\gamma_p,\gamma_{p'},\gamma_q,\gamma_{q'}$. We claim that given these four lines, there are at most 4 choices for $\rho$ with $\rho(\gamma_p)=\gamma_q$ and $\rho(\gamma_{p'})=\gamma_{q'}$. To see this, let $(x,x')\in \gamma_p\times \gamma_{p'}$ be the pair of points that minimize the distance between $x$ and $x'$. Since $\gamma_q$ and $\gamma_{q'}$ are also not parallel, we may similarly define $(y,y')\in \gamma_q\times \gamma_{q'}$. Now since $\rho(\gamma_p)=\gamma_q$ and $\rho(\gamma_{p'})=\gamma_{q'}$, we must have $\rho(x)=y$ and $\rho(x')=y'$. Now $\rho(x)=y$ and $\rho(\gamma_p)=\gamma_q$ determines $\rho|_{\gamma_p}$ up to 2 possibilities. Similarly there are 2 possibilities for $\rho|_{\gamma_{p'}}$ and these together determine $\rho$.

Given $(\gamma_p,\gamma_{p'},\gamma_q,\gamma_{q'},\rho)$, we claim that there are $(N^{1/3})^2$ possibilities for $(p,q,p',q')$. Indeed, we have $p\in (\cP\cap\gamma_p)\cap\rho^{-1}(\cP\cap\gamma_q)$. By hypothesis, one of these two sets has size at most $N^{1/3}$. Then $q=\rho(p)$ is uniquely determined. The same argument applies for $p'\in (\cP\cap\gamma_{p'})\cap\rho^{-1}(\cP\cap\gamma_{q'})$ and $q'=\rho(p')$. Thus in this case we get a bound of 
\[|\Gamma|^4(N^{1/3})^2\lesssim N^{10/3}.\]
    
Finally we may assume that $\gamma_p$ is parallel to $\gamma_{p'}$, which implies that $\gamma_q$ is parallel to $\gamma_{q'}$. We are now in a position to apply \cref{thm:parallel-line-energy} with both sets of lines equal to $\Gamma$. This shows that the number of quadruples $(\gamma_p,\gamma_q,\gamma_{p'},\gamma_{q'})\in\Gamma^4$ so that $\gamma_p,\gamma_{p'}$ are parallel lines and $\gamma_q,\gamma_{q'}$ are parallel lines and there exists $\rho$ with $\rho(\gamma_p)=\gamma_q$ and $\rho(\gamma_{p'})=\gamma_{q'}$ is bounded by $\lesssim|\Gamma|^3\log|\Gamma|\lesssim N^2\log N$.

Such a quadruple does not yet determine $(p,q,p',q')$, but we see that there are $|\gamma_p||\gamma_q|\leq N$ possibilities for $(p,q)$ as one of $\gamma_p,\gamma_q$ is less than $N^{1/3}$-rich while by hypothesis, both are at most $N^{2/3}$-rich. Now there are at most 2 possibilities for $\rho|_{\gamma_p}$ given $\rho(p)=q$ and $\rho(\gamma_p)=\gamma_q$. Given that $\rho(\gamma_{p'})=\gamma_{q'}$, all of $\rho$ is determined from $\rho|_{\gamma_p}$. Now once $\rho$ is determined, there are at most $N^{1/3}$ possibilities for $p'\in(\cP\cap \gamma_{p'})\cap\rho^{-1}(\cP\cap \gamma_{q'})$ which uniquely determines $q'$. Thus in this case we get a bound of 
\[N^2\log N\cdot N\cdot N^{1/3}=N^{10/3}\log N.\qedhere\]
\end{proof}

To complete the proof, we simply need to check that all the above claims handle all the possible cases. Consider a tuple $(p,q,p',q',\rho,\fr_1,\fr_2)\in\tilde \cE$. By \cref{claim:both-reguli-rich}, we have counted this tuple unless one of $\fr_1,\fr_2$ does not lie in $\fR^{\mathrm{pt}}\cup \fR^{\mathrm{curve}}$. Without loss of generality, suppose that $\fr_2\not\in\fR^{\mathrm{pt}}\cup \fR^{\mathrm{curve}}$. Now by \cref{claim:poor-regulus-point-off-curve}, we have counted this tuple if $q$ or $q'$ lies in $\cP_{\fr_2\setminus\Gamma}$.

Thus we may assume that $q,q'\in \cP_{\fr_2\cap\Gamma}$. By \cref{claim:one-regulus-rich-pts-on-curves}, we have counted the tuple if $\fr_1\in \fR^{\mathrm{pt}}\cup \fR^{\mathrm{curve}}$. Next if $\fr_1\not\in\fR^{\mathrm{pt}}\cup\fR^{\mathrm{curve}}$, by another application of \cref{claim:poor-regulus-point-off-curve}, we have counted this tuple if $p$ or $p'$ lies in $\cP_{\fr_1\setminus\Gamma}$. Thus we may assume that $\fr_1,\fr_2\not\in\fR^{\mathrm{curve}}$ and that $p,p'\in\cP_{\fr_1\cap\Gamma}$ and $q,q'\in\cP_{\fr_2\cap\Gamma}$. Therefore the tuple is counted by \cref{claim:few-curves-on-reguli-all-pts-on-curves}.
\end{proof}

\subsection{Fixed points}
\label{ssec:rich-point}

Finally we will handle the case where we have a list $(p_1,q_1),(p_2,q_2),\ldots$ of $N^{2/3+o(1)}$ pairs of points and we wish to bound the contribution to the distance energy from rigid motions $\rho$ so that $\rho(p_i)=q_i$ for some $i$. The following result gives us the desired bound for each fixed $i$.

\begin{proposition}\label{lem:rich-points}
Let $\cP \subset \R^3$ be a set of size $N$ such that each plane and sphere contains at most $N^{2/3}$ points of $\cP$. 
Fix a point $x^*\in E^{\Phy}(\RR)$. 
Then 
\[\abs{\set{(p,q,p',q') \in \cP^4:\exists \rho \in E^+_g\textup{ s.t. }\phi(p,q), \phi(p',q'), x^* \in F_{\rho}^{\Phy}}} \lesssim N^{8/3} \log N.\]
\end{proposition}

\begin{proof}
We have two cases depending on whether or not $x^*$ is a point at infinity or not. (The two cases are essentially the same, but we have to take more care when $x^*$ is at infinity.)

First assume that $x^*\in E^{\Phy}_o(\RR)$. In this case, there exists $(p^*,q^*)\in \RR^3\times \RR^3$ so that $x^*=\phi(p^*,q^*)$. Note that the number of pairs $(p,q, p',q') \in \cP^4$ where $(p,q) = (p^*,q^*)$ or $(p',q') = (p^*,q^*)$ is at most $2N^2$, so we can restrict to the cases when $(p,q), (p',q') \neq (p^*,q^*)$. 

For $r>0$, let $\S_r(p^*)$ be the sphere of radius $r$ with center at $p^*$. Then define $\cP_r = \cP \cap \S_r(p^*)$ and $\cQ_r = \cP \cap \S_r(q^*)$. Note that for any quadruple $(p,q,p',q')$ and genuine positive rigid motion $\rho$ with $\phi(p,q),\phi(p',q'),x^*\in F^{\Phy}_{\rho}$, we see that $\rho$ sends $(p^*,p,p')$ to $(q^*,q,q')$. This immediately shows that $\norm{p-p^*} = \norm{q-q^*}$ and $\norm{p'-p^*} = \norm{q'-q^*}$ so there exist $r,r'$ such that $p\in\S_r(p^*)$ and $q\in\S_r(q^*)$ and $p'\in\S_{r'}(p^*)$ and $q'\in\S_{r'}(q^*)$. Furthermore we must have $r>0$ else $(p,q)=(p^*,q^*)$ and $r'>0$ else $(p',q')=(p^*,q^*)$.

Hence the number of quadruples we are counting can be upper bounded by
\[2N^2+ \sum_{r, r'>0} \abs{\set{(p,q,p',q')\in\cP_r\times\cQ_{r}\times\cP_{r'}\times\cQ_{r'}:\exists \rho\in E^+_g\text{ s.t. }\phi(p,q),\phi(p',q'),x^*\in F_\rho^{\Phy}}}.\]
    
Now we project the point sets onto the unit sphere. Write $\S^2$ for the origin-centered unit sphere and define $\overline\cP_r=\{u\in\S^2: p^*+ur\in\cP_r\}$ and $\overline\cQ_r=\{v\in\S^2: q^*+vr\in\cQ_r\}$. Write $T_{p^*}$ for the translation $x\mapsto x+p^*$ and $T_{q^*}$ for the translation $x\mapsto x+q^*$. Note that if $\rho$ is a genuine rigid motion of $\R^3$ which sends $p^*+ur$ to $q^*+vr$, then $T_{q^*}^{-1}\rho T_{p^*}$ is a rigid motion of $\S^2$ which sends $u$ to $v$. Thus we can bound each of the terms in the above expression as
\[E(\overline\cP_r,\overline\cQ_{r},\overline\cP_{r'},\overline\cQ_{r'})\lesssim(|\overline\cP_r||\overline\cQ_{r'}|+|\overline\cP_{r'}||\overline\cQ_r|)N^{2/3}\log N\]
by \cref{thm:bipartite-guth-katz}. Here we use the hypothesis that $|\cP_r|,|\cP_{r'}|,|\cQ_r|,|\cQ_{r'}|\leq N^{2/3}$.

Summing over $r,r'>0$ gives a total bound of
\[\lesssim 2N^2+\sum_{r,r'>0}\paren{|\overline\cP_r||\overline\cQ_{r'}|+|\overline\cP_{r'}||\overline\cQ_r|}N^{2/3}\log N\lesssim N^2+ N^{2/3}\log N\sum_r |\cP_r|\sum_{r'}|\cQ_{r'}|\lesssim N^{8/3}\log N.\]

Next assume that $x^*\not\in E^{\Phy}_o(\RR)$; this case can be dealt with similarly. First we need a geometric interpretation of what it means for $\phi(p,q),\phi(p',q')$ and $x^*$ to all lie in some $F^{\Phy}_{\rho}$ for some $\rho\in E^+_g$.

A point $x^*=[x_0:\vec x:x_{0\ast}:\vec x_\ast]$ lies in $E^{\Phy}(\RR)\setminus  E^{\Phy}_o(\RR)$ if $x_0x_{0\ast}-\vec x\cdot \vec x_\ast=0$ and $x_0=0$. Such points can be parametrized as $x^* = [0:\frac{a+b}{2}:c:\frac{a-b}{2}]$ for some $a,b\in \RR^3$ with $\norm{a}=\norm{b}$ and $c\in \RR$. Suppose that $\phi(p,q),x^*\in F^{\Phy}_{\rho}$ for some $\rho\in E^+_g$. This implies that the line connecting $\phi(p,q),x^*$ lies completely in $F_\rho^{\Phy}$ and thus lies completely in $E^{\Phy}$. By \cref{prop:intersection-from-same-family}, we have $\langle \phi(p,q),x^*\rangle = 0$. We can expand this equation (recall the definition of the inner product from \cref{defn:zorn-split-octonion} and that by \cref{eq:iota-phys-defn}, the $\phi(p,q)$ is embedded to the point $[1:\frac{p+q}2:\frac{\norm{p}^2-\norm{q}^2}4:\frac{p-q}2]$). This gives
\[c-\frac12 a\cdot p+\frac12 b\cdot q = 0.\]
Note that if $a=0$, then $\norm{a}=\norm{b}$ implies $b=0$ and the above implies $c=0$, which is impossible. Thus we may scale $a,b,c$ so that $\norm{a}=\norm{b}=1$.

Now the above equation implies that if $\phi(p,q),x^*\in F_\rho^{\Phy}$ for some $\rho$, then there exists $r$ so that $a\cdot p= 2c+r$ and $b\cdot q=r$. We foliate $\R^3$ into planes with normal $a$ and into planes with normal $b$, defining $\cP_r=\{p\in\cP:a\cdot p=2c+r\}$ and $\cQ_r=\{q\in\cP:b\cdot q = r\}$. As in the previous case, the number of quadruples we are counting can be upper bounded by
\[
\sum_{r, r'} \abs{\set{(p,q,p',q')\in\cP_r\times\cQ_r\times\cP_{r'}\times\cQ_{r'}:\norm{p-p'}=\norm{q-q'}}}.
\]
Let $F_a,F_b$ be the planes through the origin with normal vectors $a,b$. Similar to before, define $\overline\cP_r,\overline\cQ_r$ to be the orthogonal projections of $\cP_r,\cQ_r$ to $F_a,F_b$. Suppose $(p,q,p',q')\in \cP_r\times\cQ_r\times\cP_{r'}\times\cQ_{r'}$ satisfies $\norm{p-p'}=\norm{q-q'}$. As the images of these four points under the orthogonal projections are
\begin{align*}
p&\leadsto p-(2c+r)a\in\overline\cP_r,\\
q&\leadsto q-rb\in\overline\cQ_r,\\
p'&\leadsto p'-(2c+r')a\in\overline\cP_{r'},\\
q'&\leadsto q'-r'b\in\overline\cQ_{r'},
\end{align*}
the distances between the images are thus
\[\sqrt{\norm{p-p'}^2-{(r-r')^2}}=\sqrt{\norm{q-q'}^2-{(r-r')^2}}.\]
So the images form a distance quadruple. Picking isometries $F_a,F_b\cong\R^2$, we can again apply \cref{thm:bipartite-guth-katz}, to bound the desired quantity by
\[\sum_{r,r'}E(\overline\cP_r,\overline\cQ_r,\overline\cP_{r'},\overline\cQ_{r'})\lesssim\sum_{r,r'} \paren{|\cP_r||\cQ_{r'}|+|\cP_{r'}||\cQ_r|}N^{2/3}\log N\lesssim N^{8/3}\log N.\qedhere\]
\end{proof}
\section{Very rich partial symmetries}
\label{sec:main-argument}

In this section we prove \cref{thm:main}. First in \cref{ssec:input}, we record the results proved in \cref{part:i,part:ii} as \cref{thm:input}. Then in \cref{ssec:regularization}, we prove \cref{lem:finding-uniform-quadruple} which lets us pass to some $H\subseteq\cP\times\cP$ for which the problem is extremely regular. Finally in \cref{ssec:proof-of-main-thm}, we prove \cref{thm:main-thm-2-rich-lines} which implies \cref{thm:main}.

\subsection{Bounding the non-rich contribution}
\label{ssec:input}

We combine the results of \cref{part:i,part:ii} to deduce the following theorem, stating that the distance energy is dominated by the contribution of very rich partial symmetries.

\begin{theorem}\label{thm:input}
Let $\cP \subset \R^3$ be a set of size $N$ such that $|\cP \cap \gamma| \le N^{2/3} \deg \gamma$ for every curve $\gamma\subset\CC^3$, and such that every plane contains at most $N^{2/3}$ points of $\cP$. Fix $H \subseteq \cP\times \cP$ and $k \ge N^{2/3}$. Define $\cT \subset E_g^-$ to be the set of genuine orientation-reversing rigid motions $\tau$ such that $|F_\tau^{\Phy} \cap \phi(H)| \ge k$. 

For any $C\geq 2$, let $\cL$ be the set of isotropic lines $\ell^{\Phy}\subset E^{\Phy}_o$ for which 
\[\abs{\ell^{\Phy} \cap \phi(H)}\in[2,C]\]
and such that $\ell^{\Phy}$ is not contained in $F^{\Phy}_\tau$ for any $\tau\in\cT$. Then
\[\abs{\cL}\lesssim Ck N^{2/3}  |H| \log^2 N.\]
The analogous statement holds for orientation-preserving rigid motions.
\end{theorem}

\begin{proof}
Set $\cF = \{F_{pq}^+\}_{(p,q)\in H}$. By \cref{thm:analog-elekes-sharir}, all but at most $N\abs{H}\leq kN^{2/3}\abs{H}$ of the lines $\ell^{\Phy}\in\cL$ are such that $\ell^+\in\cL_2(\cF)$.

To bound $\abs{\cL_2(\cF)}$, we first apply \cref{thm:codim-0-incidence-bound} to produce sets of (complex) irreducible surfaces $\cS_1,\ldots, \cS_r$ contained in $\cF$ so that
\[\abs{\cL_2(\cF)\setminus\bigcup_i \cL_2(\cS_i;\cF)}\lesssim N^{4/3}\abs{H}\log^2N\leq kN^{2/3}|H|\log^2N.\]
We also have $\sum_i\deg(\cS_i)\lesssim N^{2/3}|H|\log N$ and that for each $i$ and $F\in\cF$, the set of surfaces $S\in\cS_i$ lying in $F$ have total degree $\lesssim N^{2/3}$.

We will show that
\begin{equation}
\label{eq:surface-bound}
\abs{\set{\ell^{\Phy}\in\cL:\ell^+\in \cL_2(\cS_i;\cF)}}\lesssim Ck\log N\deg \cS_i,
\end{equation}
which will then imply the desired bound.

To prove this, we first bound the contribution from each plane $S\in\cS_i$. For each such $S$, we claim that
\begin{equation}
\label{eq:plane-bound}
\abs{\set{\ell^{\Phy}\in\cL:\ell^+\in \cL_2(\cS_i;\cF)\text{ and }\ell^+\subset S}}\leq k.
\end{equation}
\cref{thm:codim-0-incidence-bound} guarantees that $S\subset F_{pq}^+$ for some (unique) $(p,q)\in H$. By \cref{prop:all-2-flats}, we have $S=F_{pq}^+\cap F_\tau^+$ for some $\tau\in E^-(\CC)$. First suppose that $\tau\in E^-(\CC)\setminus E^-(\RR)$. Now as every line $\ell^+$ with $\ell^{\Phy}\in\cL$ is defined over the reals, if $\ell^+\subset F_\tau^+$, then it also is contained in the complex conjugate $\ell^+\subset F_{\bar\tau}^+$. Since $\tau\not\in E^-(\RR)$, we have $\tau\neq\bar\tau$, so $\ell^+\subseteq F_\tau^+\cap F_{\bar\tau}^+$, where the intersection is a single line by \cref{prop:intersection-from-same-family}. Thus the right-hand side of \cref{eq:plane-bound} is bounded by 1 in the case that $\tau\not\in E^-(\R)$.

Now suppose $\tau\in E^-(\R)$ and $\ell^{\Phy}$ is an element of the above set. Since $\ell^+\subset S\subset F_\tau^+$ and $\ell^+\subset E^+_o$, by \cref{prop:actual-rigid-motion}, we know that $\tau$ is a genuine orientation-reversing rigid motion. Now $\ell^+\subset S'\subset F_{p'q'}^+$ for some $S'\in\cS_i$ and some $(p',q')\in H$ distinct from $(p,q)$. We then have $\ell^+\subset F_{\tau}^+\cap F_{p'q'}^+$. By \cref{prop:intersection-from-diff-family}, this implies that this intersection is a 2-flat (as it is not a single point), showing that $\phi(p',q')\in F_{\tau}^{\Phy}$ and so $\tau(p')=q'$. Since $\ell^{\Phy}\in\cL$, but $\ell^{\Phy}\subset F_{\tau}^{\Phy}$, we know that $\tau\not\in\cT$, i.e., this rigid motion is less than $k$-rich for $H$. 

Now $\ell^+$ is uniquely identified as $F_{pq}^+\cap F_{p'q'}^+$. Since there are fewer than $k$ choices of $p',q'$ with $\tau(p')=q'$, this proves \cref{eq:plane-bound} (as $(p,q)$ and $\tau$ are both uniquely determined by $S$).

Next we point out that it is easy to bound the contribution from $\ell^+\subset S\in\cS_i$ where $\ell^+$ is not in the ruling of $S$. (Recall from \cref{defn:in-the-ruling} that if $S$ is an unruled surface, for every $\ell^+\subset S$, we say that $\ell^+$ is not in the ruling of $S$.) We claim
\begin{equation}
\label{eq:unruled-bound}
\abs{\set{\ell^{\Phy}\in\cL:\ell^+\in \cL_2(\cS_i;\cF)\text{ and }\ell^+\subset S\text{ but $\ell^+$ not in the ruling of $S$}}}\lesssim (\deg S)^2.
\end{equation}
We bound the left-hand side simply by the number of lines contained in $S$ which are not in the ruling of $S$. By \cref{cor:lines-not-in-ruling}, this number is $\lesssim(\deg S)^2$, as desired.

Define $\cS'_i\subset\cS_i$ to be the set of singly or doubly ruled surfaces in $\cS_i$. Combining \cref{eq:plane-bound,eq:unruled-bound}, we have shown that
\begin{align*}
\abs{\set{\ell^{\Phy}\in\cL:\ell^+\in \cL_2(\cS_i;\cF)}}\lesssim k\deg \cS_i + \sum_{S\in\cS_i}(\deg S)^2+\abs{\cL_{i,\textrm{ruled}}}
\end{align*}
where $\cL_{i,\textrm{ruled}}\subseteq\cL$ is defined as follows. Say $\ell^{\Phy}\in\cL$ lies in $\cL_{i,\textrm{ruled}}$ if $\ell^+\in\cL_2(\cS_i;\cF)$ and there exist distinct $F_{pq}^+,F_{p',q'}^+\in\cF$ and surfaces $S,S'\in\cS'_i$ with $S\subset F_{pq}^+$ and $S'\subset F_{p'q'}^+$ so that $\ell^+\subseteq S\cap S'$ and $\ell^+$ lies in the rulings of $S$ and $S'$. We apply \cref{thm:ruled-ruled} to bound $\abs{\cL_{i,\textrm{ruled}}}$. Indeed, for $\ell^{\Phy}$ as above, we have that $\ell^{\Phy}$ connects $\phi(p,q),\phi(p',q')$. By the definition of $\cL$, this line is at most $C$-rich for $\phi(H)$. For each $(p,q)\in H$, define $\cS_{pq}$ to be the set of $S\in\cS'_i$ which lie in $F_{pq}^+$. Thus we can apply \cref{thm:ruled-ruled} to conclude that
\[\abs{\cL_{i,\textrm{ruled}}}\lesssim CN^{2/3}\log N\sum_{(p,q)\in H}\deg\cS_{pq}\leq CN^{2/3}\log N\deg\cS_i.\]

As \cref{thm:codim-0-incidence-bound} gives $\deg S\lesssim N^{2/3}$ for each $S\in\cS_i$, we have
\[\sum_{S\in\cS_i}(\deg S)^2\lesssim N^{2/3}\sum_{S\in\cS_i}\deg S= N^{2/3}\deg\cS_i.\]
Thus we have shown \cref{eq:surface-bound}, implying the desired result.
\end{proof}

\subsection{Constructing a balanced subgraph}
\label{ssec:regularization}

We start by introducing some notation.

\begin{definition}
For $\cQ \subset E^{\Phy}_o$ and $2\leq a\le b$, define $\cL_{[a,b]}(\cQ)$ to be the set of isotropic lines $\ell^{\Phy}\subset E^{\Phy}_o$ such that $\abs{\ell^{\Phy}\cap \cQ}\in [a,b]$.
We also write $\cL_{\geq a}(\cQ)$ for the set of isotropic lines $\ell^{\Phy}$ with $\abs{\ell^{\Phy}\cap \cQ}\geq a$.
\end{definition}

\begin{definition}
\label{defn:cl-x}
For a set of lines $\cL$ in $E^{\Phy}$ and sets $\cR\subset E^+, \cT\subseteq E^-$, we define for each $x\in E^{\Phy}$ the following sets
\begin{align*}
\cL(x) &= \{\ell^{\Phy}\in \cL: x\in \ell^{\Phy}\},\\
\cR(x) &= \{\rho\in \cR: x\in F^{\Phy}_{\rho}\},\\
\cT(x) &= \{\tau\in \cT: x\in F^{\Phy}_{\tau}\}.
\end{align*}
\end{definition}

By \cref{thm:beck-trick}, the distinct distances problem reduces to bounding $\cL_{[2,C_0]}(\phi(\cP\times \cP))$ for some constant $C_0$. By \cref{thm:input}, it remains to bound the lines that are captured both by a very rich positive rigid motion and a very rich negative rigid motion.

The main goal of this subsection is to ``regularize'' the setting. We pass to a subset $H\subseteq\cP\times\cP$, a subset $\cL\subseteq\cL_{[2,C_0]}(\phi(\cP\times\cP))$, and subsets $\cR,\cT$ of the very rich rigid motions so that each $x\in \phi(H)$ lies on the approximately same number of $\ell^{\Phy}\in\cL$, the same number of $F_\rho^{\Phy}$ for $\rho\in\cR$, and the same number of $F_{\tau}^{\Phy}$ for $\tau\in\cT$. Such regularity will be essential to the main argument.

\begin{lemma}\label{lem:finding-uniform-quadruple}
There exist $\epsilon \sim \sqrt{\log \log N/\log N}$ and $K\sim \log^{4/3} N$ so that the following holds.
Let $\cP \subset \R^3$ be a set of size $N$ such that $|\cP \cap \gamma| \le N^{2/3} \deg \gamma$ for every curve $\gamma$, and such that every plane contains at most $N^{2/3}$ points of $\cP$.
Suppose that $|\cL_{[2,K]}(\phi(\cP\times \cP))| \ge N^{10/3+2\varepsilon}$. 

Then there exist a subset $H \subseteq \cP \times \cP$, a subset $\cL \subseteq \cL_{[2, K]}(\phi(H))$, and sets $\cR \subset E^+_g$ and $\cT \subset E^-_g$ with the following properties:
\begin{enumerate}[(i)]
    \item $|\cL| \ge N^{-\varepsilon} |\cL_{[2, K]}(\phi(\cP \times \cP))|$;
    \item $|\cL(x)| \le N^{\varepsilon} |H|^{-1}|\cL|$ for every $x \in \phi(H)$;
    \item there exist $k, k' \ge N^{-2/3-\varepsilon}\abs{H}^{-1}\abs{\cL_{[2,K]}(\phi(\cP\times \cP))}$ so that $|F^{\Phy}_\rho\cap \phi(H)| \in [k,2k)$ and $|F^{\Phy}_\tau \cap \phi(H)| \in [k', 2k')$ for all $\rho \in \cR$ and $\tau \in \cT$;
    \item $|\cR(x)| \le N^{\varepsilon} k |\cR| |H|^{-1}$ and $|\cT(x)| \le N^{\varepsilon} k' |\cT| |H|^{-1}$ for every $x \in \phi(H)$;
    \item for every $\ell^{\Phy} \in \cL$, there exist a unique $\rho \in \cR$ and a unique $\tau \in \cT$ so that $\ell^{\Phy} \subset F^{\Phy}_\rho \cap F^{\Phy}_\tau$ (equivalently, $\rho \in \ell^+$ and $\tau \in \ell^-$); and
    \item  $N^{-\varepsilon}k^2 |\cR| \le |\cL| \le 4k^2 |\cR|$ and $N^{-\varepsilon}k'^2 |\cT| \le |\cL| \le 4k'^2 |\cT|$.
\end{enumerate}
\end{lemma}

In other words, by passing to a subset $H$ of $\cP\times \cP$, we do not lose too many lines, while the number of lines, positive $3$-flats and negative $3$-flats through each point in $\phi(H)$ is not too much larger than the average. Each $\rho\in\cR$ is approximately the same richness for $\phi(H)$ as is each $\tau\in\cT$. In addition, we will show that each line is captured by a unique $\rho$ and a unique $\tau$. Finally each $\rho,\tau$ captures close to the maximum possible number of lines. This last property is the most important for the main argument.

We begin with a fairly standard graph-theoretic result. 

\begin{definition}
For $B\geq 1$, a bipartite graph $G=(V_1\sqcup V_2,E)$ is \emph{$B$-biregular} if there exist $d_1, d_2 \ge 0$ so that $\deg_{G} v_1 \in [d_1, Bd_1]$ for every $v_1\in V_1$ and $\deg_G v_2 \in [d_2, Bd_2]$ for every $v_2\in V_2$. We call $d_1,d_2$ the \emph{bidegrees} of $G$.
\end{definition}

The following result finds a biregular subgraph, only losing a polylogarithmic number of edges.

\begin{lemma}
\label{lem:biregular}
There exists $B\sim\log^2N$ so that the following holds. Let $G = (V_1\sqcup V_2, E)$ be a bipartite graph with $|V_1|+|V_2|\le N$. Then there exist subsets $\tilde V_1 \subseteq V_1$ and $\tilde V_2 \subseteq V_2$ so that the induced subgraph $\tilde G = (\tilde V_1\sqcup \tilde V_2, \tilde E)$ is $B$-biregular and satisfies $|\tilde E| \ge B^{-1} |E|$. Furthermore, for each $v\in \tilde V_1$, we have $\deg_{\tilde G}v\geq B^{-1}\deg_G v$.
\end{lemma}

The proof proceeds in two steps: first, we take a dyadic decomposition to find a subgraph with all degrees at most a small multiple of the average degree, second, we prune out low degree vertices to obtain a $B$-biregular subgraph.

\begin{proof}
Set $B=8(1+\log_2N)^2$. Let $V_{1j} \subseteq V_1$ be the set of vertices $v_1$ such that $\deg_{G}v_1 \in [2^j, 2^{j+1})$ for $j=0, 1,\ldots, \lfloor\log_2 N\rfloor$. By the pigeonhole principle, there exists $j$ so that $|E \cap (V_{1j}\times V_2)| \ge \frac{|E|}{1+\log_2 N}$. Fix this $j$ and let $G'$ be the induced subgraph on $V_{1j}\sqcup V_2$.

Now let $V_{2i}$ be the set of $v_2 \in V_2$ such that $\deg_{G'}v_2 \in [2^i, 2^{i+1})$ for $i=0,1 \ldots, \lfloor\log_2 N\rfloor$. By the pigeonhole principle, there exists $i$ so that $|E \cap (V_{1j}\times V_{2i})| \ge \frac{|E|}{(1+\log_2 N)^2}$. Let $G''$ be the induced subgraph on $V_{1j}\sqcup V_{2i}$. By construction, we have the following estimates on the vertex degrees
\begin{equation}
\label{eq:v1-deg-upper}
    \deg_{G''}v_1 \le \deg_{G} v_1 \le 2^{j+1} \le 2\frac{|E|}{|V_{1j}|} \le 2 (1+\log_2N)^2 \frac{|E''|}{|V_{1j}|}
\end{equation}
for all $v_1\in V_{1j}$ and 
\begin{equation}
\label{eq:v2-deg-upper}
\deg_{G''}v_2 \le \deg_{G'} v_2 \le 2^{i+1} \le 2\frac{|E'|}{|V_{2i}|} \le 2 (1+\log_2N) \frac{|E''|}{|V_{2i}|}
\end{equation}    
for all $v_2\in V_{2i}$.

Now we let $d_1 = \tfrac{1}{4} |E''|/|V_{1j}|$ and $d_2 = \tfrac{1}{4} |E''|/|V_{2i}|$. Iteratively remove from $G''$ every vertex $v_1 \in V_{1j}$ with degree less than $d_1$ and every vertex $v_2 \in V_{2i}$ with degree less than $d_2$ (after removing a vertex, we compute all degrees in the new graph, but keep the parameters $d_1,d_2$ fixed). Let $\tilde G = (\tilde V_1\sqcup\tilde V_2, \tilde E)$ be the induced subgraph of $G''$ produced by these deletions. Note that
    \[
    |\tilde E| \ge |E''| - d_1 |V_{1j}| - d_2|V_{2i}| \ge |E''|/2.
    \]
Thus we have $|\tilde E|\geq |E''|/2\geq B^{-1}|E|$.

Furthermore, every vertex $v_t \in \tilde V_t$ now has degree at least $d_t$ for $t=1, 2$. It follows from \cref{eq:v1-deg-upper,eq:v2-deg-upper} $\tilde G$ is $8 (1+\log_2N)^2$-biregular. In fact, for each $v\in\tilde V_1$, \cref{eq:v1-deg-upper} implies the stronger bound
\[d_1\leq \deg_{\tilde G}v\leq \deg_{G''}v\leq\deg_G v\leq Bd_1,\]
proving the last of the desired statements.
\end{proof}

The next lemma roughly says that if there are many rigid motions which are very rich for $\phi(H)$, then there are many lines which are $2$-rich for $\phi(H)$. Indeed, if $\rho\in E_g^+$ is $k$-rich, meaning $|\phi(H)\cap F^{\Phy}_{\rho}|\geq k$, then each line in $F^{\Phy}_\rho$ connecting two points of $\phi(H)$ is such a 2-rich line for $\phi(H)$. By Beck's theorem, there are $\gtrsim k^2$ of these lines for each $\rho$. The following lemma roughly states that these $|\cR|$ set of $\gtrsim k^2$ lines are mostly disjoint. This will be an important step in establishing \cref{lem:finding-uniform-quadruple}(vi).

\begin{lemma}
\label{lem:many-lines-from-rich-motions}
For all $N$ sufficiently large and $k$ with $k\geq N^{2/3}\log^9N$, there exists $K\sim\log^{4/3}N$ so that the following holds. Let $\cP \subset \R^3$ be a set of size $N$ such that $|\cP\cap\gamma|\leq N^{2/3}\deg\gamma$ for every curve $\gamma$. Let $H \subset \cP\times \cP$ be a subset and let $\cR\subset E_g^+$ be a set of rigid motions so that $|\phi(H) \cap F^{\Phy}_{\rho}|\geq k$ for every $\rho\in\cR$. 
    Then we have
    \[
    |\cL_{[2, K]}(\phi(H))| \gtrsim \log^{-8} N |\cR| k^2.
    \]
\end{lemma}
\begin{proof}
Consider the bipartite graph $\cG = (\cR\sqcup \phi(H), E)$ where $\rho\in \cR$ and $x\in \phi(H)$ are adjacent if $x \in F^{\Phy}_\rho$. By hypothesis, $|E| \ge k |\cR|$.

Applying \cref{lem:biregular}, we produce subsets $\cR' \subseteq \cR$ and $H' \subseteq H$ so that the induced subgraph $\cG'=(\cR'\sqcup \phi(H'), E')$ is $B$-biregular and $|E'| \ge B^{-1}|E|$, where $B \sim \log^2 N$. Furthermore, we have $\deg_{\cG'} \rho\geq B^{-1}\deg_{\cG}\rho$ for all $\rho\in\cR'$.

In particular, we have $|\phi(H')\cap F^{\Phy}_\rho|=\deg_{\cG'}\rho\geq B^{-1}\deg_{\cG}\rho=B^{-1}|\phi(H)\cap F^{\Phy}_\rho|\geq B^{-1}k$ for every $\rho\in\cR'$. Since $k\geq N^{2/3}\log^9N$, we have $|\phi(H')\cap F^{\Phy}_{\rho}|\geq 2N^{2/3}$. By hypothesis, no line contains $N^{2/3}$ points of $\cP$, implying that no line contains $N^{2/3}$ points of $\phi(H')$.

We are thus in a position to apply \cref{thm:beck-variant} to the sets $\phi(H')\cap F^{\Phy}_{\rho}$ and $\phi(H)\cap F^{\Phy}_\rho$ with parameters $|\phi(H')\cap F^{\Phy}_{\rho}|$ and $B$. Let $C,c$ be the constants in that result. Then defining $\cL_\rho$ to be the set of lines contained in $F^{\Phy}_{\rho}$ which contain at least 2 points of $\phi(H')$ and at most $CB^{2/3}$ points of $\phi(H)$, \cref{thm:beck-variant} implies that $|\cL_\rho|\geq c|\phi(H')\cap F^{\Phy}_{\rho}|^2=c(\deg_{\cG'}\rho)^2$.

Let $d_\cR$ and $d_H$ be the minimum degrees of $\cG'$ from $\cR'$ and $\phi(H')$, respectively. We have 
    \begin{equation}\label{eq:drho}
    d_\cR =\min_{\rho\in\cR'} \deg_{\cG'}\rho\geq \min_{\rho\in\cR'}B^{-1}\deg_{\cG}\rho \geq B^{-1}k.
    \end{equation}
In addition, since the maximum degree of $\cG'$ from $\cR'$ is at most $Bd_\cR$, we have
    \begin{equation}\label{eq:dx}
    d_H \le \frac{|E'|}{|H'|} \le B \frac{|\cR'| d_\cR}{|H'|}.    
    \end{equation}

For $x \in \phi(H')$, we consider $\cL_\rho(x)\subseteq\cL_{\rho}$ (defined in \cref{defn:cl-x} to be the subset of lines which contain $x$). Recall \cref{prop:intersection-from-same-family} implies that for $\rho\neq\rho' \in E^+$ the intersection $F^{\Phy}_\rho \cap F^{\Phy}_{\rho'}$ contains at most a single line. Thus for any $x$ and $\rho$, the number of lines $\ell \in \cL_\rho(x)$ which are contained in some $F^{\Phy}_{\rho'}$ for some $\rho' \in \cR'\setminus\{\rho\}$ is upper bounded by $\deg_{\cG'}x\le B d_H$. Since $\cL_\rho(x) \subseteq \cL_{[2,CB^{2/3}]}(\phi(H))(x)$, we conclude that 
\[
    |\cL_{[2, CB^{2/3}]}(\phi(H))(x)| \ge \sum_{\rho:\rho \sim_{\cG'} x} (|\cL_\rho(x)| - \deg_{\cG'}x) \ge \paren{\sum_{\rho:\rho \sim_{\cG'} x} |\cL_\rho(x)|} - (B d_H)^2.
\]
Taking the sum over $x \in \phi(H')$, since each element of $\cL_{[2, CB^{2/3}]}(\phi(H))$ is included in at most $CB^{2/3}$ of the sets $\cL_{[2, CB^{2/3}]}(\phi(H))(x)$, we conclude that
\begin{equation}
\label{eq:LH}
    \begin{split}
    |\cL_{[2, CB^{2/3}]}(\phi(H))| 
    &\ge \frac{1}{CB^{2/3}} \left( \sum_{x, \rho}|\cL_\rho(x)| -(Bd_H)^2 |H'|\right) \\
    &\ge \frac{1}{CB^{2/3}}\left(\sum_{\rho \in \cR'} |\cL_\rho|  - (Bd_H)^2 |H'|  \right)\\
    &\ge \frac{1}{CB^{2/3}}\left( |\cR'| c d_\cR^2 - (Bd_H)^2 |H'|\right)
    \end{split}
\end{equation}
    
By applying \cref{thm:very-rich-bound-main} to $H'$ with parameter $d_\cR$, we see that
    \[
    |\cR'| \lesssim \frac{N|H'|\log^{3/2} N}{d_{\cR}^{3/2}}.
    \]
By \cref{eq:dx},
    \[
    (Bd_H)^2 |H'| \lesssim B^4 \frac{|\cR'|^2 d_\cR^2}{|H'|} \lesssim  (N d_\cR^{-3/2} \log^{10} N) |\cR'| d_\cR^2
    \]
By \cref{eq:drho}, we have $d_\cR \ge B^{-1}k \gtrsim N^{2/3}\log^7 N$, so we conclude that the second term in \cref{eq:LH} is negligible.

Thus \cref{eq:LH}, together with the bounds $|E'|\geq B^{-1}|E|$ and $d_\cR\geq B^{-1}k$, gives
\[
    |\cL_{[2, CB^{2/3}]}(\phi(H))| \gtrsim |\cR'| d_\cR^2 B^{-2/3}\ge B^{-5/3} |E'| d_\cR \ge B^{-11/3}|E| k \ge B^{-11/3}|\cR| k^2.
\]
Setting $K=CB^{2/3}\sim\log^{4/3}N$ completes the proof.
\end{proof}

Next, we work towards establishing \cref{lem:finding-uniform-quadruple}(ii). To do this, we wish to make sure that each point $x\in \phi(H)$ lies in roughly the same number of the lines under consideration.

\begin{lemma}\label{lem:finding-balanced-subgraph}
There exist $K \sim \log^{4/3}N$ and $\eta = \exp(-O(\sqrt{\log N\log\log N}))$ so that the following holds. Let $\cP \subset \R^3$ be a set of size $N$ such that $|\cP \cap \gamma| \le N^{2/3} \deg \gamma$ for every curve $\gamma$, and each plane contains at most $N^{2/3}$ points of $\cP$. Let $H \subseteq \cP \times \cP$ be such that $|\cL_{[2, K]}(\phi(H))| \ge \eta^{-2}N^{10/3}$. Then there exist subsets $H' \subseteq H$ and $\cL' \subseteq \cL_{[2, K]}(\phi(H'))$ so that
    \[
    |\mc L'| \gtrsim \eta|\mc L_{[2,K]}(\phi(H))| \qquad\text{ and }\qquad |\cL'|\geq\log^{-5}N\abs{\cL_{[2,K]}(\phi(H'))}
    \]
    and for every $x \in \phi(H')$ we have
    \[
    |\mc L'(x)| \lesssim  \eta^{-1} |H'|^{-1} |\mc L'|.
    \]
\end{lemma}

This is the most involved step in the proof of \cref{lem:finding-uniform-quadruple}. For its proof we use a lemma we call the  ``edge-increment or vertex-decrement dichotomy.'' To see why we need such a result, consider the graph on $\phi(H)$ where two vertices are adjacent if they span a line in $\cL_{[2,K]}(\phi(H))$. Our goal is to regularize this graph.

First, losing at most half the edges, we can pass to a bipartite subgraph. Then, \cref{lem:biregular} lets us biregularize the graph. Unfortunately, it is possible that the two parts have very different degrees. The following lemma says that in this case, either passing to the (non-bipartite) subgraph on one of the two parts has a large increment in the number of edges, or passing to the smaller subgraph decreases the number of vertices by a large amount while not losing too many edges. In either case, the density of the graph increases.

To be more precise, when we pass from vertex set $\phi(H)$ to a subset $\phi(H')$, we consider the new graph on $\phi(H')$ where two vertices are adjacent if they span a line in $\cL_{[2,K]}(\phi(H'))$. This graph may have more edges than simply taking the induced subgraph on $\phi(H')$, as a line that was more than $K$-rich for $\phi(H)$ may now be $[2,K]$-rich for $\phi(H')$. This complication is the main reason we need the edge-increment case as well as the vertex-decrement case in the following lemma. 

\begin{lemma}[edge-increment or vertex-decrement dichotomy]
\label{lem:line-boost-size-dec-dichotomy} 
There exists a constant $C>0$ and $K\sim\log^{4/3}N$ so that the following holds. Let $\cP\subset \RR^3$ be a set of size $N$ such that $\abs{P\cap \gamma}\leq N^{2/3}\deg \gamma$ for every curve $\gamma$, and each plane contains at most $N^{2/3}$ points of $\cP$.

Given parameters $A,B\geq 1$, suppose that $H=H_1\sqcup H_2\subseteq \cP\times \cP$ is a union of two disjoint subsets of $\cP\times \cP$, and let $E\subseteq \phi(H_1)\times \phi(H_2)$ be the set of $(x_1,x_2)$ where $\aff\{x_1,x_2\}\in \cL_{[2,K]}(\phi(H))$. If $(\phi(H_1)\sqcup \phi(H_2),E)$ is $B$-biregular with bidegrees $d_1,d_2$ and $\abs{E}\geq CN^{4/3}\log^{12} N\abs{H}$, then for either $H'=H_1$ or $H'=H_2$, one of the following holds:
    \begin{enumerate}
        \item $\abs{\cL_{[2,K]}(\phi(H'))}\geq A\abs{E}$; or
        \item $\abs{\cL_{[2,K]}(\phi(H'))}\gtrsim A^{-1}\log^{-20} N\abs{E}$ and $\abs{H'}\leq B\abs{H}\min\set{\frac{d_1}{d_2},\frac{d_2}{d_1}}$.
    \end{enumerate}
\end{lemma}

\begin{proof}
Pick some sufficiently large constants $C_1,C_2$ and set $C=C_1C_2$. Take $K\sim\log^{4/3}N$ to be the parameter in \cref{lem:many-lines-from-rich-motions}. Define $k = C_1 N^{2/3}\log N$ and $X$ so that $\abs{E} =X\cdot kN^{2/3}\abs{H}$. By our choice of $C$, we have $X\geq C_2\log^{11} N$.

Let $\cR\subset E^+_g$ be the set of genuine positive rigid motions $\rho$ such that $|F^{\Phy}_{\rho}\cap \phi(H)|\geq k$. By \cref{thm:input}, there are $\lesssim kN^{2/3}|H|\log^2N$ lines in $\cL_{[2,K]}(\phi(H))$ which are not contained in some $F_{\rho}^{\Phy}$ with $\rho\in\cR$. Thus $\lesssim K^2kN^{2/3}|H|\log^2N$ edges in $E$ correspond to such a line; in particular, this number is at most $|E|/2$.

For each $i,j$, let $\cR_{i,j}$ be the set of rigid motions $\rho\in \cR$ with $|\phi(H_1)\cap F_{\rho}^{\Phy}|\in [2^i,2^{i+1})$ and $|\phi(H_2)\cap F_{\rho}^{\Phy}|\in [2^j,2^{j+1})$. Since at least half the edges of $E$ are captured by some $\rho\in\cR$, by the pigeonhole principle, we may find $i,j$ so that
\begin{equation}
\label{eq:rij-dyadic}
\abs{\set{(x_1,x_2)\in E: \exists \rho\in \cR_{i,j}\textup{ with }\aff\{x_1,x_2\}\subset F_{\rho}^{\Phy}}}\gtrsim \log^{-2}N\abs{E}.
\end{equation}

Fix such $i,j$ and set $k_1=2^i$ and $k_2=2^j$. Since every $\rho\in\cR$ is $k$-rich for $\phi(H_1)\sqcup\phi(H_2)$, we must have $2k_1+2k_2\geq k$. Without loss of generality, assume that $k_1\geq k_2$. This implies that $k_1\geq k/4=\tfrac14C_1N^{2/3}\log N$. Our next goal is to show that $k_1\geq k_2\gtrsim_{C_1} XN^{2/3}\log^{-2} N$. 

Every element of $\cR_{i,j}$ is a rigid motion that is $k_1$-rich for $H_1$. If $C_1$ is sufficiently large, \cref{thm:very-rich-bound-main} bounds this quantity by
\[|\cR_{i,j}|\lesssim \frac{N|H_1|\log^{3/2}N}{k_1^{3/2}}\lesssim_{C_1} \frac{N^{2/3}|H|\log N}{k_1}.\]
Now each pair on the left-hand side of \cref{eq:rij-dyadic} is captured by some $\rho\in \cR_{i,j}$, while each $\rho$ captures at most $2^{i+1}\cdot 2^{j+1}=4k_1k_2$ pairs. Thus \cref{eq:rij-dyadic} and the above inequality together imply that
\[\abs{E}\lesssim k_1k_2\log^2 N\abs{\cR_{i,j}}\lesssim_{C_1} k_2N^{2/3}\abs{H}\log^3 N,\]
showing that $k_2\gtrsim_{C_1} Xk\log^{-3}N\sim_{C_1} XN^{2/3}\log^{-2} N$.

We have $k_1\geq k_2\gtrsim_{C_1} XN^{2/3}\log^{-2}N\geq C_2 N^{2/3}\log^9 N$. Taking $C_2$ large in terms of $C_1$, we can ensure $k_1, k_2 \geq N^{2/3}\log^9 N$. Then applying \cref{lem:many-lines-from-rich-motions} to $H_t,\cR_{i,j},k_t$ for each $t=1,2$, we have
\begin{equation}
\label{eq:many-lines-in-ht}
\abs{\cL_{[2,K]}(\phi(H_t))}\gtrsim \log^{-8} N\abs{\cR_{i,j}}k_t^2\gtrsim \log^{-10}N\abs{E}\frac{k_t}{k_{3-t}},
\end{equation}
where the last inequality follows from $|E|\lesssim k_1k_2\log^2N|\cR_{i,j}|$.

To complete the proof, first suppose $k_1/k_2\geq c^{-1}A\log^{10}N$ where $c$ is the implicit constant in \cref{eq:many-lines-in-ht}. In this case, \cref{eq:many-lines-in-ht} with $t=1$ shows that we fall into case (1) with $H'=H_1$. Otherwise suppose that $k_1/k_2\lesssim A\log^{10}N$. In this case, \cref{eq:many-lines-in-ht} implies that
\[\abs{\cL_{[2,K]}(\phi(H_t))}\gtrsim A^{-1}\log^{-20}N\abs{E},\]
for both $t=1,2$.

Finally, using that $(\phi(H_1)\sqcup\phi(H_2),E)$ is $B$-biregular with bidegrees $d_1,d_2$, we conclude that
\[\abs{H_t}\leq \frac{\abs{E}}{d_t}\leq \frac{B|H_{3-t}|d_{3-t}}{d_t}\leq B|H|\frac{d_{3-t}}{d_t}\]
for $t=1,2$. This shows that we fall into case (2) with one of $H'=H_1$ or $H'=H_2$.
\end{proof}

Now we can use this to prove \cref{lem:finding-balanced-subgraph}.

\begin{proof}[Proof of \cref{lem:finding-balanced-subgraph}]
We will iteratively prune $H$, producing sets $H=H^{(0)}\supset H^{(1)}\supset\cdots\supset H^{(t_f)}\supset H'$.

Let $K\sim \log^{4/3}N$ be the parameter in \cref{lem:many-lines-from-rich-motions}. Suppose that we have produced $H^{(t)}$ which satisfies $\abs{\cL_{[2,K]}(H^{(t)})}\geq N^{10/3}\log^{15}N$. Define $H_*\subseteq H^{(t)}$ to be a minimal subset so that $\abs{\cL_{[2,K]}(\phi(H_*))}\geq \abs{\cL_{[2,K]}(\phi(H^{(t)}))}$.

Partition $H_*=H_1\sqcup H_2$ so that at least half of the lines $\ell^{\Phy}\in\cL_{[2,K]}(\phi(H_*))$ intersect both $\phi(H_1)$ and $\phi(H_2)$. (A uniform random partition works with positive probability.) Let $G=(\phi(H_1)\sqcup \phi(H_2),E)$ be the bipartite graph on $\phi(H_1)\sqcup \phi(H_2)$ where $x_1\in \phi(H_1)$ is adjacent to $x_2\in \phi(H_2)$ if $\aff\{x_1,x_2\}\in \cL_{[2,K]}(\phi(H_*))$. By \cref{lem:biregular}, we can find subsets $H_1'\subseteq H_1$ and $H_2'\subseteq H_2$ so that $G[H_1'\times H_2']$ is $B$-biregular with $B\sim \log^2N$. Let $d_1,d_2$ be the corresponding bidegrees. Without loss of generality, assume that $d_1\geq d_2$.

We can bound the number of edges of $G[H_1'\times H_2']$ as
\begin{align*}
\abs{(\phi(H_1')\times \phi(H_2'))\cap E}
&\geq B^{-1}|E|\geq\frac 12 B^{-1}\abs{\cL_{[2,K]}(\phi(H_*))}\geq\frac 12 B^{-1}\abs{\cL_{[2,K]}(\phi(H^{(t)}))}.
\end{align*}
Indeed the first inequality is given by \cref{lem:biregular}, the second is our choice of partition, and the third is the definition of $H_*$.

Now define $T=\exp(\sqrt{\log N\log \log N})$. If $d_1/d_2>T$, we apply \cref{lem:line-boost-size-dec-dichotomy} to $G[H_1'\times H_2']$ with parameters $A=4B$ and $B$. Note that our assumption $\abs{\cL_{[2,K]}(H^{(t)})}\geq N^{10/3}\log^{15}N$ implies that $G[H_1'\times H_2']$ has enough edges that this result applies. Let $H^{(t+1)}$ be the set produced (either $H_1'$ or $H_2'$). First we claim by our choice of $H_*$, we cannot be in case (1), the edge-increment case. Indeed, if this occurs, we would have
\[\abs{\cL_{[2,K]}(\phi(H^{(t+1)}))}\geq 4B\abs{(\phi(H_1')\times \phi(H_2'))\cap E}\geq 4\abs{E}\geq 2\abs{\cL_{[2,K]}(\phi(H_*))},\]
contradicting the choice of $H_*$.

Thus we must be in case (2), the vertex-decrement case. This gives
\[\abs{\cL_{[2,K]}(\phi(H^{(t+1)}))}\gtrsim \log^{-22} N \abs{(\phi(H_1')\times \phi(H_2'))\cap E}\gtrsim \log^{-24} N\abs{\cL_{[2,K]}(\phi(H^{(t)}))}\]
and
\[\abs{H^{(t+1)}}\leq \frac{B}{T}\abs{H_*}\leq \frac{B}{T}\abs{H^{(t)}}.\]
If $\abs{\cL_{[2,K]}(H^{(t+1)})}< N^{10/3}\log^{15}N$, we say that the iteration fails and we halt; otherwise we continue the iteration with $H^{(t+1)}$.

If $d_1/d_2\leq T$, we say that the iteration succeeds. Define \[H' = H_1'\sqcup H_2'\qquad\text{ and }\qquad \cL' = \{\aff\{x_1,x_2\}:(x_1,x_2)\in (\phi(H_1')\times\phi(H_2'))\cap E\}.\]
As each edge of $G[\phi(H_1')\times\phi(H_2')]$ produces a line in $\cL'$ and each line corresponds to at most $K^2$ edges, we conclude
\begin{align*}
|\cL'|
&\geq K^{-2}\abs{(\phi(H_1')\times \phi(H_2'))\cap E} \geq \frac12 B^{-1}K^{-2}|E|\\
&\gtrsim \log^{-5}N\abs{\cL_{[2,K]}(\phi(H_*))}\geq \log^{-5} N\abs{\cL_{[2,K]}(\phi(H^{(t)}))}.
\end{align*}
As $H'\subseteq H_*$, by the minimality of $H_*$, the above also implies $|\cL'|\gtrsim\log^{-5}N\abs{\cL_{[2,K]}(\phi(H'))}$.
Furthermore, for each $i=1,2$ and $x_i\in\phi(H_i')$, we have
\begin{align*}
|\cL'(x_i)|\leq & \deg_{G[\phi(H_1')\times \phi(H_2')]}(x_i)\leq Bd_i\leq BT \min\{d_1,d_2\}\\
\leq & BT\frac{2\abs{(\phi(H_1')\times \phi(H_2'))\cap E}}{\abs{H'}}\lesssim T\log^{5} N\frac{|\cL'|}{|H'|}.
\end{align*}

Thus we can define $\eta=\exp(-O(\sqrt{\log N\log\log N}))$ so that $|\cL'(x)|\lesssim \eta^{-1}|H'|^{-1}|\cL'|$ holds for all $x\in\phi(H')$.

Our next goal is to show that if the iteration starts with $\abs{\cL_{[2,K]}(\phi(H))}\geq \eta^{-2}N^{10/3}$, then it always succeeds with $\cL_{[2,K]}(\phi(H^{(t)}))$ and thus $\cL'$ both large. Suppose that the iteration terminates at stage $t_f$.
As $\abs{H^{(t+1)}}\leq BT^{-1}\abs{H^{(t)}}$, we have
\[1\leq\abs{H^{(t_f)}}\leq \left(\frac{B}{T}\right)^{t_f}\abs{H^{(0)}}\leq\paren{\frac BT}^{t_f} N^2,\]
implying that $t_f\leq 2\log N/(\log T-\log B)\lesssim \sqrt{\log N/\log\log N}$. Hence
\[\frac{\abs{\cL_{[2,K]}(\phi(H^{(t_f)}))}}{\abs{\cL_{[2,K]}(\phi(H^{(0)}))}}\geq (C\log^{24} N)^{-t_f}\geq \exp\paren{-O\paren{\sqrt{\log N\log\log N}}},\]
for some absolute constant $C$.

Choosing the constant in the definition of $\eta$ sufficiently large, we conclude that $\abs{\cL_{[2,K]}(\phi(H^{(t_f)}))}\geq \eta \abs{\cL_{[2,K]}(\phi(H))}$. In particular, this implies that $\abs{\cL_{[2,K]}(\phi(H^{(t)}))}\geq \eta^{-1}N^{10/3}\geq N^{10/3}\log^{15}N$ for every $t\leq t_f$. Therefore the iteration succeeds and the sets $H',\cL'$ we defined have the desired properties.
\end{proof}

We are now ready to prove \cref{lem:finding-uniform-quadruple}, the main result of this subsection.
Roughly speaking, we will apply \cref{lem:many-lines-from-rich-motions,lem:finding-balanced-subgraph} to establish \cref{lem:finding-uniform-quadruple}(i)(ii)(vi).
However, we also need to do some dyadic pigeonholing and pruning to regularize the situation for the genuine rigid motions and establish (iii)(iv).
We also need to show (v), namely that each line is uniquely captured by a genuine positive rigid motion and a genuine negative rigid motion.

\begin{proof}[Proof of \cref{lem:finding-uniform-quadruple}]
Let $K\sim \log^{4/3}N$ be the parameter produced by \cref{lem:finding-balanced-subgraph}. Define $\varepsilon \sim \sqrt{\log \log N/\log N}$ so that the parameter $\eta$ in \cref{lem:finding-balanced-subgraph} satisfies $\eta^{-1}\leq N^{\epsilon/2}$.

\vspace{3pt}\noindent\textbf{Initial pruning of $H$:} We start by applying \cref{lem:finding-balanced-subgraph} to $\cP\times\cP$. Note that$\abs{\cL_{[2,K]}(\phi(\cP\times\cP))}\geq N^{10/3+2\epsilon}\geq\eta^{-4}N^{10/3}$, so this result applies. This produces $H_0\subseteq \cP\times\cP$ and $\cL_0\subseteq \cL_{[2, K]}(\phi(H_0))$ so that
\begin{equation}
\label{eq:l0-lower-absolute}
|\cL_0|\gtrsim \eta |\cL_{[2, K]}(\phi(\cP\times \cP))|\geq N^{-\epsilon/2}|\cL_{[2, K]}(\phi(\cP\times \cP))|\geq N^{10/3+3\epsilon/2},
\end{equation}
\begin{equation}
\label{eq:l0-lower-relative}
|\cL_0|\gtrsim \log^{-5}N\abs{\cL_{[2, K]}(\phi(H_0))},
\end{equation}
and
\begin{equation}
\label{eq:l0-degree-upper}
|\cL_0(x)| \lesssim \eta^{-1}|H_0|^{-1}|\cL_0|\leq N^{\epsilon/2} |H_0|^{-1} |\cL_0|.
\end{equation}

Define $k_0= N^{2/3}\log^9N$. Let $\cR_0\subset E^+_g$ and $\cT_0\subset E^-_g$ be the sets of genuine positive and negative rigid motions, respectively, which are $k_0$-rich for $\phi(H_0)$. By \cref{thm:input}, the number of lines in $\cL_0\subseteq\cL_{[2,K]}(\phi(H_0))$ which are not contained in $F_{\tau}^{\Phy}$ for some $\tau\in\cT_0$ is bounded by
\[\lesssim Kk_0N^{2/3}|H_0|\log^2N\lesssim N^{4/3+\epsilon}|H_0|\leq N^{10/3+\epsilon}.\]
Symmetrically, the same bound holds for the number of lines in $\cL_0$ which are not contained in $F_\rho^{\Phy}$ for some $\rho\in\cR_0$.
This is a negligible fraction of $\cL_0$, so we see that almost every line in $\cL_0$ lies simultaneously in some $F_\rho^{\Phy}$ and some $F_\tau^{\Phy}$. Define $\cL_1\subseteq\cL_0$ to be this set of lines. Combining the above bound with \cref{eq:l0-lower-absolute}, we see that $|\cL_1|\sim|\cL_0|$.

\vspace{3pt}\noindent\textbf{Initial dyadic decomposition of $\cR,\cT$:} For each $j$ with $2^{j+1}>k_0$, define $\cR_j \subseteq \cR_0$ and $\cT_j\subseteq \cT_0$ to be the subsets consisting of $\rho\in\cR_0$  for which $|F^{\Phy}_\rho \cap \phi(H_0)| \in [2^j, 2^{j+1})$ and analogously for $\cT_j$.

\vspace{3pt}\noindent\textbf{Final pruning of $H$:} To define the final set $H$, we need to remove any $x\in H_0$ which has too high degree to any of the sets $\cR_j,\cT_j$.

Define $H\subseteq H_0$ by removing any $x\in H_0$ where $|\cR_j(x)| > N^{3\epsilon/4}\frac{2^j |\cR_j|}{|H_0|}$ or $|\cT_j(x)| > N^{3\epsilon/4} \frac{2^j |\cT_j|}{|H_0|}$. We remove at most $2N^{-3\epsilon/4}|H_0|$ vertices for each $\cR_j$ and each $\cT_j$, for a total of $\lesssim N^{-3\epsilon/4}\log N|H_0|$ vertices removed. In particular, we have $|H|\sim|H_0|$.

Next define $\cL_2\subseteq\cL_1$ by removing any line $\ell^{\Phy}\in\cL_1$ which contains any point of $H_0\setminus H$. Since every line in $\cL_1\subseteq\cL_0$ was $[2,K]$-rich for $H_0$, we see that every line in $\cL_2$ is $[2,K]$-rich for $H\subseteq H_0$. Since removing some $x\in H$ removes $\cL_1(x)\subseteq\cL_0(x)$ from $\cL_1$, \cref{eq:l0-degree-upper} implies that
\[|\cL_1\setminus\cL_2|\lesssim N^{\epsilon/2}|H_0|^{-1}|\cL_0|\cdot N^{-3\epsilon/4} |H_0|\log N=N^{-\epsilon/4}|\cL_0|\log N,\]
so $|\cL_2|\sim|\cL_1|\sim|\cL_0|$.

\vspace{3pt}\noindent\textbf{Final dyadic decomposition of $\cR,\cT$:} Define the covering of $\cL_2$ by the sets $\cL_{j,j'}$ where $\ell^{\Phy}\in\cL_2$ lies in $\cL_{j,j'}$ if there exist $\rho \in \cR_j$ and $\tau \in \cT_{j'}$ such that $\ell^{\Phy} \subset F^{\Phy}_\rho, F^{\Phy}_\tau$. By the pigeonhole principle, there exist $j,j'$ so that $|\cL_{j, j'}|\gtrsim \log^{-2}N |\cL_2|$. Fix such values of $j,j'$ from now on.

Note that elements $\rho \in \cR_j$ are $[2^j,2^{j+1})$-rich for $\phi(H_0)$, but not necessarily for $\phi(H)$. To fix this, we do a second dyadic decomposition. For $\tj \le j$, define $\cR'_{\tj}\subseteq\cR_j$ to be the set of $\rho$ so that $|F_\rho^{\Phy} \cap \phi(H)|\in [2^{\tj}, 2^{\tj+1})$. Similarly, define $\cT'_{\tj}\subseteq\cT_j$. Define the covering of $\cL_{j,j'}$ by the sets $\cL'_{\tj,\tj'}$ where $\ell^{\Phy}\in\cL_{j,j'}$ lies in $\cL'_{\tj,\tj'}$ if there exist $\rho \in \cR'_{\tj}$ and $\tau \in \cT'_{\tj'}$ such that $\ell^{\Phy} \subset F^{\Phy}_\rho , F^{\Phy}_\tau$. By another application of the pigeonhole principle, there exist $\tj,\tj'$ so that $|\cL'_{\tj, \tj'}|\gtrsim \log^{-2}N |\cL_{j,j'}|$. Fix such values of $\tj,\tj'$ and define $\cL_3=\cL'_{\tj,\tj'}$. In addition, define $\cR=\cR'_{\tj}$ and $\cT=\cT'_{\tj'}$ and set $k=2^{\tj}$ and $k'=2^{\tj'}$. 

Note that we have the bounds $|\cL_3|\gtrsim \log^{-2}N|\cL_{j,j'}|\gtrsim\log^{-4}N|\cL_2|\sim\log^{-4}N|\cL_0|$. In addition, $k=2^{\tj}\leq 2^j$ and $k'=2^{\tj'}\leq 2^{j'}$. With these definitions we have proved (iii) aside from the lower bound on $k,k'$.

By \cref{thm:very-rich-bound-main}, we have
\begin{equation}
\label{eq:cR-upper-bound}
|\cR|\leq|\cR_j|\lesssim \frac{N|H_0|\log^{3/2}N}{2^{3j/2}}\lesssim N^{2/3}|H|2^{-j}.
\end{equation}
The last inequality uses the bounds $|H|\sim|H_0|$ and $2^j\gtrsim k_0=N^{2/3}\log^9N$.

As for each $\rho\in\cR$, the 3-flat $F^{\Phy}_\rho$ contains at most $4k^2$ lines $\ell^{\Phy} \in \cL_3$, we deduce
\begin{equation}
\label{eq:l3-upper-Rk}
|\cL_3| \le 4k^2 |\cR|\lesssim N^{2/3}|H| k^22^{-j}\leq N^{2/3}|H|k.
\end{equation}
As \cref{eq:l0-lower-absolute} gives $|\cL_3|\gtrsim \log^{-4}N|\cL_0|\gtrsim N^{-\epsilon}\abs{\cL_{[2,K]}(\phi(\cP\times\cP))}$, we conclude that
\[k\gtrsim N^{-2/3-\epsilon}|H|^{-1}\abs{\cL_{[2,K]}(\phi(\cP\times\cP))},\]
the desired lower bound on $k$. Symmetrically, the same lower bound holds for $k'$.

To prove (iv), we first need to compare the sizes $\cR_j$ and $\cR$.
By \cref{eq:l0-lower-relative,lem:many-lines-from-rich-motions}, we have
\begin{equation}
\label{eq:l3-lower-R}
|\cL_3| \gtrsim \log^{-4} N|\cL_0|\gtrsim  \log^{-9}N|\cL_{[2,K]}(\phi(H_0))| \gtrsim \log^{-17}N |\cR_j| 2^{2j}.
\end{equation}
(Note $2^j\geq k_0=N^{2/3}\log^9 N$, so the hypothesis of this lemma is satisfied.)
Combining this with the first inequality in \cref{eq:l3-upper-Rk}, we see that $|\cR|\gtrsim \log^{-17}N|\cR_j|(2^j/k)^2$.

Now for each $x\in\phi(H)$, by definition of $H$, we have
\[|\cR(x)|\leq|\cR_j(x)|<N^{3\epsilon/4}|\cR_j|2^j|H_0|^{-1}\lesssim N^{3\epsilon/4}\log^{17}N|\cR|\paren{\frac k{2^j}}^22^j|H|^{-1}\leq N^{\epsilon}k|\cR||H|^{-1}.\]
Symmetrically, the analogous bound holds on $|\cT(x)|$, proving (iv).

\vspace{3pt}\noindent\textbf{Final pruning of $\cL$:} To establish (v), we remove every line in $\cL_3$ which is contained in two $F_{\rho}^{\Phy}$ or in two $F_{\tau}^{\Phy}$. In particular, define $\cL\subseteq\cL_3$ by removing any $\ell^{\Phy}\in\cL_3$ if there exist distinct $\rho,\rho'\in\cR$ with $\ell^{\Phy}\subset F_\rho^{\Phy}, F_{\rho'}^{\Phy}$ and similarly removing any $\ell^{\Phy}$ if there exist distinct $\tau,\tau'\in\cT$ with the same property.

With this definition of $\cL$, property (v) clearly holds. Next we bound $|\cL_3\setminus\cL|$. Suppose that $\ell^{\Phy}\in\cL_3\setminus\cL$ was removed because $\ell^{\Phy}\subset F_\rho^{\Phy}, F_{\rho'}^{\Phy}$. By \cref{prop:intersection-from-same-family}, we have $\ell^{\Phy}= F_\rho^{\Phy}\cap F_{\rho'}^{\Phy}$, so the line is uniquely determined by $\rho,\rho'$. Any such $\ell^{\Phy}$ contains at least two $x\in\phi(H)$; for any such $x$, we have $\rho,\rho'\in\cR(x)$. Thus at most $|\cR(x)|^2$ lines through $x$ are removed because of $\cR$. Applying the symmetric argument to $\cT$, we can bound
\[|\cL_3\setminus\cL|\leq\sum_{x\in\phi(H)}|\cR(x)|^2+|\cT(x)|^2.\]

Now we can bound
\begin{align*}
\sum_{x\in\phi(H)}|\cR(x)|^2
&\lesssim N^\epsilon k|\cR||H|^{-1}\cdot 2k|\cR|\leq N^{4/3+\epsilon}k^2|H|2^{-2j}\leq N^{10/3+\epsilon}
\end{align*}
by (iv) and \cref{eq:cR-upper-bound}. Using a symmetric bound for the contribution from $\cT$, we get $|\cL_3\setminus\cL|\lesssim N^{10/3+\epsilon}$. As \cref{eq:l0-lower-absolute} implies $|\cL_3|\gtrsim \log^{-4}N|\cL_0|\gtrsim N^{-3\epsilon/4}\abs{\cL_{[2,K]}(\phi(\cP\times\cP))}\geq N^{10/3+5\epsilon/4}$, we conclude that $|\cL|\sim|\cL_3|$ and that (i) holds.

To check (ii), note that \cref{eq:l0-degree-upper} implies $|\cL(x)|\leq|\cL_0(x)|\lesssim N^{\epsilon/2}|H_0|^{-1}|\cL_0|\lesssim N^{\epsilon}|H|^{-1}|\cL|$ as $|H_0|\sim|H|$ and $|\cL_0|\lesssim \log^4N|\cL|$.

We already proved (iii)(iv)(v). For (vi), we have $|\cL|\leq|\cL_3|\leq 4k^2|\cR|$ by \cref{eq:l3-upper-Rk}, while $|\cL|\sim|\cL_3|\gtrsim \log^{-17}N|\cR_j|2^{2j}\geq\log^{-17}N|\cR|k^2$ by \cref{eq:l3-lower-R}. The bound in terms of $\cT$ is analogous, proving (vi).
\end{proof}

\subsection{Proof of the main theorem}\label{ssec:proof-of-main-thm}

In this subsection we combine everything we have done so far to prove the following estimate on 2-rich lines.

\begin{theorem}\label{thm:main-thm-2-rich-lines}
Let $\cP \subset \R^3$ be a set of $N$ points such that $|\cP \cap \gamma| \le N^{2/3}\deg \gamma$ for every curve $\gamma \subset\CC^3$, and $|\cP \cap F| \le N^{2/3}$ for every plane $F\subset \R^3$, and $|\cP \cap \SS| \le N^{2/3}$ for every sphere $\SS\subset \R^3$. There exist $C\sim \log^{4/3} N$ and $\epsilon \sim \sqrt{ \log\log N/\log N}$ so that
\[
|\cL_{[2, C]}(\phi(\cP\times \cP))| \le N^{10/3+\varepsilon}.
\]
\end{theorem}

We first show that this implies the main theorem.

\begin{proof}[Proof of \cref{thm:main} assuming \cref{thm:main-thm-2-rich-lines}]
Let $\cP\subset\R^3$ be a set of size $N$. First suppose that there is an irreducible algebraic curve $\gamma\subset\CC^3$ containing more than $N^{2/3}\deg \gamma$ points of $\cP$. Fix a point $p\in \cP\cap \gamma$. For each positive real $r>0$, consider the sphere $\S_r$ of radius $r$ centered at $p$; clearly $\gamma\not\subset\S_r$ since $p\in\gamma$ but $p\not\in\S_r$. Since $\S_r$ is an irreducible degree 2 variety, we see that $\gamma\cap\S_r$ consists of at most $2\deg \gamma$ points by B\'ezout's theorem. Considering all of the distances $r$ spanned between $p$ and $(\cP\cap \gamma)\setminus\{p\}$, we see that $p$ forms at least $\tfrac{1}{2}N^{2/3}$ distinct distances.

Next suppose $\abs{\cP\cap F}\geq N^{2/3}$ for a plane $F$. By the Guth--Katz distinct distances estimate, $\cP\cap F$ spans at least $\gtrsim N^{2/3}/\log N$ distinct distances (see \cref{thm:guth-katz-distance-energy}). Similarly, if $\abs{\cP\cap \S}\geq N^{2/3}$ for a sphere $\S$, we have the same distinct distances bound by \cref{thm:distinct-distances-sphere}.

Therefore we may now assume that the hypotheses of \cref{thm:main-thm-2-rich-lines} hold, giving the bound
\[\abs{\cL_{[2,C_0]}(\phi(\cP\times\cP))}\leq \abs{\cL_{[2,C]}(\phi(\cP\times\cP))}\leq N^{10/3+\epsilon},\]
where $C\sim \log^{4/3} N$ and $\epsilon\sim \sqrt{\log \log N/\log N}$ are the constants in \cref{thm:main-thm-2-rich-lines} and $C_0\sim 1$ is the constant in \cref{thm:beck-trick}.

Finally, by \cref{thm:beck-trick} we conclude that $\cP$ spans $\gtrsim N^{2/3-\epsilon}$ distinct distances. Increasing the value of $\epsilon$ by $\lesssim 1/\log N$, we can absorb the hidden multiplicative factor, which proves the desired result.
\end{proof}

Now it remains to prove \cref{thm:main-thm-2-rich-lines}.
For the sake of contradiction, assume that there are too many $[2,C]$-rich isotropic lines for $\phi(\cP\times \cP)$.
We may thus apply \cref{lem:finding-uniform-quadruple} to obtain $(H,\cL,\cR,\cT)$ that regularizes the setting.
Recall that each $\ell^{\Phy}\in \cL$ is captured by a unique $\rho\in \cR$ and $\tau\in \cT$.
Therefore we may denote by $\rho(\ell^{\Phy}), \tau(\ell^{\Phy})$ the unique $\rho\in \cR, \tau\in \cT$ with $\rho\in \ell^+, \tau\in \ell^-$, respectively.

The first key step is to only look locally at a point $x\in \phi(H)$ and consider the incidence structure between $\cR(x)$ and $\cT(x)$.
Each $\rho\in \cR(x)$ is a point in $F^+_x$ and each $\tau\in \cT(x)$ corresponds to a $3$-flat $F^+_{\tau}$ that intersects $F^+_x$ in a $2$-flat.
Also notice that for any $\ell^{\Phy}\in \cL(x)$, we have $\rho(\ell^{\Phy})\in \ell^+\subset F^+_{\tau(\ell^{\Phy})}\cap F^+_x$.
Therefore $\cL(x)$ gives rise to many point--$2$-flat incidences in $F^+_x$ between the set of points $\cR(x)$ and set of planes $\{F^+_\tau\cap F^+_x: \tau\in \cT(x)\}$.
We may now apply \cref{thm:point-plane-bipartite-graph-decomp}, which says that the incidences between $\cR(x)$ and $\{F^+_\tau\cap F^+_x: \tau\in \cT(x)\}$ can be decomposed to some structured pieces (the vertex-disjoint complete bipartite graphs) and some sparse leftover.
The leftover is easier to control, so we first prove a lemma doing so.

\begin{lemma}\label{lem:main-argument-sparse-contribution}
Let $\cP\subset\R^3$ be a set of $N$ points, such that $|\cP\cap\gamma|\leq N^{2/3}\deg\gamma$ for every curve $\gamma$, and each plane contains at most $N^{2/3}$ points of $\cP$.
Let $H\subseteq \cP\times \cP$, $\cR\subset E^+_g$, $\cT\subset E^-_g$, $k,k'\geq 1$ be such that $|F^{\Phy}_\rho\cap \phi(H)|\in [k,2k)$ for every $\rho\in \cR$ and $|F^{\Phy}_{\tau}\cap \phi(H)|\in [k',2k')$ for every $\tau\in\cT$.
For each $x\in \phi(H)$, pick a subset $\cI_x\subseteq \cR(x)\times \cT(x)$.
Define $\cL\subseteq \cL_{\geq 2}(\phi(H))$ to be the set of isotropic lines $\ell^{\Phy}\in \cL_{\geq 2}(\phi(H))$ such that there is some $x\in \ell^{\Phy}\cap \phi(H)$ and $(\rho,\tau)\in \cI_x$ with $\ell^{\Phy}\subset F^{\Phy}_{\rho}\cap F^{\Phy}_{\tau}$.
For any $\eta\in(0,1/2)$ such that $k\geq \eta^{-1}N^{2/3}$,
\[\abs{\cL}\lesssim N^{10/3}\log^2 N+\eta k^2\abs{\cR}+\eta^{-3}k^{1/2}\sum_{x\in \phi(H)}\abs{\cI_x}.\]
\end{lemma}

\begin{proof}
For each $\rho\in \cR$, the 3-flat $F^{\Phy}_{\rho}$ contains $[k,2k)$ points of $\phi(H)$. We claim that at most $N^{2/3}$ of these points are coplanar. Indeed, since $\rho$ is a genuine rigid motion, if this were not the case, there would be a set of $N^{2/3}+1$ coplanar points in $\cP$ which are mapped by $\rho$ to another set of $N^{2/3}+1$ coplanar points in $\cP$. Since each plane contains at most $N^{2/3}$ points of $\cP$, this is impossible.

Write $\cX_\rho=F^{\Phy}_\rho\cap\phi(H)$. Picking an isomorphism $F^{\Phy}_\rho\cong\R^3$, we can view $\cX_\rho$ as a subset of $\R^3$. 
We will apply \cref{lem:pair-covering-structure-new} with $(\cP,N,k)=(\cX_\rho,|\cX_\rho|,C\eta^{-3}|\cX_\rho|^{1/2})$ where $C$ is the constant in \cref{lem:pair-covering-structure-new}. Note that $\eta|\cX_\rho|\geq \eta k\geq N^{2/3}$ by hypothesis, so at most $\eta|\cX_\rho|$ points of $\cX_\rho$ are coplanar. Thus the hypotheses of \cref{lem:pair-covering-structure-new} are satisfied.  

For any line $\ell^{\Phy}\in \cL$, we may pick $x(\ell^{\Phy})\in \ell^{\Phy}\cap \phi(H)$ and $\rho(\ell^{\Phy})\in \cR$ and $\tau(\ell^{\Phy})\in \cT$ so that $(\rho(\ell^{\Phy}), \tau(\ell^{\Phy}))\in \cI_{x(\ell^{\Phy})}$ and $\ell^{\Phy}\subset F^{\Phy}_{\rho(\ell^{\Phy})}\cap F^{\Phy}_{\tau(\ell^{\Phy})}$. We can also choose $x'(\ell^{\Phy})\in\ell^{\Phy}\cap \phi(H)$ so that $x(\ell^{\Phy})\neq x'(\ell^{\Phy})$.

Write $\cL=\bigcup_{\rho\in\cR}\cL^\rho$ where $\cL^\rho$ is the set of lines $\ell^{\Phy}\in\cL$ so that $\rho(\ell^{\Phy})=\rho$. It is clear that for any $\ell^{\Phy}\in \cL^{\rho}$, we have $x(\ell^{\Phy}), x'(\ell^{\Phy})\in \cX_\rho=F^{\Phy}_{\rho}\cap \phi(H)$. Let us write
\[\cL^{\rho}=\cL^{\rho}_{\textrm{poor-plane}}\cup\cL^{\rho}_{\textrm{rich-line}}\cup\cL^{\rho}_{\textrm{poor-line}}\]
where $\ell^{\Phy}\in\cL^{\rho}$ satisfies
\begin{itemize}
    \item $\ell^{\Phy}\in\cL^{\rho}_{\textrm{poor-plane}}$ if $\ell^{\Phy}$ is not contained in a 2-flat $F\subset F_\rho^{\Phy}$ which is $C\eta^{-3}|\cX_\rho|^{1/2}$-rich for $\cX_{\rho}$;
    \item $\ell^{\Phy}\in\cL^{\rho}_{\textrm{rich-line}}$ if there exists a 2-flat $F\subset F_\rho^{\Phy}$ containing $\ell^{\Phy}$ that is $C\eta^{-3}|\cX_\rho|^{1/2}$-rich for $\cX_{\rho}$ and a line $\ell^{\Phy}_*(\ell^{\Phy})$ that lies in $F$, is $|\cX_{\rho}|^{1/2}$-rich for $\cX_{\rho}$, and contains at least one of $x(\ell^{\Phy})$ and $x'(\ell^{\Phy})$;
    \item $\ell^{\Phy}\in\cL^{\rho}_{\textrm{poor-line}}$ if there exists a 2-flat $F\subset F_\rho^{\Phy}$ containing $\ell^{\Phy}$ that is $C\eta^{-3}|\cX_\rho|^{1/2}$-rich for $\cX_{\rho}$, but there does not exist any line that lies in $F$, is $|\cX_{\rho}|^{1/2}$-rich for $\cX_{\rho}$, and contains at least one of $x(\ell^{\Phy})$ and $x'(\ell^{\Phy})$.
\end{itemize}

By \cref{lem:pair-covering-structure-new}, we have $|\cL^{\rho}_{\textrm{poor-line}}|\lesssim \eta|\cX_{\rho}|^2\lesssim \eta k^2$, which implies
\[\abs{\bigcup_{\rho\in\cR}\cL^{\rho}_{\textrm{poor-line}}}\lesssim \eta k^2|\cR|.\]

Define $\cL_{\textrm{rich-line}}=\bigcup_{\rho\in\cR}\cL^\rho_{\textrm{rich-line}}$. For each $\ell^{\Phy}\in\cL_{\textrm{rich-line}}$, we have that there exist a line $\ell_*^{\Phy}$ and a genuine rigid motion $\rho\in\cR\subset E_g^+$ so that $\ell^{\Phy}_*\subset F_\rho^{\Phy}$, one of $x(\ell^{\Phy}),x'(\ell^{\Phy})$ lies on $\ell^{\Phy}_*$, and $\ell^{\Phy}_*$ is $|\cX_{\rho}|^{1/2}\geq k^{1/2}\geq N^{1/3}$-rich for $\cX_\rho\subset\phi(H)$. If $x(\ell^{\Phy})\in\ell^{\Phy}_*$, define $(p,q)=\phi^{-1}(x(\ell^{\Phy}))$ and $(p',q')=\phi^{-1}(x'(\ell^{\Phy}))$. If $x'(\ell^{\Phy})\in\ell^{\Phy}_*$ instead, swap the definitions of $(p,q)$ and $(p',q')$. Now define $\ell_p=\pi_1(\phi^{-1}(\ell^{\Phy}_*))$ and $\ell_q=\pi_2(\phi^{-1}(\ell^{\Phy}_*))$, where $\pi_1,\pi_2$ are the projections onto the first three and last three coordinates, respectively. Since $\ell^{\Phy}_*\subset F_\rho^{\Phy}$, we have that $\rho$ maps $\ell_p$ to $\ell_q$, each of which are $N^{1/3}$-rich for $\pi_i(\phi^{-1}(\phi(H)))\subseteq\cP$. As $\ell^{\Phy}\in\cL_{\textrm{rich-line}}$ is uniquely identified as $\aff\{\phi(p,q),\phi(p',q')\}$, we see that our bound on the distance energy contribution of rich lines, \cref{lem:rich-lines}, gives
\[|\cL_{\textrm{rich-line}}|\lesssim N^{10/3}\log^2N.\]

Finally, define $\cL_{\textrm{poor-plane}}=\bigcup_{\rho\in\cR}\cL^\rho_{\textrm{poor-plane}}$. We claim that for each $x\in\phi(H)$ and each $(\rho,\tau)\in\cI_x$, there are at most $\lesssim C\eta^{-3}k^{1/2}$ lines $\ell^{\Phy}\in\cL_{\textrm{poor-plane}}$ for which $x=x(\ell^{\Phy})$, and $\rho=\rho(\ell^{\Phy})$, and $\tau=\tau(\ell^{\Phy})$. To see this, note that such a line $\ell^{\Phy}$ is uniquely identified as $\aff\{x(\ell^{\Phy}),x'(\ell^{\Phy})\}$. Furthermore, $x'(\ell^{\Phy})\in\phi(H)\cap F_\rho^{\Phy}\cap F_\tau^{\Phy}=\cX_\rho\cap F$, where $F$ is the 2-flat $F_\rho^{\Phy}\cap F_\tau^{\Phy}$. By definition of $\cL_{\textrm{poor-plane}}$, we see that $|\cX_\rho\cap F|< C\eta^{-3}|\cX_\rho|^{1/2}\lesssim C\eta^{-3}k^{1/2}$, as desired. Summing over all $x\in\phi(H)$, and all $|\cI_x|$ choices for $(\rho,\tau)$ given $x$, we see
\[|\cL_{\textrm{poor-plane}}|\lesssim C\eta^{-3}k^{1/2}\sum_{x\in\phi(H)}|\cI_x|.\qedhere\]
\end{proof}

Having proved \cref{lem:main-argument-sparse-contribution}, it suffices to bound the ``structured'' part of the point--2-flat incidence problem in $F_x^+$ for each $x\in \phi(H)$.
Recall that for each line $\ell^{\Phy}$, it is captured by a unique $\rho(\ell^{\Phy})\in E^+_g$ and $\tau(\ell^{\Phy})\in E^-_g$. \cref{thm:point-plane-bipartite-graph-decomp} produces a small number of lines $\ell_*^+$ which capture these incidences, meaning that $\rho(\ell^{\Phy})\in \ell_*^+\subset F_x^+\cap F_{\tau(\ell^{\Phy})}^+$.
Switching back to $E^{\Phy}_o$, this means $x\in \ell_*^{\Phy}\subset F_{\rho(\ell^{\Phy})}^{\Phy}\cap F_{\tau(\ell^{\Phy})}^{\Phy}$.
If we repeat this argument for another $x'\in \ell^{\Phy}$ to get $x'\in (\ell_*')^{\Phy}\subset F_{\rho(\ell^{\Phy})}^{\Phy}\cap F_{\tau(\ell^{\Phy})}^{\Phy}$, then we get that each $\ell^{\Phy}$ in the structured part corresponds to $\ell_*^{\Phy},(\ell_*')^{\Phy}$ both lying in the $2$-flat $F_{\rho(\ell^{\Phy})}^{\Phy}\cap F_{\tau(\ell^{\Phy})}^{\Phy}$.
The key observation is that working in $E^{\Phy}$, then $F_{\rho(\ell^{\Phy})}^{\Phy}\cap F_{\tau(\ell^{\Phy})}^{\Phy}\cong \PP_{\RR}^2$ and thus $\ell_*^{\Phy},(\ell_*')^{\Phy}$ must intersect.

Repeating this argument for all $\ell^{\Phy}$ corresponding to a single $\rho$, we produce many lines in $F^{\Phy}_{\rho}$ that intersect a lot.
This will put us in position to apply \cref{lem:line-structure-new} to show that most intersections come from a few points, a few lines, a few planes or a few reguli.
For most of the cases, we will be done by what we have proven in \cref{sec:degenerate}.
However, the contribution from the set of points $\fP_\rho$ requires a significantly more involved argument.
In fact, to carry out the argument, we will also need to run the same argument with $\tau(\ell^{\Phy})=\tau$ fixed instead, which again produces a set of few points $\fP_\tau$ by \cref{lem:line-structure-new}.
We will only need to take care of $\ell^{\Phy}$ for which $\ell_*^{\Phy}, (\ell_*')^{\Phy}$ intersect at a point $x_*$ lying in both $\fP_{\rho(\ell^{\Phy})}$ and $\fP_{\tau(\ell^{\Phy})}$.
We begin by bounding the number of such lines for a given $x_*$ under some additional mild conditions.

\begin{lemma}\label{lem:rich-point-xstar}
There exists a constant $C>0$ so that the following holds.
Let $\cP\subset\R^3$ be a set of size $N$ such that every plane contains at most $N^{2/3}$ points of $\cP$.
Let $H\subseteq \cP\times \cP$ be a subset, and for each $\tau\in E^-$, define $H_\tau\subseteq H$ to be the set of $(p,q)\in H$ so that $\phi(p,q)\in F_\tau^{\Phy}$ and any isotropic line through $\phi(p,q)$ which lies in $F_\tau^{\Phy}$ contains at most $N^{1/3}$ points of $\phi(H)$.

Fix a point $x_*\in E^{\Phy}$.
For $k,k'$ positive integers, let $\cR\subset F^+_{x_*}\cap E^+_g$ be a set of genuine positive rigid motions $\rho$ with $|\phi(H)\cap F^{\Phy}_{\rho}|\in [k,2k)$ and similarly define $\cT\subset F^-_{x_*}\cap E^-_g$ to be a set of genuine negative rigid motions $\tau$ with $|\phi(H)\cap F^{\Phy}_{\tau}|\in [k',2k')$.
Finally, define $\cL\subseteq\cL_{\geq 2}(\phi(H))$ to be the set of isotropic lines $\ell^{\Phy}$ so that there exists $(\rho,\tau)\in \cR\times \cT$ with $\ell^{\Phy}\subset F^{\Phy}_{\rho}\cap F^{\Phy}_{\tau}$ such that $\ell^{\Phy}$ is 2-rich for $\phi(H_\tau)$.
Then for any parameter $A\geq CN^{1/3}$, we have
\[\abs{\cL}\lesssim \left(N^{2/3}k'+A^{-1}N^{1/3}(k')^2\right)\abs{\cT}+A^2\left(\abs{\cR}^{1/2}\abs{\cT}+\abs{\cR}\abs{\cT}^{1/2}\right).\]
\end{lemma}

Note that in this result, $x_*$ lies in $E^{\Phy}$, but not necessarily in $E^{\Phy}_o$.

\begin{proof}
We begin by studying the incidence structure between $\cR$ and $\cT$. By \cref{prop:intersection-from-diff-family}, for $(\rho,\tau)\in\cR\times\cT$, the intersection $F_{\rho}^{\Phy}\cap F_\tau^{\Phy}$ is either a point or a plane, with the latter occurring if and only if $\rho\in F_{\tau}^+$. As $\rho\in F_{x_*}^+$ and $\tau\in F_{x_*}^-$, this is a point--2-flat incidence $\rho\in F_{\tau}^+\cap F_{x_*}^+$ in $F_{x_*}^+\cong\PP_\R^3$.

We wish to apply \cref{thm:point-plane-bipartite-graph-decomp} to this incidence problem; to do this, pick a copy of $\R^3$ in $\PP_{\R}^3$ which contains all of the points $\cR\cap F_{x_*}^+$. Now applying \cref{thm:point-plane-bipartite-graph-decomp} to the point--2-flat incidence problem in this copy of $\R^3$, we conclude that there exist a set of lines $\cL_*^+$ in $F_{x_*}^+$ and a subset $\cI\subseteq \cR\times \cT$ with the following properties:
\begin{itemize}
    \item $\abs{\cL_*^+}\lesssim \min\{\abs{\cR}^{1/2}, \abs{\cT}^{1/2}\}$;
    \item $\abs{\cI}\lesssim \abs{\cR}^{1/2}\abs{T}+\abs{\cR}\abs{\cT}^{1/2}$; and
    \item if $(\rho,\tau)\in \cR\times \cT$ satisfies $\rho\in F_{\tau}^+$, then either $\rho\in \ell_*^+\subset F_{\tau}^+\cap F^+_{x_*}$ for some $\ell_*^+\in \cL_*^+$ or $(\rho,\tau)\in \cI$.
\end{itemize}

For any line $\ell^{\Phy}\in\cL$, we may define $\rho(\ell^{\Phy})\in\cR$ and $\tau(\ell^{\Phy})\in\cT$ so that $\ell^{\Phy}\subset F_{\rho(\ell^{\Phy})}^{\Phy}\cap F_{\tau(\ell^{\Phy})}^{\Phy}$. We define a partition $\cL=\cL_1\sqcup\cL_2$. For any $\ell^{\Phy}\in \cL$, we have $\ell^{\Phy}\subset F_{\rho(\ell^{\Phy})}^{\Phy}\cap F_{\tau(\ell^{\Phy})}^{\Phy}$ which, by \cref{prop:intersection-from-diff-family}, implies $\rho(\ell^{\Phy})\in F_{\tau(\ell^{\Phy})}^+$. Therefore either $\rho(\ell^{\Phy})\in \ell_*^+\subset F_{\tau(\ell^{\Phy})}^+\cap F^+_{x_*}$ for some $\ell_*^+\in \cL_*$ or $(\rho(\ell^{\Phy}),\tau(\ell^{\Phy}))\in \cI$.
Say $\ell^{\Phy}\in\cL_1$ if the former occurs, and $\ell^{\Phy}\in \cL_2$ if the latter occurs.

We first bound $\cL_1$.
Note that $\rho\in \ell_*^+\subset F_{\tau}^+\cap F^+_{x_*}$ implies $\ell_*^{\Phy}\subset F^{\Phy}_{\rho}\cap F^{\Phy}_{\tau}$ by \cref{prop:all-lines}. Set $\cL_* = \{\ell_*^{\Phy}:\ell_*^+\in \cL_*^+\}$.
Now for any line $\ell^{\Phy}\in \cL_1\setminus \cL_*$, we know that there is $\ell_*^{\Phy}\in \cL_*$ so that $\aff\{\ell^{\Phy}\cup \ell_*^{\Phy}\}=F_{\rho(\ell^{\Phy})}^{\Phy}\cap F_{\tau(\ell^{\Phy})}^{\Phy}$. By definition, $\ell^{\Phy}$ contains at least two points of $\phi(H_{\tau(\ell^{\Phy})})$, so at least one of these points must not lie on $\ell_*^{\Phy}$.

Given $\tau=\tau(\ell^{\Phy})$ and $\ell_*^{\Phy}$, there are fewer than $2k'$ options for this point in $\ell^{\Phy}\cap \phi(H_\tau)$ that does not lie on $\ell_*^{\Phy}$ as such a point must lie in $\phi(H)\cap F_\tau^{\Phy}$. Given this point, the 2-flat $\aff\{\ell^{\Phy}\cup \ell_*^{\Phy}\}$ is determined. As distinct $F_{\rho}^{\Phy}$ intersect in at most a line by \cref{prop:intersection-from-same-family}, there is at most one choice of $\rho$ so that $F_{\rho}^{\Phy}$ contains this 2-flat. Finally, to determine $\ell^{\Phy}$, we pick one other point of $\phi(H_\tau)$ lying on this line. There are at most $N^{2/3}$ options for this point as it lies in $\phi(H_\tau)\cap F_{\rho}^{\Phy}\cap F_{\tau}^{\Phy}$. (If this set were larger, we would have a plane in $F_\rho^{\Phy}$ containing more than $N^{2/3}$ points, implying that $\cP$ has a set of more than $N^{2/3}$ coplanar points). Therefore we see that for each $\tau,\ell_*^{\Phy}$, there are at most $2N^{2/3}k'$ corresponding $\ell^{\Phy}\in\cL_1\setminus\cL_*$.

Finally, note that $\tau,\ell_*^{\Phy}$ must satisfy $\ell_*^{\Phy}\subset F^{\Phy}_\tau$. As the line--3-flat incidence graph between $\cL_\ast$ and $\{F_\tau^{\Phy}\}_{\tau\in\cT}$ is $K_{2,2}$-free (by \cref{prop:intersection-from-same-family}), we see that there are at most $|\cT|+|\cL_\ast|^2$ possible pairs with $\ell_*^{\Phy}\subset F^{\Phy}_\tau$. Putting everything together, we have shown that
\[\abs{\cL_1\setminus \cL_*}\leq 2N^{2/3}k'\abs{\set{(\ell_*^{\Phy},\tau)\in \cL_*\times \cT: \ell_*^{\Phy}\subset F^{\Phy}_{\tau}}}\leq 2N^{2/3}k'\paren{|\cT|+|\cL_*|^2}.\]
As $|\cL_\ast|\lesssim|\cT|^{1/2}$, we conclude
\[|\cL_1|\lesssim N^{2/3}k'|\cT|.\]

It remains to bound $\cL_2$.
Recall that for any $\ell^{\Phy}\in \cL_2$, we have $(\rho(\ell^{\Phy}),\tau(\ell^{\Phy}))\in \cI$.
We define a partition $\cL_2=\cL_2^{\mathrm{rich}}\sqcup\cL_2^{\mathrm{poor}}$ where $\ell^{\Phy}\in\cL_2^{\mathrm{rich}}$ if $|\phi(H_{\tau(\ell^{\Phy})})\cap F^{\Phy}_{\rho(\ell^{\Phy})}\cap F^{\Phy}_{\tau(\ell^{\Phy})}|\geq A$ while $\ell^{\Phy}\in\cL_2^{\mathrm{poor}}$ otherwise.

It is clear that
\[|\cL_2^{\mathrm{poor}}|< A^2\abs{\cI}\lesssim A^2\left(\abs{\cR}^{1/2}\abs{\cT}+\abs{\cR}\abs{\cT}^{1/2}\right),\]
so it remains to bound $\cL_2^{\mathrm{rich}}$.

Fix some $\tau\in \cT$. For each $r\geq A$, let $\cR_r\subseteq \cR$ be the set of $\rho\in\cR$ with $|\phi(H_\tau)\cap F^{\Phy}_{\rho}\cap F^{\Phy}_{\tau}| \geq r$. Our first goal is to bound $|\cR_r|$. We will do this by projecting radially from $x_*$ and applying the Szemer\'edi--Trotter theorem.

Let $\psi_{x_*}\colon F^{\Phy}_{\tau}\setminus  \{x_*\}\to \PP^2_{\RR}$ be the radial projection, and let $S = \psi_{x_*}\left((F^{\Phy}_{\tau}\cap \phi(H_\tau))\setminus \{x_*\}\right)$.
For any point $p\in S$, define its weight $w(p)$ to be the size of the fiber:
\[w(p)=\abs{\psi_{x_*}^{-1}(p)\cap \left((F^{\Phy}_{\tau}\cap \phi(H_\tau))\setminus \{x_*\}\right)}.\]
Since no point of $\phi(H_\tau)$ is on an $N^{1/3}$-rich line in $F_\tau^{\Phy}$, we conclude that $w(p)< N^{1/3}$.

Define $S_i\subseteq S$ to be the subset of points $p\in S$ with $w(p)\geq i$. Since $F_\rho^{\Phy}\cap F_\tau^{\Phy}$ contains at least $r$ points of $\phi(H_\tau)$ (and thus at least $r-1$ points of $\phi(H_\tau)\setminus\{x_*\}$), we conclude that
\[\sum_{i=1}^{N^{1/3}}\abs{\psi_{x_*}\left(F^{\Phy}_{\rho}\cap F^{\Phy}_{\tau}\right)\cap S_i}\geq r-1.\]
Now the key observation is that for each $\rho\in \cR_r$, we know that $F_\rho^{\Phy}\cap F_\tau^{\Phy}$ is a 2-plane through $x_*$ in $F_\tau^{\Phy}$. Thus $\psi_{x_*}(F^{\Phy}_{\rho}\cap F^{\Phy}_{\tau})$ is a line in $\PP_{\RR}^2$. Therefore we can apply the Szemer\'erdi--Trotter theorem to bound the number of incidences between $S_i$ and the set of lines $\set{\psi_{x_*}(F^{\Phy}_{\rho}\cap F^{\Phy}_{\tau})}_{\rho\in\cR_r}$, to conclude
\[(r-1)\abs{\cR_r}\leq \sum_{i=1}^{N^{1/3}}\sum_{\rho\in \cR_r}\abs{\psi_{x_*}(F^{\Phy}_{\rho}\cap F^{\Phy}_{\tau})\cap S_i}\lesssim \sum_{i=1}^{N^{1/3}}\left(\abs{S_i}+\abs{\cR_r}+\abs{S_i}^{2/3}\abs{\cR_r}^{2/3}\right).\]

As 
\[\sum_{i=1}^{N^{1/3}}\abs{S_i}\leq |\phi(H_\tau)\cap F^{\Phy}_{\tau}|\lesssim k',\]
H\"older's inequality implies 
\[\sum_{i=1}^{N^{1/3}}\abs{S_i}^{2/3}\lesssim (N^{1/3})^{1/3}(k')^{2/3}.\]
Thus
\[r\abs{\cR_r}\lesssim k'+N^{1/3}\abs{\cR_r}+N^{1/9}(k')^{2/3}\abs{\cR_r}^{2/3}.\]

Recall that $r\geq A\geq CN^{1/3}$. Letting $C$ be a sufficiently large absolute constant, we conclude
\[\abs{\cR_r}\lesssim \frac{k'}{r}+\frac{N^{1/3}(k')^2}{r^3}.\]
Also note that $\cR_r$ is empty for $r> N^{2/3}$ as each plane contains at most $N^{2/3}$ points of $\cP$.

Recall that we earlier fixed $\tau\in\cT$. For a given $\rho\in\cR$, the number of lines $\ell^{\Phy}\in \cL_2^{\mathrm{rich}}$ with $(\rho(\ell^{\Phy}),\tau(\ell^{\Phy}))=(\rho,\tau)$ is at most $|\phi(H_\tau)\cap F^{\Phy}_{\rho}\cap F^{\Phy}_{\tau}|^2$, implying that the contribution to $\cL_2^{\mathrm{rich}}$ for this fixed $\tau$ is at most
\[A^2|\cR_A|+\sum_{r=A+1}^{N^{2/3}}r^2\abs{\cR_r\setminus\cR_{r+1}}\lesssim A^2|\cR_A|+\sum_{r=A}^{N^{2/3}}r\abs{\cR_r}\lesssim N^{2/3}k'+\frac{N^{1/3}(k')^2}{A}.\]

Summing over all $\tau\in\cT$, we conclude that
\[|\cL_2^{\mathrm{rich}}|\lesssim\paren{N^{2/3}k'+A^{-1}N^{1/3}(k')^2}|\cT|.\qedhere\]
\end{proof}

Next we sum the bound from \cref{lem:rich-point-xstar} over all $x_*$ and bound the contribution from points in $\phi(H\setminus H_\tau)$ -- points on an $N^{1/3}$-rich line lying $F_\tau^{\Phy}$ -- using \cref{lem:rich-lines}.

\begin{lemma}\label{lem:rich-points-sum}
There exists an absolute constant $C>0$ so that the following holds.
Let $\cP\subset\R^3$ be a set of size $N$ such that $|\cP\cap\gamma|\leq N^{2/3}\deg\gamma$ for every curve $\gamma$, and $|\cP\cap F|\leq N^{2/3}$ for every plane $F\subset\R^3$, and $|\cP\cap\S|\leq N^{2/3}$ for every sphere $\S\subset\R^3$.
Let $H\subseteq \cP\times \cP, \cR\subseteq E^+_g, \cT\subseteq E^-_g$, $k,k'\geq 1$ be such that $|F^{\Phy}_{\rho}\cap \phi(H)|\in [k,2k)$ for every $\rho\in \cR$ and $|F^{\Phy}_{\tau}\cap \phi(H)|\in [k',2k')$ for every $\tau\in\cT$.

Fix $M\geq1$.
For each $\rho\in \cR$ and $\tau\in \cT$, let $\fP_{\rho}\subset F^{\Phy}_{\rho}$ and $\fP_{\tau}\subset F^{\Phy}_{\tau}$ be sets of at most $M$ points.
Define $\cL\subseteq\cL_{\geq2}(\phi(H))$ to be the set of isotropic lines $\ell^{\Phy}\in \cL_{\geq 2}(\phi(H))$ so that there exists $(\rho,\tau)\in \cR\times \cT$ with $\ell^{\Phy}\subset F^{\Phy}_{\rho}\cap F^{\Phy}_{\tau}$ and $\fP_{\rho}\cap \fP_{\tau}\neq\emptyset$.
Then for any parameter $A\geq CN^{1/3}$,
\[\abs{\cL}\lesssim N^{10/3}\log^2N+\left(N^{2/3}k'+A^{-1}N^{1/3}(k')^2\right)M\abs{\cT}+A^2M^{3/2}N^{-1/3}\left(\abs{\cR}^{1/2}\abs{\cT}+\abs{\cR}\abs{\cT}^{1/2}\right).\]
\end{lemma}

Note that in this result, the sets $\fP_\rho,\fP_\tau$ lie in $E^{\Phy}$, but not necessarily in $E_o^{\Phy}$.

\begin{proof}
Let $C$ be the constant produced by \cref{lem:rich-point-xstar}. For each $\tau\in\cT$, define $H_\tau\subseteq H$ to be the set of $(p,q)\in H$ so that $\phi(p,q)\in F_\tau^{\Phy}$ and any isotropic line through $\phi(p,q)$ which lies in $F_\tau^{\Phy}$ contains at most $N^{1/3}$ points of $\phi(H)$. For each line $\ell^{\Phy}\in\cL$, fix a choice of pair $(\rho(\ell^{\Phy}),\tau(\ell^{\Phy}))\in\cR\times\cT$ so that $\ell^{\Phy}\subset F^{\Phy}_{\rho(\ell^{\Phy})}\cap F^{\Phy}_{\tau(\ell^{\Phy})}$ and $\fP_{\rho(\ell^{\Phy})}\cap \fP_{\tau(\ell^{\Phy})}\neq\emptyset$.

Define a partition $\cL=\cL_1\sqcup\cL_2$ where $\ell^{\Phy}\in\cL$ lies in $\cL_2$ if $\phi(H)\cap\ell^{\Phy}\subseteq\phi(H_{\tau(\ell^{\Phy})})$. Thus for $\ell^{\Phy}\in\cL_1$, there exists $(p,q)\not\in H_{\tau(\ell^{\Phy})}$ with $\phi(p,q)\in \phi(H)\cap\ell^{\Phy}$. This means that there is an isotropic line $\ell_*^{\Phy}$ through $\phi(p,q)$ that lies in $F_{\tau(\ell^{\Phy})}^{\Phy}$ and that is $N^{1/3}$-rich for $\phi(H)$.

As in the proof of \cref{lem:main-argument-sparse-contribution}, define $\ell_p=\pi_1(\phi^{-1}(\ell^{\Phy}_*))$ and $\ell_q=\pi_2(\phi^{-1}(\ell^{\Phy}_*))$, where $\pi_1,\pi_2$ are the projections onto the first three and last three coordinates, respectively. Since $\tau(\ell^{\Phy})$ is a genuine rigid motion, we see that $\ell_p,\ell_q$ are lines, each of which is $N^{1/3}$-rich for $\cP$, such that $\tau(\ell^{\Phy})$ maps $\ell_p$ to $\ell_q$. Let $\phi(p',q')$ be another point of $\phi(H)\cap\ell^{\Phy}$. As $\tau(\ell^{\Phy})$ maps $p$ to $q$ and $p'$ to $q'$, we see that the distance quadruple $(p,q,p',q')$ is counted by \cref{lem:rich-lines}, applied to $\cP$. Finally, as $\ell^{\Phy}$ is uniquely determined as $\ell^{\Phy}=\aff\{\phi(p,q),\phi(p',q')\}$, \cref{lem:rich-lines} implies that
\[|\cL_1|\lesssim N^{10/3}\log^2N.\]

Next we bound $\cL_2$. Define $\fP = \left(\bigcup_{\rho\in \cR}\fP_{\rho}\right)\cup \left(\bigcup_{\tau \in \cT}\fP_{\tau}\right)$ and for each $x_*\in \fP$, set $\cR_{x_*} = \{\rho\in \cR:x_*\in \fP_{\rho}\}$ and similarly define $\cT_{x_*}$.
Define $\cL_{x_*}\subseteq\cL_2$ to be the set of lines $\ell^{\Phy}\in \cL_2$ so that $(\rho(\ell^{\Phy}),\tau(\ell^{\Phy}))\in \cR_{x_*}\times \cT_{x_*}$.
It is then clear that 
\[\abs{\cL_2}\leq \sum_{x_*\in \fP}\abs{\cL_{x_*}}.\]

As each $\ell^{\Phy}\in\cL_2$ lies in $\cL_{\geq 2}(\phi(H_{\tau(\ell^{\Phy})}))$, we see that for each fixed $x_*\in \fP$, \cref{lem:rich-point-xstar} bounds $|\cL_{x_*}|$ in terms of $|\cR_{x_*}|$ and $|\cT_{x_*}|$. On the other hand, \cref{lem:rich-points} gives a uniform upper bound on $|\cL_{x_*}|$ for each fixed $x_*$. Indeed, for any $\ell^{\Phy}\in\cL_2$, we have $\ell^{\Phy}\subset F_\rho^{\Phy}$ where $x_*\in\fP_\rho\subset F_\rho^{\Phy}$. Furthermore, there are at least two points $\phi(p,q),\phi(p',q')\in \phi(H)\cap\ell^{\Phy}$. As $\ell^{\Phy}$ is uniquely identified as $\aff\{\phi(p,q),\phi(p',q')\}$ where $\phi(p,q),\phi(p',q'),x_*\in F_\rho^{\Phy}$, \cref{lem:rich-points} applies to give the bound
\[|\cL_{x_*}|\lesssim N^{8/3}\log N\]
for every $x_*\in\fP$.

Write $\fP=\fP^{\mathrm{rich}}\sqcup\fP^{\mathrm{poor}}$ where $x_*\in\fP$ lies in $\fP^{\mathrm{rich}}$ if $|\cR_{x_*}|$ is one of the $N^{2/3}$ largest of these sets or $|\cT_{x_*}|$ is one of the $N^{2/3}$ largest of these sets. If not, $x_*\in\fP^{\mathrm{poor}}$. With this definition, clearly $|\fP^{\mathrm{rich}}|\leq 2N^{2/3}$. Furthermore, as
\[\sum_{x_*\in\fP}|\cR_{x_*}|=\sum_{\rho\in\cR}|\fP_\rho|\leq M|\cR|\qquad\text{ and }\qquad\sum_{x_*\in\fP}|\cT_{x_*}|=\sum_{\tau\in\cT}|\fP_\tau|\leq M|\cT|,\]
we conclude
\[\abs{\cR_{x_*}}\leq N^{-2/3}M\abs{\cR}\qquad\text{ and }\qquad\abs{\cT_{x_*}}\leq N^{-2/3}M\abs{\cT}\]
for all $x_*\in \fP^{\mathrm{poor}}$.

\cref{lem:rich-point-xstar} then gives the bound
\begin{align*}
|\cL_{x_*}|
&\lesssim\paren{N^{2/3}k'+A^{-1}N^{1/3}(k')^2}|\cT_{x_*}|+A^2\paren{|\cR_{x_*}|^{1/2}|\cT_{x_*}|+|\cR_{x_*}||\cT_{x_*}|^{1/2}}\\
&\leq \paren{N^{2/3}k'+A^{-1}N^{1/3}(k')^2}|\cT_{x_*}|+A^2M^{1/2}N^{-1/3}\paren{|\cR|^{1/2}|\cT_{x_*}|+|\cR_{x_*}||\cT|^{1/2}},
\end{align*}
for each $x_*\in\fP^{\mathrm{poor}}$.

Combining the two, we conclude
\begin{align*}
\abs{\cL_2}\leq \sum_{x_*\in \fP}\abs{\cL_{x_*}}
&\lesssim |\fP^{\mathrm{rich}}|N^{8/3}\log N+\sum_{x_*\in\fP^{\mathrm{poor}}}\left(\paren{N^{2/3}k'+A^{-1}N^{1/3}(k')^2}|\cT_{x_*}|\right.\\
&\qquad\qquad\qquad\qquad\qquad+\left.A^2M^{1/2}N^{-1/3}\paren{|\cR|^{1/2}|\cT_{x_*}|+|\cR_{x_*}||\cT|^{1/2}}\right),\\
&\lesssim N^{10/3}\log N+\paren{N^{2/3}k'+A^{-1}N^{1/3}(k')^2}M|\cT|\\
&\qquad\qquad\qquad\qquad +A^2M^{3/2}N^{-1/3}\paren{|\cR|^{1/2}|\cT|+|\cR||\cT|^{1/2}}.\qedhere
\end{align*}
\end{proof}

We are now ready to prove \cref{thm:main-thm-2-rich-lines} and thus \cref{thm:main}.

\begin{proof}[Proof of \cref{thm:main-thm-2-rich-lines}]
Define $C\sim \log^{4/3}N$ and $\epsilon'\sim \sqrt{\log\log N/\log N}$ to be the parameters produced by \cref{lem:finding-uniform-quadruple}. Set $\epsilon = 30\epsilon'$.

\vspace{3pt}\noindent\textbf{Regularization:} 
Suppose for the sake of contradiction that
\[\abs{\cL_{[2,C]}(\phi(\cP\times \cP))}\geq N^{10/3+\varepsilon}\geq N^{10/3+2\epsilon'}.\]
Applying \cref{lem:finding-uniform-quadruple} with parameter $\varepsilon'$, we obtain $(H,\cL,\cR,\cT,k,k')$ satisfying the six properties guaranteed by that lemma.

By \cref{lem:finding-uniform-quadruple}(i) we have
\begin{equation}
\label{eq:l-lower}
\abs{\cL}\geq N^{10/3+\varepsilon-\varepsilon'}.
\end{equation}
By \cref{lem:finding-uniform-quadruple}(vi), we have
\begin{equation}
\label{eq:l-lower-r-t}
\abs{\cL}\geq N^{-\varepsilon'}k^2\abs{\cR}\qquad\text{ and }\qquad\abs{\cL}\geq N^{-\varepsilon'}(k')^2\abs{\cT}.
\end{equation}
By \cref{lem:finding-uniform-quadruple}(iv) and \cref{eq:l-lower-r-t}, for every $x\in \phi(H)$ we have
\begin{equation}
\label{eq:rx-upper-l}
\abs{\cR(x)}\leq N^{\varepsilon'}k\abs{\cR}\abs{H}^{-1}\leq N^{2\varepsilon'}\frac{\abs{\cL}}{k\abs{H}}\qquad\text{ and }\qquad 
\abs{\cT(x)}\leq N^{\varepsilon'}k'\abs{\cT}\abs{H}^{-1}\leq N^{2\varepsilon'}\frac{\abs{\cL}}{k'\abs{H}}.
\end{equation}

For each $\ell^{\Phy}\in\cL$, by \cref{lem:finding-uniform-quadruple}(v), there exist a unique $\rho(\ell^{\Phy})\in\cR$ and $\tau(\ell^{\Phy})\in\cT$ so that $\ell^{\Phy}\subset F_{\rho(\ell^{\Phy})}^{\Phy}\cap F_{\tau(\ell^{\Phy})}^{\Phy}$. In addition, fix a choice of two distinct points $x(\ell^{\Phy}),x'(\ell^{\Phy})\in \phi(H)\cap\ell^{\Phy}$.

We will reach a contradiction by removing various sets of lines from $\cL$, each of which will have size at most $N^{-\epsilon'}|\cL|$. Throughout this process, the remaining set of lines will become more structured until eventually we can reach a contradiction by showing that the final set of lines also has size at most $N^{-\epsilon'}|\cL|$. \cref{eq:l-lower-r-t} will be crucial for showing that each of these sets of lines is small relative to $|\cL|$.

\vspace{3pt}\noindent\textbf{Point-plane incidences in positive space:}
For each $x\in\phi(H)$, our goal is to study $\cL(x)$, a set of lines through $x$ in physical space. Each such line $\ell^{\Phy}$ lies in $F^{\Phy}_{\rho(\ell^{\Phy})}\cap F^{\Phy}_{\tau(\ell^{\Phy})}$ where $(\rho(\ell^{\Phy}),\tau(\ell^{\Phy}))\in\cR(x)\times\cT(x)$. Recall from \cref{prop:intersection-from-diff-family} that the intersection $F^{\Phy}_{\rho(\ell^{\Phy})}\cap F^{\Phy}_{\tau(\ell^{\Phy})}$ is either a point or a $2$-flat, the latter occurring if and only if $\rho(\ell^{\Phy})\in F^+_{\tau(\ell^{\Phy})}\cap F^+_{x}$. Since the intersection contain the line $\ell^{\Phy}$, we conclude that $\rho(\ell^{\Phy})\in F^+_{\tau(\ell^{\Phy})}\cap F^+_{x}$. Thus our first aim is to study the pairs $(\rho,\tau)\in\cR(x)\times\cT(x)$ for which $\rho\in F^+_{\tau}\cap F^+_x$. We do this using a similar technique as the proof of \cref{lem:rich-point-xstar}.

Consider the point-plane incidence problem in $F^+_x\cong \PP_{\R}^3$ between the points $\cR(x)$ and the 2-flats $\{F^+_\tau\cap F^+_x\}_{\tau\in \cT(x)}$. Pick a copy of $\R^3$ in $F^+_x\cong \PP_{\R}^3$ which contains $\cR(x)$. Applying \cref{thm:point-plane-bipartite-graph-decomp} to the point--2-flat incidence problem in this copy of $\R^3$, the incidence graph between can be decomposed into a small number of vertex-disjoint complete bipartite graphs and a small number of additional edges. For any such complete bipartite graph, there exists a line $\ell^+_*\subset F^+_x$ so that the points lie on $\ell^+_*$, and the 2-flats all contain $\ell^+_*$. Define $\cR(x,\ell^+_*)$ and $\cT(x,\ell^+_*)$ to be the two sets of vertices of this graph. Thus \cref{thm:point-plane-bipartite-graph-decomp} produces a set of lines $\cL^+_{*,x}$ in $F^+_x$, disjoint subsets $\cR(x,\ell^+_*)\subseteq\cR(x)$ and $\cT(x,\ell^+_*)\subseteq\cT(x)$ for each $\ell^+_*\in\cL^+_{*,x}$, and a subset $\cI_x\subseteq \cR(x)\times \cT(x)$ with the following properties:
\begin{itemize}
\item
    \[\abs{\cL^+_{*,x}}\lesssim \min\{\abs{\cR(x)}^{1/2},\abs{\cT(x)}^{1/2}\}\leq \left(N^{2\varepsilon'}\frac{\abs{\cL}}{k\abs{H}}\right)^{1/2};\]
\item
\begin{align*}
    \abs{\cI_x}&\lesssim \abs{\cR(x)}^{1/2}\abs{\cT(x)}+\abs{\cR(x)}\abs{\cT(x)}^{1/2}\\&\lesssim N^{\epsilon'}\abs{\cT(x)}\frac{\abs{\cL}^{1/2}}{k^{1/2}|H|^{1/2}}+N^{\epsilon'}\abs{\cR(x)}\frac{\abs{\cL}^{1/2}}{(k')^{1/2}|H|^{1/2}};
\end{align*}
\item
    if $(\rho,\tau)\in \cR(x)\times \cT(x)$ satisfies $\rho\in F_\tau^+$, then either $(\rho,\tau)\in \cI_x$, or $(\rho,\tau)\in\cR(x,\ell_*^+)\times\cT(x,\ell_*^+)$ (which implies $\rho\in\ell_*^+\subset F^{+}_{\tau}\cap F^+_x$) for some $\ell_*^+\in \cL_{*,x}^+$. 
\end{itemize}

Define $\cL_{\textrm{sparse}}\subseteq\cL$ so that $\ell^{\Phy}\in\cL$ lies in $\cL_{\textrm{sparse}}$ if either $(\rho(\ell^{\Phy}),\tau(\ell^{\Phy}))$ lies in either $\cI_{x(\ell^{\Phy})}$ or $\cI_{x'(\ell^{\Phy})}$. 

We now apply \cref{lem:main-argument-sparse-contribution} to $H,\cR,\cT,k,k'$, the sets $\cI_x$, and $\cL_{\mathrm{sparse}}$ with parameter $\eta=N^{-2\epsilon'}$. By \cref{lem:finding-uniform-quadruple}(iii) we have
\begin{equation}
\label{eq:k-k'-lower}
k,k'\geq N^{-2/3-\epsilon'}|H|^{-1}N^{10/3+\epsilon}\geq N^{2/3+\epsilon-\epsilon'}\geq N^{2/3+2\epsilon'},
\end{equation}
so the hypotheses of this lemma are satisfied. Thus we conclude that
\begin{align*}
|\cL_{\textrm{sparse}}|\lesssim N^{10/3}\log^2N + N^{-2\epsilon'}k^2|\cR|+N^{6\epsilon'}k^{1/2}\sum_{x\in\phi(H)}|\cI_x|.
\end{align*}
By \cref{eq:l-lower,eq:l-lower-r-t}, we see that each of the first two terms are bounded by $\lesssim N^{-\epsilon'}|\cL|$.
The third term can be bounded by
\begin{align*}
\lesssim N^{6\epsilon'}k^{1/2}&\sum_{x\in\phi(H)}\paren{N^{\epsilon'}\abs{\cT(x)}\frac{\abs{\cL}^{1/2}}{k^{1/2}|H|^{1/2}}+N^{\epsilon'}\abs{\cR(x)}\frac{\abs{\cL}^{1/2}}{(k')^{1/2}|H|^{1/2}}}\\
&\lesssim  N^{7\epsilon'}k^{1/2}\frac{|\cL|^{1/2}}{|H|^{1/2}}\paren{\frac{k'|\cT|}{k^{1/2}}+\frac{k|\cR|}{(k')^{1/2}}}.
\end{align*}
By \cref{thm:very-rich-bound-main}, we have
\begin{equation}
\label{eq:R-T-upper}
|\cR|\lesssim \frac{N|H|\log^{3/2}N}{k^{3/2}}\qquad\text{ and }\qquad|\cT|\lesssim \frac{N|H|\log^{3/2}N}{(k')^{3/2}}
\end{equation}
so the above expression can be bounded by
\begin{align*}
\lesssim N^{7\epsilon'}|\cL|^{1/2}|H|^{1/2}N(k')^{-1/2}\log^{3/2}N
&\lesssim N^{-1/3-\epsilon/2+8\epsilon'}|\cL|^{1/2}|H|\\
&\leq N^{5/3-\epsilon/2+8\epsilon'}|\cL|^{1/2}
\leq N^{-\epsilon+9\epsilon'}|\cL|,
\end{align*}
using \cref{eq:k-k'-lower}, the trivial bound $|H|\leq N^2$, and \cref{eq:l-lower}. Recalling $\epsilon=30\epsilon'$, we conclude
\[\abs{\cL_{\textrm{sparse}}}\leq N^{-\epsilon'}\abs{\cL}.\]

\vspace{3pt}\noindent\textbf{Line-line incidences in physical space:} Define $\cL_1=\cL\setminus\cL_{\textrm{sparse}}$. Partition $\cL_1=\bigsqcup_{\rho\in\cR}\cL_1(\rho)$ where $\ell^{\Phy}\in\cL_1(\rho)$ if $\rho=\rho(\ell^{\Phy})$.

By the definition of $\cL_1$, we know that $(\rho(\ell^{\Phy}),\tau(\ell^{\Phy}))$ does not lie in $\cI_{x(\ell^{\Phy})}$ nor in $\cI_{x'(\ell^{\Phy})}$. Thus there exists a unique choice of $\ell_*^+\in\cL_{*,x(\ell^{\Phy})}^+$ so that $(\rho(\ell^{\Phy}),\tau(\ell^{\Phy}))\in \cR(x(\ell^{\Phy}),\ell_*^+)\times\cT(x(\ell^{\Phy}),\ell_*^+)$. (The uniqueness follows since the $\cR(x(\ell^{\Phy}),\ell_*^+),\cT(x(\ell^{\Phy}),\ell_*^+)$ are disjoint.) Similarly, there is a unique choice of $(\ell'_*)^+\in\cL_{*,x'(\ell^{\Phy})}^+$ so that $(\rho(\ell^{\Phy}),\tau(\ell^{\Phy}))\in \cR(x'(\ell^{\Phy}),\ell_*^+)\times\cT(x'(\ell^{\Phy}),\ell_*^+)$. Moving to $E^{\Phy}_o$, we have $\ell_*^{\Phy},(\ell'_*)^{\Phy}\subset F_\rho^{\Phy}\cap F_\tau^{\Phy}$ by \cref{prop:all-lines}. Moving from the affine space $E^{\Phy}_o$ to the projective space $E^{\Phy}$, these two lines $\ell_*^{\Phy},(\ell'_*)^{\Phy}$ lie in the 2-plane $F_\rho^{\Phy}\cap F_\tau^{\Phy}\cong \P_{\R}^2$ and thus they intersect.

For each $\rho\in\cR$, we will perform the above procedure for all lines $\ell^{\Phy}\in\cL_1(\rho)$. This produces many lines in $F^{\Phy}_\rho\cong\P_{\R}^3$ with many pairs intersecting. We will then be in a position to apply \cref{lem:line-structure-new} to find structure in this configuration.

In detail, fix $\rho\in\cR$. Define $\cX_\rho\subseteq \phi(H)\cap F_\rho^{\Phy}$ to be the set of $x\in\phi(H)\cap F_{\rho}^{\Phy}$ for which there is a line $\ell^+_{*}\in\cL^+_{*,x}$ satisfying $\rho\in\cR(x,\ell^+_*)$. By the disjointness of the sets $\cR(x,\ell^+_*)$, there is a unique line with these properties for each $x\in\cX_\rho$. Call this line $\ell^+_{*,x}$.

Moving to physical space, consider the list of lines $\ell^{\Phy}_{*,x}$ for $x\in\cX_\rho$, all lying in $F_\rho^{\Phy}\cong\P_{\R}^3$. (Note that this list may have repetitions.) We apply \cref{lem:line-structure-new} to this list of lines with parameter $\eta=N^{-2\epsilon'}$. This produces sets $\fP_\rho,\fF_\rho,\fR_\rho,\fL_\rho$ of points, planes, reguli, and lines, respectively, each of size $\lesssim N^{16\epsilon'}$.

Define a decomposition
\[\cL_1(\rho)=\cL^{\mathrm{pt}}(\rho)\cup\cL^{\mathrm{pln}}(\rho)\cup\cL^{\mathrm{reg}}(\rho)\cup\cL^{\mathrm{line}}(\rho)\cup\cL^{\mathrm{sparse}}(\rho),\]
as follows. For each $\ell^{\Phy}\in\cL_1(\rho)$, note that the definition of $\cL_1$ implies that $x(\ell^{\Phy}),x'(\ell^{\Phy})\in\cX_\rho$. As argued above, the lines $\ell^{\Phy}_{*,x(\ell^{\Phy})},\ell_{*,x'(\ell^{\Phy})}^{\Phy}$ intersect. Now we say
\begin{itemize}
    \item $\ell^{\Phy}\in\cL^{\mathrm{pt}}(\rho)$ if $\ell^{\Phy}_{*,x(\ell^{\Phy})}\cap\ell_{*,x'(\ell^{\Phy})}^{\Phy}=\{p\}$ for some point $p\in\fP_\rho$;
    \item $\ell^{\Phy}\in\cL^{\mathrm{pln}}(\rho)$ if $\ell^{\Phy}_{*,x(\ell^{\Phy})},\ell_{*,x'(\ell^{\Phy})}^{\Phy}\subset F$ for some plane $F\in\fF_\rho$;
    \item $\ell^{\Phy}\in\cL^{\mathrm{reg}}(\rho)$ if $\ell^{\Phy}_{*,x(\ell^{\Phy})},\ell_{*,x'(\ell^{\Phy})}^{\Phy}\subset R$ for some regulus $R\in\fR_\rho$;
    \item $\ell^{\Phy}\in\cL^{\mathrm{line}}(\rho)$ if either $\ell^{\Phy}_{*,x(\ell^{\Phy})}$ or $\ell_{*,x'(\ell^{\Phy})}^{\Phy}$ lies in $\fL_\rho$; and
    \item $\ell^{\Phy}\in\cL^{\mathrm{sparse}}(\rho)$ otherwise.
\end{itemize}

First, we bound $\cL^{\mathrm{sparse}}(\rho)$. We note that $\ell^{\Phy}\in\cL_1(\rho)$ corresponds to an intersecting pair $\ell^{\Phy}_{*,x(\ell^{\Phy})},\ell_{*,x'(\ell^{\Phy})}^{\Phy}$. Then $\ell^{\Phy}$ can be uniquely identified as $\aff\{x(\ell^{\Phy}),x'(\ell^{\Phy})\}$ from this pair of lines. \cref{lem:line-structure-new} guarantees that there are $\lesssim \eta|\cX_\rho|^2$ pairs of intersecting lines which are not captured by one of $\fP_\rho,\fF_\rho,\fR_\rho,\fL_\rho$, implying that  
\[|\cL^{\mathrm{sparse}}(\rho)|\lesssim \eta|\cX_\rho|^2\leq N^{-2\epsilon'}|\phi(H)\cap F_\rho^{\Phy}|^2\lesssim N^{-2\epsilon'}k^2.\]

Second, we bound $\cL^{\mathrm{pln}}(\rho)$. Note that $x\in\ell^{\Phy}_{*,x}$ for each $x\in\cX_\rho$, as $\ell^+_{*,x}\subset F^+_x$. Thus if $\ell_{*,x},\ell_{*,x'}$ are contained in a plane $F\in\fF_\rho$, we have $x,x'\in F\cap\phi(H)$. However $|F\cap\phi(H)|\leq N^{2/3}$, as each plane in $\R^3$ contains at most $N^{2/3}$ points of $\cP$. Now each $\ell^{\Phy}\in\cL^{\mathrm{pln}}(\rho)$ can be uniquely identified as $\aff\{x(\ell^{\Phy}),x'(\ell^{\Phy})\}$ where $x(\ell^{\Phy}),x'(\ell^{\Phy})$ are both elements of one of $|\fF_\rho|$ sets of size at most $N^{2/3}$. Thus
\[|\cL^{\mathrm{pln}}(\rho)|\leq (N^{2/3})^2|\fF_\rho|\lesssim N^{4/3+16\epsilon'}.\]

Third, we study $\cL^{\mathrm{reg}}(\rho)$. Decompose $\cL^{\mathrm{reg}}(\rho)=\cL^{\textrm{reg-rich}}(\rho)\sqcup\cL^{\textrm{reg-poor}}(\rho)$ where $\ell^{\Phy}\in\cL^{\mathrm{reg}}(\rho)$ lies in $\cL^{\textrm{reg-rich}}(\rho)$ if the lines $\ell^{\Phy}_{*,x(\ell^{\Phy})},\ell_{*,x'(\ell^{\Phy})}^{\Phy}$ are both contained in some regulus $R\in\fR_\rho$ with $|R\cap\phi(H)|\geq N^{-16\epsilon'}k$. Say $\ell^{\Phy}$ lies in $\cL^{\textrm{reg-poor}}(\rho)$ otherwise.

Similar to the previous case, if $\ell^{\Phy}_{*,x(\ell^{\Phy})},\ell_{*,x'(\ell^{\Phy})}^{\Phy}$ are contained in a poor regulus $R\in\fR_\rho$, we have $x,x'\in R\cap\phi(H)$, where the final set has size less than $N^{-16\epsilon'}k$. Thus
\[|\cL^{\textrm{reg-poor}}(\rho)|\leq (N^{-16\epsilon'}k)^2|\fR_\rho|\lesssim N^{-16\epsilon'}k^2.\]

Define $\cL_{\cR}^{\textrm{reg-rich}}=\bigcup_{\rho\in\cR}\cL^{\textrm{reg-rich}}(\rho)$. We apply \cref{lem:rich-reguli} to $\cR$, the sets $H\cap\phi^{-1}(F_\rho^{\Phy})$, the sets of reguli $\phi^{-1}(\fR_\rho)$, with parameter $\eta=N^{-16\epsilon'}$. Then for each line in $\ell^{\Phy}\in \cL^{\textrm{reg-rich}}_\cR$ there is a regulus $R\in\phi^{-1}(\fR_{\rho(\ell^{\Phy})})$ so that $\phi^{-1}(x(\ell^{\Phy})),\phi^{-1}(x'(\ell^{\Phy}))\in R$. As $\ell^{\Phy}$ is uniquely identified from $(\phi^{-1}(x(\ell^{\Phy})),\phi^{-1}(x'(\ell^{\Phy})))$, we see that \cref{lem:rich-reguli} gives the bound
\[|\cL_{\cR}^{\textrm{reg-rich}}|\lesssim N^{10/3+16\epsilon'}\log^2N.\]

Fourth, we bound $\cL^{\mathrm{line}}(\rho)$. Given $\rho$ and a line $\ell\in\fL_\rho$, there are at most $N^{2/3}$ choices for $x\in\cX_\rho$ so that $\ell=\ell^{\Phy}_{*,x}$. Indeed, $x\in\ell^{\Phy}_{*,x}$ and $|\phi(H)\cap\ell|\leq N^{2/3}$. Now if $\ell^{\Phy}\in\cL^{\mathrm{line}}(\rho)$ we know that $x'(\ell^{\Phy})\in\phi(H)\cap F_\rho^{\Phy}$, a set of size at most $2k$. Thus we have shown the bound
\[|\cL^{\mathrm{line}}(\rho)|\leq|\fL_\rho|\cdot N^{2/3}\cdot 2k.\]

Define $\cL^{\mathrm{pt}}_\cR=\bigcup_{\rho\in\cR}\cL^{\mathrm{pt}}(\rho)$ and similarly for the other sets. Putting everything together and using the bounds \cref{eq:k-k'-lower,eq:l-lower-r-t,eq:l-lower}, we have shown
\begin{align*}
|\cL^{\mathrm{pln}}_\cR|&=\sum_{\rho\in\cR} |\cL^{\mathrm{pln}}(\rho)|\lesssim N^{4/3+16\epsilon'}|\cR|\lesssim N^{-2\epsilon+18\epsilon'}k^2|\cR|\leq N^{-2\epsilon+19\epsilon'}|\cL|\\
|\cL^{\mathrm{reg}}_\cR|&=|\cL^{\textrm{reg-rich}}_\cR|+\sum_{\rho\in\cR} |\cL^{\textrm{reg-poor}}(\rho)|\lesssim N^{10/3+16\epsilon'}\log^2N+N^{-16\epsilon'}k^2|\cR|\\
&\qquad\qquad\qquad\qquad\qquad\leq \paren{N^{-\epsilon+17\epsilon'}\log^2N+N^{-15\epsilon'}}|\cL|\\
|\cL^{\mathrm{line}}_\cR|&=\sum_{\rho\in\cR} |\cL^{\mathrm{line}}(\rho)|\lesssim N^{2/3+16\epsilon'}k|\cR|\lesssim N^{-\epsilon+17\epsilon'}k^2|\cR|\leq N^{-\epsilon+18\epsilon'}|\cL|\\
|\cL^{\mathrm{sparse}}_\cR|&=\sum_{\rho\in\cR} |\cL^{\mathrm{sparse}}(\rho)|\lesssim N^{-2\epsilon'}k^2|\cR|\lesssim N^{-\epsilon'}|\cL|.
\end{align*}

We can write
\[\cL_1=\cL^{\mathrm{pt}}_{\cR}\cup\cL^{\mathrm{pln}}_{\cR}\cup\cL^{\mathrm{reg}}_{\cR}\cup\cL^{\mathrm{line}}_{\cR}\cup\cL^{\mathrm{sparse}}_{\cR}.\]
By our choice of $\epsilon=30\epsilon'$, we have that the latter four sets all have size at most $N^{-\epsilon'}|\cL|$, as desired.
It remains to bound $\cL^{\mathrm{pt}}_{\cR}$.

\vspace{3pt}\noindent\textbf{The contribution of very rich points:} Define the sets $\cL^{\mathrm{pt}}_{\cT},\cL^{\mathrm{pln}}_{\cT},\cL^{\mathrm{reg}}_{\cT},\cL^{\mathrm{line}}_{\cT},\cL^{\mathrm{sparse}}_{\cT}$ symmetrically, replacing $(\cR,k)$ with $(\cT,k')$. In more detail, for each $\tau\in\cT$, define $\cX_\tau\subseteq\phi(H)\cap F_\tau^{\Phy}$ to be the set of $x\in\phi(H)\cap F_\tau^{\Phy}$ for which there is a line $\ell^+_*\in\cL_{*,x}^+$ satisfying $\tau\in\cT(x,\ell^+_*)$. As before, this line is unique for each $x\in\cX_\tau$ and it is the same line that we previously called $\ell^+_{*,x}$. Applying \cref{lem:line-structure-new} now to the list of lines $\ell^{\Phy}_{*,x}$ for $x\in\cX_\tau$ in $F_\tau^{\Phy}\cong\PP_{\R}^3$, we define these five sets analogously and prove the same bound of $\lesssim N^{-\epsilon'}|\cL|$ on the sizes of the latter four. Thus it remains to bound $\cL^{\mathrm{pt}}_{\cR}\cap\cL^{\mathrm{pt}}_{\cT}$.

To do this, we apply \cref{lem:rich-points-sum} to $H,\cR,\cT,k,k'$ with parameter $M\sim N^{16\epsilon'}$ chosen so that the $\fP_\rho,\fP_\tau$ are sets of at most $M$ points. For each line $\ell^{\Phy}\in \cL^{\mathrm{pt}}_{\cR}\cap\cL^{\mathrm{pt}}_{\cT}$, the corresponding lines $\ell^{\Phy}_{*,x(\ell^{\Phy})},\ell_{*,x'(\ell^{\Phy})}^{\Phy}\in\cL_{x,*}$ have a unique intersection lying in $\fP_\rho$. However, the unique intersection point of these lines also lies in $\fP_\tau$, showing that $\fP_\rho\cap\fP_\tau\neq\emptyset$. Thus \cref{lem:rich-points-sum} gives the bound
\begin{align*}
|\cL^{\mathrm{pt}}_{\cR}\cap\cL^{\mathrm{pt}}_{\cT}|\lesssim N^{10/3}\log^2N&+\left(N^{2/3}k'+A^{-1}N^{1/3}(k')^2\right)N^{16\epsilon'}\abs{\cT}\\&+A^2N^{24\epsilon'}N^{-1/3}\left(\abs{\cR}^{1/2}\abs{\cT}+\abs{\cR}\abs{\cT}^{1/2}\right)
\end{align*}
for any $A\geq C_0N^{1/3}$, where $C_0\sim 1$ is the constant in \cref{lem:rich-points-sum}.

We set $A=N^{1/3+25\epsilon'}$. Then using \cref{eq:R-T-upper} together with \cref{eq:l-lower-r-t,eq:k-k'-lower,eq:l-lower}, we see that
\begin{align*}
|\cL^{\mathrm{pt}}_{\cR}\cap\cL^{\mathrm{pt}}_{\cT}|
&\lesssim N^{-\epsilon'}|\cL|+N^{-\epsilon+17\epsilon'}(k')^2|\cT|+N^{-9\epsilon'}(k')^2|\cT|\\ &\qquad\qquad\qquad\qquad\qquad+N^{1/3+74\epsilon'}\left(\abs{\cR}^{1/2}\abs{\cT}+\abs{\cR}\abs{\cT}^{1/2}\right)\\
&\lesssim N^{-\epsilon'}|\cL|+N^{5/6+75\epsilon'}|H|^{1/2}\paren{k^{-3/4}|\cT|+(k')^{-3/4}|\cR|}\\
&\lesssim N^{-\epsilon'}|\cL|+N^{11/6+75\epsilon'}(N^{-2/3-\epsilon+\epsilon'})^{11/4}\paren{(k')^2|\cT|+k^{2}|\cR|}\\
&\lesssim N^{-\epsilon'}|\cL|.
\end{align*}

Combining all the bounds, we have proven that $|\cL|\lesssim N^{-\epsilon'}|\cL|$ which is a contradiction for $N$ sufficiently large. This contradiction proves the desired bound $\abs{\cL_{[2,C]}(\phi(\cP\times\cP))}\leq N^{10/3+\epsilon}$.
\end{proof}

\bibliographystyle{amsplain0}
\bibliography{ref}

\appendix
\section{Deferred proofs from \texorpdfstring{\cref{sec:approx-complete-intersection}}{Section 4}}
\label{sec:appendix-jacobian}

\jacobian*

\begin{proof}
The first part of this result is standard commutative algebra. Localizing at a prime ideal produces a local ring and the quotient of a local ring remains local. Now the (Krull) dimension of $R$ is the length of the longest chain of prime ideals $\sang0\subseteq\fp_0\subsetneq\fp_1\subsetneq\cdots\subsetneq\fp_s\subsetneq R$. These are in one-to-one correspondence with chains of prime ideas $I\subseteq \fq_0\subsetneq \fq_1\subsetneq\cdots\subsetneq\fq_s\subseteq I(V)\subset\CC[x_1,\ldots,x_n]$. By Hilbert's Nullstellensatz, these in turn are in one-to-one correspondence with chains of irreducible varieties $V\subseteq V_s\subsetneq V_{s-1}\subsetneq\cdots\subsetneq V_1\subsetneq V_0\subseteq Z(I)\subset\CC^n$. By hypothesis, the longest such chain has length $(n-r)-\dim V=\codim_{\CC^n}(V)-r$.

To prove the second part of this result, we apply \cite[Theorem 16.19]{Eis95} with $k=\CC$ and $P=I(V)$. This essentially gives the desired conclusion; we just need to check the following:
\begin{enumerate}
    \item the codimension of $I\CC[x_1,\ldots,x_n]_{I(V)}$ in $\CC[x_1,\ldots, x_n]_{I(V)}$ is equal to $r$;
    \item a matrix of polynomials in $\CC[x_1,\ldots, x_n]$ has rank $r$ on a Zariski-dense open subset of $V$ if and only if it has rank $r$ when taken modulo $I(V)$.
\end{enumerate}
The codimension of $I\CC[x_1,\ldots,x_n]_{I(V)}$ in $\CC[x_1,\ldots, x_n]_{I(V)}$ is defined (see, e.g., \cite[Chapter 9]{Eis95}) to be the minimum codimension of prime ideals $\fp\subset \CC[x_1,\ldots, x_n]_{I(V)}$ which contain $I\CC[x_1,\ldots, x_n]_{I(V)}$. The prime ideals in $\CC[x_1,\ldots, x_n]_{I(V)}$ which contain $I\CC[x_1,\ldots, x_n]_{I(V)}$ are in one-to-one correspondence with the prime ideals of $\CC[x_1,\ldots, x_n]$ which contain $I$ and are contained in $I(V)$. These are in one-to-one correspondence with irreducible varieties which contain $V$ and are contained in $Z(I)$. Consider such a $U$ with $V\subseteq U\subseteq Z(I)$ so $\fp= I(U)\CC[x_1,\ldots, x_n]_{I(V)}$.

Now the codimension of $\fp$ is defined to be the Krull dimension of $(\CC[x_1,\ldots, x_n]_{I(V)})_{\fp}=\CC[x_1,\ldots,x_n]_{I(U)}$. As before, chains of prime ideals $\sang0\subseteq\fp_0\subsetneq\fp_1\subsetneq\cdots\subsetneq\fp_s\subsetneq \CC[x_1,\ldots,x_n]_{I(U)}$ are in one-to-one correspondence with chains of irreducible varieties containing $U$ in $\CC^n$, i.e., $U\subseteq U_s\subsetneq\cdots\subsetneq U_1\subsetneq U_0\subseteq\CC^n$. By assumption, every irreducible variety $U$ which contains $V$ and is contained in $Z(I)$ has codimension at least $r$ and there exists at least one such $U$ with codimension exactly $r$, proving (1).

To see (2), note that both are equivalent to the following: the matrix has an $r\times r$ minor whose determinant does not lie in $I(V)$.
\end{proof}

\valuation*

\begin{proof}
We can write
\[\CC[x_1,\ldots, x_n]_{I(V)}/\sang{f_1,f_2}\CC[x_1,\ldots, x_n]_{I(V)} = R/\sang{[f_2]}.\]
If $[f_2]\neq 0$, then $\langle [f_2]\rangle = \langle g^{v_R(f_2)}\rangle$ and so
\[m_V(f_1,f_2)=\len_{\CC[x_1,\ldots, x_n]_{I(V)}}\paren{R/\langle g^{v_R(f_2)}\rangle}.\]
Now for any $\CC[x_1,\ldots, x_n]_{I(V)}$-submodule of $R/\sang{g^{v_R(f_2)}}$, the action of multiplication by $f_1$ is trivial since $[f_1]=0\in R$. This means that any such $\CC[x_1,\ldots, x_n]_{I(V)}$-submodule is in fact an $R$-submodule. We conclude that
\[m_V(f_1,f_2)=\len_{R}\paren{R/\langle g^{v_R(f_2)}\rangle}.\]
Thus this quantity is the length of the longest chain of ideals in $R/\sang{g^{v_R(f_2)}}$ which is equal to the length of the longest chain of ideals in $R$ which contain $\sang{g^{v_R(f_2)}}$. Since every ideal of $R$ is of the form $\sang 0$ or $\sang{g^k}$ for $k\geq 0$ (see, e.g. \cite[Proposition 11.1]{Eis95}), we conclude that the longest such chain is
\[\sang{g^{v_R(f_2)}}\subsetneq \sang{g^{v_R(f_2)-1}}\subsetneq \cdots\subseteq \sang{g^{1}}\subsetneq\sang{g^0}.\]
Thus $m_V(f_1,f_2)=v_R(f_2)$.
\end{proof}

\section{Deferred proofs from \texorpdfstring{\cref{sec:csm-for-flats}}{Section 5}}
\label{sec:appendix-grobner}

In this section we prove \cref{claim:radical-bounded,claim:grobner-saturation}. Both proofs use the theory of Gr\"obner bases, which we briefly review here. For a more thorough introduction, see the textbook \cite{BW93}.

Let $\F$ be an arbitrary field. Given $\alpha\in\Z_{\geqslant 0}^n$, we write $x^\alpha$ for the \emph{term} $x_1^{\alpha_1}\cdots x_n^{\alpha_n}$ and $|\alpha|=\alpha_1+\cdots+\alpha_n$. A \emph{term order} is a total order $\prec$ on terms in $x_1,\ldots,x_n$ such that $x^\alpha\prec x^\beta$ implies $x^{\alpha+\gamma}\prec x^{\beta+\gamma}$ and such that $1\prec x^\alpha$ for all $|\alpha|>0$. Given a nonzero polynomial $f\in \F[x_1,\ldots,x_n]$, we write $\LT(f),\LM(f),\LC(f)$ for its \emph{head term}, \emph{head monomial}, and \emph{head coefficient} with respect to a given term order. That is, the largest term which appears in $f$ with nonzero coefficient, and the corresponding monomial and coefficient of that term. A term order is \emph{graded} if $x^\alpha\prec x^\beta$ whenever $|\alpha|<|\beta|$.

For an ideal $I\subseteq \F[x_1,\ldots,x_n]$, define $\LT(I)=\sang{\LT(f):f\in I\setminus\{0\}}$. A \emph{Gr\"obner basis} for $I$ (with respect to a given term order) is a set of nonzero polynomials $g_1,\ldots,g_t\in I$ such that $\sang{\LT(g_1),\ldots,\LT(g_t)}=\LT(I)$. In other words, such that for each nonzero $f\in I$ we have $\LT(f)$ is divisible by $\LT(g_i)$ for some $i\in[t]$.

Recall that for an ideal $I\subseteq \F[x_1,\ldots,x_n]$ and an element $g\in\F[x_1,\ldots,x_n]$, the \emph{radical} $\sqrt{I}$ is the ideal consisting of those $f\in\F[x_1,\ldots,x_n]$ such that $f^k\in I$ for some $k\geq 1$ and the \emph{saturation} $(I:g^\infty)$ is the ideal consisting of those $f\in\F[x_1,\ldots,x_n]$ such that $fg^k\in I$ for some $k\geq 1$.

Let $f_1,\ldots,f_s\in \F[x_0,\ldots,x_n]$ be homogeneous polynomials of degree at most $d$. With respect to any graded order, it is known that $I=\sang{f_1,\ldots,f_s}$ has a Gr\"obner basis $g_1,\ldots,g_t$ consisting of (homogeneous) polynomials of degree at most $2((d^2)/2+d)^{2^{n-1}}$ \cite[Theorem 8.2]{Dub90}. We now explain how to construct such a Gr\"obner basis using a bounded number of operations. We say that polynomials $g_1,\ldots,g_s$ can be computed from $f_1,\ldots,f_t$ via $C$ \emph{elementary operations} if the coefficients of the $g_i$ can be expressed in terms of the coefficients of the $f_j$ as an arithmetic expression consisting of at most $C$ additions, subtractions, multiplications, and divisions.

\begin{claim}
\label{clm:grobner-basis-alg}
Let $f_1,\ldots,f_t\in \F[x_1,\ldots,x_n]$ be polynomials of degree at most $d$. Then a Gr\"obner basis $g_1,\ldots,g_s$ for $\sang{f_1,\ldots,f_t}$ (with respect to any term order) can be computed via at most $O_{d,n}(1)$ elementary operations. Furthermore, we may ensure that the $g_1,\ldots,g_s$ all have degrees $O_{d,n}(1)$, all have leading coefficient 1, and that $s=O_{d,n}(1)$.

In this Gr\"obner basis, each 
\[g_i=\sum_{j=1}^t c_{ij}f_j\]
where the $c_{ij}\in\F[x_1,\ldots,x_n]$ are polynomials of degree $O_{d,n}(1)$ whose coefficients can be computed from the coefficients of $f_1,\ldots,f_t$ via at most $O_{d,n}(1)$ elementary operations. Furthermore, in each sum, only $O_{d,n}(1)$ of the summands are nonzero.
\end{claim}

\begin{proof}
For a polynomial $f\in \F[x_1,\ldots,x_n]$, let $\tilde f\in \F[x_0,\ldots,x_n]$ denote its homogenization $\tilde f=x_0^{\deg f} f(x_1/x_0,\ldots,x_n/x_0)$. For a polynomial $f \in \F[x_0,\ldots,x_n]$, let $\ol f\in \F[x_1,\ldots,x_n]$ denote its dehomogenization $\ol f=f(1,x_1,\ldots,x_n)$.

Write $I=\sang{f_1,\ldots,f_t}\subseteq \F[x_1,\ldots,x_n]$ and $J=\sang{\tilde f_1,\ldots,\tilde f_t}\subseteq \F[x_0,\ldots,x_n]$. Note that $J$ is not necessarily the homogenization of $I$; we will see later that the homogenization $\sang{\tilde f:f\in I}$ is equal to the saturation $(J:x_0^\infty)$. Let $<$ be any term order on the terms in $x_1,\ldots,x_n$. Define the graded term order $\prec$ on the terms in $x_0,\ldots,x_n$ by $x^\alpha\prec x^\beta$ if $|\alpha|<|\beta|$ or if $|\alpha|=|\beta|$ and $\ol{x^{\alpha}}<\ol{x^\beta}$. By \cite[Theorem 8.2]{Dub90}, there exists a Gr\"obner basis $g_1,\ldots,g_s$ for $J$ with respect to $\prec$ where $g_1,\ldots,g_s$ are (homogeneous) polynomials of degree at most $D_{\max}:=2(d^2/2+d)^{2^{n-1}}$. We first show how to compute $g_1,\ldots,g_s$ in $O_{d,n}(1)$ elementary operations.

For $D\geq 1$, define the matrix $M_D$ with entries in $\F$ as follows. The columns of $M_D$ are indexed by the term $x^\alpha$ with $|\alpha|=D$, ordered in decreasing order by $\prec$. Thus any row of this matrix corresponds to a homogeneous degree $D$ polynomial in $\F[x_0,\ldots,x_n]$. Define $M_D$ so that its rows correspond to the polynomials $\tilde f_im$ for $i\in[t]$ and $m$ a term in $x_0,\ldots,x_n$ of degree $D-\deg f_i$. Now the rowspace of $M_D$ corresponds exactly to the polynomials in $J_D$: the degree $D$ homogeneous part of $J$.\footnote{Suppose one has $f=\sum_{i=1}^t a_i\tilde f_i$ where $f$ is homogeneous of degree $D$ but the $a_i$ are not necessarily homogeneous. Then setting $b_i$ be the homogeneous degree $(D-\deg(f_i))$ part of $a_i$, we still have $f=\sum_{i=1}^t b_i\tilde f_i$.}

The matrix $M_D$ can be put into reduced row echelon form via $O_{d,n}(1)$ elementary operations: in particular, the only operations we perform are taking the sum, difference, product, or quotient of entries of the matrix. At the end of this process, let $g_1,\ldots,g_{s_D}$ be the homogeneous degree $D$ polynomials corresponding to the nonzero rows of the resulting matrix. These polynomials form a basis for $J_D$ as an $\F$-vector space. Furthermore, since the columns are ordered by $\prec$, the position of the leftmost nonzero entry corresponds to the head term of the polynomial. In particular, the set of head terms of $J_D$ is the same as the set of head terms corresponding to this $\F$-vector space basis.

Now each row of the reduced row echelon form of $M_D$ is an $\F$-linear combination of the rows of $M_D$. This means that each $g_1,\ldots,g_{s_D}$ is an $\F$-linear combination of the polynomials $\tilde f_i m$ where $m$ is a term in $x_0,\ldots,x_n$. In addition, the number of nonzero summands in this linear combination is at most $\rank M_D$. We can bound $\rank M_D$ by the number of columns in the matrix which is $O_{D,n}(1)$. For $D\leq D_{\max}=O_{d,n}(1)$, thus the number of nonzero summands is $O_{d,n}(1)$. Furthermore, clearly each of these coefficients can be computed from the coefficients of $f_1,\ldots,f_t$ via at most $O_{d,n}(1)$ elementary operations. 

Apply this process for $1\leq D\leq D_{\max}$ and let $g_1,\ldots,g_s\in J$ be the collection of polynomials produced. For every polynomial in $J$ with degree at most $D_{\max}$, its head term is equal to $\LT(g_i)$ for some $i\in[s]$. Since we know that $J$ has some Gr\"obner basis consisting of polynomials of degree at most $D_{\max}$, we conclude that $g_1,\ldots,g_s$ is also a Gr\"obner basis for $J$. In addition, $s\leq\sum_{D=1}^{D_{\max}}\rank M_D=O_{d,n}(1)$.

Finally, $\ol g_1,\ldots,\ol g_s\in \ol J=I$ is a Gr\"obner basis for $I$ by \cite[Lemma 7]{MR13}. For completeness, we include a self-contained proof here. For any nonzero $f\in I$, clearly $f=\ol{\tilde f}$. In addition, we claim $\tilde f\in(J:x_0^\infty)$. Indeed, $f=\sum_i f_ih_i$ for some polynomials $h_1,\ldots,h_t\in\F[x_1,\ldots,x_n]$. Then
\[\tilde f=x_0^{\deg f}\sum_{i=1}^t f_i(x_1/x_0,\ldots,x_n/x_0)h_i(x_1/x_0,\ldots,x_n/x_0).\]
Multiplying through by $x_0^k$ where $k=\max_i(\deg f_i+\deg h_i)-\deg f$, we conclude that
\[\tilde f x_0^k=\sum_{i=1}^t \tilde f_i  x_0^{k+\deg f-\deg f_i}h_i(x_1/x_0,\ldots,x_n/x_0),\]
which implies that $\tilde f\in (\sang{\tilde f_1,\ldots,\tilde f_t}:x_0^\infty)$.

Furthermore, $\LT_{<}(f)=\ol{\LT_{\prec}(\tilde f)}$. This follows from the definition: since $\tilde f$ is homogeneous, all its terms have the same degree so $\prec$ orders them by dehomogenizing and ordering based on $<$. Then we clearly have
\[\LT_{<}(f)=\ol{\LT_{\prec}(\tilde f)}=\ol{\LT_{\prec}(x_0^k\tilde f)}\]
for every $k$. Taking $k$ such that $\tilde f x_0^k\in J$, we conclude that $\LT_{\prec}(x_0^k\tilde f)$ is divisible by $\LT_{\prec}(g_i)$ for some $i\in[s]$, implying that $\LT_{<}(f)$ is divisible by $\ol{\LT_{\prec}(g_i)}=\LT_<(\ol g_i)$. As this is true for every nonzero $f\in I$, we conclude that $\ol g_1,\ldots,\ol g_s$ is a Gr\"obner basis for $I$. As $M_D$ was put in reduced row echelon form, these polynomials have head coefficient equal to 1 and all have degree at most $D_{\max}=O_{d,n}(1)$, as desired.
\end{proof}

For $A\subseteq[n]$, we write $\F[x_I]\subseteq \F[x_1,\ldots,x_n]$ for the polynomial ring formed by adjoining the variables $x_i$ with $i\in A$.

As many operations in commutative algebra can be expressed in terms of Gr\"obner bases, \cref{clm:grobner-basis-alg} implies that these operations can be computed via a bounded number of elementary operations. We collect the results we will use here.

\begin{claim}
\label{clm:elimination-alg}
Let $f_1,\ldots,f_t\in \F[x_1,\ldots,x_n]$ be polynomials of degree at most $d$. Set $I=\sang{f_1,\ldots,f_t}$. For any $A\subseteq[n]$, a Gr\"obner basis $g_1,\ldots,g_s$ for $I\cap \F[x_A]$ consisting of polynomials of degree at most $O_{d,n}(1)$ can be computed via at most $O_{d,n}(1)$ elementary operations.
\end{claim}

\begin{proof}
Let $\prec$ be a term order on terms in $x_1,\ldots, x_n$ so that every term in the variables $\{x_i:i\in A\}$ precedes every nontrivial term in the variables $\{x_i:i\in[n]\setminus A\}$. By \cref{clm:grobner-basis-alg}, we can compute a Gr\"obner basis $g_1,\ldots,g_s$ for $I$ with respect to this term order. Consider the set of these basis elements which lie in $\F[x_A]$. By \cite[Proposition 6.15]{BW93}, this set is a Gr\"obner basis for $I\cap \F[x_A]$.
\end{proof}

\begin{claim}
\label{clm:intersection-alg}
Let $f_1,\ldots,f_t,g_1,\ldots,g_s\in \F[x_1,\ldots,x_n]$ be polynomials of degree at most $d$. Set $I=\sang{f_1,\ldots,f_t}$ and $J=\sang{g_1,\ldots,g_s}$. A Gr\"obner basis for $I\cap J$ consisting of polynomials of degree at most $O_{d,n}(1)$ can be computed via at most $O_{d,n}(1)$ elementary operations.
\end{claim}

\begin{proof}
Define $K=\sang{yf_1,\ldots,yf_t,(1-y)g_1,\ldots,(1-y)g_s}\subseteq \F[x_1,\ldots,x_n,y]$. By \cite[Proposition 6.19]{BW93}, we have $K\cap\F[x_1,\ldots,x_n]=I\cap J$. By \cref{clm:elimination-alg}, we can compute a Gr\"obner basis for $K\cap \F[x_1,\ldots,x_n]$, as desired.
\end{proof}

\begin{claim}
\label{clm:saturation-alg}
Let $f_1,\ldots,f_t,g\in \F[x_1,\ldots,x_n]$ be polynomials of degree at most $d$. Set $I=\sang{f_1,\ldots,f_t}$. A Gr\"obner basis for $(I:g^\infty)$ consisting of polynomials of degree at most $O_{d,n}(1)$ can be computed via at most $O_{d,n}(1)$ elementary operations.
\end{claim}

\begin{proof}
Define $J=\sang{f_1,\ldots,f_t,1-yg}\subseteq \F[x_1,\ldots,x_n,y]$. By \cite[Proposition 6.37]{BW93}, we have $J\cap\F[x_1,\ldots,x_n]=(I:g^\infty)$. By \cref{clm:elimination-alg}, we can compute a Gr\"obner basis for $J\cap \F[x_1,\ldots,x_n]$, as desired.
\end{proof}

\begin{claim}
\label{clm:gcd-alg}
Let $f,g\in \F[x_1,\ldots,x_n]$ be polynomials of degree at most $d$. Then $\lcm(f,g)$ can be computed via at most $O_{d,n}(1)$ elementary operations.
\end{claim}

\begin{proof}
Note that since $\F[x_1,\ldots,x_n]$ is a unique factorization domain, $\lcm(f,g)$ is well-defined (up to multiplication by a nonzero element of $\F$). Furthermore, $\sang{f}\cap\sang{g}=\sang{\lcm(f,g)}$ and by \cref{clm:grobner-basis-alg} a Gr\"obner basis for this intersection can be calculated in $O_{d,n}(1)$ elementary operations. Since this ideal is principal, the element of the Gr\"obner basis of minimal degree must be $\lcm(f,g)$.
\end{proof}

\begin{claim}
\label{clm:elimination-2-alg}
Let $f_1,\ldots,f_t\in \F(x_1,\ldots,x_k)[x_{k+1},\ldots,x_n]$ be polynomials of degree at most $d$ in $x_{k+1},\ldots,x_n$ whose coefficients are rational functions whose numerators and denominators are polynomials of degree at most $d$ in $x_1,\ldots,x_k$. Set $I=\sang{f_1,\ldots,f_t}$. A Gr\"obner basis for $I\cap \F[x_1,\ldots,x_n]$ consisting of polynomials of degree at most $O_{d,n}(1)$ can be computed via at most $O_{d,n}(1)$ elementary operations.
\end{claim}

\begin{proof}
By \cref{clm:grobner-basis-alg}, we can compute a Gr\"obner basis $g_1,\ldots,g_s$ for $I$ with respect to an arbitrary term order. This computation only involves $O_{d,n}(1)$ elementary operations, meaning that each coefficient of $g_1,\ldots,g_s$ can be expressed in terms of the coefficients of $f_1,\ldots,f_t$ via $O_{d,n}(1)$ additions, subtractions, multiplications, and divisions. Note that all of these coefficients are rational functions in $\F(x_1,\ldots,x_k)$; the coefficients of $f_1,\ldots,f_t$ have numerators and denominators which are polynomials of degree at most $d$, so the coefficients of $g_1,\ldots,g_s$ have numerators and denominators which are polynomial of degree $O_{d,n}(1)$. Furthermore, $f_1,\ldots,f_t$ have degree at most $O_{d,n}(1)$ as polynomials in $x_{k+1},\ldots,x_n$.

Clear the denominators of all the coefficients of $g_1,\ldots,g_s$. This results in a new Gr\"obner basis $h_1,\ldots,h_s$ for $I$, since all we have done is multiply by nonzero elements of $\F[x_1,\ldots,x_k]$ which are units in $\F(x_1,\ldots,x_k)[x_{k+1},\ldots,x_n]$. Next set $h=\lcm(\LC(h_1),\ldots,\LC(h_s))$. Note that here $\LC(h_1),\ldots,\LC(h_s)$ a priori lie in $\F(x_1,\ldots,x_k)$, but since we cleared the denominators, they actually lie in $\F[x_1,\ldots,x_k]$. By \cref{clm:gcd-alg}, the lcm can be computed via $O_{d,n}(1)$ elementary operations, as $\LC(h_1),\ldots,\LC(h_s)$ are polynomials of degree $O_{d,n}(1)$.

Now we view $h_1,\ldots,h_s$ as elements of $\F[x_1,\ldots,x_k][x_{k+1},\ldots,x_n]$; these are polynomials of degree $O_{d,n}(1)$ which can be computed from $f_1,\ldots,f_t$ via $O_{d,n}(1)$ elementary operations. Define $J=\sang{h_1,\ldots,h_s}\subseteq \F[x_1,\ldots,x_n]$. By \cite[Lemma 8.91]{BW93}, we have $I\cap \F[x_1,\ldots,x_n]=(J:h^\infty)$ so by \cref{clm:saturation-alg} we can compute a Gr\"obner basis for this ideal in $O_{d,n}(1)$ elementary operations.
\end{proof}

\begin{claim}
\label{clm:zero-dimensional-radical-alg}
Suppose $\F$ has characteristic zero. Let $f_1,\ldots,f_t\in \F[x_1,\ldots,x_n]$ be polynomials of degree at most $d$. Set $I=\sang{f_1,\ldots,f_t}$. If $I\cap \F[x_i]\neq\sang0$ for all $i\in[n]$, then a generating set for $\sqrt{I}$ consisting of polynomials of degree at most $O_{d,n}(1)$ can be computed via at most $O_{d,n}(1)$ elementary operations.
\end{claim}

\begin{proof}
By \cref{clm:elimination-alg}, we can compute a Gr\"obner basis for $I_i=I\cap \F[x_i]$. Since $\F[x_i]$ is a PID, $I_i=\sang{g_i}$ for some polynomial $g_i$. This polynomial must be the element of lowest degree in the Gr\"obner basis for $I_i$. The formal derivative $g'_i\in \F[x_i]$ can clearly be computed in $O_{d,n}(1)$ elementary operations. Define the squarefree part of $g_i$ to be $h_i=g_i/\gcd(g_i,g'_i)=\lcm(g_i,g_i')/g'_i$.(If $g_i$ is contstant, define $h_i=g_i$.) By \cref{clm:gcd-alg}, the lcm can be computed in $O_{d,n}(1)$ elementary operations, so clearly the division algorithm for univariate polynomials allows us to compute $h_i$ in $O_{d,n}(1)$ elementary operations. By \cite[Lemma 8.19]{BW93}, we have $\sqrt{I}=\sang{f_1,\ldots,f_t,h_1,\ldots,h_n}$.
\end{proof}

\begin{claim}
\label{clm:radical-alg}
Suppose $\F$ has characteristic zero. Let $f_1,\ldots,f_t\in \F[x_1,\ldots,x_n]$ be polynomials of degree at most $d$. Set $I=\sang{f_1,\ldots,f_t}$. A generating set for $\sqrt{I}$ consisting of polynomials of degree at most $O_{d,n}(1)$ can be computed via at most $O_{d,n}(1)$ elementary operations.
\end{claim}

\begin{proof}
We follow the algorithm given in \cite{Lap06}. Initialize $\cL_0=\{1\}$. We perform the following iteration for steps $i=1,2,3,\ldots$. Suppose we have $\cL_{i-1}$, a list of polynomials in  $\F[x_1,\ldots,x_n]$ of degree at most $d_{i-1}$. Pick $g\in\cL_{i-1}$ with $g\not\in\sqrt I$ (if such an element exists) and compute a Gr\"obner basis $g_1,\ldots,g_s$ for $J_i=(I:g^\infty)$. Let $A\subseteq[n]$ be a maximal subset so that $J\cap \F[x_A]=\sang 0$. Define $\tilde J_i=\sang{g_1,\ldots,g_s}\subseteq \F(x_A)[x_{[n]\setminus A}]$. (This is the polynomial ring with variables $\{x_j:j\in[n]\setminus A\}$ with coefficients in the field of rational functions $\F(x_A)$, i.e., rational functions in variables $\{x_j:j\in A\}$ and coefficients in $\F$.) Compute $\cL_i$, a Gr\"obner basis for $\sang{\cL_{i-1}}\cap(\sqrt{\tilde J_i}\cap \F[x_1,\ldots,x_n])$. Since $\F(x_A)$ has characteristic zero, we can compute generators for $\sqrt{\tilde J_i}$ in $O_{d_{i-1},n}(1)$ elementary operations. The polynomials involved in these generators have degree $O_{d_{i-1},n}(1)$. Then by \cref{clm:elimination-2-alg,clm:intersection-alg}, we can compute $\cL_i$ in $O_{d_{i-1},n}(1)$ elementary operations. These polynomials have degree at most $d_i=O_{d_{i-1},n}(1)$.

We repeat this process until every $g\in\cL_i$ satisfies $g\in\sqrt I$. By \cite{Lap06}, at this point we have $\sang{\cL}=\sqrt{I}$ and this process halts after at most $d^n$ iterations. Thus we have computed the desired generating set via at most $O_{d,n}(1)$ elementary operations.
\end{proof}

\grobnerRadical*

\begin{proof}
Recall the notation $[\cL]$ means that we view each element of $\cL$ as an element of $K(V)[y_1,\ldots, y_m]$. By hypothesis, these polynomials have degree at most $d$ and their coefficients are elements of $\CC[x_1,\ldots,x_n]$ with degree at most $D$.

By \cref{clm:radical-alg}, there is a generating set for $\sqrt{\sang{[\cL]}}$ consisting of polynomials of degree $O_{d,m}(1)$ which can be computed via at most $O_{d,m}(1)$ elementary operations. This means that each coefficient of each polynomial in the generating set can be obtained from the coefficients of the polynomials in $[\cL]$ via at most $O_{d,m}(1)$ additions, subtractions, multiplications, or divisions. Applying the same operations to the coefficient of $\cL$, we can represent each coefficient as a rational function where the numerator is a polynomial in $\CC[x_1,\ldots,x_n]$ of degree $O_{d,m}(D)$ and the denominator is a polynomial of the same degree that does not vanish identically on $V$. Now clearing the denominators (which is just multiplying by units in $K(V)[y_1,\ldots,y_m]$), we produce rational functions $a_1/b_1,\ldots,a_s/b_s$ where $a_1,\ldots,a_s\in\CC[x_1,\ldots,x_n,y_1,\ldots,y_m]$ and $b_1,\ldots,b_s\in\CC[x_1,\ldots,x_n]\setminus I(V)$ such that $\sang{[a_1/b_1],\ldots[a_s/b_s]}=\sqrt{\sang{[\cL]}}$. Given how the $a_i,b_i$ are produced from $\cL$, we see that each $b_i$ has degree $O_{d,m}(D)$ and each $a_i$ has degree $O_{d,m}(D)$ in $x_1,\ldots,x_n$ and degree $O_{d,m}(1)$ in $y_1,\ldots,y_m$.

Now since $[1/b_1],\ldots,[1/b_s]$ are units in $K(V)[y_1,\ldots,y_m]$, we conclude that $\sang{[a_1],\ldots,[a_s]}=\sqrt{\sang{[\cL]}}$. Defining $\tilde\cL=\set{a_1,\ldots,a_s}$ thus gives the desired conclusion.
\end{proof}

\grobnerSaturation*

\begin{proof}
Writing $f_1,\ldots,f_t$ for the elements of $\cL^{\new}$, by \cref{clm:grobner-basis-alg}, a monic Gr\"obner basis for $\sang{[f_1],\ldots,[f_t]}\subseteq K(V)[y_1,\ldots,y_m]$ can be computed in $O_{d,m}(1)$ elementary operations. In particular, the Gr\"obner basis consists of $\tilde g_1,\ldots,\tilde g_s$ of the form
\[\tilde g_i=\sum_j \tilde c_{ij} [f_j]\]
where the $\tilde c_{ij}\in K(V)[y_1,\ldots,y_m]$ are polynomials of degree $O_{d,m}(1)$ whose coefficients can be expressed in terms of the coefficients of the $[f_j]$ via $O_{d,m}(1)$ elementary operations. (Note that the coefficients in question lie in $K(V)$.) In the above, the sum is over the $O_{d,m}(1)$ terms which are nonzero. 

Viewing $f_1,\ldots,f_t$ as elements of $\CC[x_1,\ldots,x_n][y_1,\ldots,y_m]$ and applying the same operations to the coefficients of these polynomials (now the coefficients are elements of $\CC[x_1,\ldots,x_n]$ of degree at most $D$), we produce polynomials $a_{ij}\in\CC[x_1,\ldots,x_n,y_1,\ldots,y_m]$ and $b_{ij}\in\CC[x_1,\ldots,x_n]\setminus I(V)$ such that $\tilde c_{ij}=[a_{ij}/b_{ij}]$ and such that the $a_{ij},b_{ij}$ have degree $O_{d,m}(D)$ in $x_1,\ldots,x_n$ and the $a_{ij}$ have degree $O_{d,m}(1)$ in $y_1,\ldots,y_m$.

Clearing denominators, define
\[b_i=\prod_j b_{ij}\qquad\text{and}\qquad g_i=b_i\sum_j \frac{a_{ij}}{b_{ij}}f_j\]
so that $\tilde g_i=[g_i/b_i]$. Finally, set $h=b_1\cdots b_s\in\CC[x_1,\ldots,x_n]\setminus I(V)$. As each of the products are of length $O_{d,m}(1)$ and \cref{clm:grobner-basis-alg} implies that $s=O_{d,m}(1)$, we have the desired degree bound on $h$.

We next claim that $h$ satisfies \cref{eq:saturation}. Suppose that $f\in\CC[x_1,\ldots,x_n,y_1,\ldots,y_m]$ is such that $\psi(f) \in \psi(\sang{\cL^{\new}})[S^{-1}]$. (Recall that $\psi$ is the quotient map by $I(V)$ and $S=\CC[x_1,\ldots,x_n]\setminus I(V)$.) We need to show that $\psi(f)\psi(h)^k\in\psi(\sang{\cL^{\new}})$ for some $k\geq 1$. 

The assumption on $f$ implies that $[f]\in \sang{[\cL^{\new}]}$. Now consider applying the polynomial division algorithm to $[f]\in K(V)[y_1,\ldots,y_m]$ with respect to $\tilde g_1,\ldots,\tilde g_s$. As the latter set is a Gr\"obner basis, this will produce a decomposition
\[[f]=\sum_{i=1}^s \tilde d_i \tilde g_i\]
for some $\tilde d_i\in K(V)[y_1,\ldots,y_m]$. (See \cite[Definition 5.18 and Theorem 5.35]{BW93}.) The division algorithm, applied with respect to a Gr\"obner basis consisting of monic polynomials, only uses the operations of addition, subtraction, and multiplication of polynomials. Thus applying the same operations to $g_1/b_1,\ldots,g_s/b_s$ and $f$, we produce a decomposition
\[f\equiv\sum_{i=1}^s\frac{g_i}{b_i}\frac{c_i}{d_i}\pmod{I(V)}\]
where the rational functions $c_i/d_i$ are produced by applying the operations of addition, subtraction, and multiplication to $g_1/b_1,\ldots,g_s/b_s$ and $f$. In particular, this means that each denominator $d_i$ divides a product of powers of $b_1,\ldots, b_s$. Thus there exists some $k$ so that all the $d_i$ divide $(b_1\cdots b_s)^k= h^k$. Since each $g_i$ is defined to be a $\CC[x_1,\ldots,x_n,y_1,\ldots,y_m]$-linear combination of the $f_j$, we have $g_1,\ldots,g_s\in\sang{\cL^{\new}}$. Thus the above implies that $fh^k\in\sang{I(V),\cL^{\new}}$ and so $\psi(f)\psi(h)^k\in\psi(\sang{\cL^{\new}})$, as desired.
\end{proof}

\end{document}